\documentclass[11pt, a4paper]{article}
\usepackage{amsfonts,amssymb,amsmath,amsthm,amscd,epsfig}
\usepackage{stmaryrd}
\usepackage[pdfborder={0 0 0}]{hyperref}
\usepackage{pifont}
\usepackage{dsfont}
\usepackage{enumerate}
\usepackage{subcaption}
\usepackage{paralist}
\usepackage{accents}
\usepackage[all,2cell, color]{xy} \UseAllTwocells \SilentMatrices
\usepackage{verbatim}
\usepackage[normalem]{ulem}
\usepackage{tikz,pgfplots}
\usetikzlibrary{arrows.meta,calc,decorations.markings,math}

\usepackage{calc}
\usepackage{import}

\def\N                {\mathbb N}
\def\C                {\mathbb C}

\def\cala             {{\cal A}}
\def\calb             {{\cal B}}
\def\calc             {{\cal C}}

\def\cale             {{\cal E}}

\def\calm             {{\cal M}}
\def\caln             {{\cal N}}
\def \mac             {{\cal C}}
\def\mad              {{\cal B}}

\def\ii               {{\rm i}}
\def\oo               {{\otimes}}
\newcommand{\inv}     {{-1}}
\newcommand \low [2]  {{#1}_{({#2})}}
    
\def \Ob            {{\mathrm{Ob}}}
\def\End              {{\rm End}}
\def\Hom              {{\rm Hom}}

\def \Vec            {{\mathrm{Vec}}}
\def\Mod            {\mathrm{-Mod}}

\def\tr               {{\rm tr\;}}
\def\id               {{\rm id}}

\def\dim               {{\rm dim}}   
 
\def\ev               {\mathrm{ev}}

\newcounter{Def}[section]

\renewcommand{\theDef}{\arabic{section}.\arabic{Def}}
\newtheoremstyle{upright}
  {\topsep}
  {\topsep}
  {\upshape} 
  {}
  {\bfseries}
  {:}
  {1em}
  {}
\theoremstyle{upright} 

\newtheorem{Theorem}[Def]{Theorem}
\newtheorem{Example}[Def]{Example}
\newtheorem{Remark}[Def]{Remark}
\newtheorem{Corollary}[Def]{Corollary}
\newtheorem{Conjecture}[Def]{Conjecture}
\newtheorem{Definition}[Def]{Definition}
\newtheorem{Lemma}[Def]{Lemma}
\newtheorem{Proposition}[Def]{Proposition}

\newenvironment{Proof}[1][\quad]
{\vspace{-.7cm}\mbox{}\\\noindent \textbf{Proof: {
#1}\\}}{\hfill$\Box$}

\newcommand{\arrowIn}{
\tikz \draw[-stealth, line width=1pt] (-1pt,0) -- (1pt,0);
}
\newcommand{\arrowOut}{
\tikz \draw[-stealth, line width=1.5pt] (-1pt,0) -- (1pt,0);
}

\usepackage[text={17cm, 24.2cm}, centering]{geometry}
\pgfplotsset{compat=1.12}
\def\mytitle{Categorical 4-manifold invariants from  trisection diagrams}
\def\myauthors{C.~Meusburger,  V.~Mulevi\v{c}ius, F.~Torzewska}
\hypersetup{pdftitle=\mytitle, pdfauthor=\myauthors}

\begin{document}

\begin{center}
  {\LARGE\mytitle}

  \vspace{1em}

  {\large
    Catherine Meusburger\footnote{{\tt catherine.meusburger@math.uni-erlangen.de}, orcid \url{https://orcid.org/0000-0003-3380-3398}}\\
     Department Mathematik \\
  Friedrich-Alexander-Universit\"at Erlangen-N\"urnberg \\
  Cauerstra\ss e 11, 91058 Erlangen, Germany\\[+2ex]
  Vincentas Mulevi\v{c}ius\footnote{{\tt vincentas.mulevicius@univie.ac.at}, orcid \url{https://orcid.org/0000-0003-3589-9872}}\\[+1ex]
  \begin{minipage}{8cm}
  \begin{center}
  Universit\"at Wien, Fakult\"at f\"ur Physik,\\ Boltzmanngasse 5, 1090 Wien, Austria
  \end{center}
  \end{minipage} 
  \begin{minipage}{8cm}
  \begin{center}
  ITPA, Vilnius University,\\ Saul\.etekio al.\,3, 10257 Vilnius, Lithuania
  \end{center}
  \end{minipage}
  \\[+2ex]
  Fiona Torzewska\footnote{{\tt fiona.torzewska@bristol.ac.uk} \url{https://orcid.org/0000-0002-1625-0288}}\\
  {School of Mathematics\\
  University of Bristol\\
  Fry Building, Woodland Road\\
  Bristol BS8 1UG, 
  United Kingdom}\\[+4ex]
 }

{November 24, 2025}

\end{center}

\begin{abstract}
We use Gay and Kirby's  description of 4-manifolds in terms of trisections and trisection diagrams to define a new 4-manifold invariant. The algebraic data are an indecomposable finite semisimple bimodule category over a pair of spherical fusion categories, equipped with a bimodule trace, and a pivotal functor from another spherical fusion category into the spherical fusion category of its bimodule endofunctors and natural transformations between them.
The 4-manifold invariant has a simple description in terms of a diagrammatic calculus for this data, in which the three spherical fusion categories  correspond to the three colours of the trisection diagram. 

It includes the Hopf algebraic 4-manifold invariants by Chaidez, Cotler and Cui, which arise when the bimodule category is the category of finite-dimensional complex vector spaces. We also recover the 4-manifold invariants of B\"arenz and Barrett defined by a pivotal functor from a spherical fusion category into a modular fusion category. 
\end{abstract}

\section{Introduction}

In \cite{GK} Gay and Kirby introduce  a new description of 4-manifolds, \emph{trisections} and \emph{trisection diagrams}. Trisections of 4-manifolds and the associated trisection diagrams resemble Heegaard decompositions of 3-manifolds and the associated Heegaard diagrams. A trisection diagram consists of curves of three distinct colours on an orientable surface of genus $g$ such that curves of any two colours form a standard Heegaard diagram for a connected sum $\#_k (S^1\times S^2)$ with $0\leq k\leq g$.

It is shown in \cite{GK} that trisection diagrams that describe diffeomorphic 4-manifolds are related by finite sequences of trisection moves. These are  diffeomorphisms of surfaces, isotopies,  handle slides and (de)stabilisation moves. These are  analogous to the moves for Heegaard diagrams. The stabilisation move amounts to taking a connected sum with the standard trisection diagram of $S^4$.

In \cite{CCC} Chaidez, Cotler and Cui use trisections and trisection diagrams to construct 4-manifold invariants that resemble  Kuperberg invariants of 3-manifolds. 
Their algebraic datum  is a so-called Hopf triplet, 
a triple of finite-dimensional involutive Hopf algebras that satisfy certain compatibility conditions. They show that the resulting 4-manifold invariant includes certain cases of an earlier 4-manifold invariant by B\"arenz and Barrett \cite{BB}, in particular certain Crane-Yetter invariants \cite{CY,CYK}. 
They also conjecture that some of their invariants are related to Kashaev's invariants \cite{Ka}.

The notion of a Hopf triplet in \cite{CCC} encodes the required structures and properties that allow one to assign a number to each trisection diagram and to prove invariance under trisection moves. However, Hopf triplets are not a standard notion from the theory of Hopf algebras, and it is not easy to find examples. Given that the B\"arenz-Barrett invariants \cite{BB} are defined for more general categorical data, one suspects that there should be  generalised trisection invariants  based on categorical data that include the ones from \cite{CCC} as special cases.

\textbf{Content of the article}

In this article, we define such 4-manifold invariants based on trisection diagrams and categorical data and show that they generalise the 4-manifolds from \cite{CCC}. 

Our data consists of three spherical fusion categories $\cala$, $\mad$, $\mac$ over $\C$,  an indecomposable finite semisimple $(\mac,\mad)$-bimodule category $\calm$ with a bimodule trace and a pivotal functor $\Phi:\cala\to \End_{(\mac,\mad)}(\calm)$ into the spherical fusion category of $(\mac,\mad)$-bimodule endofunctors of $\calm$ and bimodule natural transformations between them.

Our 4-manifold invariants are then obtained as follows. We assign  the three spherical fusion categories $\cala$, $\mad$, $\mac$ to the three colours in a trisection diagram. We  then interpret a trisection diagram whose curves are labelled with simple objects of $\cala$, $\mad$ and $\mac$ and whose regions are labelled with simple objects of $\calm$ as a diagram in the diagrammatic calculus for bimodule categories, bimodule functors and bimodule natural transformations introduced by one of the authors in \cite{M}. This allows us to evaluate each labelled trisection diagram to a complex number. Rescaling appropriately with  dimensions  and averaging over all assignments of simple objects then eliminates the dependence on the labelling. 

Invariance under diffeomorphisms of surfaces and isotopies of curves is already built into this description, and invariance under handle slides is a consequence of the summation over simple objects. Invariance under (de)stabilisation moves requires an additional condition, namely that the evaluation of the standard trisection diagram of $S^4$ is non-zero, and a further rescaling. We then have

\textbf{Theorem A} (Theorem \ref{th:gentrisection}, Corollary \ref{cor:stabilising}) Under the  above assumptions on the data, the rescaled averaged evaluation of a trisection diagram is a 4-manifold invariant.  

We  show that over $\C$ the 4-manifold invariants constructed by Chaidez, Cotler and Cui in \cite{CCC} coincide precisely with our 4-manifold invariants for the case $\calm=\mathrm{Vect}_\C$. Each Hopf triplet $(A,B,C)$ defines spherical fusion categories  $\cala\cong A^{*cop}\Mod$, $\mad\cong B^*\Mod$ and $\mac\cong C^{*op}\Mod$ as well as a $(\mac,\mad)$-bimodule category structure on $\mathrm{Vect}_\C$ and
 a pivotal functor $\Phi:\cala\to \mathrm{End}_{(\mac,\mad)}(\mathrm{Vect}_\C)$. Conversely,  
 a $(\mac,\mad)$-bimodule category structure on $\mathrm{Vect}_\C$ and a pivotal functor $\Phi:\cala\to \mathrm{End}_{(\mac,\mad)}(\mathrm{Vect}_\C)$ define
fibre functors for $\cala$, $\mad$, $\mac$ and  a complex Hopf triplet $(A,B, C)$ as above. 

\textbf{Theorem B} (Propositions \ref{prop:hopfmod} and \ref{prop:reducetrip}, Theorem \ref{th:cccth}) The 
 4-manifold invariants from \cite{CCC}  coincide with our 4-manifold invariants from Theorem A for $\calm=\mathrm{Vect}_\C$.

We  give two examples of our 4-manifold invariants that go beyond the Hopf algebraic ones from \cite{CCC}. We first show that our 4-manifold invariants recover the B\"arenz-Barrett invariants \cite{BB} for a  pivotal functor $\phi:\cala\to\mac$ into a modular fusion category $\mac$. This shows in particular that the Crane-Yetter invariants \cite{CY,CYK} defined by a modular fusion category are special cases of our invariants.

\textbf{Theorem C} (Theorem \ref{th:bbtric})
For a modular fusion category $\mac$ as a $(\mac,\mac^{rev})$-bimodule category and $\Phi$ induced by 
a pivotal functor $\phi:\cala\to\mac$ the 4-manifold invariants from Theorem A coincide with the B\"arenz-Barrett invariants from \cite{BB}, up to rescaling.

 Our second example  arises from  vector spaces graded by  finite sets with  transitive group actions.  In this case, the 4-manifold invariants from Theorem A  admit a  description in terms of generalised Hopf triplets, where all three Hopf algebras in the  triplet  are replaced by  \emph{weak Hopf algebras}.   Although the resulting 4-manifold invariants are not stronger than the invariants from  \cite{CCC}, they extend the invariants from \cite{CCC} to weak Hopf algebras, as anticipated in \cite[Rem 5.10]{CCC}.
The resulting 4-manifold invariants admit a  simple  combinatorial description.  In this case, the underlying diagrammatic calculus encodes the structure maps of the weak Hopf algebras in the triplet.

\textbf{Open questions}

The categorical framework and the diagrammatic description of our 4-manifold invariants relates them to the framework of topological quantum field theories (TQFTs), in particular defect TQFTs. 

The diagrams that define our 4-manifold invariants were originally developed to treat defects in 3d Turaev-Viro-Barrett-Westbury TQFTs \cite{TV,BW}. From this viewpoint, the trisection diagrams labelled with categorical data  appear as defect surfaces with defect lines in Turaev-Viro-Barrett-Westbury state sums.  We expect that they are also related 2d defect TQFTs of state sum type introduced by Davydov, Kong and Runkel in \cite{DKR} and to the defect TQFTs treated by Carqueville and M\"uller \cite{CM}. In this framework the  
evaluation of a trisection diagram should correspond to a 2d defect TQFT evaluation of a surface with defect lines and defect points. 

We  expect that other examples of defect TQFTs could also define trisection invariants, if they allow for a suitable averaging of the data on the defect lines. For instance, it would be interesting to apply the functorial formalism for 2d defect TQFTs 
  introduced in \cite{DKR} to evaluate trisection diagrams for 2d defect TQFTs  that are not of state sum type. We refer to \cite{Ca} and references therein for a summary on 2d defect TQFTs and a list of examples.

On the other hand, the categorical data of our trisection invariants suggests that  constructions of TQFTs from non-semisimple tensor categories such as the ones using modified traces could be applied to generalise them to non-semisimple settings, see for instance  Costantino, Geer, Patureau-Mirand, Virelizier \cite{GPMV} or De Renzi, Gainutdinov,  Geer, Patureau-Mirand, and Runkel \cite{DGGPMR}. It would be interesting to see, if this allows one to make contact with the skein-theoretic 4-manifold invariants by   Costantino, Geer, Haïoun and Patureau-Mirand in \cite{CGHPM}.

Finally, it is natural to ask, if the 4-manifold invariants in this article arise from a 4d TQFT 
or are related to specific 4d TQFTs in the literature. We relate them to Crane-Yetter TQFTs \cite{CY, CYK} for specific choices of data, but we do not know, if they are related to Douglas-Reutter TQFTs \cite{DR}.  To our knowledge, this remains an open question even for B\"arenz-Barrett invariants \cite{BB}.

\textbf{Structure of the article}

Our article is structured as follows. In Section \ref{sec:trisection} we summarise trisections of 4-manifolds and the associated trisection diagrams. Section \ref{sec:hopftrisecccc} introduces the Hopf algebraic trisection invariants from \cite{CCC} and their underlying algebraic data.

Section \ref{sec:cattrisec} contains our core results. It starts with a summary of the categorical data used in our construction in Section \ref{sec:bimodbackground}. Section \ref{sec:diagcalc} contains a summary of the diagrammatic calculus from \cite{M} for the  setting of this article. In Section \ref{sec:triseccchord} we then define our categorical 4-manifold invariants based on trisection diagrams and show that they are indeed 4-manifold invariants (Theorem \ref{th:gentrisection}). 

Section \ref{sec:trisecbimod} identifies the  Hopf algebraic 4-manifold invariants from \cite{CCC} as special cases of our construction.  In Section \ref{subsec:reptheoret} we   reformulate the invariants from \cite{CCC} in terms of Hopf algebra  representations.  In Section \ref{sec:tripletbimod} we show that every Hopf triplet defines a bimodule category structure on $\mathrm{Vect}_\C$ and a pivotal functor into its spherical fusion category of bimodule endofunctors. We also show the converse of this statement, namely that every pair of a $(\mac,\mad)$-bimodule category structure on $\mathrm{Vect}_\C$ and a pivotal functor $\Phi:\cala\to \End_{(\mac,\mad)}(\mathrm{Vect}_\C)$ defines a Hopf triplet. In Section \ref{subsec:trisecbimod} we then apply these results to prove that the 4-manifold invariants from \cite{CCC} correspond to  our invariants for bimodule category structures on $\calm=\mathrm{Vect}_\C$.

Section \ref{sec:newexamples}  treats examples that go beyond the construction in \cite{CCC}. In Section \ref{subsec:bbinv} we show how the B\"arenz-Barrett invariants defined by a modular fusion category $\mac$ and a pivotal functor $\phi:\cala\to\mac$ can be recovered from our construction. Section \ref{subsec:DWinvariant} investigates our invariants for vector spaces graded by finite sets with transitive group actions. It establishes the weak Hopf algebra symmetries of these examples and relates them to the diagrammatic calculus. The appendix contains some basic results on Hopf algebras, on their representations and on weak Hopf algebras.

\tableofcontents

\section{Trisections of 4-manifolds and trisection diagrams}
\label{sec:trisection}
Trisections of closed 4-manifolds and their description by trisection diagrams were introduced by Gay and Kirby in \cite{GK}. 
They can be viewed as a 4-manifold counterpart of Heegaard splittings of 3-manifolds and their description by Heegaard diagrams. 
In this section we summarise the required results on trisections and trisection diagrams from \cite{GK}.

All manifolds considered in the following are compact smooth oriented manifolds and all diffeomorphisms are assumed to be orientation preserving. We denote by $M\# N$ the connected sum of connected manifolds $M$ and $N$ and by $\#_k M$ for $k\in\N_0$ the $k$-fold connected sum of  $M$ with itself. For $k=0,1$ this is defined as $\#_1 M=M$ and  $\#_0 M=S^n$, where $n\geq 1$ is the dimension of $M$. 

\begin{Definition}\cite[Def 1]{GK} \label{def:trisection}Let $0\leq k \leq g$. A \textbf{$(g,k)$-trisection} of a closed connected 4-manifold $X$ is a decomposition $X=X_1\cup X_2\cup X_3$ into submanifolds $X_1,X_2,X_3\subset X$ (with corners) 
such that 
\begin{compactenum}[(i)]
\item there are orientation preserving diffeomorphisms $\phi_i: X_i\to \#_k (S^1\times B^3)$ for $i=1,2,3$,
\item for $i=1,2,3$ taken mod 3, $\phi_i(X_i\cap X_{i+1})$ and $\phi_i(X_i\cap X_{i-1})$ form the standard Heegaard splitting of $\#_k(S^1\times S^2)$ obtained by $(g-k)$-fold stabilisation of the standard genus $k$ splitting.
\end{compactenum}
\end{Definition}

\begin{figure}
\centering

	\begin{tikzpicture}[baseline={([yshift=-.5ex]current bounding box.center)}]
	\node at (0,0) {\def\svgscale{.5} 
\begingroup%
  \makeatletter%
  \providecommand\color[2][]{%
    \errmessage{(Inkscape) Color is used for the text in Inkscape, but the package 'color.sty' is not loaded}%
    \renewcommand\color[2][]{}%
  }%
  \providecommand\transparent[1]{%
    \errmessage{(Inkscape) Transparency is used (non-zero) for the text in Inkscape, but the package 'transparent.sty' is not loaded}%
    \renewcommand\transparent[1]{}%
  }%
  \providecommand\rotatebox[2]{#2}%
  \newcommand*\fsize{\dimexpr\f@size pt\relax}%
  \newcommand*\lineheight[1]{\fontsize{\fsize}{#1\fsize}\selectfont}%
  \ifx\svgwidth\undefined%
    \setlength{\unitlength}{841.88976378bp}%
    \ifx\svgscale\undefined%
      \relax%
    \else%
      \setlength{\unitlength}{\unitlength * \real{\svgscale}}%
    \fi%
  \else%
    \setlength{\unitlength}{\svgwidth}%
  \fi%
  \global\let\svgwidth\undefined%
  \global\let\svgscale\undefined%
  \makeatother%
  \begin{picture}(1,0.70707071)%
    \lineheight{1}%
    \setlength\tabcolsep{0pt}%
    \put(0,0){\includegraphics[width=\unitlength,page=1]{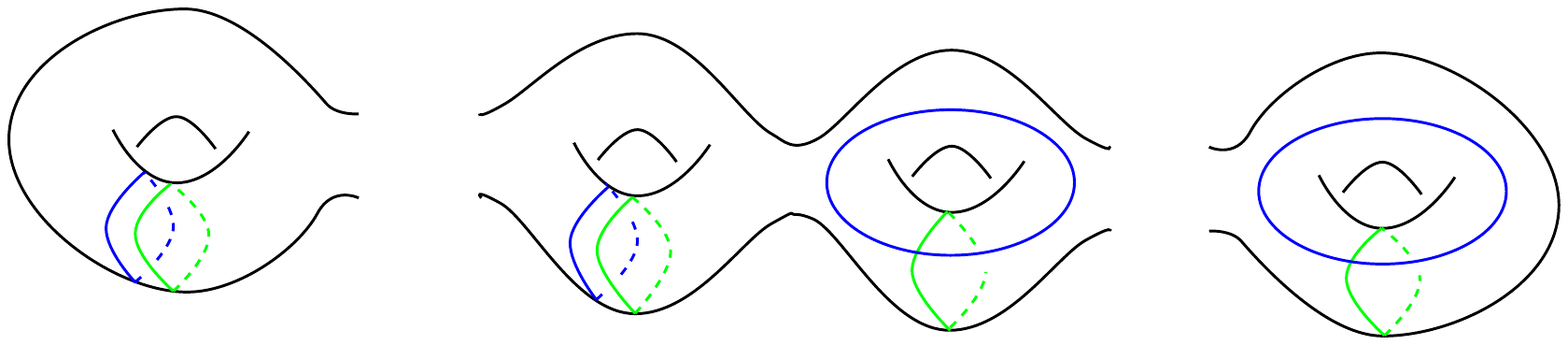}}%
    \put(0.69072166,0.53579116){\makebox(0,0)[lt]{\lineheight{1.25}\smash{\begin{tabular}[t]{l}$\ldots$\end{tabular}}}}%
    \put(0.2382196,0.55141463){\makebox(0,0)[lt]{\lineheight{1.25}\smash{\begin{tabular}[t]{l}$\ldots$\end{tabular}}}}%
    \put(0,0){\includegraphics[width=\unitlength,page=2]{standardheegard.pdf}}%
    \put(0.22629683,0.3679835){\makebox(0,0)[lt]{\lineheight{1.25}\smash{\begin{tabular}[t]{l}$k\times$\end{tabular}}}}%
    \put(0.72958261,0.36971943){\makebox(0,0)[lt]{\lineheight{1.25}\smash{\begin{tabular}[t]{l}$(g-k)\times$\end{tabular}}}}%
  \end{picture}%
\endgroup%
 };
	\end{tikzpicture}

$\quad$\\[-30ex]
\caption{Standard genus $g$ Heegaard diagram for $\#_k(S^1\times S^2)$.}
\label{fig:standardheegard}
\end{figure}

The standard genus $g$ Heegaard diagram for the splitting of $\#_k(S^1\times S^2)$ in (ii)  is given in Figure \ref{fig:standardheegard}. 

It is shown in \cite[Th 4]{GK} that every connected closed 4-manifold admits a $(g,k)$-trisection for suitable $0\leq k\leq g$. 
As explained in \cite{GK}, the numbers $k$ and $g$ are related by the Euler characteristic of $X$
\begin{align}\label{eq:euler}
\chi(X)=2+g-3k.
\end{align}

Just as Heegaard splittings of 3-manifolds are given by Heegaard diagrams, trisections of 4-manifolds are given by trisection diagrams. Each trisected 4-manifold $X=X_1\cup X_2\cup X_3$ defines a trisection diagram. Conversely, a trisection diagram allows one to construct a trisected 4-manifold  $X=X_1\cup X_2\cup X_3$.
This construction is described in \cite{GK} and also in more detail in \cite[Def 2.29]{CCC}.

\begin{Definition}\label{def:trisec} \cite[Def 2]{GK} 
A \textbf{$(g,k)$-trisection diagram} 
is an oriented  surface $\Sigma$ of genus $g\geq 1$ with three sets $\alpha=\{\alpha_i\}_{i=1}^g$, $\beta=\{\beta_i\}_{i=1}^g$, $\gamma=\{\gamma_i\}_{i=1}^g$ of curves on $\Sigma$ such that
each of the triples $(\Sigma,\alpha,\beta)$, $(\Sigma,\beta,\gamma)$ and $(\Sigma,\gamma,\alpha)$ is a genus $g$ Heegaard diagram for $\#_k (S^1\times S^2)$.
\end{Definition}

In the following,  we call the curves in $\alpha$, $\beta$ and $\gamma$ the red, blue and green curves of the trisection diagram, respectively. Note that Definition \ref{def:trisec} implies that each red, green and blue curve  is non-separating and intersects only curves of different colours. 
The Heegaard diagrams $(\Sigma,\alpha,\beta)$, $(\Sigma,\beta,\gamma)$ and $(\Sigma,\gamma,\alpha)$ represent 
the three Heegaard splittings in Definition \ref{def:trisection} (ii). In particular, by applying diffeomorphisms of surfaces, isotopies of curves,  handle slides, it is  possible to transform them  into the standard diagram in Figure \ref{fig:standardheegard}, up to permutations of the colours.

\begin{Example} \label{ex:trisecstandard}The following trisection diagram $\Sigma_{st}$ of genus $g=3$  with $k=1$ describes the 4-sphere $S^4$ with Euler characteristic $\chi(S^4)=2+g-3k=2$.
\begin{align}\label{eq:stabilisation}

	\begin{tikzpicture}[baseline={([yshift=-.5ex]current bounding box.center)}]
	\node at (0,0) {\def\svgscale{.2} 
\begingroup%
  \makeatletter%
  \providecommand\color[2][]{%
    \errmessage{(Inkscape) Color is used for the text in Inkscape, but the package 'color.sty' is not loaded}%
    \renewcommand\color[2][]{}%
  }%
  \providecommand\transparent[1]{%
    \errmessage{(Inkscape) Transparency is used (non-zero) for the text in Inkscape, but the package 'transparent.sty' is not loaded}%
    \renewcommand\transparent[1]{}%
  }%
  \providecommand\rotatebox[2]{#2}%
  \newcommand*\fsize{\dimexpr\f@size pt\relax}%
  \newcommand*\lineheight[1]{\fontsize{\fsize}{#1\fsize}\selectfont}%
  \ifx\svgwidth\undefined%
    \setlength{\unitlength}{841.88976378bp}%
    \ifx\svgscale\undefined%
      \relax%
    \else%
      \setlength{\unitlength}{\unitlength * \real{\svgscale}}%
    \fi%
  \else%
    \setlength{\unitlength}{\svgwidth}%
  \fi%
  \global\let\svgwidth\undefined%
  \global\let\svgscale\undefined%
  \makeatother%
  \begin{picture}(1,0.70707071)%
    \lineheight{1}%
    \setlength\tabcolsep{0pt}%
    \put(0,0){\includegraphics[width=\unitlength,page=1]{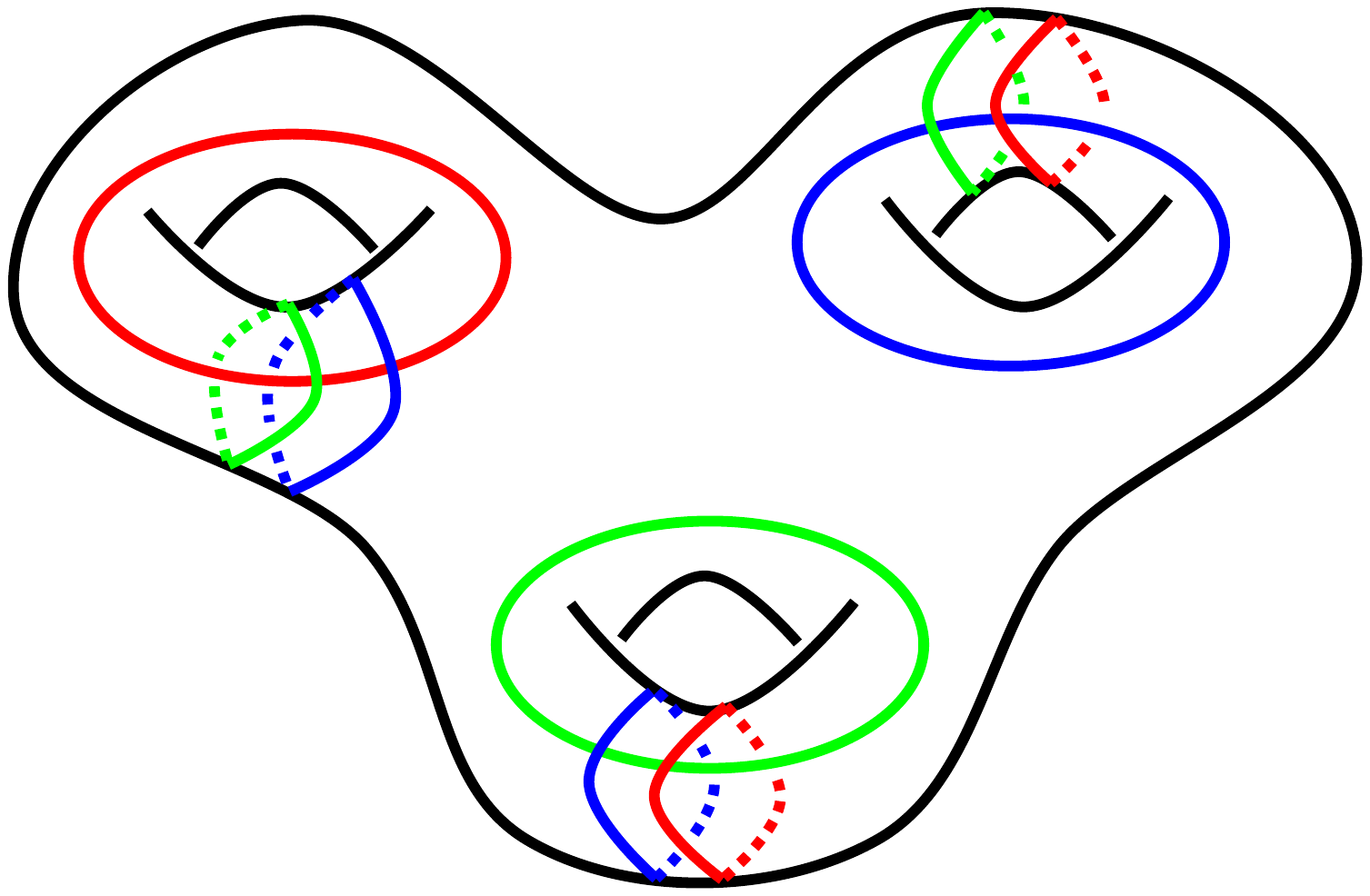}}%
  \end{picture}%
\endgroup%
 };
	\end{tikzpicture}

\end{align}
\end{Example}

\begin{Example} \label{ex:triseccp2}The following trisection diagram of genus $g=1$  with $k=0$ describes complex projective space $\C\mathrm{P}^2$  with Euler characteristic $\chi(\C\mathrm{P}^2)=2+g-3k=3$.
\begin{align}\label{eq:cp2}

	\begin{tikzpicture}[baseline={([yshift=-.5ex]current bounding box.center)}]
	\node at (0,0) {\def\svgscale{.18} 
\begingroup%
  \makeatletter%
  \providecommand\color[2][]{%
    \errmessage{(Inkscape) Color is used for the text in Inkscape, but the package 'color.sty' is not loaded}%
    \renewcommand\color[2][]{}%
  }%
  \providecommand\transparent[1]{%
    \errmessage{(Inkscape) Transparency is used (non-zero) for the text in Inkscape, but the package 'transparent.sty' is not loaded}%
    \renewcommand\transparent[1]{}%
  }%
  \providecommand\rotatebox[2]{#2}%
  \newcommand*\fsize{\dimexpr\f@size pt\relax}%
  \newcommand*\lineheight[1]{\fontsize{\fsize}{#1\fsize}\selectfont}%
  \ifx\svgwidth\undefined%
    \setlength{\unitlength}{841.88976378bp}%
    \ifx\svgscale\undefined%
      \relax%
    \else%
      \setlength{\unitlength}{\unitlength * \real{\svgscale}}%
    \fi%
  \else%
    \setlength{\unitlength}{\svgwidth}%
  \fi%
  \global\let\svgwidth\undefined%
  \global\let\svgscale\undefined%
  \makeatother%
  \begin{picture}(1,0.70707071)%
    \lineheight{1}%
    \setlength\tabcolsep{0pt}%
    \put(0,0){\includegraphics[width=\unitlength,page=1]{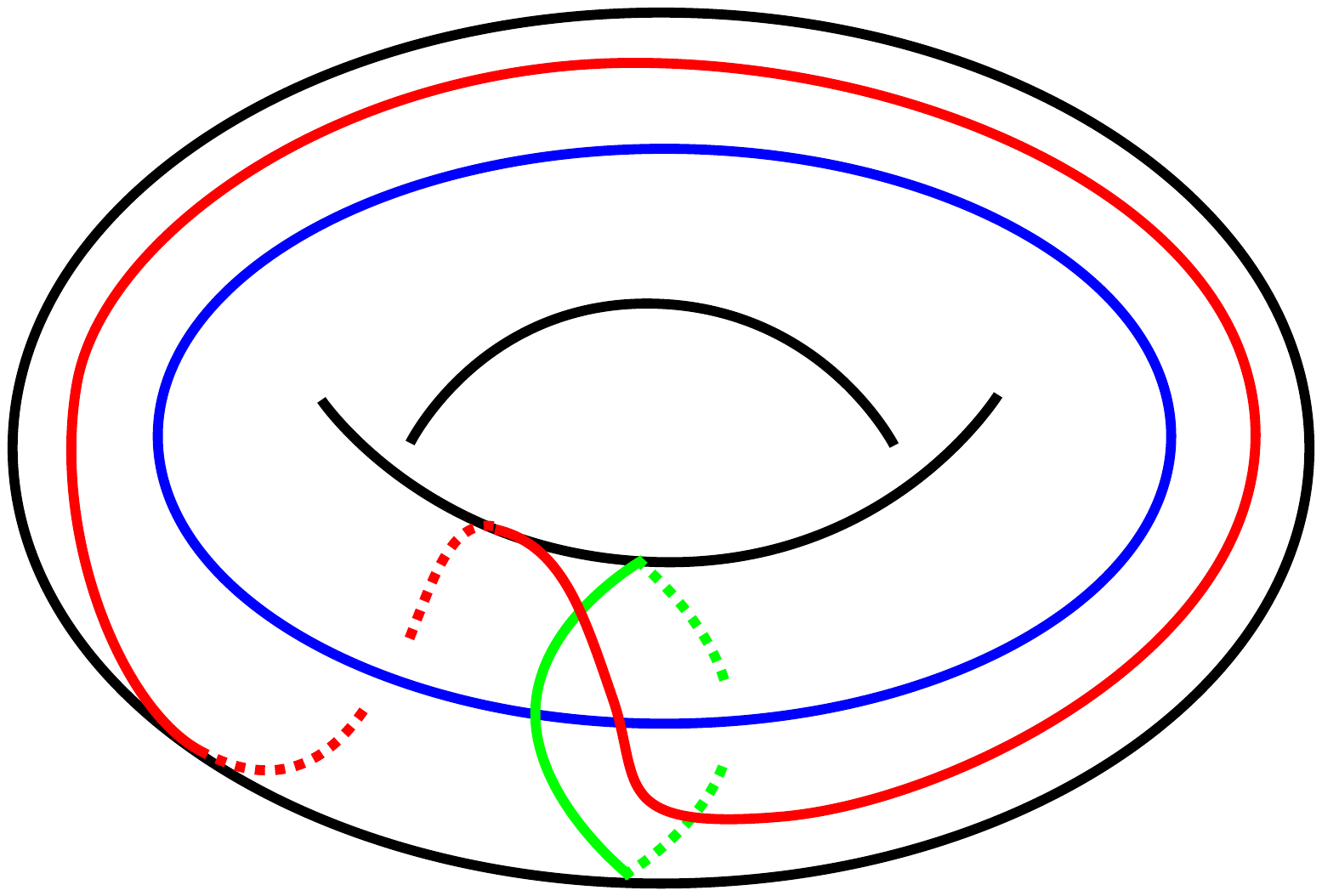}}%
  \end{picture}%
\endgroup%
 };
	\end{tikzpicture}

\end{align}
\end{Example}

Each trisection diagram defines a unique closed connected 4-manifold. However, distinct trisection diagrams may yield diffeomorphic  4-manifolds.  Corollary 12 in \cite{GK} gives sufficient and necessary conditions, under which two trisection diagrams define diffeomorphic 4-manifolds. The following theorem summarises the results from \cite[Th 11, Cor 12]{GK}, but with the more explicit description of the isotopies of curves in terms of two-point and three-point moves from 
Definition 2.25, Lemma 2.26 and Theorem 2.32 in \cite{CCC}.

\begin{Theorem}\cite[Th 11, Cor 12]{GK}\label{th:intmoves}
 The 4-manifolds defined by trisection diagrams $T,T'$  are  diffeomorphic  if and only if $T$ and $T'$ are related by a finite sequence of the following moves and their counterparts with permuted colours:
\begin{compactenum}[(i)]
\item \textbf{Diffeomorphisms of surfaces}\\[-2ex]
\item \textbf{Two-point moves}
\begin{align}\label{eq:twopoint}
\begin{tikzpicture}[scale=.4]
\begin{scope}[shift={(0,0)}]
\draw[dashed, line width=1pt] (0,0) circle (2);
\draw[blue, line width=1pt] (70:2)..controls (-1,1) and (-1,-1)..(-70:2);
\draw[red, line width=1pt] (110:2).. controls (1,1) and (1,-1)..(250:2);
\end{scope}
\node at (3,0) {$\leftrightarrow$};
\begin{scope}[shift={(6,0)}]
\draw[dashed, line width=1pt] (0,0) circle (2);
\draw[blue, line width=1pt] (70:2)--(-70:2);
\draw[red, line width=1pt] (110:2)--(250:2);
\end{scope}
\end{tikzpicture}
\end{align}

\item \textbf{Three-point moves}
\begin{align}\label{eq:threepoint}
\begin{tikzpicture}[scale=.4]
\begin{scope}[shift={(0,0)}]
\draw[dashed, line width=1pt] (0,0) circle(2);
\draw[red, line width=1pt] (60:2)--(240:2);
\draw[blue, line width=1pt] (120:2)--(300:2);
\draw[green, line width=1pt] (30:2)--(150:2);
\end{scope}
\node at (3,0) {$\leftrightarrow$};
\begin{scope}[shift={(6,0)}]
\draw[dashed, line width=1pt] (0,0) circle(2);
\draw[red, line width=1pt] (60:2).. controls (150:1) .. (240:2);
\draw[blue, line width=1pt] (120:2).. controls  (30:1).. (300:2);
\draw[green, line width=1pt] (30:2).. controls (270:1) ..(150:2);
\end{scope}
\end{tikzpicture}
\end{align}
\item \textbf{Handle slides}
\begin{align}\label{eq:handleslide}
\begin{tikzpicture}[scale=.4]
\begin{scope}[shift={(20,0)}]
\draw[dashed, line width=1pt] (5,0) circle (3);
\draw[dashed, line width=1pt] (-5,0) circle (3);
\draw[dashed, line width=1pt] (5,0) circle (1);
\draw[dashed, line width=1pt] (-5,0) circle (1);
\draw[fill=white, draw=none] (-2.5,1)--(2.5,1)--(2.5,-1)--(-2.5,-1)--(-2.5,1);
\draw[red, line width=1pt] (5,0) circle (2.5);
\draw[red, line width=1pt] (5,0) circle (1.5);
\draw[red, line width=1pt] (-5,0) circle (1.5);
\draw[fill=white, draw=none] (-3.7,.7)--(3,.7)--(3,-.7)--(-3.7,-.7)--(-3.7,.7);
\draw[dashed, line width=1pt] (-2,1)--(2,1);
\draw[dashed, line width=1pt] (-2,-1)--(2,-1);
\draw[red, line width=1pt] (-3.6,.7)--(2.55,.7);
\draw[red, line width=1pt] (-3.6,-.7)--(2.55,-.7);
\draw[blue, line width=1pt] (5,1)--(5,3);
\draw[blue, line width=1pt] (6,0)--(8,0);
\draw[green, line width=1pt] (5,-1)--(5,-3);
\draw[blue, line width=1pt] (-5,-1)--(-5,-3);
\draw[green, line width=1pt] (-6,0)--(-8,0);
\end{scope}
\node at (10,0){$\leftrightarrow$};
\begin{scope}[shift={(0,0)}]
\draw[dashed, line width=1pt] (5,0) circle (3);
\draw[dashed, line width=1pt] (-5,0) circle (3);
\draw[dashed, line width=1pt] (5,0) circle (1);
\draw[dashed, line width=1pt] (-5,0) circle (1);
\draw[fill=white, draw=none] (-2.5,1)--(2.5,1)--(2.5,-1)--(-2.5,-1)--(-2.5,1);
\draw[red, line width=1pt] (5,0) circle (1.5);
\draw[red, line width=1pt] (-5,0) circle (1.5);
\draw[dashed, line width=1pt] (-2,1)--(2,1);
\draw[dashed, line width=1pt] (-2,-1)--(2,-1);
\draw[blue, line width=1pt] (5,1)--(5,3);
\draw[blue, line width=1pt] (6,0)--(8,0);
\draw[green, line width=1pt] (5,-1)--(5,-3);
\draw[blue, line width=1pt] (-5,-1)--(-5,-3);
\draw[green, line width=1pt] (-6,0)--(-8,0);
\end{scope}
\end{tikzpicture}
\end{align}

\item \textbf{Stabilisation and destabilisation} \\
Replacing a trisection diagram $T$ by the connected sum $T\# \Sigma_{st}$,  where $\Sigma_{st}$ is the trisection diagram for $S^4$ from \eqref{eq:stabilisation}.
\end{compactenum}
\end{Theorem}

The diffeomorphisms of  surfaces  in Theorem \ref{th:intmoves} must send curves of one diagram to curves of the same colour of the other. 
For the two-point and three-point moves and the handle slides, the two diagrams are supposed to agree outside the region bordered by the dashed lines, and this region can be embedded anywhere in the surface.
The number of green and blue line segments in Diagram \eqref{eq:handleslide} for the handle slide is arbitrary.
For the stabilisation and destabilisation move, the connected sum of the diagrams is taken in such a way that no curves on the surfaces intersect the removed discs.

Note that diffeomorphisms of surfaces,  two-point moves and the handle slides of trisection diagrams correspond to diffeomorphisms, handle slides and two-point moves of the Heegaard diagrams obtained by forgetting the curves of a given colour. By applying these moves it is possible to bring the curves of the other two colours into standard position, as in Figure \ref{fig:standardheegard}.
As explained in \cite[Sec 2]{GK},  this  can be used to define a framed link diagram for a handlebody decomposition of the underlying 4-manifold, see for instance the textbooks by Kirby \cite[Ch I]{Ki}
or by Akbulut \cite[Ch 1]{A}.

\begin{Remark} \label{rem:linkdiagram}Let $(\Sigma,\alpha,\beta,\gamma)$ be a $(g,k)$-trisection diagram for a 4-manifold $X$, whose green and blue curves are in standard position, as in Figure \ref{fig:standardheegard}.  

This defines a framed link diagram for a handlebody decomposition of $X$ as follows:
\begin{compactenum}
\item Erase the last $(g-k)$ green and blue curves that form longitude-meridian pairs.
\item Replace each of the first $k$ pairs of green and blue curves by a circle that is  lifted slightly from $\Sigma$ such that it surrounds the segments of red curves intersecting this pair.
\item  Keep the red curves and equip them with the framing induced by $\Sigma$.
\end{compactenum}
The $g$ red components of the resulting link diagram  specify the  attaching maps for the 2-handles and its $k$ circles specify the attaching maps for the 1-handles.
\end{Remark}

\begin{Example}\label{ex:trisecstablink}
The trisection diagram for $S^4$ from Example \ref{ex:trisecstandard} defines the  framed link diagram 
\begin{align*}
\begin{tikzpicture}[scale=.6]
\draw[line width=.6pt, color=red] (0,0) ellipse (2cm and 1cm);
\draw[line width=.6pt, color=red] (5,0) ellipse (2cm and 1cm);
\draw[line width=.6pt, color=red] (10,0) ellipse (2cm and 1cm);
\draw[draw=none, fill=white] (-.9,-.8) circle (.3);
\draw[line width=.6pt, color=black] (0,-.7)  arc(10:320:.5cm and 1cm);
\end{tikzpicture}
\end{align*}
\end{Example}

\begin{Example} The trisection diagram for $\C\mathrm{P}^2$ from Example \ref{ex:triseccp2} defines the framed link diagram 
\begin{align*}

	\begin{tikzpicture}[baseline={([yshift=-.5ex]current bounding box.center)}]
	\node at (0,0) {\def\svgscale{.11} 
\begingroup%
  \makeatletter%
  \providecommand\color[2][]{%
    \errmessage{(Inkscape) Color is used for the text in Inkscape, but the package 'color.sty' is not loaded}%
    \renewcommand\color[2][]{}%
  }%
  \providecommand\transparent[1]{%
    \errmessage{(Inkscape) Transparency is used (non-zero) for the text in Inkscape, but the package 'transparent.sty' is not loaded}%
    \renewcommand\transparent[1]{}%
  }%
  \providecommand\rotatebox[2]{#2}%
  \newcommand*\fsize{\dimexpr\f@size pt\relax}%
  \newcommand*\lineheight[1]{\fontsize{\fsize}{#1\fsize}\selectfont}%
  \ifx\svgwidth\undefined%
    \setlength{\unitlength}{841.88976378bp}%
    \ifx\svgscale\undefined%
      \relax%
    \else%
      \setlength{\unitlength}{\unitlength * \real{\svgscale}}%
    \fi%
  \else%
    \setlength{\unitlength}{\svgwidth}%
  \fi%
  \global\let\svgwidth\undefined%
  \global\let\svgscale\undefined%
  \makeatother%
  \begin{picture}(1,0.70707071)%
    \lineheight{1}%
    \setlength\tabcolsep{0pt}%
    \put(0,0){\includegraphics[width=\unitlength,page=1]{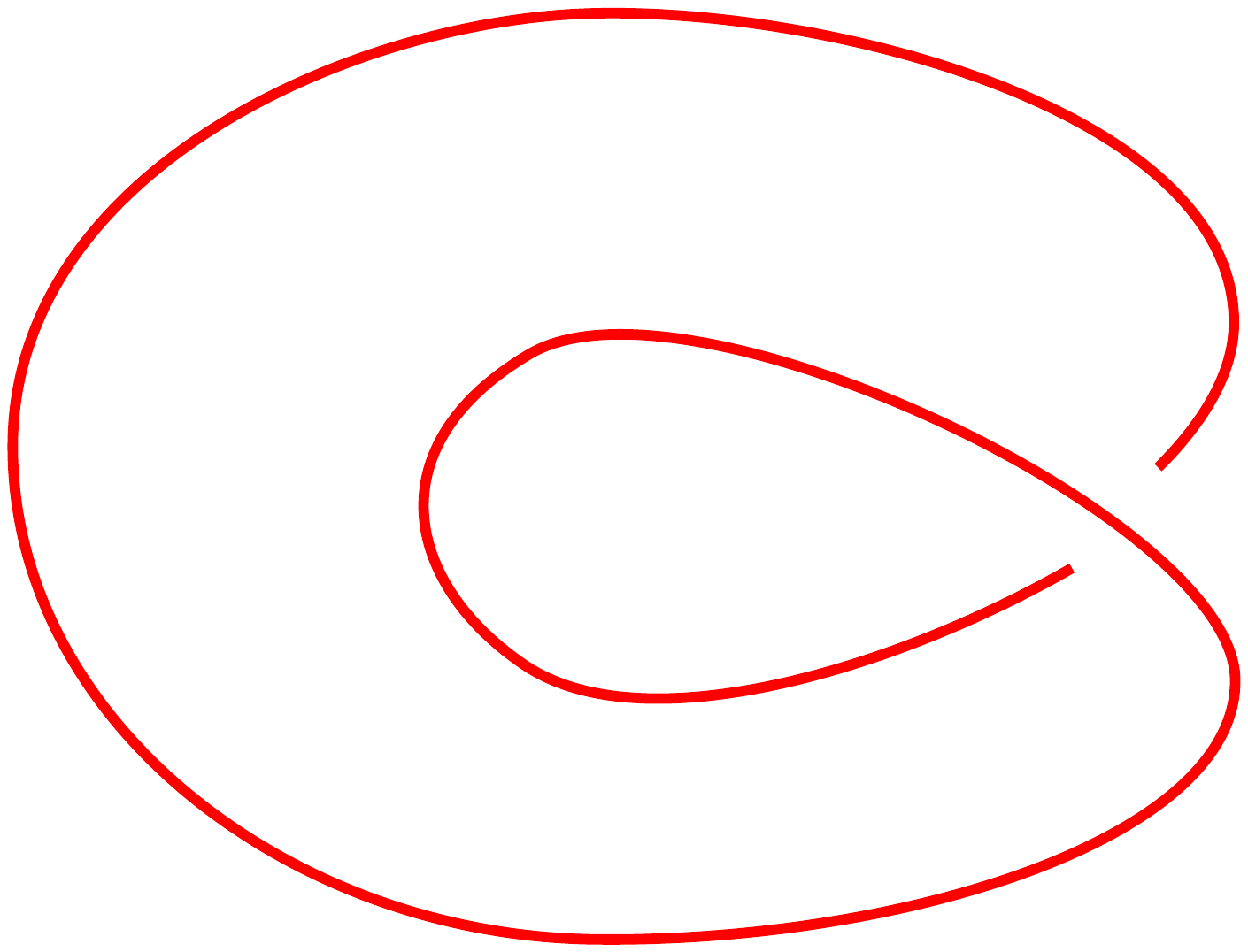}}%
  \end{picture}%
\endgroup%
 };
	\end{tikzpicture}

\end{align*}
\end{Example}

As each closed connected 4-manifold $X$ admits a  trisection diagram and different trisection diagrams of $X$ are related by the moves from Theorem \ref{th:intmoves}, they allow one to construct 4-manifold invariants.
One may define a  \emph{trisection invariant} as a map from the set of all trisection diagrams that is invariant under the moves from Theorem \ref{th:intmoves}. The first invariant of this type was constructed by Chaidez, Cotler and Cui in \cite{CCC} and is described in the next section.

\section{Hopf algebraic trisection invariants}
\label{sec:hopftrisecccc}
In this section, we summarise the construction of trisection invariants in \cite{CCC}. While the construction in \cite{CCC} is given over more general fields, we restrict attention to finite-dimensional semisimple complex Hopf algebras. The reason for this is that  some results we use in our construction are only formulated for this case, although we expect that they generalise to algebraically closed fields of characteristic zero. 
In Section \ref{sec:hopfdoubtrip} we  introduce the Hopf algebraic data for the trisection invariants in \cite{CCC}. In Section \ref{sec:cccinv} we then describe the construction of the  trisection invariants. 

We use the standard notation for the structure maps of a Hopf algebra. Coproducts are denoted in Sweedler notation without summation signs.
For a finite-dimensional Hopf algebra $H$ we denote by  $H^*$ its dual and by  $H^{op}$, $H^{cop}$ and $H^{op,cop}$ the Hopf algebras with the opposite multiplication, comultiplication and opposite multiplication and comultiplication. 
For more details, see Appendix \ref{sec:hopfback}.

\subsection{Hopf doublets and Hopf triplets}
\label{sec:hopfdoubtrip}
The algebraic data for the construction of the trisection invariants in \cite{CCC} are \emph{Hopf doublets} and \emph{Hopf triplets}. The former are related to \emph{skew pairings}, see for instance the textbook \cite{R} by Radford. A skew pairing generalises the canonical pairing between a  Hopf algebra and its dual.

\begin{Definition}\label{def:skewpair}\cite[Def 2.10]{CCC}, \cite[Def 7.7.7]{R}\\ Let $A,B$ be finite-dimensional complex Hopf algebras.
\begin{compactenum}
\item A {\bf skew pairing} between $A$ and $B$ is a $\C$-bilinear map $\tau: A\times B\to \C$ that induces a Hopf algebra homomorphism $f_\tau: A\to B^{*cop}$, $a\mapsto \tau(a,-)$. Explicitly,  for all $a,a'\in A$, $b, b'\in B$:
\begin{align}\label{eq:tauforms}
&\tau(a a', b)=\tau(a, \low b 1)\tau(a', \low b 2), \;
\tau(a, b b')=\tau(\low a 2, b)\tau(\low a 1, b'),\;
\tau(a, 1)=\epsilon(a), \; \tau(1, b)=\epsilon(b). 
\end{align}
\item If $\tau: A\times B\to\C$ is a skew pairing, then the pair $(A,B)$ is called a \textbf{Hopf doublet}.
\end{compactenum}
\end{Definition}

\begin{Example} \label{ex:pairexamples} Let $A,B$ finite-dimensional complex Hopf algebras.
\begin{compactenum}
\item The \textbf{trivial skew pairing} between  $A$ and $B$ is  
$$\tau: A\times B\to \C, \quad (a,b)\mapsto \epsilon(a)\epsilon(b).$$

\item \label{ex:ddouble}  The \textbf{canonical skew pairing} for  $A$  is 
$$\langle\;,\;\rangle: A^{*cop}\times A\to \C, \quad (\alpha,a)\mapsto \langle \alpha,a\rangle=\alpha(a).$$

\item \label{ex:Hopfhompair} Hopf algebra homomorphisms $f: A\to B^{*cop}$ or $g: B\to A^{*op}$ define  skew pairings  
$$
\tau_f: A\times B\to \C,\; (a,b)\mapsto
f(a)(b)\qquad \tau_g: A\times B\to \C,\; (a,b)\mapsto  g(b)(a).
$$
\item If $A$ is coquasitriangular, any universal $R$-matrix  $R\in A^*\oo A^*$ defines a skew pairing 
$$
\tau_R: A\times A\to\C, \quad (a,a') \mapsto \langle R, a\oo a'\rangle.
$$
\item A skew pairing $\tau: A\times B\to \C$  defines the following eight, generally distinct, Hopf doublets
\begin{align*}
&(A,B) & &(A^{op}, B^{cop}) & &(A^{cop}, B^{op}) & &(A^{op,cop}, B^{op,cop})\\
&(B, A^{op,cop}) & &(B^{op}, A^{op}) & &(B^{cop}, A^{cop}) & &(B^{op,cop},A).
\end{align*}
\end{compactenum}
\end{Example}

Any skew pairing induces a non-degenerate skew pairing between suitable quotient Hopf algebras. This relates any skew pairing to the canonical one in Example \ref{ex:pairexamples}, 2.~and to the skew pairing in 
Example \ref{ex:pairexamples}, 3.~for a Hopf algebra isomorphism. 
The radical $A_\tau=\{a\in A\mid \tau(a,b)=0\;\forall b\in B\}$
is the kernel of the Hopf algebra homomorphism $f_\tau: A\to B^{*cop}$ in Definition \ref{def:skewpair} and hence an ideal and a coideal. It defines a quotient Hopf algebra $A'=A/A_\tau$. 
Likewise, the ideal and coideal  $B_\tau=\{b\in B\mid \tau(a,b)=0\;\forall a\in A\}$ defines a Hopf algebra $B'=B/B_\tau$.

\begin{Remark}\label{rem:skewpair}
A skew pairing $\tau: A\times B\to\C$ induces a Hopf algebra isomorphism between $A/A_\tau$ and $(B/B_\tau)^{*cop}$ and a 
non-degenerate skew pairing 
$\tau': A/A_\tau\times B/B_\tau\to \C$.
\end{Remark}

The interest of a skew pairing $\tau: A\times B\to\C$ is that it defines a convolution invertible 2-cocycle for the tensor product $A\oo B$ that can be used to twist the multiplication of $A\oo B$, see Definition \ref{def:cocyc} and Proposition \ref{prop:twist}. 
This yields a new Hopf algebra, called the  \emph{double}  in \cite[Def 2.14]{CCC} and the \emph{generalised double}  in \cite[Def 13.8.1]{R}.

\begin{Proposition} \label{prop:2cocyc} \cite[Prop 7.7.10 and Def 13.8.1]{R}\\ Let $\tau: A\times B\to \C$ be a skew pairing and $\tau^\inv: A\times B\to \C$, $(a,b)\mapsto \tau(S(a),b)$ its convolution inverse. 
\begin{compactenum}
\item 
The following is a convolution invertible normalised 2-cocycle for $A\oo B$ 
\begin{align}\label{eq:sigmadef}
&\sigma: (A\oo B)\times(A\oo B)\to \C, \quad (a\oo b, a'\oo b')\mapsto\epsilon(a)\epsilon(b')\tau(a', b).
\end{align}
\item This defines a  Hopf algebra $D(A,B)=(A\oo B)_\sigma$, the \textbf{generalised double},   with the comultiplication, counit and unit  of the tensor product $A\oo B$, and  twisted multiplication and antipode
\begin{align}\label{eq:doubletmult}
&(a\oo b)\cdot_\sigma(a'\oo b')=\tau(a'_{(1)}, \low b 1)\tau^\inv(a'_{(3)}, \low b 3)\; aa'_{(2)}\oo \low b 2 b'\\
&S(a\oo b)= (1\oo S(b))\cdot_\sigma (S(a)\oo 1).\nonumber
\end{align}
\item The Hopf algebras $A$ and $B$ are Hopf subalgebras of $(A\oo B)_\sigma$ via the inclusions
$$
\iota_A: A\to (A\oo B)_\sigma,\; a\mapsto a\oo 1,\qquad \iota_B: B\to (A\oo B)_\sigma,\;b\mapsto 1\oo b.
$$ 
\end{compactenum}
\end{Proposition}

\begin{Proof}
The convolution inverse of $\sigma$ is $\sigma^\inv: (A\oo B)\times (A\oo B)\to \C$, $(a\oo b, a'\oo b')\mapsto \epsilon(a)\epsilon(b')\, \tau^\inv(a',b)$. This follows, because $\tau$ is  convolution invertible  with  inverse $\tau^\inv(a,b)=\tau(S(a), b)$, i.~e.~it satisfies
 \begin{align}\label{eq:convinv}
\tau(\low a 1,\low b 1)\, \tau^\inv(\low a 2, \low b 2)=\epsilon(a)\epsilon(b)=\tau^\inv(\low a 1,\low b 1)\, \tau(\low a 2,\low b 2)\qquad\forall a\in A,b\in B.
 \end{align}
The remaining claims  follow by direct, but lengthy computations, for details see \cite[Prop 7.7.10]{R}. 
\end{Proof}

\begin{Example} \label{ex:doubles} $\quad$
\begin{compactenum}
\item The generalised double $D(A,B)$ for the trivial  skew pairing $\tau: A\times B\to \C$, $(a,b)\mapsto \epsilon(a)\epsilon(b)$ is the tensor product $D(A,B)=A\oo B$.\\[-2ex]

\item For the canonical skew pairing $\langle\;,\;\rangle: A^{*cop}\times A\to \C$, $(\alpha,a)\mapsto \alpha(a)$ from Example \ref{ex:pairexamples}, 2.~the generalised double is the Drinfeld double $D(A^{*cop}, A)=D(A)$ with  
\begin{align*}
&(\alpha\oo a)\cdot (\beta\oo b)=\low \beta 3(\low a 1) \low \beta 1(S^\inv(\low a 3))\; \alpha\low\beta 2\oo \low a 2 b\\
&S(\alpha\oo a)=(1\oo S(a))\cdot  (S^\inv(\alpha)\oo 1).
\end{align*}
\end{compactenum}
\end{Example}

The second example shows that the generalised double generalises the construction of the Drinfeld double by twisting with the universal $R$-matrix. In fact,  Remark \ref{rem:skewpair}
and  Example \ref{ex:pairexamples}, 3.~relate every  double $D(A,B)$ to a Drinfeld double. This is discussed in detail after \cite[Th 3.13]{CCC}.

Just as for a Drinfeld double, integrals of a generalised double $D(A,B)$ that is built from finite-dimensional semisimple complex Hopf algebras $A$ and $B$ are just tensor products of integrals of  $A$ and $B$. For basic definitions on integrals, see Appendix \ref{sec:hopfback} and the references therein.

\begin{Lemma}\label{lem:intdoub}
Let $A,B$ be finite-dimensional semisimple complex Hopf algebras and $\tau: A\times B\to\C$ a skew pairing. 
Then for all normalisable integrals $\ell_A\in A$, $\ell_B\in B$ the element $\ell_A\oo \ell_B$ is a normalisable integral of the generalised double $D(A,B)$.
\end{Lemma}

\begin{Proof} This follows by a direct computation with  Formula \eqref{eq:doubletmult} for the multiplication of $D(A,B)$, the involutivity of the antipode, Definition \ref{def:inthopf} of an integral and Lemma \ref{lem:haarprops}, 4.
\end{Proof}

The main ingredient in the construction of trisection invariants in \cite{CCC} is  a \emph{Hopf triplet}.  This is 
 essentially a skew pairing of a double $D(A,B)$ and another Hopf algebra $C$, such that the three Hopf algebras are paired in a cyclically symmetric way. We summarise  this in an adapted notation and over the field $\C$, where semisimplicity and involutivity of finite-dimensional Hopf algebras are equivalent.

\begin{Definition} \cite[Def 2.16]{CCC}\label{def:triplet} 
A {\bf Hopf triplet}  consists of finite-dimensional semisimple complex Hopf algebras $A,B,C$ and skew pairings $\tau_{AB}: A\times B\to \C$, $\tau_{BC}: B\times C\to \C$ and $\tau_{CA}: C\times A\to \C$ that induce Hopf algebra homomorphisms 
\begin{align}\label{eq:triplets}
&\phi_{ABC}: D(A^{op}, B^{cop})\to C^{*} & &\phi_{ABC}(a\oo b)(c)=\tau_{CA}(\low c 1, a)\tau_{BC}(b,\low c 2) \\
&\phi_{BCA}: D(B^{op}, C^{cop})\to A^{*} & &\phi_{BCA}(b\oo c)(a)=\tau_{AB}(\low a 1, b)\tau_{CA}(c,\low a 2) \nonumber\\
&\phi_{CAB}: D(C^{op}, A^{cop})\to B^{*} & &\phi_{CAB}(c\oo a)(b)=\tau_{BC}(\low b 1, c)\tau_{AB}(a,\low b 2).\nonumber 
\end{align}
\end{Definition}

\begin{Lemma} \cite[Lemma 2.19]{CCC}\label{lem:hopfconc} 
Finite-dimensional semisimple complex Hopf algebras $A,B,C$ with skew pairings $\tau_{AB}: A\times B\to \C$, $\tau_{BC}: B\times C\to \C$ and $\tau_{CA}: C\times A\to \C$
form a Hopf triplet iff one of the following equivalent conditions holds:
\begin{compactenum}[(i)]
\item One of the three maps in \eqref{eq:triplets} is a Hopf algebra morphism.
\item The pairings satisfy for all $a\in A$, $b\in B$ and $c\in C$ 
\begin{align}\label{eq:idd}
\tau_{AB}(\low a 1, \low b 2) \tau_{BC}(\low b 1, \low c 2)\tau_{CA}(\low c 1, \low a 2)=\tau_{AB}(\low a 2, \low b 1) \tau_{BC}(\low b 2, \low c 1)\tau_{CA}(\low c 2, \low a 1).
\end{align}
\end{compactenum}
\end{Lemma}

It follows  
that every Hopf triplet defines three generalised doubles,  
$D(D(A^{op}, B^{cop}), C^{op})$, $D(D(B^{op}, C^{cop}), A^{op})$ and $D(D(C^{op}, A^{cop}), B^{op})$, by using Example \ref{ex:pairexamples}, 3.~to obtain a skew pairing, and then taking the generalised double.
Their underlying vector space is
$A\oo B\oo C$, and the comultiplication, counit and unit are the ones of the  tensor products
$A\oo B^{cop}\oo C$, $B\oo C^{cop}\oo A$ and $C\oo A^{cop}\oo B$, whereas the multiplication is twisted by a  cocycle. 
Their duals are the algebras
$A^{*}\oo B^{*op}\oo C^{*}$, $B^*\oo C^{*op}\oo A^{*}$ and $C\oo A^{*op}\oo B^{*}$, but with twisted comultiplications.

\begin{Example}\label{ex:BBtriplet} \cite[Sec 5.2]{CCC} \\
 Let  $H,K$ be finite-dimensional semisimple complex Hopf algebras  and  $\phi: D(H)\to K$ a Hopf algebra map from the Drinfeld double of $H$. 
This defines a Hopf triplet $(A,B,C)$ with
\begin{align}\label{eq:hopftripletbb}
A=H^{*op,cop}\qquad B=H^{cop}\qquad C=K^*.
\end{align}
The skew pairings are given by the canonical pairing of $H$ and $H^*$ and the restrictions of $\phi$ to the Hopf subalgebras $H\subset D(H)$ and $H^{*cop}\subset D(H)$
\begin{align}\label{eq:pairingsbb}
\qquad\qquad\qquad\qquad &\tau_{AB}: H^{*op,cop}\times H^{cop}\to \C, & &(\alpha, b)\mapsto \alpha(b) \qquad\qquad \qquad\qquad\\
&\tau_{BC}: H^{cop}\times K^*\to \C, &  &(b,\kappa)\mapsto \kappa\circ \phi(1\oo b)\nonumber\\
&\tau_{CA}: K^*\times H^{*op,cop}\to \C, & &(\kappa,\alpha)\mapsto \kappa\circ \phi(\alpha\oo 1).\nonumber 
\end{align}
\end{Example}

\begin{Example}\label{ex:Ktriplet} \cite[Def 4.11]{CCC}\\
Let $\C[\mathbb Z/n\mathbb Z]$ be the group algebra of $\mathbb Z/n\mathbb Z$ with the universal $R$-matrix
\begin{align}\label{eq:rmatkash}
R=\frac 1 n\sum_{k,l=0}^{n-1} e^{2\pi\ii kl/n}\; \bar k\oo\bar l.
\end{align}
By Example \ref{ex:pairexamples}, 4.~this defines a skew pairing of the dual $\C[\mathbb Z/n\mathbb Z]^*$ with itself.  Together with the canonical skew pairing between $\C[\mathbb Z/n\mathbb Z]$ and $\C[\mathbb Z/n\mathbb Z]^*$ this defines a Hopf triplet $(A,B,C)$ with
\begin{align}\label{eq:kashaev}
A=\C[\mathbb Z/n\mathbb Z]^*\qquad B=\C[\mathbb Z/n\mathbb Z]\qquad C=\C[\mathbb Z/n\mathbb Z]^*.
\end{align}
\end{Example}

\subsection{Trisection invariants from Hopf triplets}
\label{sec:cccinv}

In this section, we summarise the construction of trisection invariants from \cite{CCC} for the case where the underlying field is $\C$.
In the following, all Hopf algebras in a Hopf triplet are assumed to be complex finite-dimensional semisimple or, equivalently, complex finite-dimensional involutive. Note that this implies that each of these Hopf algebras has a normalisable integral $\ell$ and all other normalisable integrals  are  non-zero multiples of $\ell$, see  Definition \ref{def:inthopf}.

The construction in \cite{CCC} proceeds in two steps.  The first is the assignment of a \emph{trisection bracket} to each trisection diagram. 
This trisection bracket is already invariant under  all moves from Theorem \ref{th:intmoves} except the  (de)stabilisation moves (v).  In the second step, this trisection bracket is normalised with a number associated to the trisection diagram \eqref{eq:stabilisation} to  achieve invariance under  (de)stabilisation moves.

\begin{Definition} \cite[Def 3.1]{CCC}  \label{def:trisecbr} \\ Let $(A,B,C)$ be a Hopf triplet  and  $\ell_A\in A$, $\ell_B\in B$ and $\ell_C\in C$ normalisable integrals.
The \textbf{trisection bracket} $\langle T\rangle_{A,B,C}\in\C$ of a trisection diagram $T$ is  obtained as follows:

\begin{enumerate}
\item Denote by
$\{\alpha_i\}_{i\in I}$, $\{\beta_i\}_{i\in I}$ and $\{\gamma_i\}_{i\in I}$ for $I=\{1,\ldots, g\}$  the red, blue and green curves of $T$.  Equip them with  orientations and  basepoints. Colour each  curve $\lambda$ with a Hopf algebra $X_\lambda\in\{A,B,C\}$, where $X_\lambda=A$ for red, $X_\lambda=B$ for blue and $X_\lambda=C$ for green curves.

\item Assign to each intersection point  $p$ a sign $\epsilon_p\in\{1,-1\}$ as follows
\begin{align}
\label{eq:intsign}
&\epsilon_p=1 &
&\begin{tikzpicture}[scale=.5]
\draw[line width=1pt, color=blue, ->,>=stealth] (-1,1)--(1,-1);
\draw[line width=1pt, color=red, ->,>=stealth] (1,1)--(-1,-1);
\end{tikzpicture}
& 
&\begin{tikzpicture}[scale=.5]
\draw[line width=1pt, color=green, ->,>=stealth] (-1,1)--(1,-1);
\draw[line width=1pt, color=blue, ->,>=stealth] (1,1)--(-1,-1);
\end{tikzpicture}
&
&\begin{tikzpicture}[scale=.5]
\draw[line width=1pt, color=red, ->,>=stealth] (-1,1)--(1,-1);
\draw[line width=1pt, color=green, ->,>=stealth] (1,1)--(-1,-1);
\end{tikzpicture}\\
&\epsilon_p=-1 &
&\begin{tikzpicture}[scale=.5]
\draw[line width=1pt, color=red, ->,>=stealth] (-1,1)--(1,-1);
\draw[line width=1pt, color=blue, ->,>=stealth] (1,1)--(-1,-1);
\end{tikzpicture}
& 
&\begin{tikzpicture}[scale=.5]
\draw[line width=1pt, color=blue, ->,>=stealth] (-1,1)--(1,-1);
\draw[line width=1pt, color=green, ->,>=stealth] (1,1)--(-1,-1);
\end{tikzpicture}
&
&\begin{tikzpicture}[scale=.5]
\draw[line width=1pt, color=green, ->,>=stealth] (-1,1)--(1,-1);
\draw[line width=1pt, color=red, ->,>=stealth] (1,1)--(-1,-1);
\end{tikzpicture}\nonumber
\end{align}

\item Choose integrals $\ell_A\in A$, $\ell_B\in B$, $\ell_C\in C$. Associate to each red, blue or green curve the integral $\ell_A$, $\ell_B$, $\ell_C$, respectively,  and consider the element $\ell_A^{\oo g}\oo \ell_B^{\oo g}\oo \ell_C^{\oo g}\in A^{\oo g}\oo B^{\oo g}\oo C^{\oo g}$.

\item If a curve $\lambda$ of colour $X_\lambda$ has $n_\lambda$ intersection points with curves of different colours, apply the map
$\Delta^{(n_\lambda-1)}:X_\lambda\to X_\lambda^{\oo n_\lambda}$ to the integral $\ell_{X_\lambda}\in X_\lambda$ corresponding to $\lambda$ in the tensor product $\ell_A^{\oo g}\oo \ell_B^{\oo g}\oo \ell_C^{\oo g}$, where the order of the tensor factors in $\Delta^{(n_\lambda-1)}(\ell_{X_\lambda})$ corresponds to the  order of the intersection points on $\lambda$, defined by the orientation and basepoint.

Applying this to all curves of $T$ yields an element of $A^{\oo n_A}\oo B^{\oo n_B}\oo C^{\oo n_C}$, where $n_A$, $n_B$ and $n_C$ are the number of intersection points on red, blue and green curves, respectively.

\item Apply a symmetric braiding to the resulting element of $A^{\oo n_A}\oo B^{\oo n_B}\oo C^{\oo n_C}$
such that for each intersection point $p\in \lambda\cap \mu$  
 with colours $(X_\lambda,X_\mu)\in \{(A,B), (B,C), (C,A)\}$
 the associated  factor of $\Delta^{(n_\lambda-1)}(\ell_{X_\lambda})$ is directly to the left to the one of $\Delta^{(n_\mu-1)}(\ell_{X_\mu})$.  Apply the 
  skew pairing $\tau_{X_\lambda X_\mu}: X_\lambda\oo X_\mu\to\C$  to these consecutive factors, if $\epsilon_p=1$, or its convolution inverse, if $\epsilon_p=-1$.
\end{enumerate}
\end{Definition}

\begin{Remark}\label{rem:rescale} While \cite{CCC} uses a specific  normalisable integral, the cotrace, we admit arbitrary normalisable integrals, which are all related by non-zero rescalings. A rescaling $\ell_X\mapsto z_X\cdot \ell_X$ with $z_X\in\C^\times$ for $X=A,B,C$ rescales  the trisection bracket by $\langle T\rangle_{A,B,C}\mapsto z_A^g\cdot z_B^g\cdot z_C^g\cdot \langle T\rangle_{A,B,C}$.
\end{Remark}

By \cite[Lemma 3.2]{CCC}  the trisection bracket is independent of the orientations of the curves in the trisection diagram. This is essentially a consequence of the identity $S(\ell)=\ell$ for a normalisable integral of a finite-dimensional complex semisimple Hopf algebra. Independence of the choice of the basepoint follows from  Lemma \ref{lem:haarprops}, 4. 
It is  shown in \cite{CCC} that the trisection bracket also has the following properties that ensure invariance under the moves (i) to (iv) in Theorem \ref{th:intmoves}.

\begin{Proposition} \cite[Prop 3.3]{CCC} 
The trisection bracket satisfies:\\[-3ex]
\begin{compactenum}
\item \textbf{Diffeomorphism invariance:} If trisection diagrams $T,T'$ are related by  diffeomorphisms, then  their trisection brackets coincide. \\[-2ex]

\item \textbf{Isotopy invariance:} If trisection diagrams $T,T'$ are related by two-point and three-point moves, then their trisection brackets coincide. \\[-2ex]

\item \textbf{Handle slides:} If trisection diagrams $T,T'$ are related by handle slides, 
then their trisection brackets coincide.\\[-2ex]

\item \textbf{Connected sums:}  The trisection bracket of a connected sum of trisection diagrams $T$ and $T'$, taken in such a way that no curves intersect the removed discs, is given by
$$\langle T\# T'\rangle_{A,B,C}=\langle T\rangle_{A,B,C}\cdot \langle T'\rangle_{A,B,C}.$$

\end{compactenum}
\end{Proposition}

To obtain a trisection invariant, one rescales the trisection bracket  by a complex number to ensure invariance under the (de)stabilisation moves from  Theorem \ref{th:intmoves} (v). This requires a non-zero  trisection bracket for the trisection diagram $\Sigma_{st}$ of $S^4$ from \eqref{eq:stabilisation}.

\begin{Definition} \label{def:cccinv}\cite[Def 3.6]{CCC}  Let $(A,B,C)$ be a Hopf triplet 
such that $\langle \Sigma_{st}\rangle_{A,B,C}\neq 0$ and $\xi\in\C^\times$ with $\xi^3=\langle \Sigma_{st}\rangle_{A,B,C}$. 
The
 \textbf{trisection invariant} of a trisection diagram $T$ is
$$
I_{A,B,C}(X)=\xi^{-g}\cdot \langle T\rangle_{A,B,C}.
$$
\end{Definition}

\begin{Remark} The trisection invariant is invariant under a rescaling of the normalisable integrals of $A,B,C$ as in Remark
    \ref{rem:rescale}. Such a rescaling yields  $\langle \Sigma_{st}\rangle_{A,B,C}\mapsto z_A^3\cdot z_B^3\cdot z_C^3\cdot \langle \Sigma_{st}\rangle_{A,B,C}$ and can be absorbed by rescaling $\xi\mapsto z_A\cdot z_B\cdot z_C\cdot \xi$. 
\end{Remark}

The rescaling in Definition \ref{def:cccinv} ensures invariance under  \emph{all} moves from Theorem \ref{th:intmoves},  including the (de)stabilisation move. It defines a 4-manifold invariant, as outlined at the end of Section \ref{sec:trisection}.

\begin{Theorem}\cite[Th 3.7]{CCC} \label{th:ccc} Let $(A,B,C)$ be a Hopf triplet and $\xi\in\C^\times$ with $\xi^3=\langle \Sigma_{st}\rangle_{A,B,C}$. Then the associated trisection invariant is an invariant of closed connected 4-manifolds.
\end{Theorem}

The article  \cite{CCC} investigates several examples of this trisection invariant. The main example takes as algebraic input 
the Hopf triplet in Example \ref{ex:BBtriplet}.
The associated  trisection invariant is related to the  4-manifold invariant introduced  by B\"arenz and Barrett \cite{BB}, which in turn generalises the Crane-Yetter invariant \cite{CY, CYK}, Broda's invariant \cite{Br}, its refinement by Roberts \cite{Ro1,Ro2} and Petit's dichromatic invariant \cite{P}. 

The algebraic data for the B\"arenz-Barrett invariant is a pivotal functor $F:\mac\to \mathcal D$ from a spherical fusion category $\mac$ to a   ribbon fusion category $\mathcal D$ with a trivial twist on the transparent objects. This includes in particular pivotal functors into a modular fusion category $\mathcal D$.

The Hopf triplet in Example \ref{ex:BBtriplet} defines a spherical fusion category $\mac=K\Mod$, a modular fusion category $\mathcal D=D(H)\Mod$, and a pivotal functor $F=\phi^*: K\Mod\to D(H)\Mod$.
It is shown in \cite{CCC} that the associated trisection invariant coincides with the B\"arenz-Barrett invariant, up to a factor given by the algebraic data and the Euler characteristic.
We  show in Section \ref{subsec:bbinv} that an analogous, but more general relation to B\"arenz-Barrett invariants holds for our  trisection invariants.

\begin{Theorem}\cite[Th 5.7]{CCC}\label{th:cccbb} \\
For any closed connected 4-manifold $X$, the trisection invariant $I_{A,B,C}(X)$ defined by the Hopf triplet $(A,B,C)$ in Example \ref{ex:BBtriplet}
and the B\"arenz-Barrett invariant $I_{\phi^*}(X)$ for $\phi^*$ are related by
$$
I_{A,B,C}(X)=m_{\phi}^{2/3-\chi(X)/3} \cdot I_{\phi^*}(X).
$$
\end{Theorem}

The factor $m_{\phi}\in\C^\times$ depends only on the  Hopf algebra map $\phi: H\to K$ in Example \ref{ex:BBtriplet}. A concrete expression can be found in \cite[Th 5.7]{CCC}, but is not needed in the following. 

It is also shown in \cite{CCC} that for a quasitriangular finite-dimensional semisimple Hopf algebra $H=K$, which is automatically ribbon, and the Hopf algebra  map  
\begin{align}\label{eq:rmap}
\phi: D(H)\to H, \quad \alpha\oo h\mapsto \langle \alpha, \low R 1\rangle\, \low R 2\cdot h
\end{align}
the trisection invariant for Example \ref{ex:BBtriplet}  gives the Crane-Yetter invariant \cite{CY, CYK} for the ribbon fusion category $H\Mod$.

\begin{Corollary}\cite[Cor 5.8]{CCC} If $H=K$ in Example \ref{ex:BBtriplet}  is quasitriangular and $\phi$ as in \eqref{eq:rmap}, then for any closed connected 4-manifold $X$ the trisection invariant $I_{A,B,C}(X)$ and the Crane-Yetter invariant $CY_{H\Mod}(X)$ of $X$ for $H\Mod$ are related by
$$
I_{A,B,C}(X)=\dim_\C(H)^{1-\chi(X)}\cdot CY_{H\Mod}(X).
$$
\end{Corollary}

Other examples of trisection invariants investigated in \cite{CCC} are trisection invariants for Hopf triplets associated to the eight-dimensional semisimple Hopf algebra $H_8$ and trisection invariants for the Hopf triplet $(A,B,C)$ in Example \ref{ex:Ktriplet}.
Trisection invariants for the latter are computed for a number of 4-manifolds and  conjecturally  related to Kashaev's invariant  \cite{Ka}.

\begin{Conjecture}\cite[Conj 4.3]{CCC} For any $n\geq 2$ and closed connected 4-manifold $X$ the trisection invariant $I_{A,B,C}(X)$ for the Hopf triplet \eqref{eq:kashaev} is related to Kashaev's invariant $K_n(X)$ by
$$
I_{A,B,C}(X)=n^{1+\chi(X)}\cdot K_n(X).
$$    
\end{Conjecture}

\section{Categorical trisection invariants}
\label{sec:cattrisec}
In this section we construct 4-manifold invariants from  trisection diagrams and  different, more general algebraic data than the triples of Hopf algebras considered in the previous section. The input datum is an indecomposable finite semisimple bimodule category over a pair of spherical fusion categories and a pivotal functor into the   spherical fusion category of bimodule endofunctors and bimodule natural transformations between them.  The only additional requirements are that the bimodule categories have bimodule traces in the sense of Schaumann \cite[Def 3.7]{S} and 
that the bimodule categories are \emph{stabilising} with respect to the pivotal functor.
The latter corresponds to the  condition on the trisection bracket of $\Sigma_{st}$ in Definition \ref{def:cccinv}.

The 4-manifold invariants are obtained by labelling trisection diagrams with these algebraic data and interpreting them as  special cases of the diagrammatic calculus in \cite{M}. 
This diagrammatic calculus   associates to each labelled trisection diagram  a complex number, its evaluation. 
After summation over simple labels, the  evaluation of a diagram is  invariant under diffeomorphisms of the surface, two-point and three-point moves and under handle slides. The stabilising condition allows one to rescale the evaluation to achieve invariance under  (de)stabilisation moves. 

We  show  in Section \ref{sec:trisecbimod} that the trisection invariants from \cite{CCC} are  special cases of our construction. In Section \ref{sec:newexamples} we give new examples that go beyond  trisection invariants from Hopf triplets.

\subsection{Bimodule categories,  functors and  natural transformations}
\label{sec:bimodbackground}

We start by summarising the required background on (bi)module categories over pivotal fusion categories, (bi)module endofunctors and (bi)module natural transformations. All  (bi)module categories are assumed to be finite and semisimple. For more background, see  Chapter 7 of the textbook \cite{EGNO}   by Etingof, Gelaki, Nikshych and Ostrik.

\subsubsection{Basic definitions and examples}
\label{subsec:catbasic}

A \textbf{fusion category} over $\C$ is a  finite semisimple  $\C$-linear abelian rigid monoidal category  $\mac$ with a tensor product $\oo:\mac\times\mac\to \mac$  that is bilinear on the morphisms and with $\End_\mac(e)\cong \C$ for the tensor unit $e\in \Ob\mac$. 
A \textbf{pivotal structure} on $\mac$ is a monoidal natural isomorphism  $\omega:**\Rightarrow\id_\mac$, where $**:\mac\to\mac$ is induced by the (left) duals.  
A pivotal fusion category is called \textbf{spherical} if the left and right traces of all endomorphisms agree: 
$\mathrm{tr}^L(\alpha)=\mathrm{tr}^R(\alpha)=:\mathrm{tr}(\alpha)$ for all $\alpha\in \End_\mac(c)$ and $c\in \Ob\mac$. 

The \textbf{dimension} of an object $c$  in a spherical fusion category $\mac$ is $\dim(c)=\mathrm{tr}(1_c)$, and $\dim \mac=\Sigma_{c\in I_\mac} \dim(c)^2$, where $I_\mac$ is a set of representatives of the isomorphism classes of simple objects.

A \textbf{tensor functor} from a fusion category $\mac$ to a fusion category $\mathcal D$ is a $\C$-linear exact monoidal functor $F:\mac\to\mathcal D$. Note that any tensor functor preserves duals and is automatically faithful, see for instance Natale \cite[Sec 2.1]{N}. It is called a \textbf{tensor equivalence}, if it is an equivalence of categories. In this case, the fusion categories are called \textbf{equivalent}.
A tensor functor $F:\mac\to\mathcal D$ between pivotal fusion categories $(\mac,\omega^\mac)$ and $(\mathcal D,\omega^{\mathcal D})$ is called a \textbf{pivotal functor}, if $F(\omega^\mac_c)=\omega^{\mathcal D}_{F(c)}$ for all $c\in \Ob\mac$. Note that any pivotal functor preserves traces and hence dimensions.

\begin{Example} Let $H$ be a finite-dimensional semisimple complex Hopf algebra. The category $H\Mod$ of finite-dimensional $H$-modules and $H$-linear maps is a spherical fusion category with the tensor product, duals and pivotal structure from the category $\mathrm{Vect}_\C$ of finite-dimensional vector spaces.
\end{Example}

\begin{Example}\label{ex:vectgomega}\cite[Ex 2.3.8]{EGNO}\\ Let $G$ be a finite group,  $\omega: G\times G\times G\to \C^\times$ a normalised  3-cocycle. 
 The fusion category
$\mathrm{Vec}_G^\omega$ has
\begin{compactitem}
\item  as objects finite-dimensional $G$-graded vector spaces $V=\oplus_{g\in G} V_g$,
\item  as morphisms linear maps $f: V\to W$ with $f(V_g)\subset W_g$ for all $g\in G$,
\item  simple objects  $\delta^g$ for $g\in G$ with $(\delta^g)_g=\C$ and $(\delta^g)_h=0$ for $g\neq h$, 
\item  the  tensor product $V\oo W=\oplus_{g\in G}(\oplus_{hk=g} V_h\oo W_k)$, 
\item  the associator  given by $a_{\delta^g,\delta^h,\delta^k}=\omega(g,h,k) \id_{\delta^{ghk}}$ on the simple objects,  

\item dual objects $V^*=\oplus_{g\in G} V_{g^\inv}^*$ with the  usual (co)evaluation for  vector spaces,
\item pivotal structures  in bijection with characters $\kappa: G\to\C^\times$. 
 \end{compactitem}
A pivotal structure is spherical iff $\kappa(g)\in \{1,-1\}$ for all $g\in G$.
In particular, the trivial character $\kappa\equiv 1$ determines a spherical structure on  $\mathrm{Vec}_G^\omega$, its {\bf standard spherical structure}.
\end{Example}

The fusion categories $\mathrm{Vec}_G^\omega$ and $\mathrm{Vec}_G^{\omega'}$ are equivalent if and only if $\omega,\omega'$ are related by a 3-coboundary \cite[Prop 2.6.1]{EGNO}. Note  that the fusion category $\mathrm{Vec}^\omega_G$ is equivalent to the  category $H\Mod$ for a Hopf algebra $H$ if and only if $\omega$ is a 3-coboundary.  In this case $H=\C^G=\C[G]^*$ is the Hopf algebra of complex functions on $G$. 

Throughout the article we details  on finite semisimple (bi)module categories over fusion categories as well as (bi)module functors and (bi)module natural transformations between them. To simplify the terminology we include finiteness and semisimplicity in the definition. The coherence isomorphisms of a fusion category are denoted   $a:\oo (\oo\times\id)\Rightarrow \oo(\id\times\oo)$,  $L: e\oo - \Rightarrow \id$ and $R:-\oo e\Rightarrow \id$.

\pagebreak

\begin{Definition}\label{def:modulecat} \cite[Defs 7.1.1, 7.1.7]{EGNO}
Let $\mac, \mad$ be fusion categories.
\begin{compactenum}
\item  A {\bf $\mac$-(left) module category} is a finite semisimple $\C$-linear abelian category $\calm$  with
\begin{compactitem}
\item a functor $\rhd:\mac\times\calm\to\calm$,  $\C$-bilinear on the morphisms and exact in the first variable, 
\item  natural isomorphisms 
$l: \rhd(\oo\times\id_\calm )\Rightarrow \rhd(\id_\mac\times\rhd)$ and $\lambda: e\rhd -\Rightarrow \id_\calm$
\end{compactitem}
such that  the following diagrams commute for all  $x,y,z\in\Ob \mac$ and $m\in\Ob \calm$: 
\begin{align}\label{eq:pentagoncdef}
&\xymatrix{
((x\oo y)\oo z)\rhd m \ar[r]^{l_{x\oo y, z,m}} \ar[d]_{a_{x,y,z}\rhd 1_m}& \;(x\oo y)\rhd (z\rhd m) \ar[r]^{l_{x,y,z\rhd m}} & x\rhd(y\rhd(z\rhd m))\\
(x\oo(y\oo z)) \rhd m \ar[r]_{l_{x,y\oo z, m}}& x\rhd ((y\oo z)\rhd m)\ar[ru]_{\quad1_x\rhd l_{y,z,m}}
}\\
&\xymatrix{ (x\oo e)\rhd m \ar[rr]^{l_{x,e,m}} \ar[rd]_{R_X\rhd 1_m}  & & x\rhd (e\rhd m) \ar[ld]^{1_x\rhd \lambda_m}\\
&x\rhd m.
}\label{eq:trianglec}
\end{align}

\item A {\bf $\mad$-right module category}
is  a finite semisimple $\C$-linear abelian category $\calm$  with 
\begin{compactitem}
\item a functor $\lhd:\calm\times \mad\to\calm$, $\C$-bilinear on the morphisms and exact in the second variable, 
\item  natural isomorphisms 
$r: \lhd(\lhd \times\id_\mad)\Rightarrow \lhd(\id_\calm\times\oo)$ and $\rho: -\lhd e\Rightarrow \id_\calm$
\end{compactitem}
that make the following diagrams commute for all  $x,y,z\in\Ob \mad$ and $m\in\Ob \calm$
\begin{align}\label{eq:pentagonddef}
&\xymatrix{
((m\lhd x)\lhd y)\lhd z \ar[r]^{r_{m\lhd x, y,z}} \ar[d]_{r_{m,x,y}\lhd 1_z}& (m\lhd x)\lhd (y\oo z)\; \ar[r]^{r_{m,x,y\oo z}} & m\lhd(x\oo(y\oo z))\\
(m\lhd(x\oo y)) \lhd z \ar[r]_{r_{m,x\oo y,z}}& m\lhd ((x\oo y)\oo z)\ar[ru]_{\quad 1_m\lhd a_{x,y,z}}
}\\
&\xymatrix{ m\lhd (e\oo x) \ar[rr]^{r_{m,e,x}} \ar[rd]_{1_m\lhd L_X}  & & (m\lhd e)\lhd x \ar[ld]^{\rho_m\lhd 1_x}\\
&m\lhd x.
}\label{eq:triangled}
\end{align}

\item A  {\bf $(\mac,\mad)$-bimodule category}  is a  finite semisimple $\C$-linear abelian category $\calm$  with
\begin{compactitem}
\item a $\mac$-left module category structure $(\rhd,l,\lambda)$, 
\item a $\mad$-right module category structure $(\lhd, r,\rho)$, 
\item a natural isomorphism
 $q: \lhd(\rhd\times\id_\mad)\Rightarrow \rhd(\id_\mac\times\lhd)$
 \end{compactitem}
such that  the following diagrams commute for all  $x,y\in\Ob \mac$, $u,v\in\Ob \mad$ and $m\in\Ob \calm$
\begin{align}\label{eq:pentagonbdef}
&\xymatrix{
((x\oo y)\rhd m)\lhd u \ar[r]^{q_{x\oo y, m,u}} \ar[d]_{l_{x,y,m}\lhd 1_u}& (x\oo y)\rhd (m\lhd u) \ar[r]^{l_{x,y,m\lhd u}} & x\rhd(y\rhd(m\lhd u))\\
 (x\rhd(y\rhd m))\lhd u\; \ar[r]_{q_{x, y\rhd m, u}}& \; x\rhd ((y\rhd m)\lhd u)\ar[ru]_{\quad 1_x\rhd q_{y,m,u}}
}\\
&\xymatrix{
((x\rhd m)\lhd u)\lhd v \ar[r]^{r_{x\rhd m,u,v}} \ar[d]_{q_{x,m,u}\lhd 1_v}& (x\rhd m)\lhd (u\oo v) \ar[r]^{q_{x,m,u\oo v}} & x\rhd(m\lhd(u\oo v))\\
 (x\rhd(m\lhd u))\lhd v\; \ar[r]_{q_{x, m\lhd u, v}}& \;x\rhd ((m\lhd u)\lhd v).\ar[ru]_{\quad 1_x\rhd r_{y,m,u}}
}\nonumber
\end{align}
\end{compactenum}
A $\mac$-module, $\mad$-right module  or $(\mac,\mad)$-bimodule category $\calm$ is called \textbf{indecomposable}, if it is not a  non-trivial direct sum of  $\mac$-module, $\mad$-right module or $(\mac,\mad)$-bimodule categories, respectively.
\end{Definition}

Note that a $\mad$-right module category is simply a $\mad^{rev}$-module category, where $\mad^{rev}$ is obtained from $\mad$ by reversing the tensor product. 
A $(\mac,\mad)$-bimodule category is a $\mac\boxtimes \mad^{rev}$-module category, where $\boxtimes$ denotes the Deligne product.  

\begin{Example}\label{ex:vecmodcat} For any  finite-dimensional semisimple complex Hopf algebra $H$ the forgetful functor $F:H\Mod\to \mathrm{Vect}_\C$ defines a $H\Mod$-module category structure on $\mathrm{Vect}_\C$. It is given by the functor
$
\rhd=\oo (F\times \id): H\Mod\times\mathrm{Vect}_\C\to\mathrm{Vect}_\C$ and the associator of $\mathrm{Vect}_\C$ as a coherence datum.
\end{Example}

\begin{Example} \label{ex:modvectgomega} \cite[Ex 7.4.10, Cor 7.12.20]{EGNO} \cite[Ex 3.13]{S} \cite[Prop 4.19 and 4.23]{L}
\\ Finite semisimple module categories $\calm$  over  the fusion category $\mathrm{Vec}_G^\omega$ from Example \ref{ex:vectgomega} are given by pairs 
 $(X, \psi)$ of a finite $G$-set $X$ and a 2-cochain $\psi\in C^2(G, (\C^{\times})^X)$ with $d\psi=\omega^\inv$.
\begin{compactitem}
\item The objects of $\calm$ are $X$-graded vector spaces.
\item The simple objects of $\calm$ are of the form $\delta^x$ for $x\in X$ . 

\item The  functor $\rhd: \mathrm{Vec}_G^\omega\times\calm\to\calm$  is given by $\delta^g\rhd \delta^x=\delta^{g\rhd x}$  on the simple objects. 

\item The natural isomorphism $l:\rhd(\oo\times \id)\Rightarrow \rhd(\id\times\rhd)$ is given by   $\psi$.
\end{compactitem}
Here, the $G$-action on $(\C^\times)^X$ is given by $(g\rhd f)(x)=f(g^\inv \rhd x)$ for  $f: G\to \C^\times$, $g\in G$ and $x\in X$.

Indecomposable module categories $\calm$  over   $\mathrm{Vec}_G^\omega$ are classified by pairs $(H,\chi)$ of a subgroup $H\subset G$ 
with $\omega\vert_{H\times H\times H}\in B^3(H,\C^\times)$ and a cohomology class $\chi\in H^2(H,\C^\times)$.
\end{Example}

We now review (bi)module functors between (bi)module categories and (bi)module  natural transformations between them. In summary,
(bi)module functors are functors between (bi)module categories,  with a coherence datum that describes their compatibility with the (bi)module category structure. (Bi)module natural transformations are natural transformations that respect this coherence datum.

\begin{Definition}\label{def:modulefunc}
Let $\mac,\mad$ be fusion categories.  
\begin{compactenum}
\item  A {\bf $\mac$-module functor} $F:\calm\to\caln$ between $\mac$-module categories $\calm,\caln$ is a $\C$-linear right exact functor $F:\calm\to\caln$ together with a natural isomorphism
$s: F\rhd \Rightarrow \rhd(\id_\mac\times F)$ that makes the following diagrams commute  for all $x,y\in\Ob\mac$ and $m\in\Ob\calm$
\begin{align}\label{eq:pentagoncfunc}
&\xymatrix{ F((x\oo y)\rhd m) \ar[d]_{F(l_{x,y,m})}\ar[r]^{s_{x\oo y, m}} & (x\oo y)\rhd F(m) \ar[r]^{l_{x,y, F(m)}} & x\rhd (y\rhd F(m))\\
F(x\rhd(y\rhd m)) \ar[r]_{s_{x, y\rhd m}} & x\rhd F(y\rhd m)\ar[ru]_{\;\; 1_x\rhd s_{y,m}}
}\\
&\xymatrix{ F(e\rhd m)\ar[rd]_{F(\lambda_m)} \ar[rr]^{s_{e,m}} & & e\rhd F(m)\ar[ld]^{\lambda_{F(m)}}\\
& F(m).
}\label{eq:trianglecfunc}
\end{align}

\item A {\bf $\mad$-right module functor} $F:\calm\to\caln$ between $\mad$-right module categories $\calm,\caln$ is a $\C$-linear right exact functor $F:\calm\to\caln$ together with a natural isomorphism
$t: \lhd(F\times\id_\mad)\Rightarrow F\lhd$ that 
makes the following diagrams commute for all $x,y\in\Ob\mad$ and $m\in\Ob\calm$ 
\begin{align}\label{eq:pentagondfunc}
&\xymatrix{ (F(m)\lhd x)\lhd y \ar[d]_{t_{m,x}\lhd 1_y}\ar[r]^{r_{F(m),x, y}} &  F(m)\lhd(x\oo y) \ar[r]^{t_{m, x\oo y}} & F(m\lhd (x\oo y))\\
F(m\lhd x)\lhd y \ar[r]_{t_{m\lhd x, y}} & F((m\lhd x)\lhd y) \ar[ru]_{\;\; F(r_{m,x,y})}
}\\
&\xymatrix{ F(m)\lhd e \ar[rd]_{\rho_{F(m)}} \ar[rr]^{t_{m,e}} & & F(m\lhd e)\ar[ld]^{F(\rho_m)}\\
& F(m).
}\label{eq:triangledfunc}
\end{align}

\item A  {\bf $(\mac,\mad)$-bimodule functor}  from  $\calm$ to  $\caln$ is a $\C$-linear right exact functor $F:\calm\to\caln$ with a $\mac$-module functor structure $s$ and a $\mad$-right module functor structure $t$ that make the following diagram commute  for all $x\in\Ob\mac$, $y\in\Ob\mad$ and $m\in\Ob\calm$ 
\begin{align}
\label{eq:hexa}
\xymatrix{
F(x\rhd m)\lhd y\ar[d]_{s_{x,m}\lhd 1_y} \ar[r]^{t_{x\rhd m, y}} & F((x\rhd m)\lhd y)\; \ar[r]^{F(q_{x,m,y})}& \;F(x\rhd (m\lhd y))\ar[d]^{s_{x, m\lhd y}}\\
(x\rhd F(m))\lhd y\; \ar[r]_{\;q_{x, F(m),y}} & \; x\rhd (F(m)\lhd y)  \ar[r]_{1_x\rhd t_{m,y}} & x\rhd F(m\lhd y).
}
\end{align}
\end{compactenum}
An {\bf equivalence} of (bi)module categories is a (bi)module functor that is an equivalence of categories.
\end{Definition}

\begin{Definition}\label{def:modnat} Let $\mac,\mad$ be fusion categories.
\begin{compactenum}
\item A {\bf $\mac$-module natural transformation} $\nu: F\Rightarrow G$ for $\mac$-module functors $F,G: \calm\to\caln$ is a natural transformation  such that the
following diagram commutes for all $x\in\Ob\mac$,   $m\in\Ob\calm$
\begin{align}
\label{eq:cnat}
\xymatrix{ F(x\rhd m)\ar[d]_{s^F_{x,m}} \ar[r]^{\nu_{x\rhd m}} & G(x\rhd m)\ar[d]^{s^G_{x,m}}\\
x\rhd F(m) \ar[r]_{1_x\rhd \nu_m} & x\rhd G(m).
}
\end{align}
\item A {\bf $\mad$-right module natural transformation} $\nu: F\Rightarrow G$ between $\mad$-right module functors $F,G: \calm\to\caln$ is a natural transformation  such that the
following diagram commutes for all $x\in\Ob\mad$,  $m\in\Ob\calm$
\begin{align}
\label{eq:dnat}
\xymatrix{ F(m)\lhd x \ar[d]_{t^F_{m,x}} \ar[r]^{\nu_{m}\lhd 1_x} & G( m)\lhd x\ar[d]^{t^G_{m,x}}\\
F(m\lhd x) \ar[r]_{\nu_{m\lhd x}} & G(m\lhd x).
}
\end{align}
\item A {\bf $(\mac,\mad)$-bimodule natural transformation} $\nu: F\Rightarrow G$ between two $(\mac,\mad)$-bimodule functors $F,G: \calm\to\caln$ is a
$\mac$-module and a $\mad$-right module natural transformation $\nu: F\Rightarrow G$.
\end{compactenum}
\end{Definition}

For fixed $\mac$-module categories $\calm$ and $\caln$ the $\mac$-module functors from $\calm$ to $\caln$ and $\mac$-module natural transformations between them form a category $\Hom_\mac(\calm,\caln)$. Likewise, for  $(\mac,\mad)$-bimodule categories $\calm$, $\caln$  the $(\mac,\mad)$-bimodule functors from $\calm$ to $\caln$ and $(\mac,\mad)$-bimodule natural transformations between them form a category  $\Hom_{(\mac,\mad)}(\calm,\caln)$. 
The categories $\Hom_\mac(\calm,\caln)$ and $\Hom_{(\mac,\mad)}(\calm,\caln)$ are finite semisimple $\C$-linear abelian, see  \cite[Prop 7.11.6]{EGNO} and  \cite[Th 2.16]{ENO}. 

In fact,  module categories, module functors and module natural transformations over  a spherical fusion category $\mac$ form a 2-category $\mathrm{Mod}(\mac)$. Likewise, $(\mac,\mad)$-bimodule categories, functors and natural transformations form a 2-category $\mathrm{BiMod}(\mac,\mad)=\mathrm{Mod}(\mac\boxtimes \mad^{rev})$. 

The  2-categories $\mathrm{Mod}(\mac)$ and  $\mathrm{BiMod}(\mac,\mad)$  are equipped with duals. Any $\mac$-module  functor $F:\calm\to\caln$ between finite semisimple $\mac$-module categories $\calm$, $\caln$ has a left adjoint $F^l:\caln\to\calm$ and a right adjoint $F^r:\caln\to\calm$. This follows, because all module categories are assumed to be finite semisimple and hence exact, which implies that  module functors between them are exact \cite[Prop 7.6.9]{EGNO}.  As  these adjoints are also $\mac$-module functors \cite[Sec 7.12]{EGNO}, the  duals are  given by these adjunctions. 

In particular,  all categories  $\End_\mac(\calm):=\Hom_\mac(\calm,\calm)$ and $\End_{(\mac,\mad)}(\calm):=\Hom_{(\mac,\mad)}(\calm,\calm)$ are left rigid $\C$-linear monoidal categories. The monoidal structure is given by the composition of (bi)module functors.  The left dual of a $\mac$-module endofunctor $F:\calm\to\calm$ is its left adjoint $F^l:\calm \to\calm$. The associated coevaluation morphism is the unit $\mathrm{coev}^L_F=\eta^F: \id_\calm\Rightarrow FF^l$ and the evaluation morphism is the counit  $\ev^L_F=\epsilon^F:F^lF\Rightarrow \id_\calm$
 of the adjunction. 

The categories $\End_\mac(\calm)$ and $\End_{(\mac,\mad)}(\calm)$ are  fusion categories if and only if $\calm$
is indecomposable,   see \cite[Sec 7.3.6]{EGNO} and  the remark after \cite[Def 7.2.12]{EGNO}. In this case  $\calm$ is also a module category over the fusion category $\End_\mac(\calm)$ or $\End_{(\mac,\mad)}(\calm)$ with the action that applies endofunctors of $\calm$ and natural transformations between them to objects and morphisms in $\calm$. 
The fusion category $\End_\mac(\calm)$ is  denoted $\mac^*_\calm$ and called the dual fusion category to $\mac$ with respect to $\calm$ in \cite{EGNO}.

\begin{Example}\label{ex:center} \cite[Prop 7.13.8]{EGNO} For a spherical fusion category $\mac$ considered as a 
$(\mac,\mac)$-bimodule category  $\End_\mac(\mac)$ is  equivalent  to the Drinfeld centre $ Z(\mac)$.
\end{Example}

\begin{Example}\label{ex:hmodmodule} \cite[Ex 7.12.26]{EGNO}\\
For a finite-dimensional semisimple complex Hopf algebra $H$ and the $H\Mod$-module category structure on $\mathrm{Vect}_\C$ from Example \ref{ex:vecmodcat},
the fusion category $\End_{H\Mod}(\mathrm{Vect}_\C)$ is  equivalent to  $H^{*cop}\Mod$. 

The equivalence assigns to a $H^{*cop}$-module $M$ the  endofunctor $F=M\oo -:\mathrm{Vect}_\C\to \mathrm{Vect}_\C$. Its $H\Mod$-module functor structure is given by the coherence isomorphism
$
s: F\rhd\Rightarrow \rhd(\id\times F)
$, whose components for a finite-dimensional $H$-module $N$ and vector space $V$ are
\begin{align*}
s^F_{N,V}: M\oo N\oo V\to N\oo M\oo V, \quad m\oo n\oo v\mapsto \Sigma_i (b_i\rhd n)\oo (\beta^i\rhd m)\oo v,
\end{align*}
where $\{b_i\}$ is a basis of $H$ and $\{\beta^i\}$ the dual basis of $H^*$. This corresponds to the action of the universal $R$-matrix $R=\Sigma_i (1\oo b_i)\oo (\beta^i\oo 1)$ of the Drinfeld double $D(H)$.
\end{Example}

\subsubsection{(Bi)module categories with (bi)module traces}
\label{subsec:bimoduletr}

(Bi)module traces on (bi)module categories over  fusion categories were   introduced by Schaumann in \cite{S, G14}.
We summarise the main results from \cite{S,G14}, restricting attention to pivotal fusion categories. For ease of notation we focus on left module categories and identify $(\mac,\mad)$-bimodule categories with $\mac\boxtimes\mad^{rev}$-left module categories.

\begin{Definition}\label{def:moduletrace} \cite[Def 3.7]{S} Let $\mac$ be a pivotal fusion category and $\calm$  a  $\mac$-module category.\\
A {\bf trace} on $\calm$ is a collection of maps $\mathrm{tr}_m: \End_\calm(m)\to \C$ that satisfy
\begin{compactenum}
\item {\bf cyclicity:} $\mathrm{tr}_m(\beta\circ \alpha)=\mathrm{tr}_{m'}(\alpha\circ \beta)$ for all $\alpha\in \Hom_\calm(m,m')$ and $\beta\in \Hom_\calm(m',m)$.
\item {\bf non-degeneracy:} for all $m,m'\in\Ob\calm$ the following bilinear map is non-degenerate
$$
\Hom_\calm(m',m)\times \Hom_\calm(m,m')\to \C,\quad (\beta,\alpha)\mapsto \mathrm{tr}_m(\beta\circ \alpha).
$$
\end{compactenum}
A trace  on $\calm$ is called a {\bf $\mac$-module trace} if $\mathrm{tr}_{c\rhd m}(\alpha)=\mathrm{tr}_{m}(\mathrm{tr}^\mac(\alpha))$ for all $\alpha\in \End_\calm(c\rhd m)$, where $\mathrm{tr}^\mac(\alpha)$  is the partial trace of $\alpha$ with respect to $\mac$, in the sense of \cite[(3.24)]{S}).
\end{Definition}

A {\bf $\mad$-right module trace} on a $\mad$-right module category $\calm$ is defined as a $\mad^{rev}$-module trace. A
{\bf $(\mac,\mad)$-bimodule trace} on a $(\mac,\mad)$-bimodule category $\calm$ is a $\mac\boxtimes \mad^{rev}$-module trace.

The \textbf{dimension} of an object $m$ in a $\mac$-module category $\calm$ with a $\mac$-module trace is defined as $\dim(m)=\mathrm{tr}_m(1_m)$. This implies $\dim(c\rhd m)=\dim(c)\dim(m)$ for all $c\in \Ob\mac$ and $m\in\Ob\calm$.

Existence of a $\mac$-module trace on a $\mac$-module category $\calm$ is a condition on  the dimensions of simple objects of $\mac$ and $\calm$, see \cite[Sec 5]{S}.
In particular, each spherical fusion category $\mac$ with its canonical $(\mac,\mac)$-bimodule category structure has a bimodule trace. Uniqueness of a $\mac$-module
trace is easy to establish for indecomposable $\mac$-module categories.

\begin{Proposition}\label{prop:uniquetrace}
\cite[Prop 4.4]{S}  A $\mac$-module trace on an indecomposable $\mac$-module category $\calm$ is unique up to a  rescaling $\mathrm{tr}_m\to z\cdot \mathrm{tr}_m$ for all $m\in\Ob\calm$ with  a global constant $z\in\C^\times$. 
\end{Proposition}

\begin{Example}\cite[Ex 3.13]{S} \label{ex:tracevect} 
An indecomposable $\mathrm{Vec}_G^\omega$-module category $\calm$ given by a pair $(G/H,\psi)$ as in Example \ref{ex:modvectgomega}
admits a module trace if and only if the character $\kappa: G\to\C^\times$  from Example \ref{ex:vectgomega} that defines the pivotal structure on $\mathrm{Vec}_G^\omega$  satisfies $\kappa\vert_H\equiv 1$. 
In particular,
any indecomposable semisimple module category over $\mathrm{Vec}^{\omega}_G$ with the standard spherical structure admits a module trace.
\end{Example}

It is shown in \cite[Th 4.5]{S} and \cite[Th 5.8]{G14} 
that for any pivotal fusion category $\mac$ the $\mac$-module categories with $\mac$-module traces, $\mac$-module functors between them and $\mac$-module transformations form a \emph{pivotal} 2-category $\mathrm{Mod}^\tr(\mac)$ in the sense of  \cite[Def A.12]{G14}.
This implies the following result for the category $\End_\mac(\calm)$ of $\mac$-module endofunctors of $\calm$ and natural transformations between them.

\begin{Proposition}\cite[Th 5.16]{Mu} 
\cite[Prop 5.10]{S}\\
For any spherical fusion category $\mac$ and indecomposable $\mac$-module category $\calm$ with a $\mac$-module trace the fusion category $\End_\mac(\calm)$ is spherical.
\end{Proposition}

\begin{Corollary}\label{lem:macmodule}\cite[Cor 4.6]{S} Let $\mac$ be a spherical fusion category and  $\calm$  a $\mac$-module category with a $\mac$-module trace. Then the $\mac$-module trace is also an $\End_\mac(\calm)$-module trace.  
\end{Corollary}

More concretely, the pivotal structure of $\End_\mac(\calm)$ is given as follows. 
The component morphisms
 $\omega^F: F^{ll}\Rightarrow F$ of the pivotal structure are characterised by the condition
\begin{align}
\label{eq:pivcond}
\mathrm{tr}_{F^l(n)}(\beta\circ\epsilon^F_m\circ F^l(\alpha))=\mathrm{tr}_n (\epsilon^{F^l}_n\circ F^{ll}(\beta)\circ(\omega^{F}_m)^{\inv} \circ \alpha)
\end{align}
for all $\alpha\in \Hom_\calm(n, F(m))$ and $\beta\in \Hom_\calm(m, F^l(n))$, see  \cite[Eq (4.14)]{S}.
They define the right (co)evaluations  
$\mathrm{coev}^R_F=\eta'^{F}=F^l\omega^F\circ \eta^{F^l}:\id_\calm\Rightarrow F^lF$ and $\mathrm{ev}^R_F=\epsilon'^F= \epsilon^{F^l}\circ (\omega^{F} )^{\inv}F^l:FF^l \Rightarrow \id_\calm.$
Corollary \ref{lem:macmodule} then states  that for all $m,n\in\Ob\calm$ and  $\alpha\in \End_\calm(F(m))$
\begin{align}\label{eq:traceconds}
&\mathrm{tr}_m(\epsilon^{F}_m\circ F^l(\alpha)\circ \eta'^F_m)=\mathrm{tr}_{F(m)}(\alpha).
\end{align}

\subsection{Diagrammatic calculus}
\label{sec:diagcalc}

In this section we summarise the diagrammatic calculus for bimodule categories, bimodule functors and bimodule natural transformations from \cite{M} for fixed spherical fusion categories $\mac$ and $\mad$. 
We only require a special case involving a single  
$(\mac,\mad)$-bimodule category $\calm$ with a bimodule trace and bimodule natural transformations  that are identities or (co)evaluations  for duals.

The diagrams from \cite{M} are constructed in several steps. The first step, summarised in Section \ref{subsec:simplediags}, involves diagrams that are superpositions of diagrams for the spherical fusion categories $\mac$, $\mad$ and $\cale:=\End_{(\mac,\mad)}(\calm)$. In a second step, summarised in Section \ref{subsec:bimodeval}, these diagrams are evaluated on  
 objects and morphisms of $\calm$. In
Section \ref{subsec:tracediags} we include
 bimodule traces into the diagrammatic calculus. Section \ref{subsec:tricchord} then defines the  evaluations of  chord diagrams, whose chords, boundary segments and chord endpoints are labelled by objects in $\mad$, $\mac$, $\cale$, objects of $\calm$ and morphisms of $\calm$ respectively.
 Section \ref{subsec:surfdiag} generalises this to  diagrams involving curves on oriented surfaces

\subsubsection{Diagrams for bimodule categories,  functors and natural transformations}
\label{subsec:simplediags}

The starting point for this section is the usual graphical calculus for pivotal monoidal categories. Briefly, diagrams are an embedding of lines and vertices into a rectangle, where lines may not overlap except at the vertices. The lines are oriented and labelled by objects, where a line pointing upwards and labelled with $c\in \Ob\mac$ corresponds to a line that points downwards and is labelled with the dual $c^*$.  Maxima and minima correspond to evaluations and coevaluations. Each vertex is labelled by a morphism from the tensor product of the objects labelling the lines above the vertex  to the tensor product of the objects labelling the lines below.
Each diagram can then be split into a series of horizontal slices.   Tensoring the morphisms in each horizontal slice and then composing the resulting morphisms gives a morphism from the tensor product of the objects labelling the lines at the top of the diagram, to the tensor product of the objects labelling the lines at the bottom of the diagram. See for example Turaev and Virelizier \cite{TVBook} for more details.

The  diagrams in \cite{M} are  obtained by stacking 
three such diagrams, one for the spherical fusion category $\mac$, one for $\mad$ and one for $\cale:=\End_{(\mac,\mad)}(\calm)$, where the labelling of the diagrams in $\mac$ and $\mad$ are reinterpreted as objects and morphisms in $\cale$, as explained below.
These three diagrams are required to overlap in such a way that all crossing points are transversal double crossings. All diagrams considered in the following are as in \cite[Sec 2.2.3]{M}, but the black, grey and dashed lines in \cite[Sec 2.2.3]{M} are replaced by green, blue and red lines in this article.

In more detail, the diagrams involve lines labelled with
\begin{compactitem}
\item  objects  $c\in \Ob\mac$,  representing endofunctors $c\rhd -: \mathcal M\to\mathcal M$,  

\item  objects  $b\in \Ob\mad$,  representing endofunctors $ -\lhd b: \mathcal M\to\mathcal M$,

\item  $(\mac,\mad)$-bimodule endofunctors $F:\mathcal M\to\mathcal M$,
\end{compactitem}
and vertices labelled with 
\begin{compactitem}
\item  morphisms  $\gamma: c\to c'$ in $\mac$,  representing natural transformations $\gamma\rhd -: c\rhd -\Rightarrow c'\rhd -$,  

\item  morphisms  $\beta: b\to b'$ in $\mad$,  representing natural transformations $-\lhd \beta:  -\lhd b\Rightarrow -\lhd b'$,  

\item  $(\mac,\mad)$-bimodule natural transformations $\nu:F\Rightarrow F'$.
\end{compactitem}

Maxima and minima in the diagram are labelled with the evaluations and coevaluations of the associated spherical fusion category. They represent the following natural transformations of endofunctors
\begin{compactitem}
\item  for each $c\in\Ob\mac$
\begin{align*}
&\mathrm{ev}^L_c \rhd -: (c^*\oo c)\rhd -\Rightarrow \id_\calm & &\mathrm{coev}^L_c \rhd -: \id_\calm \Rightarrow (c\oo c^*)\rhd -\\
&\mathrm{ev}^R_c \rhd -: (c\oo c^*)\rhd -\Rightarrow \id_\calm  & &\mathrm{coev}^R_c \rhd -: \id_\calm \Rightarrow (c^*\oo c)\rhd -  
\end{align*}

\item  for each $b\in\Ob\mad$
\begin{align*}
&-\lhd \mathrm{ev}^L_b: -\lhd (b^*\oo b)\Rightarrow \id_\calm & 
&-\lhd \mathrm{coev}^L_b : \id_\calm \Rightarrow -\lhd (b\oo b^*)\\
&-\lhd \mathrm{ev}^R_b: -\lhd (b\oo b^*)\Rightarrow \id_\calm  & &-\lhd \mathrm{coev}^R_b: \id_\calm \Rightarrow -\lhd (b^*\oo b) 
\end{align*}
\item for each $F\in \Ob \cale$
\begin{align*}
&\mathrm{ev}^L_F: F^lF\Rightarrow \id_\calm & & \mathrm{coev}^L_F : \id_\calm \Rightarrow FF^l\\
&\mathrm{ev}^R_F: FF^l\Rightarrow \id_\calm  & & \mathrm{coev}^R_F: \id_\calm \Rightarrow F^lF .
\end{align*}
\end{compactitem}

Crossings only occur between lines 
labelled with distinct spherical fusion categories 
and are labelled with the coherence isomorphisms from Definitions \ref{def:modulecat} and \ref{def:modulefunc}, in particular crossings between
\begin{compactitem}
\item  $c\in\Ob\mac$ and  $b\in\Ob\mad$ with 
\begin{align}\label{pic:bimodule}
&\begin{tikzpicture}[scale=.3, baseline=(current bounding box.center)]
\draw[line width=1pt, color=blue] (-2,2) node[anchor=south]{$b$}--(2,-2)node[sloped, pos=0.3, allow upside down]{\arrowOut};
\draw[line width=1pt,color=green,] (2,2) node[anchor=south]{$c$}--(-2,-2)node[sloped, pos=0.3, allow upside down]{\arrowOut};
\end{tikzpicture}
&
&\begin{tikzpicture}[scale=.3,baseline=(current bounding box.center)]
\draw[line width=1pt, color=blue] (2,2) node[anchor=south]{$b$}--(-2,-2)node[sloped, pos=0.3, allow upside down]{\arrowOut};
\draw[line width=1pt,color=green,] (-2,2) node[anchor=south]{$c$}--(2,-2) node[sloped, pos=0.3, allow upside down]{\arrowOut};
\end{tikzpicture}
\\
&q_{c,-,b}: (c\rhd -)\lhd b\Rightarrow c\rhd(-\lhd b)
&
&q^\inv_{c,-,b}:  c\rhd(-\lhd b)\Rightarrow (c\rhd -)\lhd b,
 \nonumber
\end{align}
\item   $c\in\Ob\mac$ and  $F\in\Ob\cale$ with 
\begin{align}\label{pic:modulefunc}
&\begin{tikzpicture}[scale=.3, baseline=(current bounding box.center)]
\draw[line width=1pt, color=red] (-2,2) node[anchor=south]{$F$}--(2,-2) node[sloped, pos=0.3, allow upside down]{\arrowOut};
\draw[line width=1pt, color=green,] (2,2) node[anchor=south]{$c$}--(-2,-2) node[sloped, pos=0.3, allow upside down]{\arrowOut};
\end{tikzpicture}
&
&\begin{tikzpicture}[scale=.3, baseline=(current bounding box.center)]
\draw[line width=1pt,  color=red] (2,2) node[anchor=south]{$F$}--(-2,-2) node[sloped, pos=0.3, allow upside down]{\arrowOut};
\draw[line width=.5pt, color=green,] (-2,2) node[anchor=south]{$c$}--(2,-2) node[sloped, pos=0.3, allow upside down]{\arrowOut};
\end{tikzpicture}
\\
 &s^F_{c,-}: F(c\rhd -)\Rightarrow c\rhd F(-)
& & (s^{F}_{c,-})^{\inv}:  c\rhd F(-)\Rightarrow F(c\rhd -),\nonumber
\end{align}
\item  $b\in\Ob\mad$ and  $F\in\Ob\cale$ with 
\begin{align}\label{pic:modulefuncinv}
&\begin{tikzpicture}[scale=.3, baseline=(current bounding box.center)]
\draw[line width=1pt, color=blue] (-2,2) node[anchor=south]{$b$}--(2,-2) node[sloped, pos=0.3, allow upside down]{\arrowOut};
\draw[line width=1pt,  color=red,] (2,2) node[anchor=south]{$F$}--(-2,-2)node[sloped, pos=0.3, allow upside down]{\arrowOut};
\end{tikzpicture}
&
&\begin{tikzpicture}[scale=.3, baseline=(current bounding box.center)]
\draw[line width=1pt, color=blue] (2,2) node[anchor=south]{$b$}--(-2,-2)node[sloped, pos=0.3, allow upside down]{\arrowOut};
\draw[line width=1pt, color=red] (-2,2) node[anchor=south]{$F$}--(2,-2)node[sloped, pos=0.3, allow upside down]{\arrowOut};
\end{tikzpicture}\\
 &t^F_{-,b}: F(-)\lhd b\Rightarrow F(-\lhd b) 
& & (t^{F}_{-,b})^{\inv}:  F(-\lhd b)\Rightarrow F(-)\lhd b.\nonumber
\end{align}
\end{compactitem}

\medskip
Each such diagram $D$ represents a natural transformation $\mathrm{ev}(D): T\Rightarrow B$ between  endofunctors $T:\calm\to\calm$ and $B:\calm \to\calm$ associated with the top and bottom of the diagram. The endofunctors $T,B$ are obtained by composing the endofunctors for the lines at the top and bottom of the diagram. The natural transformation $\mathrm{ev}(D)$ is obtained by 
\begin{compactitem}
\item using the horizontal composition of natural transformations to compose
for each horizontal slice  the natural transformations for the crossings, maxima and minima on the slice with the functors (treated as identity natural transformations) labelling lines to their left and right,
\item composing the resulting natural transformations for each slice from  top to  bottom of the diagram.
\end{compactitem}
The conventions for the horizontal composition  are as in Diagrams \eqref{pic:bimodule} to \eqref{pic:modulefuncinv}.

Due to the coherence theorem for (bi)module categories and the strict associativity of the composition of functors, the crossings are compatible with the tensor products in the three spherical fusion categories, see \cite[Sec 2.3.3]{M}. This ensures that the diagrammatic calculus is consistent.

The axioms in Definitions \ref{def:modulefunc} and \ref{def:modnat}  ensure that crossings are natural with respect to maxima and minima of lines. 
This gives
the following identity and analogous identities  obtained by reflecting it on a vertical or horizontal axis, by reversing line orientations and by permuting line colours. 
\begin{align} \label{ex:maxmin}
\begin{tikzpicture}[scale=.5]
\begin{scope}[shift={(0,0)}]
\draw[red, line width=1pt] (-1,2) node[anchor=south]{$F$}--(-3,-2) node[sloped, pos=0.3, allow upside down]{\arrowOut};
\draw[blue, line width=1pt] (-4,1) .. controls (-4,2) and (-3,2) ..(-3,1);
\draw[blue, line width=1pt] (-3,1).. controls (-3,0)..(0,-2) node[sloped, pos=0.8, allow upside down]{\arrowOut} node[anchor=north]{$b$};
\draw[blue, line width=1pt] (-4,1).. controls (-4,0)..(-1,-2);
\end{scope}
\node at (1,0){$=$};
\begin{scope}[shift={(5,0)}]
\draw[red, line width=1pt] (-1,2)node[anchor=south]{$F$} --(-3,-2) node[sloped, pos=0.3, allow upside down]{\arrowOut};
\draw[blue, line width=1pt] (0,-1).. controls (0,0) and (-1,0)..(-1,-1);
\draw[blue, line width=1pt] (-1,-1)--(-1,-2);
\draw[blue, line width=1pt] (0,-1)--(0,-2) node[sloped, pos=0.5, allow upside down]{\arrowOut} node[anchor=north]{$b$};
\end{scope}
\end{tikzpicture}
\end{align}

It is also shown in \cite[Sec 2.3.3]{M} that the evaluations of these diagrams is invariant under  coloured variants of the Reidemeister 2 move (R2) and Reidemeister 3 move (R3), also called two-point and three-point move.
This follows from the invertibility and naturality of the coherence isomorphisms at the crossings and from  Condition \eqref{eq:hexa} for the coherence isomorphisms for bimodule functors. 
The diagrammatic representations of these identities are the following,
with additional identities being given by reversing line orientations and permuting line colours and labels.
\begin{align}\label{pic:rm2}
&\begin{tikzpicture}[scale=.15]
\begin{scope}[shift={(-4,0)}]
\draw[line width=1pt, color=blue] (-2,2) node[anchor=south]{$b$}--(2,-2)--(-2,-6) node[sloped, pos=0.3, allow upside down]{\arrowOut};
\draw[line width=1pt, color=green] (2,2) node[anchor=south]{$c$}--(-2,-2)--(2,-6)node[sloped, pos=0.3, allow upside down]{\arrowOut};
\end{scope}
\node at (0,-2){$=$};
\begin{scope}[shift={(4,0)}]
\draw[line width=1pt, color=blue] (-2,2) node[anchor=south]{$b$}--(-2,-6) node[sloped, pos=0.5, allow upside down]{\arrowOut};
\draw[line width=1pt, color=green] (2,2) node[anchor=south]{$c$}--(2,-6) node[sloped, pos=0.5, allow upside down]{\arrowOut};
\end{scope}
\end{tikzpicture}\quad
\begin{tikzpicture}[scale=.15]
\begin{scope}[shift={(-4,0)}]
\draw[line width=1pt, color=blue] (2,2) node[anchor=south]{$b$}--(-2,-2)--(2,-6)node[sloped, pos=0.3, allow upside down]{\arrowOut};
\draw[line width=1pt, color=green] (-2,2) node[anchor=south]{$c$}--(2,-2)--(-2,-6)node[sloped, pos=0.3, allow upside down]{\arrowOut};
\end{scope}
\node at (0,-2){$=$};
\begin{scope}[shift={(4,0)}]
\draw[line width=1pt, color=blue] (2,2) node[anchor=south]{$b$}--(2,-6) node[sloped, pos=0.5, allow upside down]{\arrowOut};
\draw[line width=1pt, color=green] (-2,2) node[anchor=south]{$c$}--(-2,-6)node[sloped, pos=0.5, allow upside down]{\arrowOut};
\end{scope}
\end{tikzpicture}\quad
\begin{tikzpicture}[scale=.15]
\begin{scope}[shift={(-4,0)}]
\draw[line width=1pt, red] (-2,2) node[anchor=south]{$F$}--(2,-2)--(-2,-6)node[sloped, pos=0.3, allow upside down]{\arrowOut};
\draw[line width=1pt, color=green] (2,2) node[anchor=south]{$c$}--(-2,-2)--(2,-6)node[sloped, pos=0.3, allow upside down]{\arrowOut};
\end{scope}
\node at (0,-2){$=$};
\begin{scope}[shift={(4,0)}]
\draw[line width=1pt, red] (-2,2) node[anchor=south]{$F$}--(-2,-6)node[sloped, pos=0.5, allow upside down]{\arrowOut};
\draw[line width=.5pt, color=green] (2,2) node[anchor=south]{$c$}--(2,-6)node[sloped, pos=0.5, allow upside down]{\arrowOut};
\end{scope}
\end{tikzpicture}\quad
\begin{tikzpicture}[scale=.15]
\begin{scope}[shift={(-4,0)}]
\draw[line width=1pt, red] (2,2) node[anchor=south]{$F$}--(-2,-2)--(2,-6)node[sloped, pos=0.3, allow upside down]{\arrowOut};
\draw[line width=.5pt, color=green] (-2,2) node[anchor=south]{$c$}--(2,-2)--(-2,-6) node[sloped, pos=0.3, allow upside down]{\arrowOut};
\end{scope}
\node at (0,-2){$=$};
\begin{scope}[shift={(4,0)}]
\draw[line width=1pt, red] (2,2) node[anchor=south]{$F$}--(2,-6)node[sloped, pos=0.5, allow upside down]{\arrowOut};
\draw[line width=1pt, color=green] (-2,2) node[anchor=south]{$c$}--(-2,-6)node[sloped, pos=0.5, allow upside down]{\arrowOut};
\end{scope}
\end{tikzpicture}
\quad \begin{tikzpicture}[scale=.15]
\begin{scope}[shift={(-4,0)}]
\draw[line width=1pt, color=blue] (-2,2) node[anchor=south]{$b$}--(2,-2)--(-2,-6) node[sloped, pos=0.3, allow upside down]{\arrowOut};
\draw[line width=1pt, red] (2,2) node[anchor=south]{$F$}--(-2,-2)--(2,-6) node[sloped, pos=0.3, allow upside down]{\arrowOut};
\end{scope}
\node at (0,-2){$=$};
\begin{scope}[shift={(4,0)}]
\draw[line width=1pt, color=blue] (-2,2) node[anchor=south]{$b$}--(-2,-6) node[sloped, pos=0.5, allow upside down]{\arrowOut};
\draw[line width=1pt, red] (2,2) node[anchor=south]{$F$}--(2,-6) node[sloped, pos=0.5, allow upside down]{\arrowOut};
\end{scope}
\end{tikzpicture}
\quad
\begin{tikzpicture}[scale=.15]
\begin{scope}[shift={(-4,0)}]
\draw[line width=1pt, color=blue] (2,2) node[anchor=south]{$b$}--(-2,-2)--(2,-6)node[sloped, pos=0.3, allow upside down]{\arrowOut};
\draw[line width=1pt, red] (-2,2) node[anchor=south]{$F$}--(2,-2)--(-2,-6)node[sloped, pos=0.3, allow upside down]{\arrowOut};
\end{scope}
\node at (0,-2){$=$};
\begin{scope}[shift={(4,0)}]
\draw[line width=1pt, color=blue] (2,2) node[anchor=south]{$b$}--(2,-6) node[sloped, pos=0.5, allow upside down]{\arrowOut};
\draw[line width=1pt, red] (-2,2) node[anchor=south]{$F$}--(-2,-6)node[sloped, pos=0.5, allow upside down]{\arrowOut};
\end{scope}
\end{tikzpicture}
\end{align}
  \begin{align}\label{pic:hexagon}
 &\begin{tikzpicture}[scale=.3, baseline=(current bounding box.center)]
 \begin{scope}[shift={(-4,0)}]
\draw[blue, line width=1pt] (0,0) node[anchor=south]{$b$}--(2,-2)--(0,-6) node[sloped, pos=0.2, allow upside down]{\arrowOut};
\draw[green, line width=1pt] (2,0) node[anchor=south]{$c$}--(0,-2)--(-2,-6) node[sloped, pos=0.2, allow upside down]{\arrowOut};
\draw[red, line width=1pt] (-2,0) node[anchor=south]{$F$}--(-2,-2)--(2,-6) node[sloped, pos=0.2, allow upside down]{\arrowOut};
    \end{scope}
    \node at (0,-3)[anchor=west]{$=$};
     \begin{scope}[shift={(6,0)}]
     \draw[blue, line width=1pt] (0,0) node[anchor=south]{$b$}--(-2,-4) node[sloped, pos=0.2, allow upside down]{\arrowOut}--(0,-6) ;
\draw[green, line width=1pt] (2,0) node[anchor=south]{$c$}--(0,-4) node[sloped, pos=0.4, allow upside down]{\arrowOut} --(-2,-6);
\draw[red, line width=1pt] (-2,0) node[anchor=south]{$F$}--(2,-4)node[sloped, pos=0.2, allow upside down]{\arrowOut} --(2,-6) ;
    \end{scope}
 \end{tikzpicture}
\end{align}

\subsubsection{Evaluation on  bimodule categories}
\label{subsec:bimodeval}
We now consider the diagrams from \cite[Sec 2.3.1]{M} that describe the interaction of the diagrams in Section \ref{subsec:simplediags}
with objects and morphisms of the underlying bimodule category $\calm$. As the diagrams from  Section \ref{subsec:simplediags} represent endofunctors of $\calm$ and natural transformations between them, they can be evaluated on objects and morphisms of $\calm$. This is represented diagrammatically as follows. 

Objects  and morphisms of $\calm$ label segments and vertices  of a line to the right of the diagram.
Horizontal slices are then chosen in such a way that
each horizontal slice contains at most one vertex on the line, crossing point or maximum or minimum of the diagram.

Lines of the diagram labelled with objects of $\mac$, $\mad$ and $\cale$ may end on a vertex of the vertical line from above or below.
The functors which end at the vertex from above (below) are evaluated on the element labelling the line above (below) the vertex. The vertex must then be labelled with a morphism in $\calm$ from the evaluation above the vertex to the evaluation below, as in the following example 
\begin{align}
&\begin{tikzpicture}[scale=.6]
\draw[line width=1pt] (0,2)--(0,-2) node[pos=.7, anchor=west]{$m$} node[pos=.3, anchor=west]{$m'$};
\draw[line width=1pt, red] (0,0)--(-1,2) node[anchor=south]{$F$} node[sloped, pos=0.7, allow upside down]{\arrowOut};
\draw[fill=black] (0,0) circle (.15) node[anchor=west]{$\;\;\alpha$};
\draw[line width=1pt, blue] (-2,2) node[anchor=south]{$b$}--(0,0) node[sloped, pos=0.5, allow upside down]{\arrowOut} ;
\draw[line width=1pt, green] (-1,-2) node[anchor=north]{$c$}--(0,0) node[sloped, pos=0.5, allow upside down]{\arrowOut} ;
\draw[line width=1pt, red] (0,0)--(-2,-2) node[anchor=north]{$G$} node[sloped, pos=0.5, allow upside down]{\arrowOut} ;
\draw[fill=black] (0,0) circle (.15) node[anchor=west]{$\;\;\alpha$};
\end{tikzpicture}\nonumber\\
&\alpha: F^l(m)\lhd b\to G(c^*\rhd m').
\end{align}

The diagrammatic calculus assigns to each horizontal slice of the diagram a morphism in $\calm$.
If the slice contains a vertex on the vertical line, 
it is obtained by composing the functors on diagram lines intersecting the horizontal slice to obtain an endofunctor of $\calm$. This endofunctor is then applied to  the morphism labelling the vertex, as in the following example
\begin{align}
&\begin{tikzpicture}[scale=.6]
\draw[line width=1pt] (0,2)--(0,-2) node[pos=.3, anchor=west]{$m$} node[pos=.7, anchor=west]{$m'$};
\draw[line width=1pt, red] (0,0)--(-1,2) node[anchor=south]{$F$} node[sloped, pos=0.7, allow upside down]{\arrowOut};
\draw[fill=black] (0,0) circle (.15) node[anchor=west]{$\;\;\alpha$};
\draw[line width=1pt, blue] (-2,2) node[anchor=south]{$b$}--(0,0) node[sloped, pos=0.5, allow upside down]{\arrowOut} ;
\draw[line width=1pt, green] (-1,-2) node[anchor=north]{$c$}--(0,0) node[sloped, pos=0.5, allow upside down]{\arrowOut} ;
\draw[line width=1pt, red] (0,0)--(-2,-2) node[anchor=north]{$G$} node[sloped, pos=0.5, allow upside down]{\arrowOut} ;
\draw[fill=black] (0,0) circle (.15) node[anchor=west]{$\;\;\alpha$};
\draw[line width=1pt, red] (-3,2) node[anchor=south]{$H$}--(-3,-2) node[sloped, pos=0.5, allow upside down]{\arrowOut} ;
\draw[line width=1pt, green] (-4,2) node[anchor=south]{$c'$}--(-4,-2) node[sloped, pos=0.5, allow upside down]{\arrowOut} ;
\end{tikzpicture}\nonumber\\
&c'\rhd H(\alpha): c'\rhd H(F^l(m)\lhd b)\to c'\rhd HG(c^*\rhd m').
\end{align}

Vertices, crossing points, maxima and minima of the diagram only occur on horizontal slices that do not intersect a vertex on the vertical line.  Each of them represents a natural transformation between certain endofunctors of $\calm$ as above, obtained by composing it with all functors on diagram lines to its left and  right. One
then evaluates it on the object of $\calm$ at the intersection of the slice with the vertical line as in the following example
\begin{align}
&\begin{tikzpicture}[scale=.5]
\draw[line width=1pt] (0,2)--(0,-2) node[pos=.5, anchor=west]{$m$};
\draw[line width=1pt, blue] (-1,2) node[anchor=south]{$b$}--(-1,-2) node[sloped, pos=0.5, allow upside down]{\arrowOut};
\draw[line width=1pt, red] (-2,2) node[anchor=south]{$F$}--(-4,-2) node[sloped, pos=0.3, allow upside down]{\arrowOut};
\draw[line width=1pt, red] (-5,2) node[anchor=south]{$G$}--(-5,-2) node[sloped, pos=0.3, allow upside down]{\arrowOut};
\draw[line width=1pt, green] (-4,2) node[anchor=south]{$c$}--(-2,-2) node[sloped, pos=0.3, allow upside down]{\arrowOut};
\end{tikzpicture}\nonumber\\
&G((s^{F}_{c,  m\lhd b})^{\inv}): G(c\rhd F(m\lhd b))\to GF(c\rhd (m\lhd b)).
\end{align}

The construction assigns  an identity morphism, if the horizontal slice does not contain any vertices, crossing points, maxima or minima of the diagram nor vertices on the line at the right.

Composing the morphisms for horizontal slices from the top to the bottom of the diagram yields a morphism in $\calm$. Its source (target) is the image of the object of $\calm$ at the top (bottom) of the vertical line under all functors labelling diagram lines to its left. 

\begin{Example}\label{ex:evaluationmor} The vertices of the following  diagram are labelled with morphisms
\begin{align*}
&\alpha: m\to F(n) & &\beta:n\lhd b\to p & &\gamma: p\to c\rhd q & &\delta: F(q)\to r & &\epsilon: c\rhd r\to s
\end{align*}
and the two crossing points with the natural transformations 
$$t^F_{-,b}: F(-)\lhd b\Rightarrow F(-\lhd b)\qquad s^F_{c,-}: F(c\rhd -)\Rightarrow c\rhd F(-).$$
\begin{align*}
\begin{tikzpicture}[scale=.7]
\draw[line width=1pt] (0,3)--(0,-5);
\draw[line width=1pt, color=red] (0,2.5)..controls (-2,.25)..(0,-2) node[sloped, pos=0.5, allow upside down] {\arrowOut}node[pos=.5, anchor=east] {$F$};
\draw[line width=1pt, color=green] (0,-.5)..controls (-2,-2.25)..(0,-4) node[sloped, pos=0.7, allow upside down] {\arrowOut} node[pos=.7, anchor=north east]{$c$};
\draw[line width=1pt, color=blue] (-2,3) node[anchor=south]{$b$}--(0,1) node[sloped, pos=0.3, allow upside down]{\arrowOut};
\draw[fill=black](0,2.5) circle (.1) node[anchor=west]{$\alpha$};
\draw[fill=black] (0,1) circle (.1) node[anchor=west]{$\beta$};
\draw[fill=black] (0,-.5) circle (.1) node[anchor=west]{$\gamma$};
\draw[fill=black] (0,-2) circle (.1) node[anchor=west]{$\delta$};
\draw[fill=black] (0,-4) circle (.1) node[anchor=west]{$\epsilon$};
\node at (0,3) [anchor=west]{$m$};
\node at (0,1.75) [anchor=west]{$n$};
\node at (0,.25) [anchor=west]{$p$};
\node at (0,-1.25) [anchor=west]{$q$};
\node at (0,-3) [anchor=west]{$r$};
\node at (0,-5) [anchor=west]{$s$};
\end{tikzpicture}
\end{align*}
The diagram evaluates to the following morphism in $\calm$
\begin{align*}
m\lhd b\xrightarrow{\alpha\lhd b} F(n)\lhd b\xrightarrow{t^F_{n,b}} F(n\lhd b) \xrightarrow{F(\beta)} F(p)\xrightarrow{F(\gamma)} F(c\rhd q) \xrightarrow{s^F_{c,q}} c\rhd F(q) \xrightarrow{c\rhd\delta} c\rhd r\xrightarrow{\epsilon} s.
\end{align*}
\end{Example}

\subsubsection{Diagrams for bimodule traces}
\label{subsec:tracediags}
Bimodule traces are integrated into the diagrammatical calculus as in \cite[Sec 2.3.3]{M}. 
The bimodule trace of an endomorphism in $\Hom_\calm(m,m)$, depicted graphically as in the previous section, is given by adding horizontal lines labelled by $ m$. 
 The cyclicity condition from Definition \ref{def:moduletrace} 
takes the form
\begin{align}\label{pic:cyclic trace}
\begin{tikzpicture}[scale=.3]
\node at (-6,0)[anchor=east]{$\mathrm{tr}_m(\beta\circ\alpha)=$};
\begin{scope}[shift={(-4,0)}]
\draw[line width=1pt, color=black] (0,3) --(0,-3);
\draw[line width=1pt, color=black] (-1,3)--(1,3) node[anchor=west]{$m$};
\draw[line width=1pt, color=black] (-1,-3)--(1,-3)node[anchor=west]{$m$};
\draw[fill=black, color=black] (0,1.5) circle (.2) node[anchor=east]{$\alpha\;$};
\draw[fill=black, color=black] (0,-1.5) circle (.2) node[anchor=east]{$\beta\;$};
\node (0,0)[anchor=west, color=black]{$m'$};
\end{scope}
\node at (0,0){$=$};
\begin{scope}[shift={(4,0)}]
\draw[line width=1pt, color=black] (0,3) --(0,-3);
\draw[line width=1pt, color=black] (-1,3)--(1,3) node[anchor=west]{$m'$};
\draw[line width=1pt, color=black] (-1,-3)--(1,-3)node[anchor=west]{$m'$};
\draw[fill=black, color=black] (0,1.5) circle (.2) node[anchor=east]{$\beta\;$};
\draw[fill=black, color=black] (0,-1.5) circle (.2) node[anchor=east]{$\alpha\;$};
\node (0,0)[anchor=west, color=black]{$m$};
\end{scope}
\node at (6,0)[anchor=west]{$=\mathrm{tr}_{m'}(\alpha\circ\beta).$};
\end{tikzpicture}
\end{align}
Compatibility  between  bimodule traces and the partial traces in $\mac$ and $\mad$ from Definition \ref{def:moduletrace}  reads
\begin{align}\label{pic:cmoduletrace}
&\begin{tikzpicture}[scale=.35]
\begin{scope}[shift={(-3,0)}]
\draw[line width=1pt, color=black] (0,3) node[anchor=north west]{$m$} --(0,-3) node[anchor=south west]{$m$};
\draw[line width=1pt, color=black] (-2.5,3)--(1,3);
\draw[line width=1pt, color=black] (-2.5,-3)--(1,-3);
\draw[line width=1pt, color=green] (-2,3)node[anchor=north east]{$c$} --(0,0) node[sloped, pos=0.3, allow upside down]{\arrowOut};
\draw[line width=1pt, color=green] (-2,-3) node[anchor=south east]{$c$} --(0,0) node[sloped, pos=0.7, allow upside down]{\arrowOut};
\draw[fill=black, color=black] (0,0) circle (.2) node[anchor=west]{$\;\alpha$};
\end{scope}
\node at (0,0) {$=$};
\begin{scope}[shift={(5,0)}]
\draw[line width=1pt, color=black] (0,3) node[anchor=north west]{$m$} --(0,-3) node[anchor=south west]{$m$};
\draw[line width=1pt, color=black] (-1,3)--(1,3);
\draw[line width=1pt, color=black] (-1,-3)--(1,-3);
\draw[line width=1pt, green] plot [smooth, tension=0.6] coordinates 
      {(0,0)(-1,2)(-2,1.5) (-2.2,0)(-2,-1.5) (-1,-2)(0,0)};
\draw[line width=1pt, ->,>=stealth, green]  (-2.2,-.1)--(-2.2,.1) node[anchor=east]{$c$};    
\draw[fill=black, color=black] (0,0) circle (.2) node[anchor=west]{$\;\alpha$};
\end{scope}
\end{tikzpicture}
&
&\begin{tikzpicture}[scale=.35]
\begin{scope}[shift={(-3,0)}]
\draw[line width=1pt, color=black] (0,3) node[anchor=north west]{$m$} --(0,-3) node[anchor=south west]{$m$};
\draw[line width=1pt, color=black] (-2.5,3)--(1,3);
\draw[line width=1pt, color=black] (-2.5,-3)--(1,-3);
\draw[line width=1pt, color=blue] (-2,3)node[anchor=north east]{$b$} --(0,0) node[sloped, pos=0.3, allow upside down]{\arrowOut};
\draw[line width=1pt, color=blue] (-2,-3) node[anchor=south east]{$b$} --(0,0) node[sloped, pos=0.3, allow upside down]{\arrowOut};
\draw[fill=black, color=black] (0,0) circle (.2) node[anchor=west]{$\;\alpha$};
\end{scope}
\node at (0,0) {$=$};
\begin{scope}[shift={(5,0)}]
\draw[line width=1pt, color=black] (0,3) node[anchor=north west]{$m$} --(0,-3) node[anchor=south west]{$m$};
\draw[line width=1pt, color=black] (-1,3)--(1,3);
\draw[line width=1pt, color=black] (-1,-3)--(1,-3);
\draw[line width=1pt, color=blue] plot [smooth, tension=0.6] coordinates 
      {(0,0)(-1,2)(-2,1.5) (-2.2,0)(-2,-1.5) (-1,-2)(0,0)};
\draw[line width=1pt, ->,>=stealth, color=blue]  (-2.2,-.1)--(-2.2,.1) node[anchor=east]{$b$};    
\draw[fill=black, color=black] (0,0) circle (.2) node[anchor=west]{$\;\alpha$};
\end{scope}
\end{tikzpicture}
\end{align}
The statement that a $(\mac,\mad)$-bimodule trace is also an $\End_{(\mac,\mad)}(\calm)$-module trace from Corollary \ref{lem:macmodule} and formula \eqref{eq:traceconds} is depicted as follows
\begin{align*}
\begin{tikzpicture}[scale=.35]
\begin{scope}[shift={(-2.5,0)}]
\draw[line width=1pt, color=black] (0,3) node[anchor=north west]{$m$} --(0,-3) node[anchor=south west]{$m$};
\draw[line width=1pt] (-2.5,3)--(1,3);
\draw[line width=1pt] (-2.5,-3)--(1,-3);
\draw[line width=1pt, red] (-2,3)node[anchor=north east]{$F$} --(0,0)node[sloped, pos=0.3, allow upside down]{\arrowOut} ;
\draw[line width=1pt, red] (-2,-3) node[anchor=south east]{$F$} --(0,0) node[sloped, pos=0.7, allow upside down]{\arrowOut};
\draw[fill=black] (0,0) circle (.15) node[anchor=west]{$\;\alpha$};
\end{scope}
\node at (0,0) {$=$};
\begin{scope}[shift={(4.5,0)}]
\draw[line width=1pt, color=black] (0,3) node[anchor=north west]{$m$} --(0,-3) node[anchor=south west]{$m$};
\draw[line width=1pt] (-1,3)--(1,3);
\draw[line width=1pt] (-1,-3)--(1,-3);
\draw[line width=1pt, red] plot [smooth, tension=0.6] coordinates 
      {(0,0)(-1,2)(-2,1.5) (-2.2,0)(-2,-1.5) (-1,-2)(0,0)}node[sloped, pos=0.3, allow upside down]{\arrowOut};
\draw[line width=1pt, ->,>=stealth, red]  (-2.2,-.1)--(-2.2,.1) node[anchor=east]{$F$};    
\draw[fill=black] (0,0) circle (.15) node[anchor=west]{$\;\alpha$};
\end{scope}
\end{tikzpicture}
\end{align*}

\subsubsection{Tricoloured chord diagrams}
\label{subsec:tricchord}

We now consider diagrams analogous to the ones from Section \ref{subsec:bimodeval}, but with the vertical line  replaced by a circle and with additional restrictions on the vertices. 
These diagrams are special cases of the polygon diagrams in
\cite[Sec 2.3.3]{M} and are
obtained by placing the diagrams from Section \ref{subsec:simplediags}  in a disc, such that  lines labelled with objects of $\mac$, $\mad$ and $\cale$ may end on its boundary. 

These diagrams are subject to the following additional restrictions. All  vertices in the interior of the diagram must be crossings, maxima or minima, as opposed to the general natural transformations used to label vertices in Section \ref{subsec:simplediags}. We also impose that each vertex at the boundary of the disc 
either (i) involves no red, blue or green lines or (ii) is the endpoint of a single segment of red, blue or green line or (iii)
is the common endpoint of the two ends of a red, blue or green loop. 
The resulting diagram is a chord diagram, up to the loops in (iii). In the following, we use the term chord diagram for such diagrams, which may include loops as in (iii).

As in Sections \ref{subsec:simplediags} and \ref{subsec:bimodeval}, blue, green and red lines in the interior of the diagram are labelled with objects of the spherical fusion categories $\mad$, $\mac$ and $\cale=\End_{(\mac,\mad)}(\calm)$. Segments of the circle between vertices are labelled with objects of $\calm$ and vertices on the circle  with morphisms in $\calm$. 
In case (i) a vertex between segments labelled with objects
$m$ and $m'$, clockwise, is labelled with a morphism in $\Hom_\calm(m,m')$. In case (ii) it is labelled with
a morphism in 
$\Hom_\calm(m, F(m'))$ or $\Hom_\calm(F(m'),m)$, where $F:\calm \to\calm$ is either a $(\mac,\mad)$-bimodule endofunctor or a functor of the form $F=c\rhd -$
or $F=-\lhd b$ associated with the red, green or blue line. In case (iii) such a vertex is labelled with a morphism in $\Hom_\calm(F(m), F(m'))$ or in $\Hom_\calm(F^l(m), F^l(m'))$. 
\begin{align}\label{eq:boundlabelmor}
\begin{tikzpicture}[scale=.7]
\draw[line width=1pt] (-1,0) node[anchor=east]{$m$} -- (1,0) node[anchor=west]{$m'$};
\draw[line width=1pt] (0,0)--(0,-1.5) node[sloped, pos=0.5, allow upside down]{\arrowOut} node[pos=.5, anchor=west]{$F$};
\node at (0,.5)[anchor=south]{$\Hom_\calm(m,F(m'))$};
\end{tikzpicture}
\qquad
\begin{tikzpicture}[scale=.7]
\draw[line width=1pt] (-1,0) node[anchor=east]{$m$}--(1,0) node[anchor=west]{$m'$};
\draw[line width=1pt] (0,-1.5)--(0,0) node[sloped, pos=0.5, allow upside down]{\arrowOut} node[pos=.5, anchor=west]{$F$};
\node at (0,.5)[anchor=south]{$\Hom_\calm(F(m),m')$};
\end{tikzpicture}
\qquad 
\begin{tikzpicture}[scale=.7]
\draw[line width=1pt] (-1,0) node[anchor=east]{$m$}--(1,0) node[anchor=west]{$m'$};
\draw[line width=1pt] (0,0).. controls (1,-1.5) and (-1,-1.5) .. (0,0) node[sloped, pos=0.9, allow upside down]{\arrowOut} node[pos=.7, anchor=east]{$F$};
\node at (0,.5)[anchor=south]{$\Hom_\calm(F(m),F(m'))$};
\end{tikzpicture}
\qquad
\begin{tikzpicture}[scale=.7]
\draw[line width=1pt] (-1,0) node[anchor=east]{$m$}--(1,0) node[anchor=west]{$m'$};
\draw[line width=1pt] (0,0).. controls (-1,-1.5) and (1,-1.5) .. (0,0) node[sloped, pos=0.2, allow upside down]{\arrowOut} node[pos=.3, anchor=east]{$F$};
\node at (0,.5)[anchor=south]{$\Hom_\calm(F^l(m),F^l(m'))$};
\end{tikzpicture}
\end{align}
We use the natural isomorphisms
\begin{align}\label{eq:natiso}
&\Hom_\calm(m, F(m'))\cong \Hom_\calm(F^l(m), m')\cong \Hom_\calm(F(m'),m)^*\cong \Hom_\calm(m', F^l(m))^*\\
&\Hom_\calm(F(m),m')\cong \Hom_\calm(m, F^l(m'))\cong \Hom_\calm(m', F(m))^*\cong \Hom_\calm(F^l(m'), m)^*,
\end{align}
induced by the adjunctions and bimodule traces to identify different morphism spaces at the boundary.
These identifications can be used to choose morphism spaces for vertices on the boundary circle to contain only $F$ and not $F^l$, whenever there is only one chord endpoint at the vertex. Thus we can depict chord endpoints as meeting the boundary at right angles, as in \eqref{eq:boundlabelmor}.

\begin{Definition} \label{def:chord}$\quad$\\[-5ex]
\begin{enumerate}
    \item A \textbf{tricoloured  chord diagram} is a    chord diagram with oriented chords of colours green, blue and red such that crossings occur only between chords of different colours. 
\item A \textbf{labelling} of a tricoloured chord diagram is an assignment
\begin{compactitem}
\item of objects of $\mac$, $\mad$ and $\cale$ to the green, blue and red chords, respectively,
\item of objects of $\calm$ to the segments of the boundary circle,
\item of morphisms of $\calm$ to the chord endpoints as in \eqref{eq:boundlabelmor}.
    \end{compactitem}
    \end{enumerate}
\end{Definition}

The diagrams in \cite[Sec 2.3.3]{M} and hence the labelled
tricoloured chord diagrams considered in this article evaluate to a complex number. This a special case of Definitions 2.9 and 2.13 in   \cite{M}.

\begin{Definition}\label{def:tricev}
  The \textbf{evaluation} $\mathrm{ev}_l(C)$ of a labelled  tricoloured  chord diagram $C$ is the complex number obtained as follows:\\[-4ex]
  \begin{enumerate}

\item Assign to each crossing point of the chords the associated isomorphism from \eqref{pic:bimodule}, \eqref{pic:modulefunc} or \eqref{pic:modulefuncinv}.

\item Cut the  chord diagram at a boundary segment and straighten the boundary circle to a vertical line. 
Apply an isotopy to ensure that 
each horizontal slice of the resulting diagram intersects at most one chord endpoint, maximum, minimum or crossing point.

\item Evaluate the resulting diagram as in Example \ref{ex:evaluationmor}. This yields an endomorphism of the object $m\in \Ob\calm$ assigned to the top and bottom of the vertical line.

\item Take the bimodule trace of this endomorphism to obtain a complex number.  
\end{enumerate}
\end{Definition}

It is shown in \cite[Sec 2.2, 2.3]{M} that the evaluation of the chord diagram is independent of the isotopy  that moves  crossing points and chord endpoints to different slices of the diagram. This is a consequence of the naturality of the  isomorphisms at the crossing points and at  the maxima and minima.
It is also shown in \cite[Sec 2.3]{M} that all  diagrams obtained by cutting a  given  chord diagram have the same evaluation.  This is a consequence of the properties of a bimodule trace and shows that the evaluation of a labelled tricoloured chord diagram is well-defined.

\subsubsection{Labelled surface diagrams}
\label{subsec:surfdiag}

The evaluation of diagrams on a disc is used in \cite{M}  to define the evaluation of analogous diagrams on oriented surfaces with boundaries. We specialise this to the following diagrams. 

\begin{Definition}\label{def:surfdiag} A \textbf{surface diagram} is
 an oriented surface with at most one boundary component
and finitely many oriented curves of colours red, blue and green such that
\begin{compactenum}[(i)]
\item curves of the same colour do not intersect,
\item all intersections are transversal crossings in the interior of the surface,
\item curves are either closed and do not intersect 
the boundary or have both endpoints on the boundary, both at the same point or at two distinct points. 
\end{compactenum}
The \textbf{regions} of a surface diagram are the connected components of the complement of its curves. A region is called \textbf{outer}, if it contains a boundary point, and \textbf{inner}, if it does not.
\end{Definition}

Just as the chords of tricoloured chord diagrams, the green, blue and red curves of a surface diagram 
are labelled with objects of  $\mac$, $\mad$ and $\cale$. Boundary segments between curve endpoints are labelled with objects of $\calm$ and curve endpoints at the boundary with morphisms in $\calm$, as in Definition \ref{def:chord}. 

\begin{Definition}\label{def:labsurf}
 A \textbf{labelling} $l$ of a surface diagram $X$ is an assignment
 \begin{compactitem}
\item of an object  of $\mac$ to every green curve on $\calm$,
\item an object  of $\mad$ to every blue curve on $\calm$,
\item an object of $\cale$ to every red curve on $\calm$,
\item an object of $\calm$ to each  segment of $\partial X$, if present, 
\item a morphism in $\calm$ to each curve endpoint on $\partial X$, if present, as in Definition \ref{def:chord}.
 \end{compactitem}
\end{Definition}

To evaluate a labelled surface diagram $X$, we
choose a basepoint $b\in X$ that does not lie on any curve of $X$ and with $b\in \partial X$, if $X$ has a boundary circle. We choose
a set of curves based at $b\in X$ that represent a set of generators of the
fundamental group $\pi_1(X,b)$  such that
\begin{compactenum}[(i)]
\item the complement of these curves is a disc, 

\item green, blue or red curves of $X$ intersect these  curves transversally.
\end{compactenum}
Such a set of curves is called a set of \textbf{cutting generators} in the following.  

Cutting a labelled surface diagram $X$ along a set of cutting generators yields a  tricoloured  chord diagram in the sense of Definition \ref{def:chord}. This chord diagram has two types of segments on its boundary circle, \textbf{outer segments} and \textbf{inner segments}, 
for which the corresponding points on $X$ are contained in outer and inner regions, respectively. 
It also has two types of chord endpoints, \textbf{outer chord endpoints}
that arise from endpoints of curves on $\partial X$ and \textbf{inner chord endpoints} that arise from cutting the surface. Chord endpoints with two chord ends are always outer chord endpoints. 

The labelling of the surface diagram $X$ induces a labelling of the chords of the tricoloured chord diagram  with objects of $\mad$, $\mac$, $\cale$. It also induces a labelling of the outer segments of its boundary circle with objects of $\calm$ and of the outer chord endpoints with morphisms of $\calm$. 

To obtain a labelled tricoloured chord diagram in the sense of Definition \ref{def:chord} we need to extend this labelling to the inner segments and the inner chord endpoints.

We fix a set $I_\calm$ of representatives of the isomorphism classes of simple objects in $\calm$. For $m,n\in \mathrm{Ob}\calm$ we call bases $B=\{p^\alpha\mid \alpha=1,\ldots, k\}$ of $\Hom_\calm(m,n)$ and  $B^*=\{j^\alpha\mid \alpha=1,\ldots,k\}$ of $\Hom_\calm(n,m)\cong \Hom_\calm(m,n)^*$ \textbf{trace normalised dual bases} if 
\begin{align}\label{eq:tracenormalised}
\tr(p^\alpha\circ j^\beta)=\tr(j^\beta\circ p^\alpha)=\delta_{\alpha\beta}\qquad \qquad \forall \alpha,\beta\in\{1,\ldots,k\}.
\end{align}
We label the inner boundary segments of the tricoloured chord diagram by elements of $I_\calm$ and inner chord endpoints with elements of fixed trace normalised dual bases. 
The evaluation of the surface diagram is then obtained from the evaluation of this labelled chord diagram: One  rescales with the dimensions of the objects at the inner boundary segments, sums over the dual bases of morphism spaces for inner chord endpoints and sums over  labellings of inner segments with elements of $I_\calm$:

\bigskip
\begin{Definition}\label{def:surfdiagev}
 The \textbf{evaluation} $\mathrm{ev}_l(X)$ of a labelled surface diagram $X$  is obtained as follows.
\begin{enumerate}
    \item Cut $X$ along the set of cutting generators. If $X$ is a closed surface, this yields a polygon whose sides are identified pairwise. 
    If $X$ has a boundary circle, this yields a polygon whose sides are identified pairwise, up to an additional side that corresponds to the boundary circle. 
    
    This polygon defines a tricoloured  chord diagram on a disc $D$  with inner chord 
endpoints on $\partial D$  identified pairwise. 
The chords inherit a labelling  with elements of $\Ob\mac$, $\Ob\mad$ and $\Ob\cale$. 

\item If $X$ has a boundary circle, the outer boundary segments of $D$ inherit a labelling with objects of $\calm$ from $X$. 
Label the inner boundary segments of $D$ with elements of $I_\calm$, such that the labels coincide, whenever  any of the corresponding points in $\partial D$ are identified in $X$.

\item Assign morphism spaces in $\calm$ to the chord endpoints, as in Definition \ref{def:chord}. The outer chord endpoints inherit a labelling by morphisms in these spaces from $X$.
Label pairs of inner chord endpoints that are identified on $X$ by dual elements of  trace normalised dual bases, using the isomorphisms from \eqref{eq:natiso}.

\item
This yields a labelled tricoloured chord diagram $C$.   
Sum its evaluation $\mathrm{ev}_l(C)$  over the  trace normalised dual bases of the morphism spaces at the  inner chord endpoints. 

For each pair of sides of the polygon that are identified in $X$, i.~e.~correspond to a single cutting generator, but do not contain a chord endpoint, add a factor $\dim(m)^\inv$, where $m\in I_\calm$ is the label of the segment. This is summation over trace normalised dual bases of $\End_\calm(m)\cong \C$.

\item For each inner boundary segment of the chord diagram, multiply by the dimension $\dim(m)$ of its label $m\in I_\calm$ and sum over all assignments of elements of  $I_\calm$ to inner boundary segments.

\end{enumerate}
\end{Definition}

\bigskip
\begin{Example}\label{ex:evaluation}
Consider the  following surface diagram $X'$ on the  torus with a disc removed 
\begin{align*}
\begin{tikzpicture}[scale=.5]
\draw[line width=1pt] (0,0) ellipse (4 and 3);
\draw[line width=1pt] (-1.2,0).. controls (-1,-.7) and (1,-.7)..(1.2,0);
\draw[line width=1pt] (-.8,0).. controls (-.8,.5) and (.8,.5)..(.8,0);
\draw[fill=gray, line width=1pt, gray] (-2,0) ellipse (.5 and .3);
\node at (-2,.3)[anchor=south]{$m$};
\draw[line width=1pt, red] (0,0) ellipse (3 and 2) ;
\draw[line width=1pt, red, ->,>=stealth] (-.2,2)--(.2,2) node[pos=.5, anchor=south]{$F$};
\draw[line width=1pt, blue] (0, -.5).. controls (-.7,-1) and (-1,-2)..(0,-3)node[sloped, pos=0.5, allow upside down]{\arrowOut} node[pos=.8, anchor=east]{$b$};
\draw[line width=1pt, green] (.5, -.5).. controls (-.2,-1) and (-.5,-2)..(.5,-3) node[sloped, pos=0.5, allow upside down]{\arrowOut} node[pos=.8, anchor=west]{$c$};
\end{tikzpicture}
\end{align*}
labelled with $m\in \Ob\calm$, $b\in \Ob\mad$, $c\in \Ob\mac$ and a bimodule endofunctor $F:\calm\to\calm$.
Cutting it along a set of cutting generators and labelling the inner boundary segments of the resulting disc  with $m'\in I_\calm$ yields the following labelled chord diagram with chord endpoints identified pairwise in $X'$ 
\begin{align}\label{pic:circeval}
&\begin{tikzpicture}[scale=.8, baseline=(current bounding box.center)]
\draw[line width=1pt, color=black] (0,0) circle (3);
\draw[line width=1pt, red] (45:3)--(225:3) node[sloped, pos=0.5, allow upside down]{\arrowOut} node[pos=.5, anchor=north west]{$F$};
\node at (45:3)[anchor=south west]{$\Hom_\calm(m, F(m))$};
\node at (225:3)[anchor=north east]{$\Hom_\calm(F(m), m)$};
\draw[line width=1pt, green] (110:3)--(330:3) 
node[sloped, pos=0.3, allow upside down]{\arrowOut} node[pos=.3, anchor=south west]{$c$};
\node at (110:3)[anchor=south east]{$\Hom_\calm(m', c\rhd m)$};
\node at (330:3)[anchor=north west]{$\Hom_\calm(c\rhd m, m')$};
\draw[line width=1pt, color=blue] (155:3)--(290:3) node[sloped, pos=0.3, allow upside down]{\arrowOut} node[pos=.3, anchor=south west]{$b$};
\node at (155:3)[anchor=east]{$\Hom_\calm(m, m'\lhd b)$};
\node at (290:3)[anchor=north west]{$\Hom_\calm(m'\lhd b, m)$};
\node at (90:3) [color=black, anchor=south west]{$m$};
\node at (130:3) [color=black, anchor=south east]{$m'$};
\node at (190:3) [color=black, anchor=north east]{$m$};
\node at (260:3) [color=black, anchor=north east]{$m$};
\node at (310:3) [color=black, anchor=north west]{$m'$};
\node at (10:3) [color=black, anchor=south west]{$m$};
\end{tikzpicture}
\qquad
&\begin{tikzpicture}[scale=.8, baseline=(current bounding box.center)]
\draw[line width=1pt, color=black] (0,4)--(0,-4);
\draw[line width=1pt, color=black] (-1,4)--(1,4);
\draw[line width=1pt, color=black] (-1,-4)--(1,-4);
\draw[color=blue, line width=1pt] (0,3).. controls (-3,2) and (-3,-1)..(0,-2)
node[sloped, pos=0.5, allow upside down]{\arrowOut} node[pos=.5, anchor=east]{$b$};
\node at (0,3)[anchor=west, color=blue]{$\beta$};
\node at (0,-2)[anchor=west, color=blue]{$\beta$};
\draw[color=black, line width=1pt, green] (0,2).. controls (-2,1) and (-2,0)..(0,-1)
node[sloped, pos=0.7, allow upside down]{\arrowOut} node[pos=.5, anchor=west]{$c$};
\node at (0,2)[anchor=west, color=green]{$\gamma$};
\node at (0,-1)[anchor=west, color=green]{$\gamma$};
\draw[line width=1pt, red] (0,0).. controls (-2.5,-1) and (-2.5,-2)..(0,-3)
node[sloped, pos=0.7, allow upside down]{\arrowOut} node[pos=.7, anchor=north east]{$F$};
\node at (0,0)[anchor=west, color=red]{$\alpha$};
\node at (0,-3)[anchor=west, color=red]{$\alpha$};
\node at (0,3.5)[anchor=west, color=black]{$m$};
\node at (0,2.5)[anchor=west, color=black]{$m'$};
\node at (0,1)[anchor=west, color=black]{$m$};
\node at (0,-.5)[anchor=west, color=black]{$m$};
\node at (0,-1.5)[anchor=west, color=black]{$m'$};
\node at (0,-2.5)[anchor=west, color=black]{$m$};
\node at (0,-3.5)[anchor=west, color=black]{$m$};
\end{tikzpicture}
\end{align}
If chord endpoints at the boundary  are labelled with pairs of dual morphisms
\begin{align*}
&\beta\in \Hom_\calm(m, m'\lhd b) & &\gamma\in \Hom_\calm(m', c\rhd m) & &\alpha\in\Hom_\calm(m, F(m))\\
&\overline \beta\in \Hom_\calm(m'\lhd b,m) & &\overline\gamma\in \Hom_\calm(c\rhd m,m') & &\overline\alpha\in\Hom_\calm(F(m),m),
\end{align*}
with the bars of the duals  suppressed in \eqref{pic:circeval} for better legibility,
then the evaluation  of the chord diagram is the bimodule trace of the endomorphism
\begin{align*}
\phi_{\alpha,\beta,\gamma, m,m'}: 
m&\xrightarrow{\beta} m'\lhd b
\xrightarrow{\gamma\lhd b} (c\rhd m)\lhd b
\xrightarrow{(c\rhd \alpha)\lhd b} (c\rhd F(m))\lhd b
\xrightarrow{(s^{F}_{c, m})^{\inv}\lhd b} F(c\rhd m)\lhd b\\
&
\xrightarrow{F(\overline\gamma)\lhd b} F(m')\lhd b\xrightarrow{t^F_{m',b}} F(m'\lhd b)\xrightarrow{F(\overline\beta)} F(m)\xrightarrow{\overline\alpha} m.
\end{align*}
Summing over dual bases of the morphism spaces for the chord endpoints and over $m'\in I_\calm$ yields
\begin{align*}
\mathrm{ev}_l(X')&=\sum_{m'\in I_\calm}\dim(m')\, \sum_{\alpha,\beta,\gamma}  \mathrm{tr}_m (\phi_{\alpha,\beta,\gamma, m, m'})
\end{align*}
For the associated closed surface diagram $X$ on the torus we obtain
\begin{align*}
\mathrm{ev}_l(X)&=\sum_{m,m'\in I_\calm}\dim(m)\cdot  \dim(m')\, \sum_{\alpha,\beta,\gamma}  \mathrm{tr}_m (\phi_{\alpha,\beta,\gamma, m, m'}).
\end{align*}
\end{Example}

\bigskip
\begin{Example}\label{ex:stabilisationsurf_gen} Consider the  surface diagram $\Sigma'_{st}$ obtained from the trisection diagram $\Sigma_{st}$ in  \eqref{eq:stabilisation} by removing  a disc. The green, blue and red curves are oriented and labelled with objects $c_1,c_2,c_3\in \Ob\mac$,  $b_1,b_2,b_3\in \Ob\mad$ and  $F_1,F_2,F_3\in \Ob\cale$. The boundary  is labelled with a simple object $m\in I_\calm$. 

\begin{figure}
\centering

	\begin{tikzpicture}[baseline={([yshift=-.5ex]current bounding box.center)}]
	\node at (0,0) {\def\svgscale{.35} 
\begingroup%
  \makeatletter%
  \providecommand\color[2][]{%
    \errmessage{(Inkscape) Color is used for the text in Inkscape, but the package 'color.sty' is not loaded}%
    \renewcommand\color[2][]{}%
  }%
  \providecommand\transparent[1]{%
    \errmessage{(Inkscape) Transparency is used (non-zero) for the text in Inkscape, but the package 'transparent.sty' is not loaded}%
    \renewcommand\transparent[1]{}%
  }%
  \providecommand\rotatebox[2]{#2}%
  \newcommand*\fsize{\dimexpr\f@size pt\relax}%
  \newcommand*\lineheight[1]{\fontsize{\fsize}{#1\fsize}\selectfont}%
  \ifx\svgwidth\undefined%
    \setlength{\unitlength}{841.88976378bp}%
    \ifx\svgscale\undefined%
      \relax%
    \else%
      \setlength{\unitlength}{\unitlength * \real{\svgscale}}%
    \fi%
  \else%
    \setlength{\unitlength}{\svgwidth}%
  \fi%
  \global\let\svgwidth\undefined%
  \global\let\svgscale\undefined%
  \makeatother%
  \begin{picture}(1,0.70707071)%
    \lineheight{1}%
    \setlength\tabcolsep{0pt}%
    \put(0,0){\includegraphics[width=\unitlength,page=1]{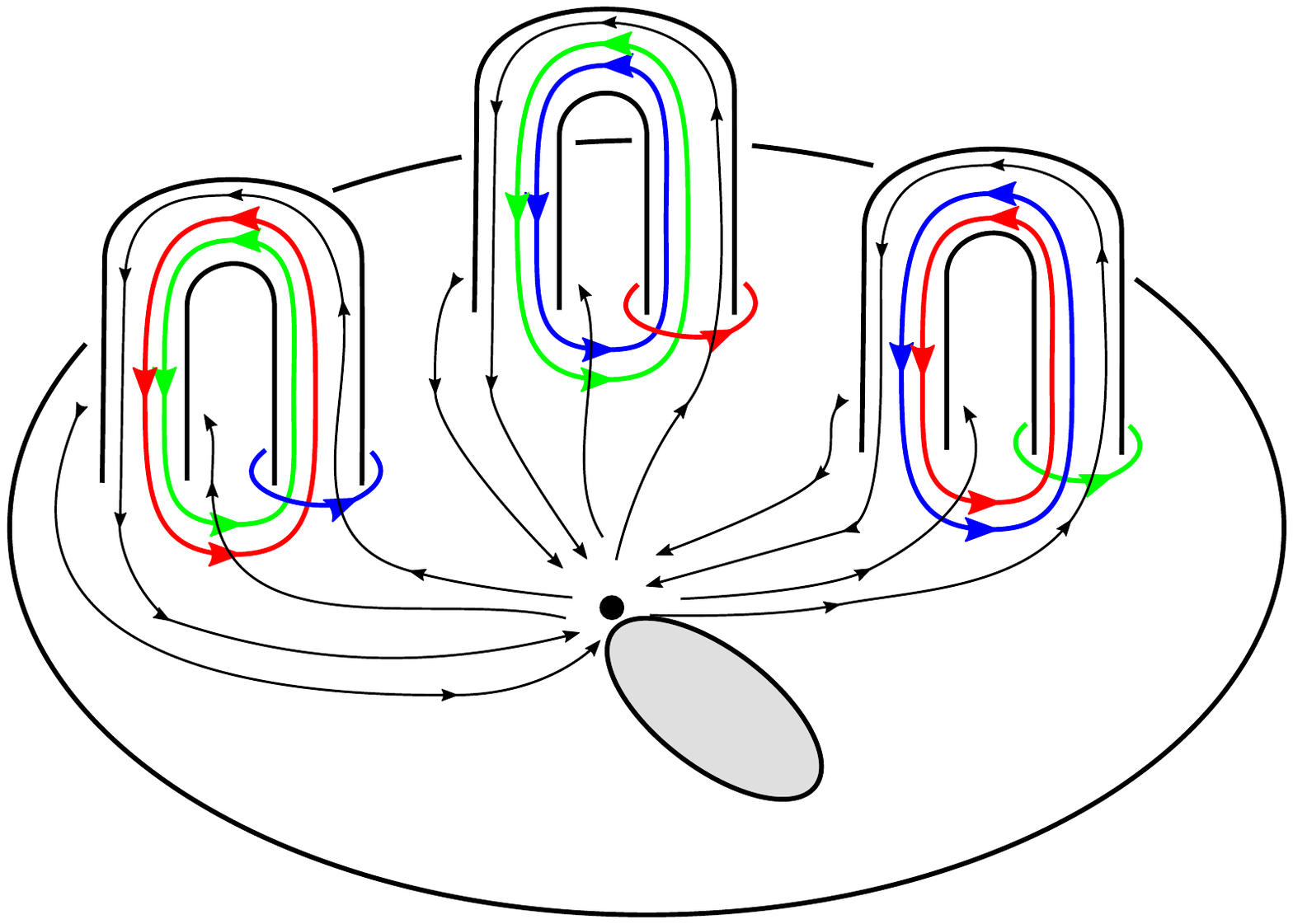}}%
    \put(0.19859403,0.56637622){\makebox(0,0)[lt]{\lineheight{1.25}\smash{\begin{tabular}[t]{l}3\end{tabular}}}}%
    \put(0.47471043,0.68324792){\makebox(0,0)[lt]{\lineheight{1.25}\smash{\begin{tabular}[t]{l}2\end{tabular}}}}%
    \put(0.74963328,0.5857326){\makebox(0,0)[lt]{\lineheight{1.25}\smash{\begin{tabular}[t]{l}1\end{tabular}}}}%
  \end{picture}%
\endgroup%
 };
	\end{tikzpicture}

\caption{A set of cutting generators for the trisection diagram $\Sigma'_{st}$ of $S^4$ with a disc removed}
\label{fig:stab_surfg}
\end{figure}

Cutting the surface along a set of cutting generators as in Figure \ref{fig:stab_surfg} and labelling the inner boundary segments with objects $m_1,m_2,m_3\in I_\calm$  yields the following chord diagram $D$, where we omit the morphism spaces for better legibility 
\begin{align}\label{eq:chordstablab}
\begin{tikzpicture}[scale=.7]
\draw[line width=1pt] (0,0) circle (5);
\draw[line width=1pt, red] (-330:5) ..controls (-340:3) and (-350:3) ..(0:5) node[sloped, pos=0.3, allow upside down]{\arrowOut} node[pos=.3, anchor=west]{$F_1$};
\draw[line width=1pt, blue] (-320:5) ..controls (-340:2) and (-350:2) ..(-10:5) node[sloped, pos=0.4, allow upside down]{\arrowOut} node[pos=.4, anchor=east]{$b_1$};
\draw[line width=1pt, green] (-345:5) ..controls (0:2) and (-20:2) ..(-35:5) node[sloped, pos=0.6, allow upside down]{\arrowOut} node[pos=.6, anchor=south east]{$c_1$};
\draw[line width=1pt, green] (-90:5) ..controls (-100:3) and (-110:3) ..(-120:5) node[sloped, pos=0.3, allow upside down]{\arrowOut} node[pos=.3, anchor=north east]{$c_2$};
\draw[line width=1pt, red] (-80:5) ..controls (-100:2) and (-110:2) ..(-130:5) node[sloped, pos=0.4, allow upside down]{\arrowOut} node[pos=.4, anchor=south west]{$F_2$};
\draw[line width=1pt, blue] (-105:5) ..controls (-120:2) and (-140:2) ..(-155:5) node[sloped, pos=0.6, allow upside down]{\arrowOut} node[pos=.6, anchor=south]{$b_2$};
\draw[line width=1pt, blue] (-200:5) ..controls (-220:3) and (-230:3) ..(-240:5) node[sloped, pos=0.3, allow upside down]{\arrowOut} node[pos=.3, anchor=south]{$b_3$};
\draw[line width=1pt, green] (-190:5) ..controls (-220:2) and (-230:2) ..(-250:5) node[sloped, pos=0.4, allow upside down]{\arrowOut} node[pos=.4, anchor=north west]{$c_3$};
\draw[line width=1pt, red] (-225:5) ..controls (-240:2) and (-260:2) ..(-275:5) node[sloped, pos=0.6, allow upside down]{\arrowOut} node[pos=.6, anchor=north west]{$F_3$};
\node at (70:5)[anchor=south west]{$m$};
\node at (35:5)[anchor=south west]{$m_1$};
\node at (20:5)[anchor=south west]{$m$};
\node at (5:5)[anchor=south west]{$m$};
\node at (-5:5)[anchor=west]{$m_1$};
\node at (-20:5)[anchor=north west]{$m$};
\node at (-60:5)[anchor=north west]{$m$};
\node at (90:5)[anchor=south east]{$m$};
\node at (110:5)[anchor=south east]{$m_3$};
\node at (125:5)[anchor=south east]{$m$};
\node at (145:5)[anchor=south east]{$m$};
\node at (165:5)[anchor=east]{$m_3$};
\node at (185:5)[anchor=north east]{$m$};
\node at (220:5)[anchor=north east]{$m$};
\node at (235:5)[anchor=north east]{$m_2$};
\node at (250:5)[anchor=north east]{$m$};
\node at (265:5)[anchor=north east]{$m$};
\node at (275:5)[anchor=north]{$m_2$};
\end{tikzpicture}
\end{align}
As $m\in I_\calm$ is  simple, the evaluation of $D$ is given by the following three bimodule traces
\begin{align}\label{eq:chordtri}
&\begin{tikzpicture}[scale=.8, baseline=(current bounding box.center)]
\draw[line width=1pt, color=black] (0,4)--(0,-4);
\draw[line width=1pt, color=black] (-1,4)--(1,4);
\draw[line width=1pt, color=black] (-1,-4)--(1,-4);
\draw[color=blue, line width=1pt] (0,3).. controls (-3,2) and (-3,-1)..(0,-2)
node[sloped, pos=0.5, allow upside down]{\arrowOut} node[pos=.5, anchor=east]{$b_1$};
\node at (0,3)[anchor=west, color=blue]{$\beta_1$};
\node at (0,-2)[anchor=west, color=blue]{$\beta_1$};
\draw[color=black, line width=1pt, red] (0,2).. controls (-2,1) and (-2,0)..(0,-1)
node[sloped, pos=0.7, allow upside down]{\arrowOut} node[pos=.5, anchor=west]{$F_1$};
\node at (0,2)[anchor=west, color=red]{$\alpha_1$};
\node at (0,-1)[anchor=west, color=red]{$\alpha_1$};
\draw[line width=1pt, green] (0,0).. controls (-2.5,-1) and (-2.5,-2)..(0,-3)
node[sloped, pos=0.7, allow upside down]{\arrowOut} node[pos=.7, anchor=north east]{$c_1$};
\node at (0,0)[anchor=west, color=green]{$\gamma_1$};
\node at (0,-3)[anchor=west, color=green]{$\gamma_1$};
\node at (0,3.5)[anchor=west, color=black]{$m$};
\node at (0,2.5)[anchor=west, color=black]{$m_1$};
\node at (0,1)[anchor=west, color=black]{$m$};
\node at (0,-.5)[anchor=west, color=black]{$m$};
\node at (0,-1.5)[anchor=west, color=black]{$m_1$};
\node at (0,-2.5)[anchor=west, color=black]{$m$};
\node at (0,-3.5)[anchor=west, color=black]{$m$};
\end{tikzpicture}
& 
&\begin{tikzpicture}[scale=.8, baseline=(current bounding box.center)]
\draw[line width=1pt, color=black] (0,4)--(0,-4);
\draw[line width=1pt, color=black] (-1,4)--(1,4);
\draw[line width=1pt, color=black] (-1,-4)--(1,-4);
\draw[color=red, line width=1pt] (0,3).. controls (-3,2) and (-3,-1)..(0,-2)
node[sloped, pos=0.5, allow upside down]{\arrowOut} node[pos=.5, anchor=east]{$F_2$};
\node at (0,3)[anchor=west, color=red]{$\alpha_2$};
\node at (0,-2)[anchor=west, color=red]{$\alpha_2$};
\draw[color=black, line width=1pt, green] (0,2).. controls (-2,1) and (-2,0)..(0,-1)
node[sloped, pos=0.7, allow upside down]{\arrowOut} node[pos=.5, anchor=west]{$c_2$};
\node at (0,2)[anchor=west, color=green]{$\gamma_2$};
\node at (0,-1)[anchor=west, color=green]{$\gamma_2$};
\draw[line width=1pt, blue] (0,0).. controls (-2.5,-1) and (-2.5,-2)..(0,-3)
node[sloped, pos=0.7, allow upside down]{\arrowOut} node[pos=.7, anchor=north east]{$b_2$};
\node at (0,0)[anchor=west, color=blue]{$\beta_2$};
\node at (0,-3)[anchor=west, color=blue]{$\beta_2$};
\node at (0,3.5)[anchor=west, color=black]{$m$};
\node at (0,2.5)[anchor=west, color=black]{$m_2$};
\node at (0,1)[anchor=west, color=black]{$m$};
\node at (0,-.5)[anchor=west, color=black]{$m$};
\node at (0,-1.5)[anchor=west, color=black]{$m_2$};
\node at (0,-2.5)[anchor=west, color=black]{$m$};
\node at (0,-3.5)[anchor=west, color=black]{$m$};
\end{tikzpicture}
& 
&\begin{tikzpicture}[scale=.8, baseline=(current bounding box.center)]
\draw[line width=1pt, color=black] (0,4)--(0,-4);
\draw[line width=1pt, color=black] (-1,4)--(1,4);
\draw[line width=1pt, color=black] (-1,-4)--(1,-4);
\draw[color=green, line width=1pt] (0,3).. controls (-3,2) and (-3,-1)..(0,-2)
node[sloped, pos=0.5, allow upside down]{\arrowOut} node[pos=.5, anchor=east]{$c_3$};
\node at (0,3)[anchor=west, color=green]{$\gamma_3$};
\node at (0,-2)[anchor=west, color=green]{$\gamma_3$};
\draw[line width=1pt, blue] (0,2).. controls (-2,1) and (-2,0)..(0,-1)
node[sloped, pos=0.7, allow upside down]{\arrowOut} node[pos=.5, anchor=west]{$b_3$};
\node at (0,2)[anchor=west, color=blue]{$\beta_3$};
\node at (0,-1)[anchor=west, color=blue]{$\beta_3$};
\draw[line width=1pt, red] (0,0).. controls (-2.5,-1) and (-2.5,-2)..(0,-3)
node[sloped, pos=0.7, allow upside down]{\arrowOut} node[pos=.7, anchor=north east]{$F_3$};
\node at (0,0)[anchor=west, color=red]{$\alpha_3$};
\node at (0,-3)[anchor=west, color=red]{$\alpha_3$};
\node at (0,3.5)[anchor=west, color=black]{$m$};
\node at (0,2.5)[anchor=west, color=black]{$m_3$};
\node at (0,1)[anchor=west, color=black]{$m$};
\node at (0,-.5)[anchor=west, color=black]{$m$};
\node at (0,-1.5)[anchor=west, color=black]{$m_3$};
\node at (0,-2.5)[anchor=west, color=black]{$m$};
\node at (0,-3.5)[anchor=west, color=black]{$m$};
\end{tikzpicture}\nonumber\\[+2ex]
&\mathrm{tr}_m(\phi_{\alpha_1,\beta_1,\gamma_1, m,m_1})
& &\mathrm{tr}_m(\phi_{\alpha_2,\beta_2,\gamma_2, m,m_2})
& &\mathrm{tr}_m(\phi_{\alpha_3,\beta_3,\gamma_3,m,m_3}),
\end{align}
where
\begin{align*}
&\phi_{\alpha_1,\beta_1,\gamma_1,m,m_1}:\\
&m \xrightarrow{\beta_1} m_1\lhd b_1
\xrightarrow{\alpha_1\lhd b_1} F_1(m)\lhd b_1 \xrightarrow{F_1(\gamma_1)\lhd b_1} F_1(c_1\rhd m)\lhd b_1
\xrightarrow{s_{c_1,m}^{F_1}\lhd b_1} (c_1\rhd F_1(m))\lhd b_1\\
&
\xrightarrow{(c_1\rhd \overline\alpha_1)\lhd b_1}
(c_1\rhd m_1)\lhd b_1
\xrightarrow{q_{c_1,m_1, b_1}} c_1\rhd (m_1\lhd b_1) \xrightarrow{c_1\rhd \overline \beta_1} c_1\rhd m\xrightarrow{\overline\gamma_1} m
\\[+2ex]
&\phi_{\alpha_2,\beta_2,\gamma_2,m,m_2}:\\
&m \xrightarrow{\alpha_2} F(m_2)
\xrightarrow{F(\gamma_2)} F_2(c_2\rhd m) 
\xrightarrow{F_2(c_2\rhd \beta_2)} F_2(c_2\rhd (m\lhd b_2))
\xrightarrow{F_2(q_{c_2,m,b_2}^\inv)} F_2((c_2\rhd m)\lhd b_2)\\
&
\xrightarrow{F_2(\overline\gamma_2\lhd b_2)}
F_2(m_2\lhd b_2)
\xrightarrow{(t_{m_2,b_2}^{F_2})^{\inv}} F_2(m_2)\lhd b_2 \xrightarrow{\overline\alpha_2\lhd b_2} m\lhd b_2\xrightarrow{\overline\beta_2} m
\\[+2ex]
&\phi_{\alpha_3,\beta_3,\gamma_3,m,m_3}:\\
&m \xrightarrow{\gamma_3} c_3\rhd m_3
\xrightarrow{c_3\rhd \beta_3} c_3\rhd(m\lhd b_3)
\xrightarrow{c_3\rhd(\alpha_3\lhd b_3)} 
c_3\rhd(F_3(m)\lhd b_3)
\xrightarrow{c_3\rhd t_{m,b_3}^{F_3}} 
c_3\rhd F_3(m\lhd b_3)\\
&
\xrightarrow{c_3\rhd F_3(\overline\beta_3)}
c_3\rhd F_3(m_3) 
\xrightarrow{(s_{c_3,m_3}^{F_3})^{\inv}} 
F_3(c_3\rhd m_3) \xrightarrow{F_3(\overline\gamma_3)}
F_3(m)\xrightarrow{\overline\alpha_3} m.
\end{align*}
The evaluation of the surface diagram $\Sigma'_{st}$ with boundary label $m\in I_\calm$ reads 
\begin{align*}
  &\mathrm{ev}_l(\Sigma'_{st})= \dim(m)^{-2} \!\!\!\!\!\!\!\!\!\!\sum_{m_1,m_2,m_3\in I_\calm}  \!\!\!\!\!\dim(m_1)\dim(m_2)\dim(m_3)\\
  &\sum_{\alpha_1,\alpha_2,\alpha_3}\sum_{\beta_1,\beta_2,\beta_3}\sum_{\gamma_1,\gamma_2,\gamma_3} 
  \mathrm{tr}_m(\phi_{\alpha_1,\beta_1,\gamma_1,m,m_1})
\,\mathrm{tr}_m(\phi_{\alpha_2,\beta_2,\gamma_2,m,m_2})
\,\mathrm{tr}_m(\phi_{\alpha_3,\beta_3,\gamma_3,m,m_3}).
\end{align*}
The evaluation of the  labelled trisection diagram $\Sigma_{st}$ from \eqref{eq:stabilisation}   is
\begin{align*}
  &\mathrm{ev}_l(\Sigma_{st})= \sum_{m\in I_\calm}\dim(m)^{-1} \!\!\!\!\!\!\!\!\!\!\sum_{m_1,m_2,m_3\in I_\calm} \!\!\!\!\! \dim(m_1)\dim(m_2)\dim(m_3)\\
  &\sum_{\alpha_1,\alpha_2,\alpha_3}\sum_{\beta_1,\beta_2,\beta_3}\sum_{\gamma_1,\gamma_2,\gamma_3} 
  \mathrm{tr}_m(\phi_{\alpha_1,\beta_1,\gamma_1,m,m_1})
\,\mathrm{tr}_m(\phi_{\alpha_2,\beta_2,\gamma_2,m,m_2})
\,\mathrm{tr}_m(\phi_{\alpha_3,\beta_3,\gamma_3,m,m_3}).
\end{align*}
\end{Example}

It is shown in \cite{M} and follows by a direct computation that the evaluation of a surface diagram $X$  does not depend on the choice of the dual bases of the morphism spaces at inner  chord endpoints. 

Moreover,  the evaluation of a surface diagram is
independent of the choice of cutting generators. This is a consequence of \cite[Lemma 5.15]{M}, which follows from the fact that the evaluations of the diagrams satisfy the cutting and gluing identities from \cite[Th 3.1]{M}.

We summarise the argument for the convenience of the reader:
The cutting and gluing identities allow one to compute the evaluation of a surface diagram by cutting it into pieces. Its evaluation is obtained by summing over dual bases of the morphism spaces at the chord ends on the cuts and over simple objects, whenever gluing moves a whole boundary segment into the interior of the diagram.

\begin{Theorem} \cite[Th 3.1]{M} \label{th:cutdiag}\\
Tricoloured chord diagrams satisfy the cutting and gluing identities
\begin{enumerate}
\item {\bf Gluing sides:} for all  $m,n\in I_\calm$ 
\begin{align}\label{pic:gluepoly}
&\begin{tikzpicture}[scale=.32]
\node at (-19,0) {$\sum_{\alpha}$};
\begin{scope}[shift={(-15,0)}]
\draw[draw=none, fill=gray, fill opacity=.2] (-90:5)--(-30:5)--(30:5)--(90:5);
\draw[line width=1pt, color=black] (0:4.35)--(30:5) node[anchor=south]{$\;m$};
\draw[line width=1pt, color=black] (30:5)--(60:4.35);
\draw[line width=1pt, color=black] (-30:5)--(-60:4.35);
\draw[line width=1pt, color=black] (0:4.35)--(-30:5) node[anchor=north]{$\;n$};
\draw[line width=1pt, red] (.2,0) node[anchor=east]{$F$}--(0:4.35) node[sloped, pos=0.5, allow upside down]{\arrowIn};
\draw[color=black, fill=black] (0:4.35) circle (.2) node[anchor=west]{$\alpha$};
\end{scope}
\begin{scope}[shift={(-3,0)}]
\draw[draw=none, fill=gray, fill opacity=.2] (90:5)--(150:5)--(210:5)--(270:5);
\draw[line width=1pt, color=black] (180:4.35)--(150:5) node[anchor=south]{$m\;$};
\draw[line width=1pt, color=black] (150:5)--(120:4.35);
\draw[line width=1pt, color=black] (210:5)--(240:4.35);
\draw[line width=1pt, color=black] (180:4.35)--(210:5) node[anchor=north]{$n\;$};
\draw[line width=1pt, red] (180:4.35)--(0,0) node[anchor=west]{$F$} node[sloped, pos=0.5, allow upside down]{\arrowIn};
\draw[color=black, fill=black] (180:4.35) circle (.2) node[anchor=east]{$\alpha$};
\end{scope}
\node at (2,0){$=$};
\begin{scope}[shift={(7,0)}]
\draw[draw=none, fill=gray, fill opacity=.2] (-90:5)--(-30:5)--(30:5)--(90:5);
\begin{scope}[shift={(8.6,0)}]
\draw[draw=none, fill=gray, fill opacity=.2] ((90:5)--(150:5)--(210:5)--(270:5);
\draw[line width=1pt] (150: 5)--(120:4.35);
\draw[line width=1pt] (210:5)--(240:4.35);
\end{scope}
\draw[line width=1pt] (30:5) --(60:4.35);
\node at (4.2,3)[anchor=south]{$m$};
\node at (4.2,-3)[anchor=north]{$n$};
\draw[line width=1pt] (-30:5)--(-60:4.35);
\draw[line width=1pt, red] (.2,0) node[anchor=east]{$F$}--(8.8,0) node[sloped, pos=0.5, allow upside down]{\arrowIn};
\end{scope}
\end{tikzpicture}
\end{align}

\item {\bf Gluing around a corner:} for all  $m\in I_\calm$ 
\begin{align}
\label{pic:closepoly}
&\begin{tikzpicture}[scale=.32]
\node at (-22,0) {$\sum_{n\in I_\calm}\sum_{\alpha} \dim(n)$};
\begin{scope}[shift={(-7,0)}]
\draw[draw=none, fill=gray, fill opacity=.2] (0,0)--(30:5)--(90:5)--(150:10)--(210:10)--(270:5)--(330:5)--(0,0);
\draw[line width=1pt] (30:2.5)--(0,0) node[anchor=west]{$\;n$}--(-30:2.5);
\draw[line width=1pt, ] (30:2.5)--(30:5) node[anchor=west]{$m$}--(60:4.35);
\draw[line width=1pt, ] (-30:2.5)--(-30:5)node[anchor=west]{$m$} --(-60:4.35);
\draw[line width=1pt, red] (0,2.5) node[anchor=east]{$F$}--(30:2.5) node[sloped, pos=0.5, allow upside down]{\arrowIn};
\draw[line width=1pt, red] (-30:2.5)--(0,-2.5) node[anchor=east]{$F$} node[sloped, pos=0.5, allow upside down]{\arrowIn};
\draw[fill=black] (30:2.5) circle (.2) node[anchor=west]{$\;\alpha$};
\draw[fill=black] (-30:2.5) circle (.2) node[anchor=west]{$\;\alpha$};
\end{scope}
\node at (0,0){$=$};
\begin{scope}[shift={(9,0)}]
\draw[draw=none, fill=gray, fill opacity=.2] (0:5)--(72:5)--(144:8)--(216:8)--(288:5)--(0:5);
\draw[line width=1pt, ] (-36:4)--(0:5) node[anchor=west]{$m$}--(36:4);
\draw[line width=1pt, red] (50:2)  .. controls (0:2)..(-50:2) node[anchor=north east]{$F$} node[sloped, pos=1, allow upside down]{\arrowIn};
\end{scope}
\end{tikzpicture}
\end{align}
\end{enumerate}
\end{Theorem}
Here, concatenation of diagrams, as in the left-hand side of \eqref{pic:gluepoly}, stands for the product of their evaluations. All data that is not drawn is assumed to coincide on the  left-hand side and right-hand side of the equation. 
The summations indexed by $\alpha$ are over trace normalised dual bases of  the morphism spaces.  The functor $F$ in \eqref{pic:gluepoly} and \eqref{pic:closepoly} may  be a functor $F\in \Ob\cale$, a functor of the form  $F=c\rhd -$ or $F=-\lhd b$ with $c\in \Ob\mac$ or $b\in \Ob\mad$, or a composite of such functors. 

That any two sets of cutting generators yield the same evaluation may be seen as follows. Given two sets of cutting generators, one may apply an isotopy to ensure that (i) the curves in the two sets have only transversal intersection points and that (ii) the basepoint of one set does not lie on a curve or the basepoint of the other. One may then cut the underlying surface along both sets of cutting generators at once.  Gluing the pieces via identities \eqref{pic:gluepoly} and \eqref{pic:closepoly} along one set of cutting generators
then yields the diagram associated with the other. Identities \eqref{pic:gluepoly} and \eqref{pic:closepoly} then ensure that the evaluations of the two diagrams obtained this way agree.

As a direct consequence of Definition \ref{def:surfdiagev} and Theorem \ref{th:cutdiag} we also obtain an identity that allows us to compute the evaluation of a connected sum of surface diagrams from the evaluation of its summands. The \textbf{connected sum} of closed  labelled surface  diagrams  is
defined analogously to the connected sum of trisection diagrams in Sections \ref{sec:trisection} and \ref{sec:cccinv}. It 
is  obtained by
removing a  disc   that does not intersect any curves from each diagram, gluing the resulting surface diagrams along their boundary circles and equipping them with the induced labelling.

\begin{Lemma}\label{lem:gluing} $\quad$
\begin{compactenum}
\item If a surface diagram $Y'$ with labelling $l_{Y'}$ is obtained from a closed  surface diagram $Y$ with labelling $l_Y$  by removing a disc that does not intersect any curves and adding a boundary label $m\in I_\calm$, then
\begin{align}\label{eq:discsurf}\ev_{l_Y}(Y)=\sum_{m\in I_\calm} \dim(m)\cdot \mathrm{ev}_{l_{Y'}}(Y').\end{align}

\item Let  $(X, l_X)$ be  the connected sum of closed labelled surface diagrams  $(Y,l_Y)$ and $(Z,l_Z)$.
Then
\begin{align}
\ev_{l_X}(X)=\sum_{m\in I_\calm} \ev_{l_{Y'}}(Y')\cdot \ev_{l_{Z'}}(Z'),
\end{align}
where $Y'$ and $Z'$ are the associated surface diagrams with discs removed and their labellings  $l_{Y'}$, $l_{Z'}$  are obtained from $l_Y$, $l_Z$ by adding the same boundary label $m\in I_\calm$. 
\end{compactenum}
\end{Lemma}
\begin{Proof}
The first claim is a direct consequence 
of Definition \ref{def:surfdiagev},  5. 
For the second claim, 
note that the surface diagram $X$ is obtained by gluing $Y'$ and $Z'$ along their boundaries.
Choose basepoints on $Y'$ and $Z'$ that are glued to the same point in $X$ and sets 
 of cutting generators for $Y'$ and $Z'$ with respect to  these basepoints.  These cutting generators define a set of cutting generators for $X$.

Cutting the tricoloured chord diagram for $X$ along a generator that goes once around this gluing boundary and labelling the cut with $m\in I_\calm$ yields the tricoloured chord diagrams for $Y'$ and $Z'$ with the boundary components labelled by $m\in I_\calm$.
 Definition \ref{def:surfdiagev} and Equations
\eqref{pic:gluepoly} and \eqref{pic:closepoly}  imply
\begin{align*}
\ev_{l_X}(X)=\sum_{m\in I_\calm} \sum_\alpha \dim(m)\cdot \ev_{l_{Y'}}(Y')\cdot \ev_{l_{Z'}}(Z'),
\end{align*}
where the sum over $\alpha$ is over trace normalised dual bases of the morphism space $\End_\calm(m)\cong \C$ associated with the gluing boundary. By definition of the trace normalised dual bases, this yields a factor $\dim(m)^\inv$ and gives the result.
\end{Proof}

\subsection{Categorical trisection invariants}
\label{sec:triseccchord}

We now extend the evaluation of labelled surface diagrams to surface diagrams with labels assigned only to chords that end at the boundary, to boundary segments and to chord endpoints at the boundary. 

For this, we require again that $\mac$ and $\mad$ are spherical fusion categories and $\calm$ is an indecomposable $(\mac,\mad)$-bimodule category with a bimodule trace. We fix sets $I_\mac$, $I_\mad$ and $I_\calm$ of representatives of the isomorphism classes of simple objects in $\mac$, $\mad$ and $\calm$ and set $\cale=\End_{(\mac,\mad)}(\calm)$.

We now generalise the evaluation of labelled surface diagrams by considering another spherical fusion category $\cala$ and  a  pivotal functor $\Phi:\cala\to\cale$, cf.~Section \ref{subsec:catbasic}. This mirrors the construction in \cite{BB}.
We  fix a set of representatives $I_\cala$ of the isomorphism classes of simple objects in $\cala$. Note that the condition that $\calm$ is indecomposable does not allow one to modify the spherical fusion categories $\mac$ and $\mad$ freely in a similar way. 

To evaluate the surface diagram, we label its  blue and green curves by elements of  $I_\mad$ and $I_\mac$ and its red curves by images of elements of $I_\cala$ under $\Phi$. We evaluate the resulting labelled  surface diagram, rescale with the dimensions of these simple objects  and sum the evaluation over $I_\cala$, $I_\calb$ and $I_\mac$.

To express this in a formula, we  denote for a surface diagram
 $X$   by $G$, $B$, $R$  the sets of green, blue and red curves of $X$ and by $G_i$, $B_i$, $R_i$ the sets of closed green, blue and red curves that do not intersect $\partial X$. For a  labelling $l$ we denote by  $c_\gamma\in\Ob\mac$, $b_\beta\in \Ob\mad$ and $F_\rho=\Phi(a_\rho)\in \Ob\cale$ the   objects associated to curves $\gamma\in G$, $\beta\in B$ and $\rho\in R$.

\begin{Definition}\label{def:diagevalav} Let $\cala$, $\mad$, $\mac$ be spherical fusion categories, $\calm$ an indecomposable $(\mac,\mad)$-bimodule category with a bimodule trace and $\Phi:\cala\to \End_{(\mac,\mad)}(\calm)$ a pivotal functor. 

  Let $X$ be a surface diagram  with a labelling $l_\partial$ that assigns 
  \begin{compactitem}
\item objects in $\calm$ to  boundary segments of $X$,
\item objects of $\mac$, $\mad$ and $\cala$  to green, blue and red curves with endpoints on $\partial X$, 
\item morphisms in $\calm$ to curve endpoints on $\partial X$.
\end{compactitem}
 The \textbf{averaged evaluation} of $X$ is  obtained by:
\begin{compactenum}
\item Extending the boundary labelling $l_\partial$ to a labelling $l$ of $X$ 
that assigns elements of $I_\mac$ and $I_\mad$ to closed green and blue curves and images $\Phi(a)$ of elements $a\in I_\cala$ to closed red curves. \\[-2ex]
\item Multiplying  
$\mathrm{ev}_l(X)$ 
with the dimensions of the  objects assigned to the closed curves of $X$.\\[-2ex]
\item Summing over all the labellings $l$ in 1.~that extend the boundary labelling $l_\partial$.
\end{compactenum}
\begin{align*}
\mathrm{av}_{l_\partial}(X)=\sum_{l\vert_{\partial \Sigma}=l_\partial} \prod_{\gamma\in G_i}\prod_{\beta\in B_i}\prod_{\rho\in R_i} \dim(c_{\gamma})\dim(b_\beta)\dim(F_{\rho})\, \mathrm{ev}_l(X).
\end{align*}
\end{Definition}

If $X$ has  a boundary, its averaged evaluation  is a complex number that depends on the labels of the boundary segments, the labels of the chord endpoints at the boundary and the labels of the chords that extend to the boundary. In particular, this holds for the surface $\Sigma'_{st}$ from  Example \ref{ex:stabilisationsurf_gen}, which has a single object of $\calm$ as a boundary label.
We call the $(\mac,\mad)$-bimodule category $\calm$ \emph{stabilising}, if the averaged evaluation of  $\Sigma'_{st}$ with a simple boundary label  is given by a nonzero constant and the dimension of the boundary label.

\begin{Definition}\label{def:stabilisingcond} 
Let $\cala$, $\mad$, $\mac$ be spherical fusion categories, $\calm$ an indecomposable $(\mac,\mad)$-bimodule category with a bimodule trace and $\Phi:\cala \to \End_{(\mac,\mad)}(\calm)$ a pivotal functor. Then $\calm$ is called  \textbf{stabilising} with respect to $\Phi$, if there is a constant $C_{st}\in\C^\times$ with
$$\mathrm{av}_m(\Sigma'_{st})=C_{st}\cdot \dim(m)\qquad 
\forall m\in I_\calm.$$
\end{Definition}

We will look at the stabilisation condition in more depth at the end of this section. In particular, we will show  in Proposition \ref{prop:stabilising} that a $(\mac,\mad)$-bimodule category $\calm$ is stabilising with respect to $\Phi:\cala\to \End_{(\mac,\mad)}(\calm)$ if and only if $\mathrm{av}_m(\Sigma'_{st})\neq 0$ for all $m\in I_\calm$.  
This identifies the property \emph{stabilising} as the counterpart of the  condition in Definition \ref{def:cccinv} of the trisection invariant.

\begin{Example} For finite groups $B,C$  any indecomposable $(\Vec_C,\Vec_B)$-bimodule category $\calm$ is stabilising with respect to the identity functor on $\End_{(\mac,\mad)}(\calm)$. This is shown in Section \ref{subsec:dwcomputation}.
\end{Example}

Given spherical fusion categories $\cala$, $\mad$, $\mac$, an indecomposable  $(\mac,\mad)$-bimodule category $\calm$ and a pivotal functor  $\Phi:\cala\to  \End_{(\mac,\mad)}(\calm)$ such that $\calm$ is stabilising with respect to $\Phi$, we can use the averaged evaluation of surface diagrams to construct a trisection invariant. For this we interpret a trisection diagram as a surface diagram and rescale its averaged evaluation to ensure invariance under the stabilisation and destabilisation move.

\begin{Theorem}\label{th:gentrisection}
Let $\cala, \mad,\mac$ be spherical fusion categories. Suppose that  $\calm$ is a  $(\mac,\mad)$-bimodule category with a bimodule trace,  stabilising with respect to a fixed pivotal functor $\Phi:\cala\to \End_{(\mac,\mad)}(\calm)$, and let $\xi\in\C^\times$ with $\xi^3=C_{st}$.
Then the averaged evaluation of a trisection diagram $T$ on $\Sigma$, whose green, blue and red curves are labelled with $\mac$, $\mad$ and $\cala$ defines a 4-manifold invariant 
\begin{align}\label{eq:trisecresc}
|T|=\xi^{-g} \cdot \mathrm{av}(T).
\end{align}
\end{Theorem}

\pagebreak
\begin{Proof} We need to check invariance under the moves from Theorem \ref{th:intmoves}.
Diffeomorphism invariance is clear:  Applying the same diffeomorphism to the curves of $\Sigma$ and to a cutting system yields the same chord diagram.
Invariance under the two-point and three-point move follows from   \eqref{pic:rm2} and \eqref{pic:hexagon}.

Invariance under handle slides follows from the diagrammatic computation in Figure \ref{fig:handleslide}, which uses the naturality of the crossings with respect to morphisms in $\mac$, $\mad$ and $\cale=\End_{(\mac,\mad)}(\calm)$, the pivotality of $\Phi$ and that the following diagrammatic identities that hold in any spherical fusion category $\mac$.
\begin{align}\label{eq:sphids}
\begin{tikzpicture}[scale=.3]
\node at (-2,0)[anchor=east] {$\sum_{k}\sum_\alpha$ \dim(k)};
\draw[line width=1pt](-2,4) node[anchor=south]{$i$}--(0,2);
\draw[line width=1pt](2,4) node[anchor=south]{$j$}--(0,2);
\draw[line width=1pt](0,2)--(0,-2) node[pos=.5, anchor=west]{$k$};
\draw[line width=1pt](-2,-4) node[anchor=north]{$i$}--(0,-2);
\draw[line width=1pt](2,-4) node[anchor=north]{$j$}--(0,-2);
\draw[fill=black] (0,-2) circle (.2) node[anchor=west]{$\;\alpha$};
\draw[fill=black] (0,2) circle (.2) node[anchor=west]{$\;\alpha$};
\node at (2,0)[anchor=west]{$=$};
\draw[line width=1pt](4,4) node[anchor=south]{$i$}--(4,-4) node[anchor=north]{$i$};
\draw[line width=1pt](6,4) node[anchor=south]{$j$}--(6,-4) node[anchor=north]{$j$};
\end{tikzpicture}
\qquad\qquad\qquad
\begin{tikzpicture}[scale=.3]
\node at (-3,0)[anchor=east]{$\sum_{j}\sum_\alpha \dim(j)$};
\draw[line width=1pt] (-2,4) node[anchor=south]{$k$}--(-2,2).. controls (-2,1) and (0,1).. (0,2);
\draw[line width=1pt](-1,4) node[anchor=south]{$i$}--(0,2);
\draw[line width=1pt](0,2)--(0,3).. controls (0,4) and (2,4)..(2,3)--(2,-3).. controls (2,-4) and (0,-4)..(0,-2);
\draw[line width=1pt](-1,-4) node[anchor=north]{$i$}--(0,-2);
\draw[line width=1pt] (-2,-4) node[anchor=north]{$k$}--(-2,-2).. controls (-2,-1) and (0,-1).. (0,-2);
\draw[fill=black] (0,2) circle (.2) node[anchor=west]{$\;\alpha$};
\draw[fill=black] (0,-2) circle (.2) node[anchor=west]{$\;\alpha$};
\node at (.5,4)[anchor=south]{$j$};
\node at (.5,-4)[anchor=north]{$j$};
\node at (3,0)[anchor=west]{$=$};
\draw[line width=1pt] (5,4) node[anchor=south]{$k$}--(5,-4)node[anchor=north]{$k$};
\draw[line width=1pt] (7,4) node[anchor=south]{$i$}--(7,-4)node[anchor=north]{$i$};
\end{tikzpicture}
\end{align}
Here, the sums run over simple objects in a spherical fusion category and trace normalised dual bases $\{p^\alpha_{ijk}\}$ and $\{j^\alpha_{ijk}\}$ of the morphism spaces $\Hom_\mac(i\oo j,k)$ and $\Hom_\mac(k,i\oo j)$  with
\begin{align*}
p^\beta_{ijl}\circ j^\alpha_{ijk}=\frac{\delta_{\alpha\beta}\delta_{kl}}{\dim(k)} 1_k\qquad \sum_k\sum_\alpha \dim(k)\; j^\alpha_{ijk}\circ p_{ijk}^\alpha=1_{i\oo j}.
\end{align*}

\begin{figure}

	\begin{tikzpicture}[baseline={([yshift=-.5ex]current bounding box.center)}]
	\node at (0,0) {\def\svgscale{.75} 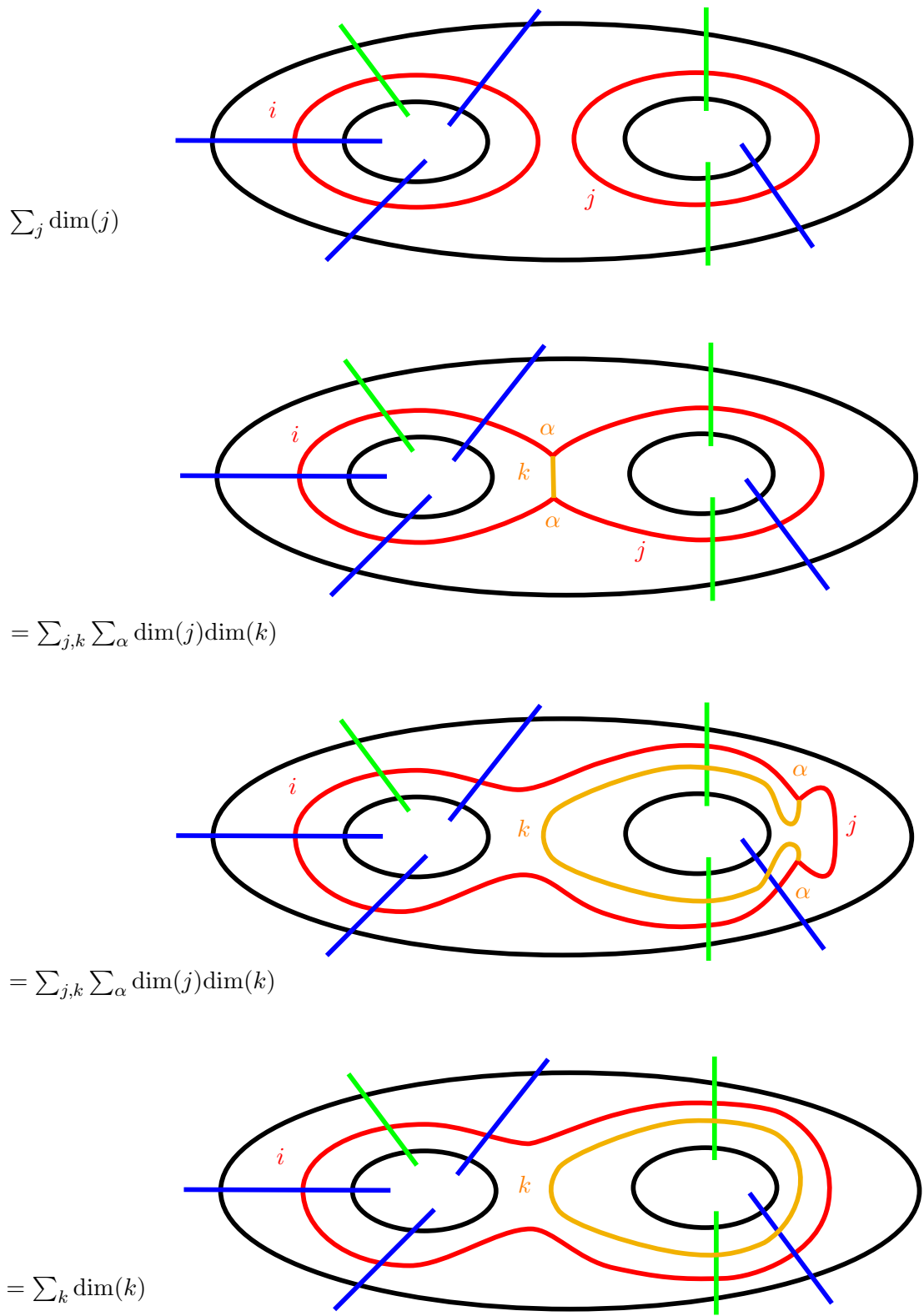 };
	\end{tikzpicture}

\caption{Diagrammatic computation: invariance under a handle slide. The labels for the blue and green curves are omitted. To pass from the first to the second diagram and from the third to the fourth we apply identities \eqref{eq:sphids}. The third diagram is obtained from the second via an isotopy. }
\label{fig:handleslide}
\end{figure}

To show invariance under the stabilisation move we consider a trisection diagram $X$  of genus $g$. We denote by $X'$ the associated diagram  with a disc removed, by $\Sigma_{st}$ the trisection diagram of $S^4$  in  \eqref{eq:stabilisation} and by  $\Sigma'_{st}$ the associated surface diagram with a disc removed from   Example \ref{ex:stabilisationsurf_gen}. 
By Lemma \ref{lem:gluing}
the evaluation of  $Y=X\#\Sigma_{st}$ with the labelling $l$  induced by  labellings $l_X$, $l_{st}$ of $X$,  $\Sigma_{st}$ is given by
\begin{align*}
\mathrm{ev}_l(Y)=\sum_{m\in I_\calm} \mathrm{ev}_{l'}(X')\cdot \mathrm{ev}_{l''}(\Sigma'_{st}),
\end{align*}
where $l'$ and $l''$ are the labellings of $X'$ and $\Sigma'_{st}$ obtained from $l_X$ and $l_{st}$ by labelling the boundary circle with $m\in I_\calm$.
As the sets of red, blue and green curves of $Y=X\#\Sigma_{st}$ are  the disjoint unions of the corresponding sets for $X'$ and $\Sigma'_{st}$, we obtain with Lemma \ref{lem:gluing} and
Definitions \ref{def:diagevalav} and \ref{def:stabilisingcond}
\begin{align*}
\mathrm{av}(Y)&=\sum_{m\in I_\calm}\sum_{l'\vert_\partial=m}\sum_{l''\vert_\partial=m} \prod_{\gamma\in G(Y)}\prod_{\beta\in B(Y)}\prod_{\rho\in R(Y)} \dim(c_\gamma)\dim(b_\beta)\dim(F_\rho) \,\mathrm{ev}_{l'}(X')\cdot \mathrm{ev}_{l''}(\Sigma'_{st})\\
&=\sum_{m\in I_\calm} \mathrm{av}_m(X')\cdot \mathrm{av}_m(\Sigma'_{st})=C_{st}\sum_{m\in I_\calm}  \dim(m)\cdot  \mathrm{av}_m(X')
=C_{st} \cdot \mathrm{av}(X).
\end{align*}
As $Y=X\#\Sigma_{st}$ is of genus $g+3$, this implies
\begin{align*}
|Y|=\xi^{-g-3}\cdot \mathrm{av}(Y)=\xi^{-g-3}\cdot C_{st}\cdot \mathrm{av}(X)=\xi^{-g}\cdot \mathrm{av}(X)=|X|.
\end{align*}
\end{Proof}

\begin{Remark} \label{rem:trisecinvrem}  The condition that $\calm$ is stabilising is required  only for invariance under the (de)stabilisation move. Invariance under diffeomorphisms of the surface, the two-point and three-point move and under handle slides holds for any indecomposable $(\mac,\mad)$-bimodule category with a bimodule trace.
This is analogous to the proof of Theorem \cite[Th 3.7]{CCC}.
\end{Remark}

We now look more closely at the stabilising condition on the datum $(\cala,\mad,\mac,\calm,\Phi)$ in Definition \ref{def:stabilisingcond}.
It turns out that it is always satisfied, as long as the averaged evaluation of the surface diagram $\Sigma'_{st}$ in Example \ref{ex:stabilisationsurf_gen}
 does not vanish for simple boundary labels. This is the counterpart of the corresponding condition in Definition \ref{def:cccinv} and Theorem \ref{th:ccc}.  To prove this, we require the following two lemmas.

\begin{Lemma}
\label{lem:stablemma} Let $(\cala,\mad,\mac,\calm,\Phi)$  be as in Definition \ref{def:stabilisingcond}, but not necessarily stabilising.
Let $X$ be obtained from a closed surface  diagram by removing a disc that intersects no curves of the diagram and then  placing a green and a blue curve around the boundary circle.
If $\mathrm{av}_{m}(X) \neq 0$ for all boundary labels $m\in I_\calm$, one has for all $m,m' \in I_{\calm}$
\begin{equation}
\label{eq:stablema}
    \frac{\mathrm{av}_{m}(X)}{\dim(m)} = \frac{\mathrm{av}_{m'}(X)}{\dim(m')}.
\end{equation}
\end{Lemma}

\begin{Proof}
The claim follows from the fact that any two $(\mac,\mad)$-bimodule traces $\mathrm{tr}$, $\mathrm{tr}'$ on an indecomposable $(\mac,\mad)$-bimodule category $\calm$ satisfy $\mathrm{tr}=z\cdot \mathrm{tr}'$ for some $z\in\C^\times$, see Proposition \ref{prop:uniquetrace}.

We consider a surface diagram $X$ as in the lemma, place  $n$  vertices on its boundary circle and view them as chord endpoints of chords labelled by the identity functor $\id_\calm:\calm\to\calm$. labelling the $n$ boundary segments  by objects  $m_1,\ldots, m_n$ in $\calm$ and the vertices with morphisms $\alpha_i: m_{i}\to m_{i+1}$ with $m_{n+1}=m_1$ yields a diagram with a boundary labelling $l_\partial$ as in Definition \ref{def:diagevalav}.

The morphisms at the boundary compose to an endomorphism $\alpha_n\circ\ldots\circ \alpha_1: m_1\to m_1$, and we define
\begin{align*}
\Theta(\alpha_n\circ \ldots\circ \alpha_1)=\mathrm{av}_{l_\partial}(X).
\end{align*}
In particular, we can consider the case $n=1$ with a single boundary vertex labelled by an endomorphism $\alpha:m\to m$ in $\calm$. This defines a family of maps $\Theta_m : \End_\calm(m) \to \C$ for each $m\in\Ob\calm$.
\begin{equation}
\Theta_m(\alpha) := \mathrm{av}_m\left(

	\begin{tikzpicture}[baseline={([yshift=-.5ex]current bounding box.center)}]
	\node at (0,0) {\def\svgscale{0.2} 
\begingroup%
  \makeatletter%
  \providecommand\color[2][]{%
    \errmessage{(Inkscape) Color is used for the text in Inkscape, but the package 'color.sty' is not loaded}%
    \renewcommand\color[2][]{}%
  }%
  \providecommand\transparent[1]{%
    \errmessage{(Inkscape) Transparency is used (non-zero) for the text in Inkscape, but the package 'transparent.sty' is not loaded}%
    \renewcommand\transparent[1]{}%
  }%
  \providecommand\rotatebox[2]{#2}%
  \newcommand*\fsize{\dimexpr\f@size pt\relax}%
  \newcommand*\lineheight[1]{\fontsize{\fsize}{#1\fsize}\selectfont}%
  \ifx\svgwidth\undefined%
    \setlength{\unitlength}{841.88976378bp}%
    \ifx\svgscale\undefined%
      \relax%
    \else%
      \setlength{\unitlength}{\unitlength * \real{\svgscale}}%
    \fi%
  \else%
    \setlength{\unitlength}{\svgwidth}%
  \fi%
  \global\let\svgwidth\undefined%
  \global\let\svgscale\undefined%
  \makeatother%
  \begin{picture}(1,0.70707071)%
    \lineheight{1}%
    \setlength\tabcolsep{0pt}%
    \put(0,0){\includegraphics[width=\unitlength,page=1]{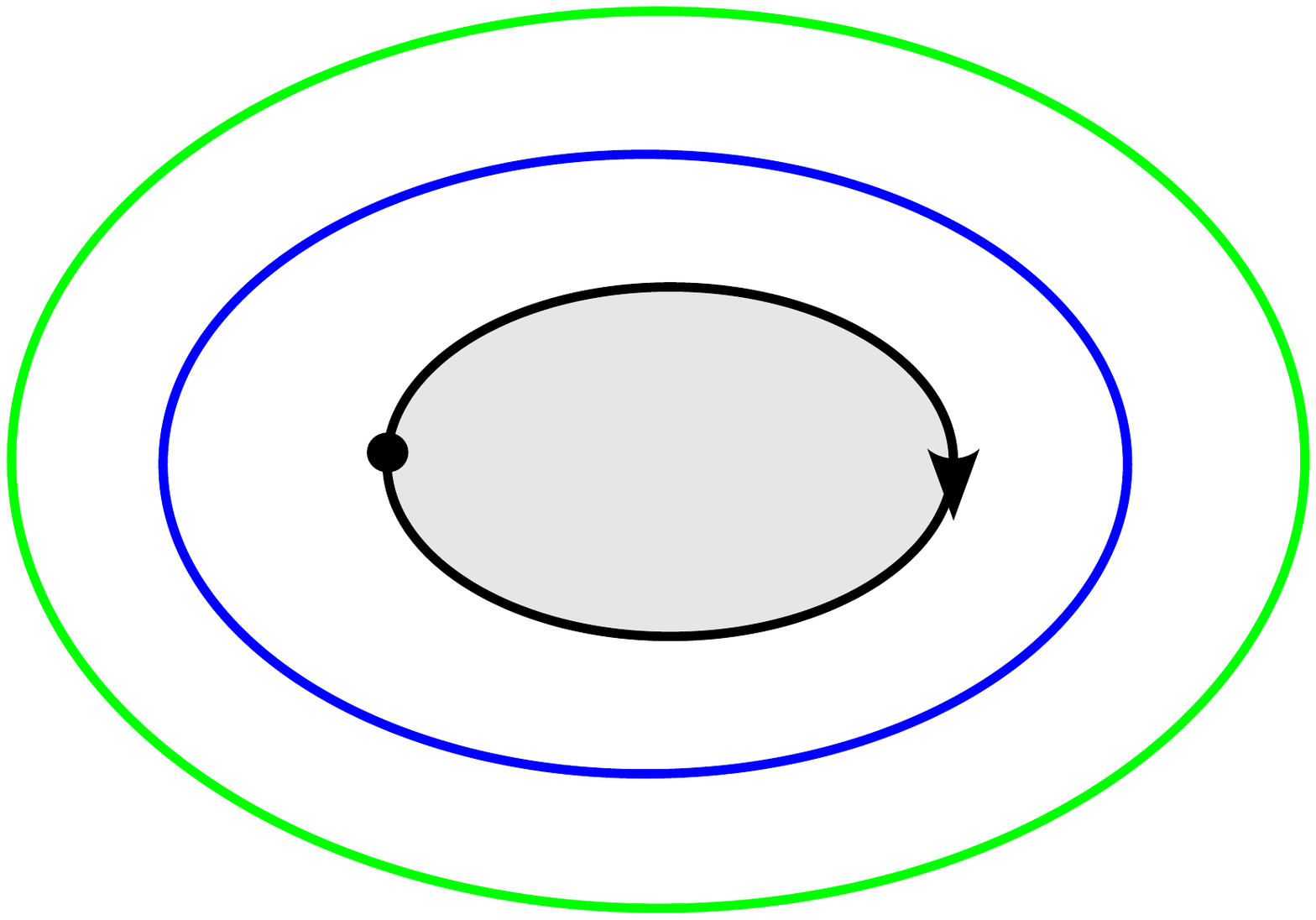}}%
    \put(0.33252182,0.36500886){\makebox(0,0)[lt]{\lineheight{1.25}\smash{\begin{tabular}[t]{l}$\alpha$\end{tabular}}}}%
    \put(0.60791688,0.35871493){\makebox(0,0)[lt]{\lineheight{1.25}\smash{\begin{tabular}[t]{l}$m$\end{tabular}}}}%
  \end{picture}%
\endgroup%
 };
	\end{tikzpicture}

\right) \, .
\end{equation}
We claim that the family $\{\Theta_m\}_{m\in \Ob \calm}$ is a $(\mac,\mad)$-bimodule trace on $\calm$.

It is easy to see that it is a trace on $\calm$, i.~e.~satisfies the cyclicity and the non-degeneracy axioms in Definition \ref{def:moduletrace}. The cyclicity  follows, because the two boundary vertices labelled with morphisms $\alpha: m\to n$ and $\beta:n\to m$ can be combined into a single boundary vertex labelled with $\alpha\circ \beta: n\to n$ or $\beta\circ \alpha: m\to m$ and because the evaluation  does not depend on the choice of basepoint at the boundary.
Non-degeneracy follows, because the dimensions of simple objects $\Theta_m(1_m) = \mathrm{av}_{m}(X)$ for $m\in I_\calm$ are non-zero by assumption, see \cite[Lemma 5.2]{S}.

It remains to prove   the module trace condition from Definition \ref{def:moduletrace}, which corresponds to the graphical identities in \eqref{pic:cmoduletrace}. As we show below, this  
follows with the diagrammatic calculus, more specifically, the invariance of the diagram evaluation under isotopies, the two-point or Reidemeister 2 move and the invariance under handle slides. Note that this does not require that the bimodule category is stabilising, see Remark \ref{rem:trisecinvrem} and the diagrammatic computation in Figure \ref{fig:handleslide}.

The first identity in \eqref{pic:cmoduletrace} can be shown as follows, and the proof of the second identity is analogous. Invariance under the handle slide follows analogously to the diagrammatic computation in Figure \ref{fig:handleslide}.

\begin{align} \nonumber
		& \Theta_{c\,\triangleright \, m}(\alpha) = 
		\mathrm{av}_m\left( 
	\begin{tikzpicture}[baseline={([yshift=-.5ex]current bounding box.center)}]
	\node at (0,0) {\def\svgscale{0.18} 
\begingroup%
  \makeatletter%
  \providecommand\color[2][]{%
    \errmessage{(Inkscape) Color is used for the text in Inkscape, but the package 'color.sty' is not loaded}%
    \renewcommand\color[2][]{}%
  }%
  \providecommand\transparent[1]{%
    \errmessage{(Inkscape) Transparency is used (non-zero) for the text in Inkscape, but the package 'transparent.sty' is not loaded}%
    \renewcommand\transparent[1]{}%
  }%
  \providecommand\rotatebox[2]{#2}%
  \newcommand*\fsize{\dimexpr\f@size pt\relax}%
  \newcommand*\lineheight[1]{\fontsize{\fsize}{#1\fsize}\selectfont}%
  \ifx\svgwidth\undefined%
    \setlength{\unitlength}{841.88976378bp}%
    \ifx\svgscale\undefined%
      \relax%
    \else%
      \setlength{\unitlength}{\unitlength * \real{\svgscale}}%
    \fi%
  \else%
    \setlength{\unitlength}{\svgwidth}%
  \fi%
  \global\let\svgwidth\undefined%
  \global\let\svgscale\undefined%
  \makeatother%
  \begin{picture}(1,0.70707071)%
    \lineheight{1}%
    \setlength\tabcolsep{0pt}%
    \put(0,0){\includegraphics[width=\unitlength,page=1]{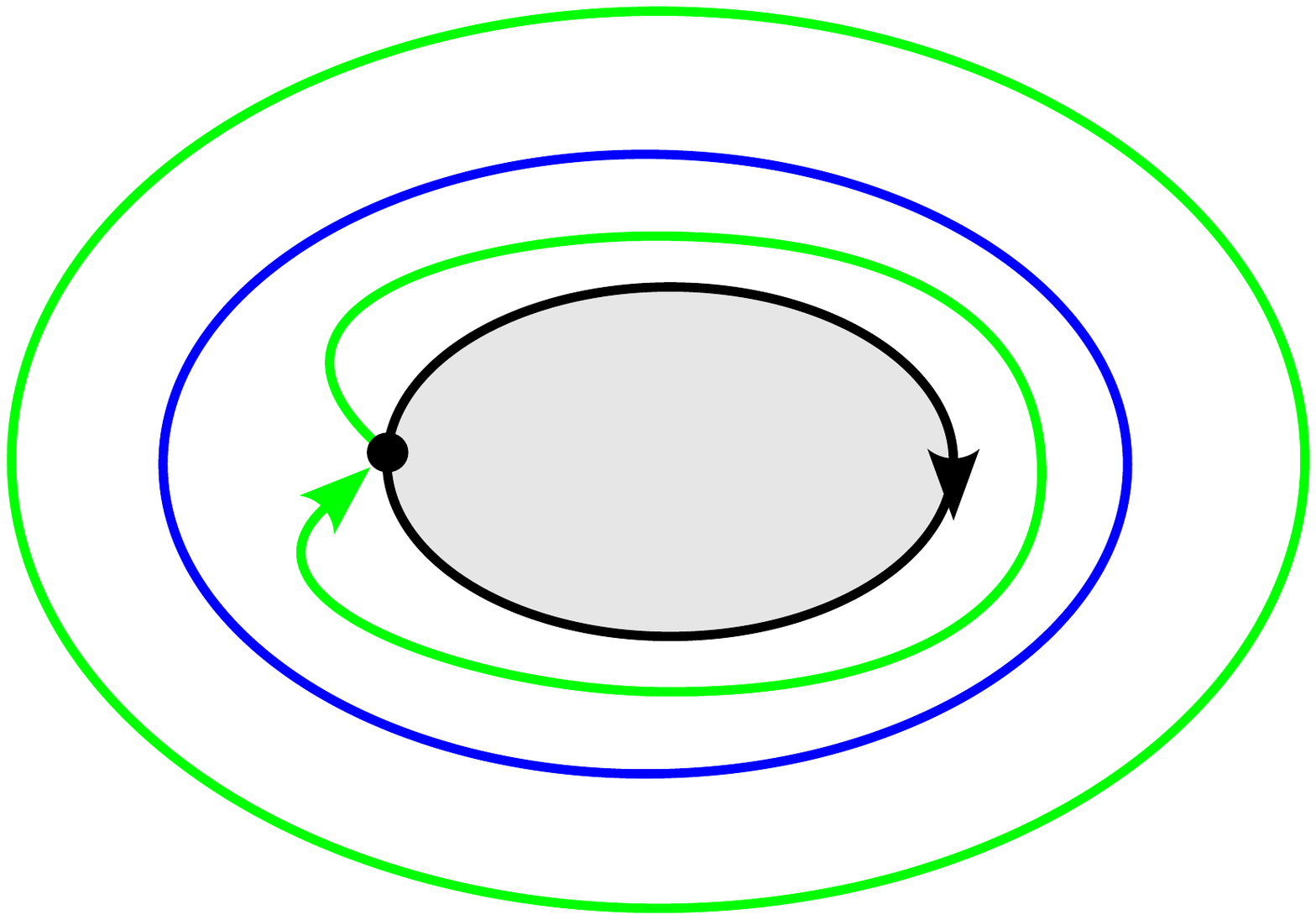}}%
    \put(0.33252182,0.36500886){\makebox(0,0)[lt]{\lineheight{1.25}\smash{\begin{tabular}[t]{l}$\alpha$\end{tabular}}}}%
    \put(0.60791688,0.35871493){\makebox(0,0)[lt]{\lineheight{1.25}\smash{\begin{tabular}[t]{l}$m$\end{tabular}}}}%
    \put(0.17096294,0.34206195){\makebox(0,0)[lt]{\lineheight{1.25}\smash{\begin{tabular}[t]{l}${\color{green}c}$\end{tabular}}}}%
  \end{picture}%
\endgroup%
 };
	\end{tikzpicture}
 \right) \stackrel{\text{R2}}{=}
		\mathrm{av}_m\left( 
	\begin{tikzpicture}[baseline={([yshift=-.5ex]current bounding box.center)}]
	\node at (0,0) {\def\svgscale{0.18} 
\begingroup%
  \makeatletter%
  \providecommand\color[2][]{%
    \errmessage{(Inkscape) Color is used for the text in Inkscape, but the package 'color.sty' is not loaded}%
    \renewcommand\color[2][]{}%
  }%
  \providecommand\transparent[1]{%
    \errmessage{(Inkscape) Transparency is used (non-zero) for the text in Inkscape, but the package 'transparent.sty' is not loaded}%
    \renewcommand\transparent[1]{}%
  }%
  \providecommand\rotatebox[2]{#2}%
  \newcommand*\fsize{\dimexpr\f@size pt\relax}%
  \newcommand*\lineheight[1]{\fontsize{\fsize}{#1\fsize}\selectfont}%
  \ifx\svgwidth\undefined%
    \setlength{\unitlength}{841.88976378bp}%
    \ifx\svgscale\undefined%
      \relax%
    \else%
      \setlength{\unitlength}{\unitlength * \real{\svgscale}}%
    \fi%
  \else%
    \setlength{\unitlength}{\svgwidth}%
  \fi%
  \global\let\svgwidth\undefined%
  \global\let\svgscale\undefined%
  \makeatother%
  \begin{picture}(1,0.70707071)%
    \lineheight{1}%
    \setlength\tabcolsep{0pt}%
    \put(0,0){\includegraphics[width=\unitlength,page=1]{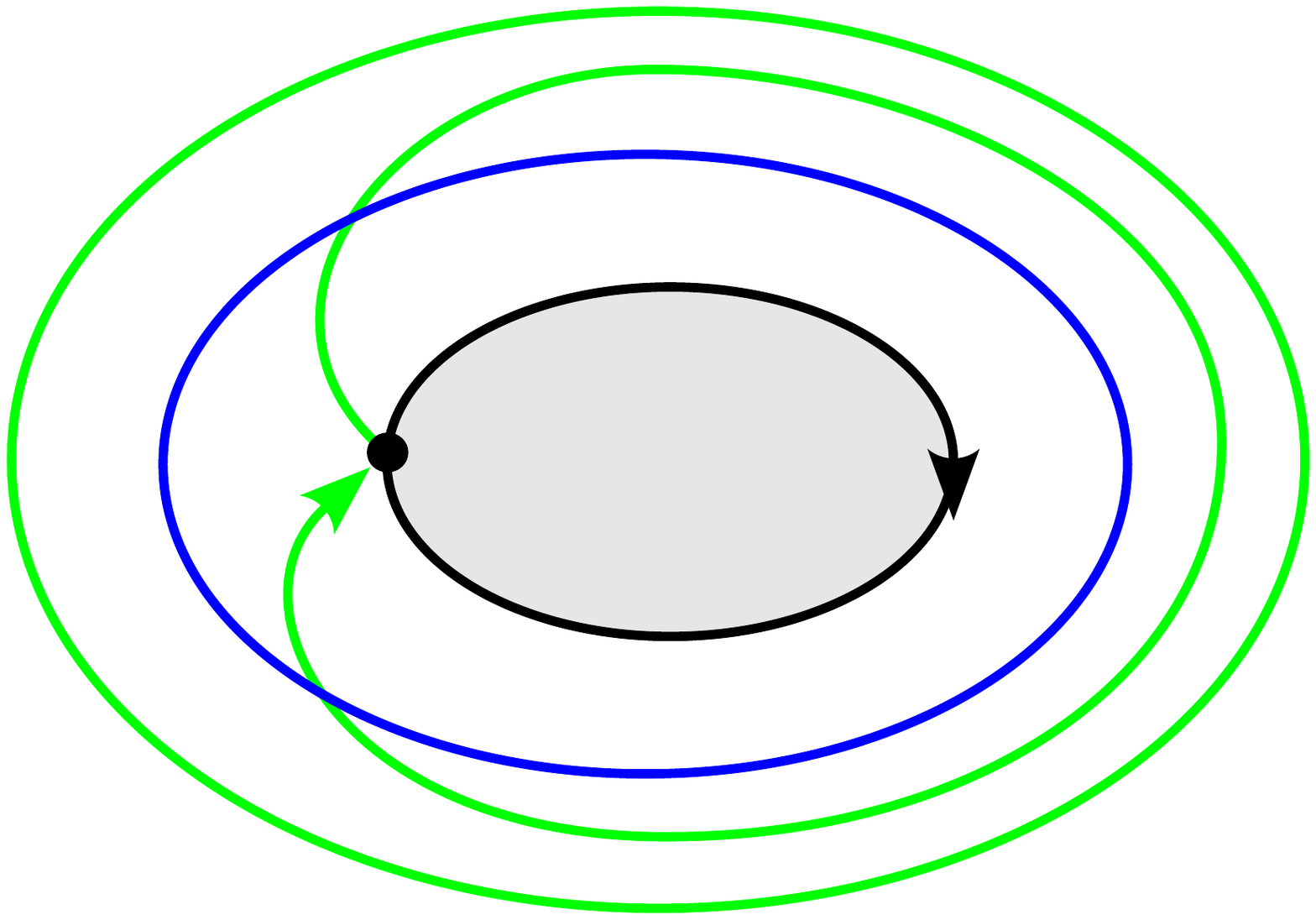}}%
    \put(0.33252182,0.36500886){\makebox(0,0)[lt]{\lineheight{1.25}\smash{\begin{tabular}[t]{l}$\alpha$\end{tabular}}}}%
    \put(0.60173001,0.35077022){\makebox(0,0)[lt]{\lineheight{1.25}\smash{\begin{tabular}[t]{l}$m$\end{tabular}}}}%
    \put(0.17597572,0.33754972){\makebox(0,0)[lt]{\lineheight{1.25}\smash{\begin{tabular}[t]{l}${\color{green}c}$\end{tabular}}}}%
  \end{picture}%
\endgroup%
 };
	\end{tikzpicture}
 \right)\\ \nonumber
        & \;\;\stackrel{\text{handle slide}}{=}
		\mathrm{av}_m\left( 
	\begin{tikzpicture}[baseline={([yshift=-.5ex]current bounding box.center)}]
	\node at (0,0) {\def\svgscale{0.18} 
\begingroup%
  \makeatletter%
  \providecommand\color[2][]{%
    \errmessage{(Inkscape) Color is used for the text in Inkscape, but the package 'color.sty' is not loaded}%
    \renewcommand\color[2][]{}%
  }%
  \providecommand\transparent[1]{%
    \errmessage{(Inkscape) Transparency is used (non-zero) for the text in Inkscape, but the package 'transparent.sty' is not loaded}%
    \renewcommand\transparent[1]{}%
  }%
  \providecommand\rotatebox[2]{#2}%
  \newcommand*\fsize{\dimexpr\f@size pt\relax}%
  \newcommand*\lineheight[1]{\fontsize{\fsize}{#1\fsize}\selectfont}%
  \ifx\svgwidth\undefined%
    \setlength{\unitlength}{841.88976378bp}%
    \ifx\svgscale\undefined%
      \relax%
    \else%
      \setlength{\unitlength}{\unitlength * \real{\svgscale}}%
    \fi%
  \else%
    \setlength{\unitlength}{\svgwidth}%
  \fi%
  \global\let\svgwidth\undefined%
  \global\let\svgscale\undefined%
  \makeatother%
  \begin{picture}(1,0.70707071)%
    \lineheight{1}%
    \setlength\tabcolsep{0pt}%
    \put(0,0){\includegraphics[width=\unitlength,page=1]{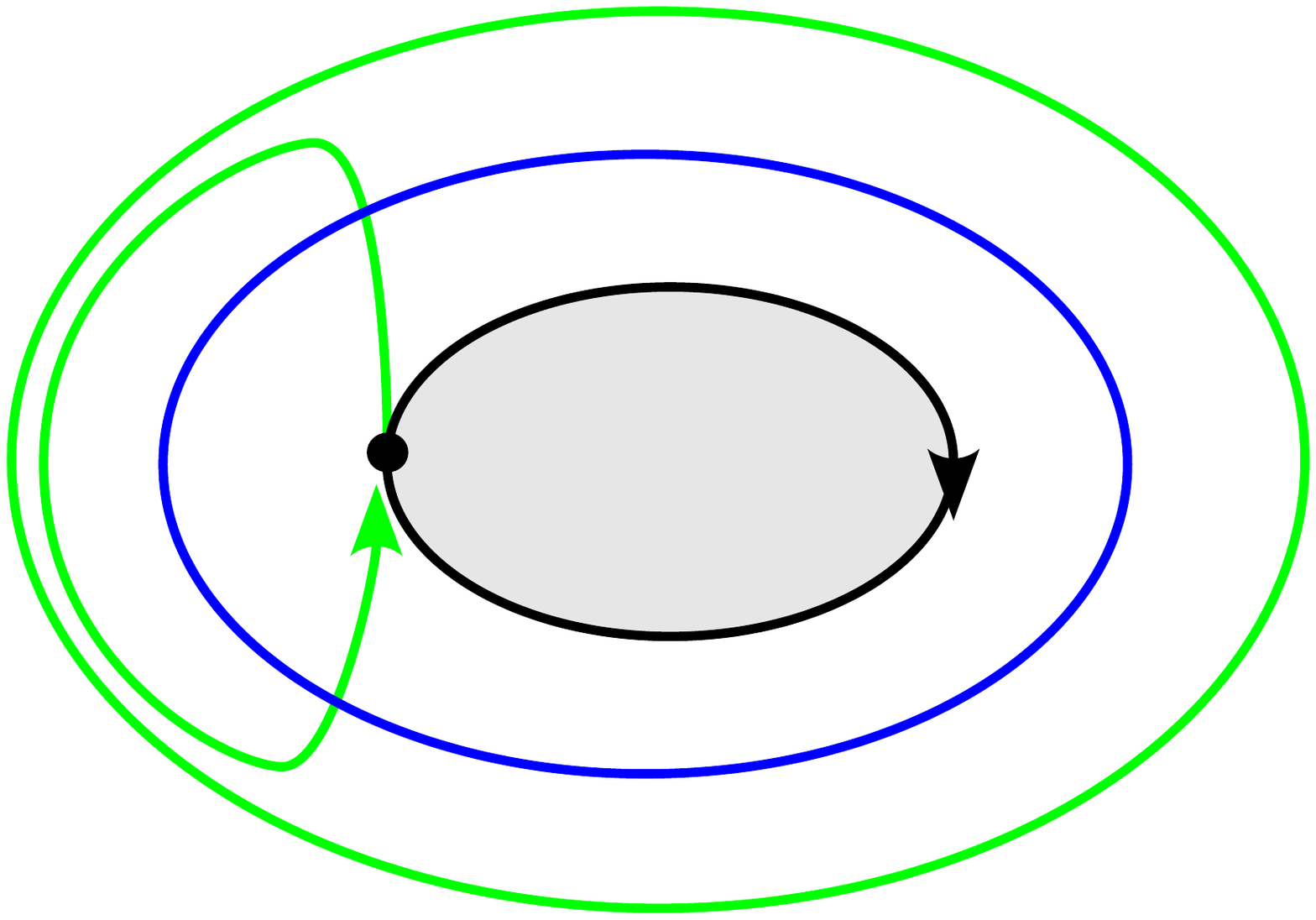}}%
    \put(0.33252182,0.36500886){\makebox(0,0)[lt]{\lineheight{1.25}\smash{\begin{tabular}[t]{l}$\alpha$\end{tabular}}}}%
    \put(0.59733392,0.35586229){\makebox(0,0)[lt]{\lineheight{1.25}\smash{\begin{tabular}[t]{l}$m$\end{tabular}}}}%
    \put(0.20746153,0.31621142){\makebox(0,0)[lt]{\lineheight{1.25}\smash{\begin{tabular}[t]{l}${\color{green}c}$\end{tabular}}}}%
  \end{picture}%
\endgroup%
 };
	\end{tikzpicture}
 \right)
		  \stackrel{\text{R2}}{=}
		\mathrm{av}_m\left( 
	\begin{tikzpicture}[baseline={([yshift=-.5ex]current bounding box.center)}]
	\node at (0,0) {\def\svgscale{0.18} 
\begingroup%
  \makeatletter%
  \providecommand\color[2][]{%
    \errmessage{(Inkscape) Color is used for the text in Inkscape, but the package 'color.sty' is not loaded}%
    \renewcommand\color[2][]{}%
  }%
  \providecommand\transparent[1]{%
    \errmessage{(Inkscape) Transparency is used (non-zero) for the text in Inkscape, but the package 'transparent.sty' is not loaded}%
    \renewcommand\transparent[1]{}%
  }%
  \providecommand\rotatebox[2]{#2}%
  \newcommand*\fsize{\dimexpr\f@size pt\relax}%
  \newcommand*\lineheight[1]{\fontsize{\fsize}{#1\fsize}\selectfont}%
  \ifx\svgwidth\undefined%
    \setlength{\unitlength}{841.88976378bp}%
    \ifx\svgscale\undefined%
      \relax%
    \else%
      \setlength{\unitlength}{\unitlength * \real{\svgscale}}%
    \fi%
  \else%
    \setlength{\unitlength}{\svgwidth}%
  \fi%
  \global\let\svgwidth\undefined%
  \global\let\svgscale\undefined%
  \makeatother%
  \begin{picture}(1,0.70707071)%
    \lineheight{1}%
    \setlength\tabcolsep{0pt}%
    \put(0,0){\includegraphics[width=\unitlength,page=1]{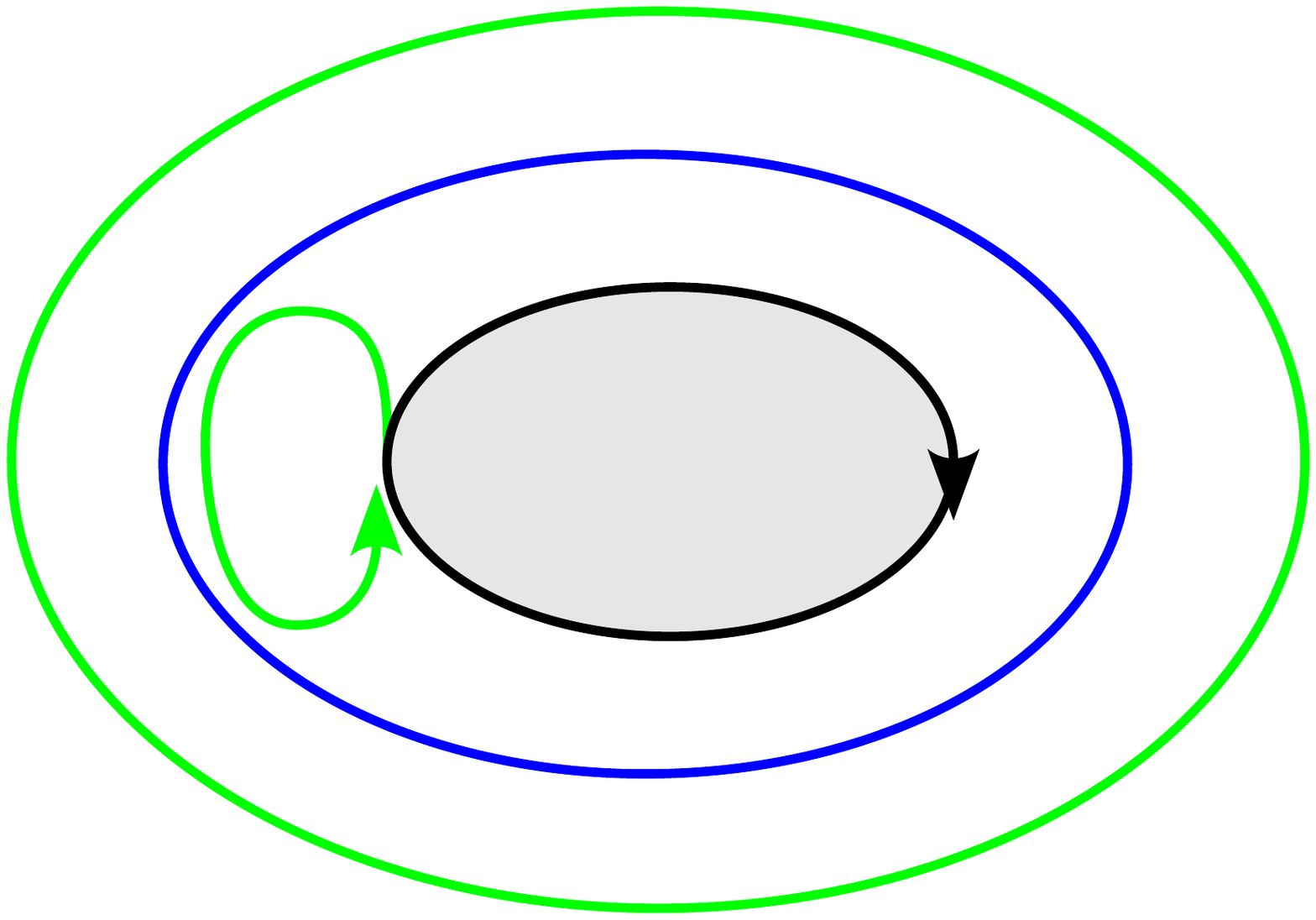}}%
    \put(0.5981666,0.3486511){\makebox(0,0)[lt]{\lineheight{1.25}\smash{\begin{tabular}[t]{l}$m$\end{tabular}}}}%
    \put(0,0){\includegraphics[width=\unitlength,page=2]{stablemma_tr4.pdf}}%
    \put(0.33252182,0.36500886){\makebox(0,0)[lt]{\lineheight{1.25}\smash{\begin{tabular}[t]{l}$\alpha$\end{tabular}}}}%
    \put(0.20566403,0.32103466){\makebox(0,0)[lt]{\lineheight{1.25}\smash{\begin{tabular}[t]{l}${\color{green}c}$\end{tabular}}}}%
  \end{picture}%
\endgroup%
 };
	\end{tikzpicture}
 \right)\nonumber\\
		& =
		\Theta_{m}(\mathrm{tr}_{\mac}(\alpha)) \, .
\end{align}

This establishes that the family of maps $\{\Theta_m\}_{m\in\calm}$ defines a $(\mac,\mad)$-bimodule trace on $\calm$. Proposition \ref{prop:uniquetrace} then implies  $\Theta_m = z \cdot \mathrm{tr}_m$ for the bimodule trace $\mathrm{tr}$ on $\calm$ and fixed $z\in\C^\times$. The claim follows. 
\end{Proof}

\begin{Lemma}\label{prop:stabilising}
Let $(\cala,\mad,\mac,\calm,\Phi)$  be as in Definition \ref{def:stabilisingcond}, but not necessarily stabilising. Let $T$ be a 
trisection diagram with the green and blue curves in standard position, as in Figure \ref{fig:standardheegard}, and let  $X$ be the surface diagram obtained by removing  a disc that intersects no curves of $T$. 
If $\mathrm{av}_m(X)\neq 0$ for all boundary labels $m\in I_\calm$, then there is a constant $C\in\C^\times$, independent of $m\in I_\calm$, such that
$$
\mathrm{av}_m(X)=\dim(m)\cdot C.
$$
\end{Lemma}

\begin{Proof}
It is sufficient to show that  for all boundary labels $m\in I_\calm$ one has
\begin{equation}
\label{eq:stabpropid}
\dim\mad \cdot \dim\mac \cdot
\mathrm{av}_m(X) = \mathrm{av}_m(X'), 
\end{equation}
where $X'$ is obtained from $X$ by placing a green and blue circle around the boundary circle. The claim then follows immediately from Lemma \ref{lem:stablemma}.

We prove Identity \eqref{eq:stabpropid} via diagrammatic calculus. We start with the observation that adding a contractible green or blue loop to $X$ changes its evaluation by $\dim\mac$ and $\dim \mad$, respectively.

Using this observation together with the invariance of the diagram evaluation under isotopies, the generalised two-point move (R2) and three-point move (R3) and  handle slides, we obtain
\begin{align} \nonumber
&
\dim\mad \cdot
\mathrm{av}_m \left( 
	\begin{tikzpicture}[baseline={([yshift=-.5ex]current bounding box.center)}]
	\node at (0,0) {\def\svgscale{0.18} 
\begingroup%
  \makeatletter%
  \providecommand\color[2][]{%
    \errmessage{(Inkscape) Color is used for the text in Inkscape, but the package 'color.sty' is not loaded}%
    \renewcommand\color[2][]{}%
  }%
  \providecommand\transparent[1]{%
    \errmessage{(Inkscape) Transparency is used (non-zero) for the text in Inkscape, but the package 'transparent.sty' is not loaded}%
    \renewcommand\transparent[1]{}%
  }%
  \providecommand\rotatebox[2]{#2}%
  \newcommand*\fsize{\dimexpr\f@size pt\relax}%
  \newcommand*\lineheight[1]{\fontsize{\fsize}{#1\fsize}\selectfont}%
  \ifx\svgwidth\undefined%
    \setlength{\unitlength}{841.88976378bp}%
    \ifx\svgscale\undefined%
      \relax%
    \else%
      \setlength{\unitlength}{\unitlength * \real{\svgscale}}%
    \fi%
  \else%
    \setlength{\unitlength}{\svgwidth}%
  \fi%
  \global\let\svgwidth\undefined%
  \global\let\svgscale\undefined%
  \makeatother%
  \begin{picture}(1,0.70707071)%
    \lineheight{1}%
    \setlength\tabcolsep{0pt}%
    \put(0,0){\includegraphics[width=\unitlength,page=1]{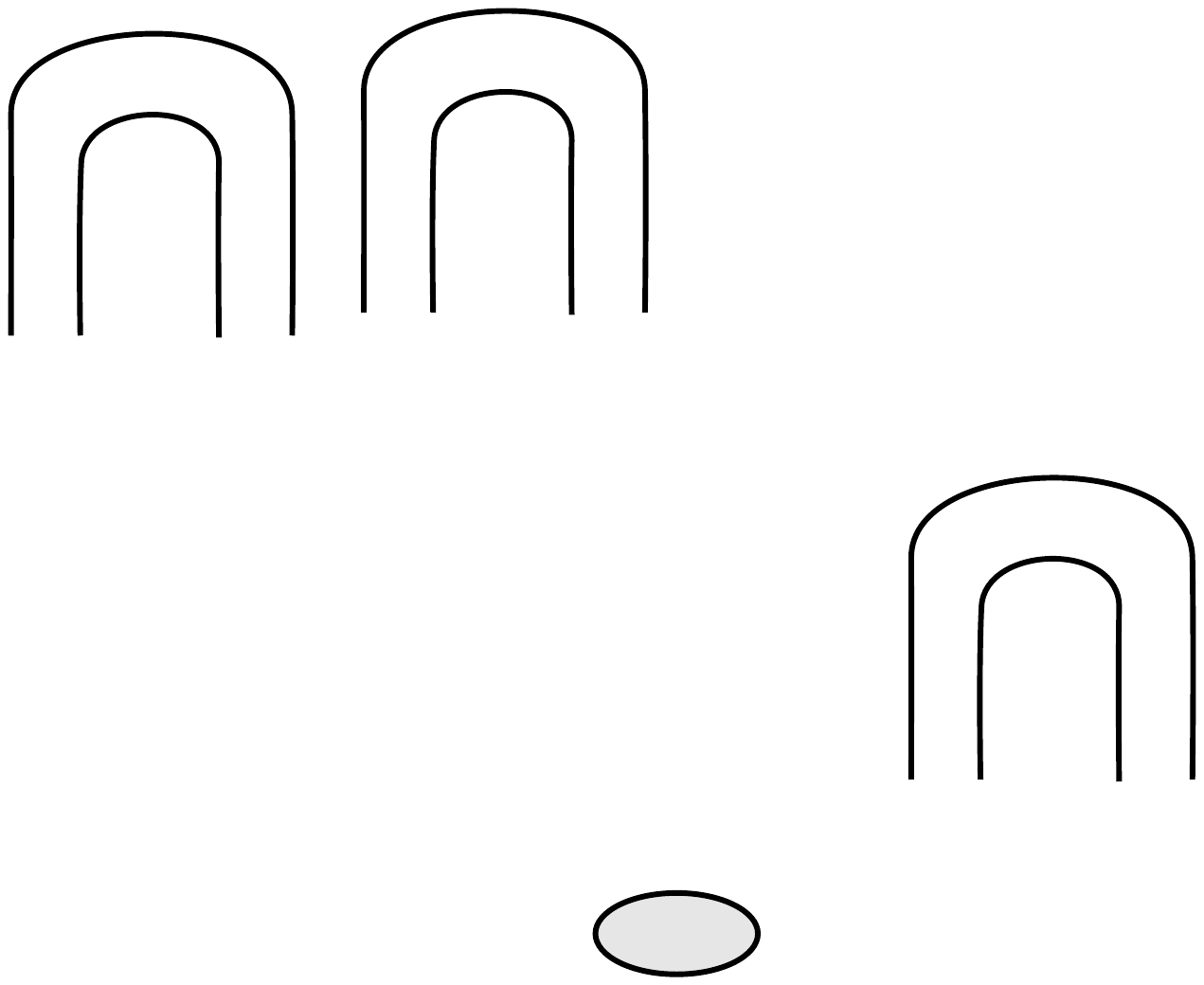}}%
    \put(0.18146344,0.44480575){\makebox(0,0)[lt]{\lineheight{1.25}\smash{\begin{tabular}[t]{l}$\ldots$\end{tabular}}}}%
    \put(0.62238519,0.44971417){\makebox(0,0)[lt]{\lineheight{1.25}\smash{\begin{tabular}[t]{l}$\ldots$\end{tabular}}}}%
    \put(0,0){\includegraphics[width=\unitlength,page=2]{stabproof1.pdf}}%
  \end{picture}%
\endgroup%
 };
	\end{tikzpicture}
 \right) =
\mathrm{av}_m\left( 
	\begin{tikzpicture}[baseline={([yshift=-.5ex]current bounding box.center)}]
	\node at (0,0) {\def\svgscale{0.18} 
\begingroup%
  \makeatletter%
  \providecommand\color[2][]{%
    \errmessage{(Inkscape) Color is used for the text in Inkscape, but the package 'color.sty' is not loaded}%
    \renewcommand\color[2][]{}%
  }%
  \providecommand\transparent[1]{%
    \errmessage{(Inkscape) Transparency is used (non-zero) for the text in Inkscape, but the package 'transparent.sty' is not loaded}%
    \renewcommand\transparent[1]{}%
  }%
  \providecommand\rotatebox[2]{#2}%
  \newcommand*\fsize{\dimexpr\f@size pt\relax}%
  \newcommand*\lineheight[1]{\fontsize{\fsize}{#1\fsize}\selectfont}%
  \ifx\svgwidth\undefined%
    \setlength{\unitlength}{841.88976378bp}%
    \ifx\svgscale\undefined%
      \relax%
    \else%
      \setlength{\unitlength}{\unitlength * \real{\svgscale}}%
    \fi%
  \else%
    \setlength{\unitlength}{\svgwidth}%
  \fi%
  \global\let\svgwidth\undefined%
  \global\let\svgscale\undefined%
  \makeatother%
  \begin{picture}(1,0.70707071)%
    \lineheight{1}%
    \setlength\tabcolsep{0pt}%
    \put(0,0){\includegraphics[width=\unitlength,page=1]{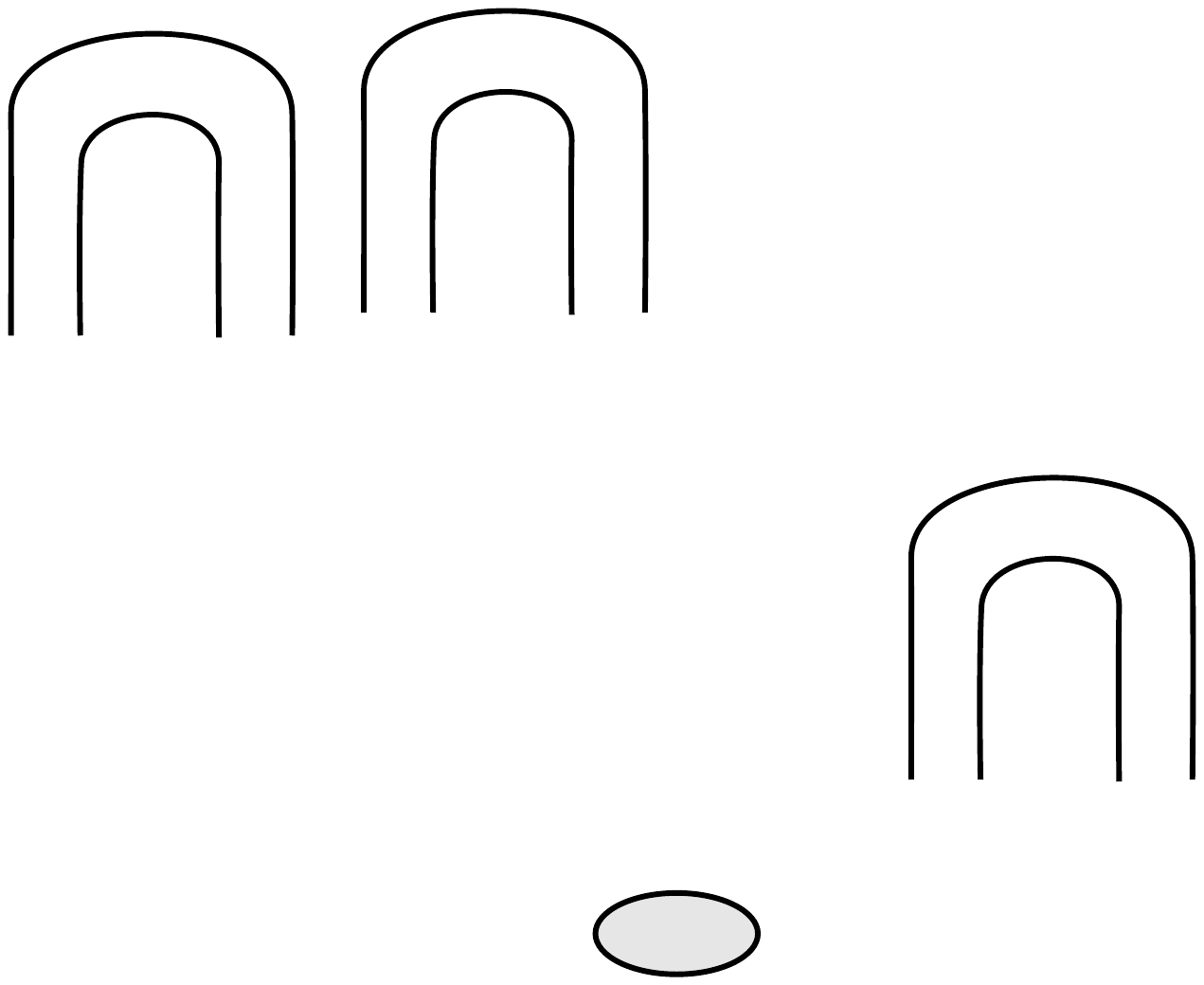}}%
    \put(0.18146344,0.44480575){\makebox(0,0)[lt]{\lineheight{1.25}\smash{\begin{tabular}[t]{l}$\ldots$\end{tabular}}}}%
    \put(0.62238519,0.44971417){\makebox(0,0)[lt]{\lineheight{1.25}\smash{\begin{tabular}[t]{l}$\ldots$\end{tabular}}}}%
    \put(0,0){\includegraphics[width=\unitlength,page=2]{stabproof2.pdf}}%
  \end{picture}%
\endgroup%
 };
	\end{tikzpicture}
 \right)
\intertext{} \nonumber
&\stackrel{\text{isotopy}}{=}
\mathrm{av}_m\left( 
	\begin{tikzpicture}[baseline={([yshift=-.5ex]current bounding box.center)}]
	\node at (0,0) {\def\svgscale{0.18} 
\begingroup%
  \makeatletter%
  \providecommand\color[2][]{%
    \errmessage{(Inkscape) Color is used for the text in Inkscape, but the package 'color.sty' is not loaded}%
    \renewcommand\color[2][]{}%
  }%
  \providecommand\transparent[1]{%
    \errmessage{(Inkscape) Transparency is used (non-zero) for the text in Inkscape, but the package 'transparent.sty' is not loaded}%
    \renewcommand\transparent[1]{}%
  }%
  \providecommand\rotatebox[2]{#2}%
  \newcommand*\fsize{\dimexpr\f@size pt\relax}%
  \newcommand*\lineheight[1]{\fontsize{\fsize}{#1\fsize}\selectfont}%
  \ifx\svgwidth\undefined%
    \setlength{\unitlength}{841.88976378bp}%
    \ifx\svgscale\undefined%
      \relax%
    \else%
      \setlength{\unitlength}{\unitlength * \real{\svgscale}}%
    \fi%
  \else%
    \setlength{\unitlength}{\svgwidth}%
  \fi%
  \global\let\svgwidth\undefined%
  \global\let\svgscale\undefined%
  \makeatother%
  \begin{picture}(1,0.70707071)%
    \lineheight{1}%
    \setlength\tabcolsep{0pt}%
    \put(0,0){\includegraphics[width=\unitlength,page=1]{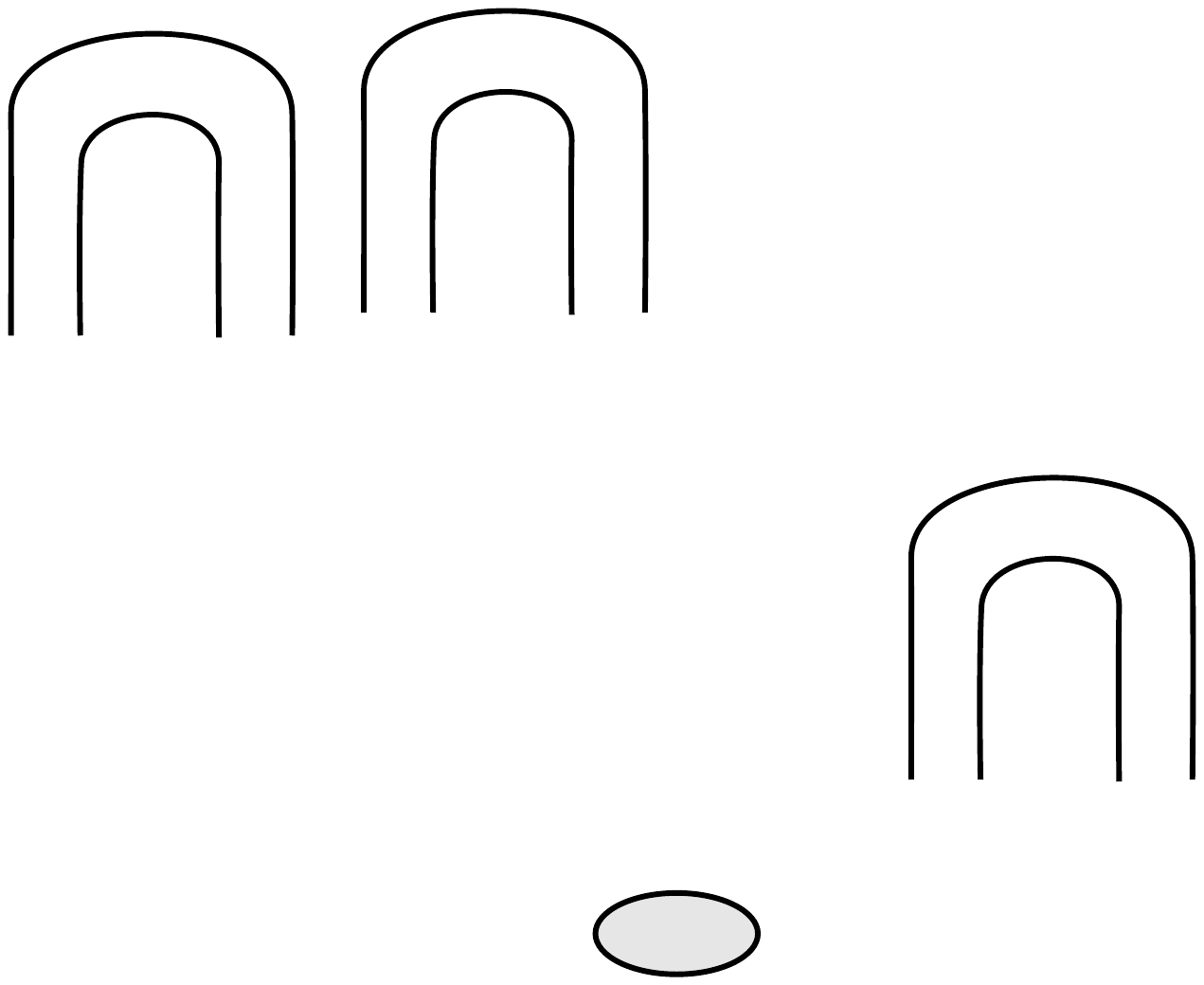}}%
    \put(0.18146344,0.44480575){\makebox(0,0)[lt]{\lineheight{1.25}\smash{\begin{tabular}[t]{l}$\ldots$\end{tabular}}}}%
    \put(0.62238519,0.44971417){\makebox(0,0)[lt]{\lineheight{1.25}\smash{\begin{tabular}[t]{l}$\ldots$\end{tabular}}}}%
    \put(0,0){\includegraphics[width=\unitlength,page=2]{stabproof3.pdf}}%
  \end{picture}%
\endgroup%
 };
	\end{tikzpicture}
 \right) \stackrel{\eqref{eq:hardslide}}{=}
\mathrm{av}_m\left( 
	\begin{tikzpicture}[baseline={([yshift=-.5ex]current bounding box.center)}]
	\node at (0,0) {\def\svgscale{0.18} 
\begingroup%
  \makeatletter%
  \providecommand\color[2][]{%
    \errmessage{(Inkscape) Color is used for the text in Inkscape, but the package 'color.sty' is not loaded}%
    \renewcommand\color[2][]{}%
  }%
  \providecommand\transparent[1]{%
    \errmessage{(Inkscape) Transparency is used (non-zero) for the text in Inkscape, but the package 'transparent.sty' is not loaded}%
    \renewcommand\transparent[1]{}%
  }%
  \providecommand\rotatebox[2]{#2}%
  \newcommand*\fsize{\dimexpr\f@size pt\relax}%
  \newcommand*\lineheight[1]{\fontsize{\fsize}{#1\fsize}\selectfont}%
  \ifx\svgwidth\undefined%
    \setlength{\unitlength}{841.88976378bp}%
    \ifx\svgscale\undefined%
      \relax%
    \else%
      \setlength{\unitlength}{\unitlength * \real{\svgscale}}%
    \fi%
  \else%
    \setlength{\unitlength}{\svgwidth}%
  \fi%
  \global\let\svgwidth\undefined%
  \global\let\svgscale\undefined%
  \makeatother%
  \begin{picture}(1,0.70707071)%
    \lineheight{1}%
    \setlength\tabcolsep{0pt}%
    \put(0,0){\includegraphics[width=\unitlength,page=1]{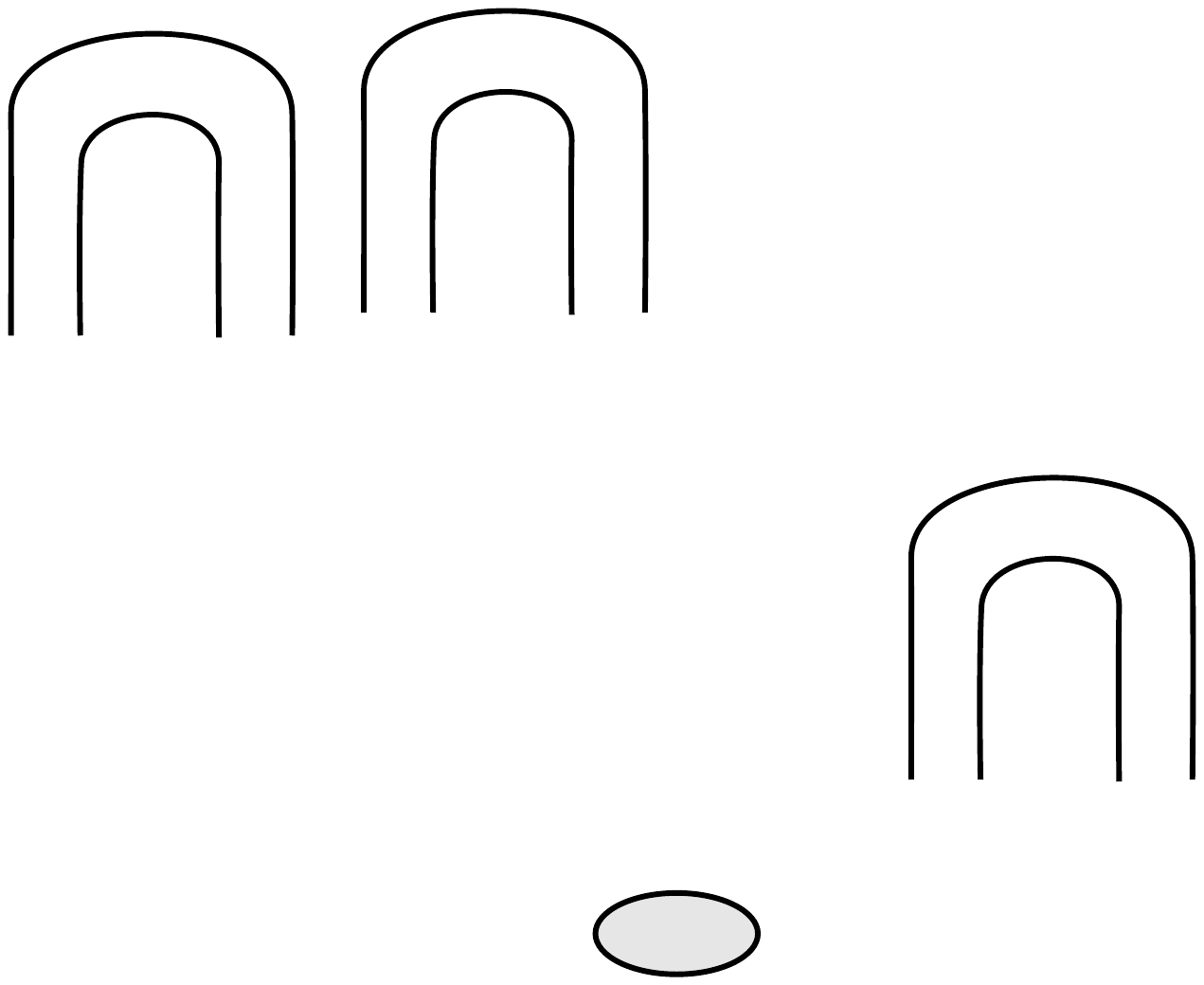}}%
    \put(0.18146344,0.44480575){\makebox(0,0)[lt]{\lineheight{1.25}\smash{\begin{tabular}[t]{l}$\ldots$\end{tabular}}}}%
    \put(0.62238519,0.44971417){\makebox(0,0)[lt]{\lineheight{1.25}\smash{\begin{tabular}[t]{l}$\ldots$\end{tabular}}}}%
    \put(0,0){\includegraphics[width=\unitlength,page=2]{stabproof4.pdf}}%
  \end{picture}%
\endgroup%
 };
	\end{tikzpicture}
 \right)
\stackrel{\eqref{eq:hardslide}}{=}\ldots
\\ 
&\;\;\;\;\stackrel{\eqref{eq:hardslide}}{=}
\mathrm{av}_m\left( 
	\begin{tikzpicture}[baseline={([yshift=-.5ex]current bounding box.center)}]
	\node at (0,0) {\def\svgscale{0.18} 
\begingroup%
  \makeatletter%
  \providecommand\color[2][]{%
    \errmessage{(Inkscape) Color is used for the text in Inkscape, but the package 'color.sty' is not loaded}%
    \renewcommand\color[2][]{}%
  }%
  \providecommand\transparent[1]{%
    \errmessage{(Inkscape) Transparency is used (non-zero) for the text in Inkscape, but the package 'transparent.sty' is not loaded}%
    \renewcommand\transparent[1]{}%
  }%
  \providecommand\rotatebox[2]{#2}%
  \newcommand*\fsize{\dimexpr\f@size pt\relax}%
  \newcommand*\lineheight[1]{\fontsize{\fsize}{#1\fsize}\selectfont}%
  \ifx\svgwidth\undefined%
    \setlength{\unitlength}{841.88976378bp}%
    \ifx\svgscale\undefined%
      \relax%
    \else%
      \setlength{\unitlength}{\unitlength * \real{\svgscale}}%
    \fi%
  \else%
    \setlength{\unitlength}{\svgwidth}%
  \fi%
  \global\let\svgwidth\undefined%
  \global\let\svgscale\undefined%
  \makeatother%
  \begin{picture}(1,0.70707071)%
    \lineheight{1}%
    \setlength\tabcolsep{0pt}%
    \put(0,0){\includegraphics[width=\unitlength,page=1]{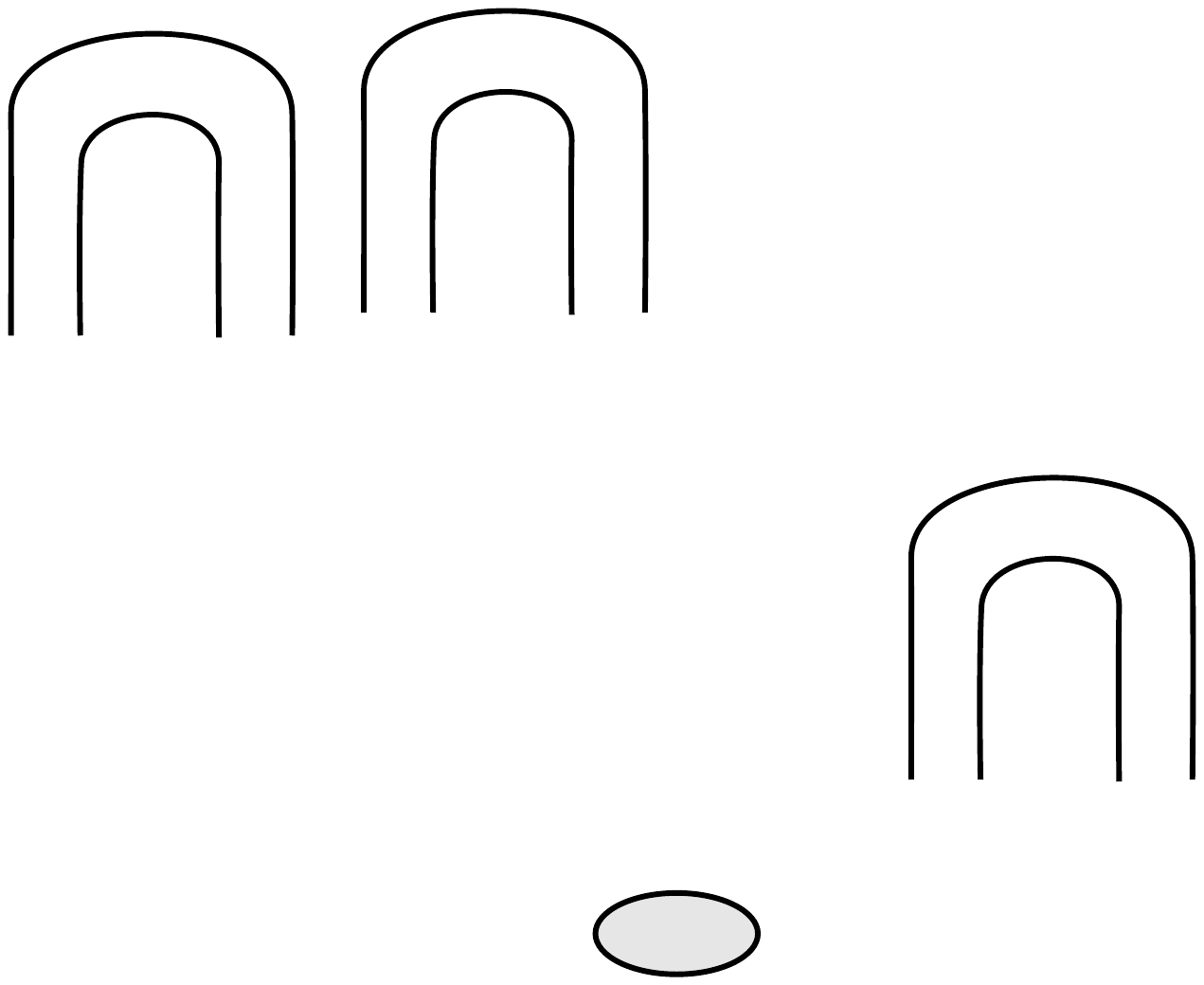}}%
    \put(0.18146344,0.44480575){\makebox(0,0)[lt]{\lineheight{1.25}\smash{\begin{tabular}[t]{l}$\ldots$\end{tabular}}}}%
    \put(0.62238519,0.44971417){\makebox(0,0)[lt]{\lineheight{1.25}\smash{\begin{tabular}[t]{l}$\ldots$\end{tabular}}}}%
    \put(0,0){\includegraphics[width=\unitlength,page=2]{stabproof6.pdf}}%
  \end{picture}%
\endgroup%
 };
	\end{tikzpicture}
 \right) \stackrel{\eqref{eq:easyslide}}{=}
\mathrm{av}_m\left( 
	\begin{tikzpicture}[baseline={([yshift=-.5ex]current bounding box.center)}]
	\node at (0,0) {\def\svgscale{0.18} 
\begingroup%
  \makeatletter%
  \providecommand\color[2][]{%
    \errmessage{(Inkscape) Color is used for the text in Inkscape, but the package 'color.sty' is not loaded}%
    \renewcommand\color[2][]{}%
  }%
  \providecommand\transparent[1]{%
    \errmessage{(Inkscape) Transparency is used (non-zero) for the text in Inkscape, but the package 'transparent.sty' is not loaded}%
    \renewcommand\transparent[1]{}%
  }%
  \providecommand\rotatebox[2]{#2}%
  \newcommand*\fsize{\dimexpr\f@size pt\relax}%
  \newcommand*\lineheight[1]{\fontsize{\fsize}{#1\fsize}\selectfont}%
  \ifx\svgwidth\undefined%
    \setlength{\unitlength}{841.88976378bp}%
    \ifx\svgscale\undefined%
      \relax%
    \else%
      \setlength{\unitlength}{\unitlength * \real{\svgscale}}%
    \fi%
  \else%
    \setlength{\unitlength}{\svgwidth}%
  \fi%
  \global\let\svgwidth\undefined%
  \global\let\svgscale\undefined%
  \makeatother%
  \begin{picture}(1,0.70707071)%
    \lineheight{1}%
    \setlength\tabcolsep{0pt}%
    \put(0,0){\includegraphics[width=\unitlength,page=1]{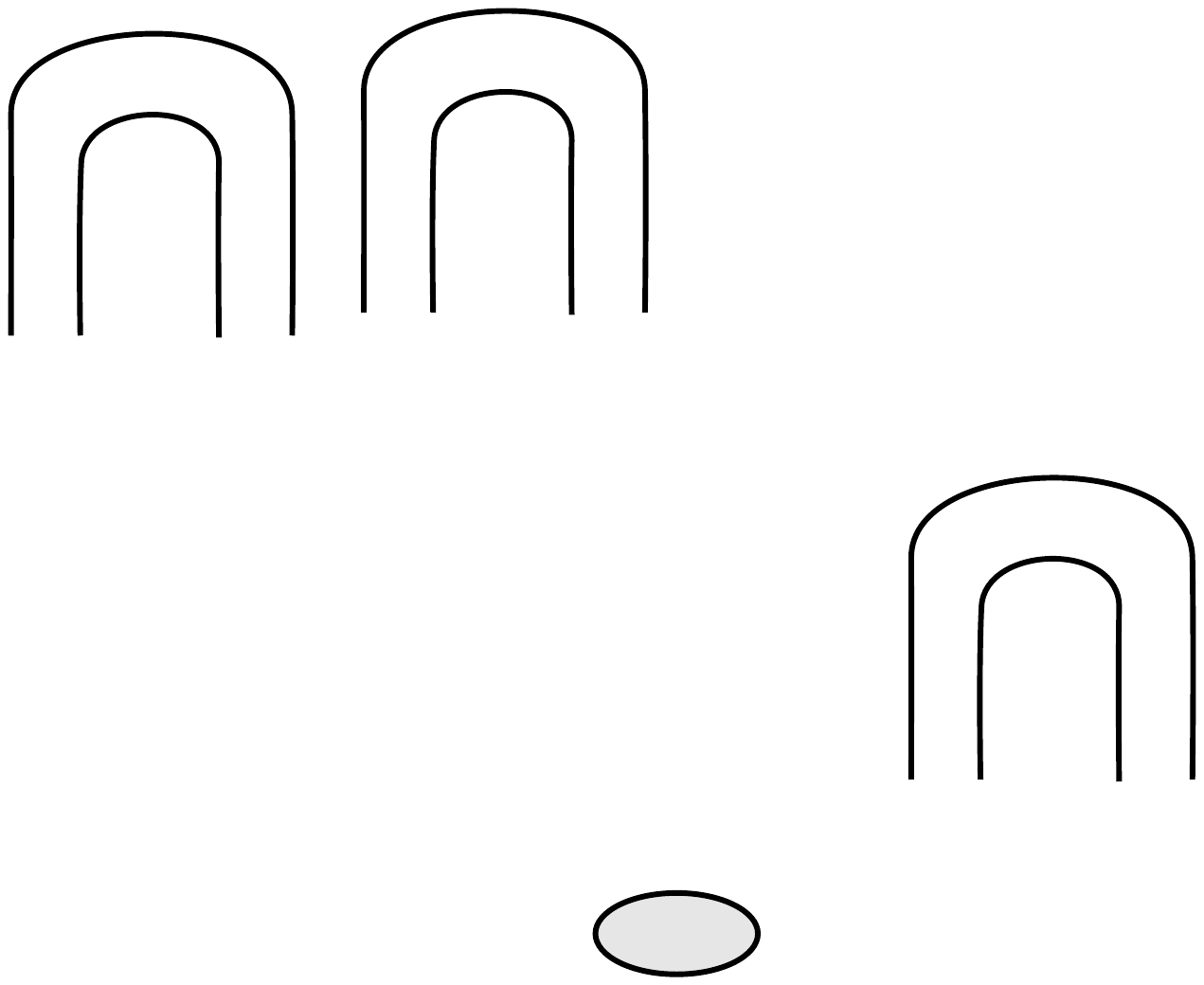}}%
    \put(0.18146344,0.44480575){\makebox(0,0)[lt]{\lineheight{1.25}\smash{\begin{tabular}[t]{l}$\ldots$\end{tabular}}}}%
    \put(0.62238519,0.44971417){\makebox(0,0)[lt]{\lineheight{1.25}\smash{\begin{tabular}[t]{l}$\ldots$\end{tabular}}}}%
    \put(0,0){\includegraphics[width=\unitlength,page=2]{stabproof7.pdf}}%
  \end{picture}%
\endgroup%
 };
	\end{tikzpicture}
 \right) 
\stackrel{\eqref{eq:easyslide}}{=}\ldots
\nonumber \\
\label{eq:stabproproof}
& 
\;\;\;\;\stackrel{\eqref{eq:easyslide}}{=}
\mathrm{av}_m\left( 
	\begin{tikzpicture}[baseline={([yshift=-.5ex]current bounding box.center)}]
	\node at (0,0) {\def\svgscale{0.18} 
\begingroup%
  \makeatletter%
  \providecommand\color[2][]{%
    \errmessage{(Inkscape) Color is used for the text in Inkscape, but the package 'color.sty' is not loaded}%
    \renewcommand\color[2][]{}%
  }%
  \providecommand\transparent[1]{%
    \errmessage{(Inkscape) Transparency is used (non-zero) for the text in Inkscape, but the package 'transparent.sty' is not loaded}%
    \renewcommand\transparent[1]{}%
  }%
  \providecommand\rotatebox[2]{#2}%
  \newcommand*\fsize{\dimexpr\f@size pt\relax}%
  \newcommand*\lineheight[1]{\fontsize{\fsize}{#1\fsize}\selectfont}%
  \ifx\svgwidth\undefined%
    \setlength{\unitlength}{841.88976378bp}%
    \ifx\svgscale\undefined%
      \relax%
    \else%
      \setlength{\unitlength}{\unitlength * \real{\svgscale}}%
    \fi%
  \else%
    \setlength{\unitlength}{\svgwidth}%
  \fi%
  \global\let\svgwidth\undefined%
  \global\let\svgscale\undefined%
  \makeatother%
  \begin{picture}(1,0.70707071)%
    \lineheight{1}%
    \setlength\tabcolsep{0pt}%
    \put(0,0){\includegraphics[width=\unitlength,page=1]{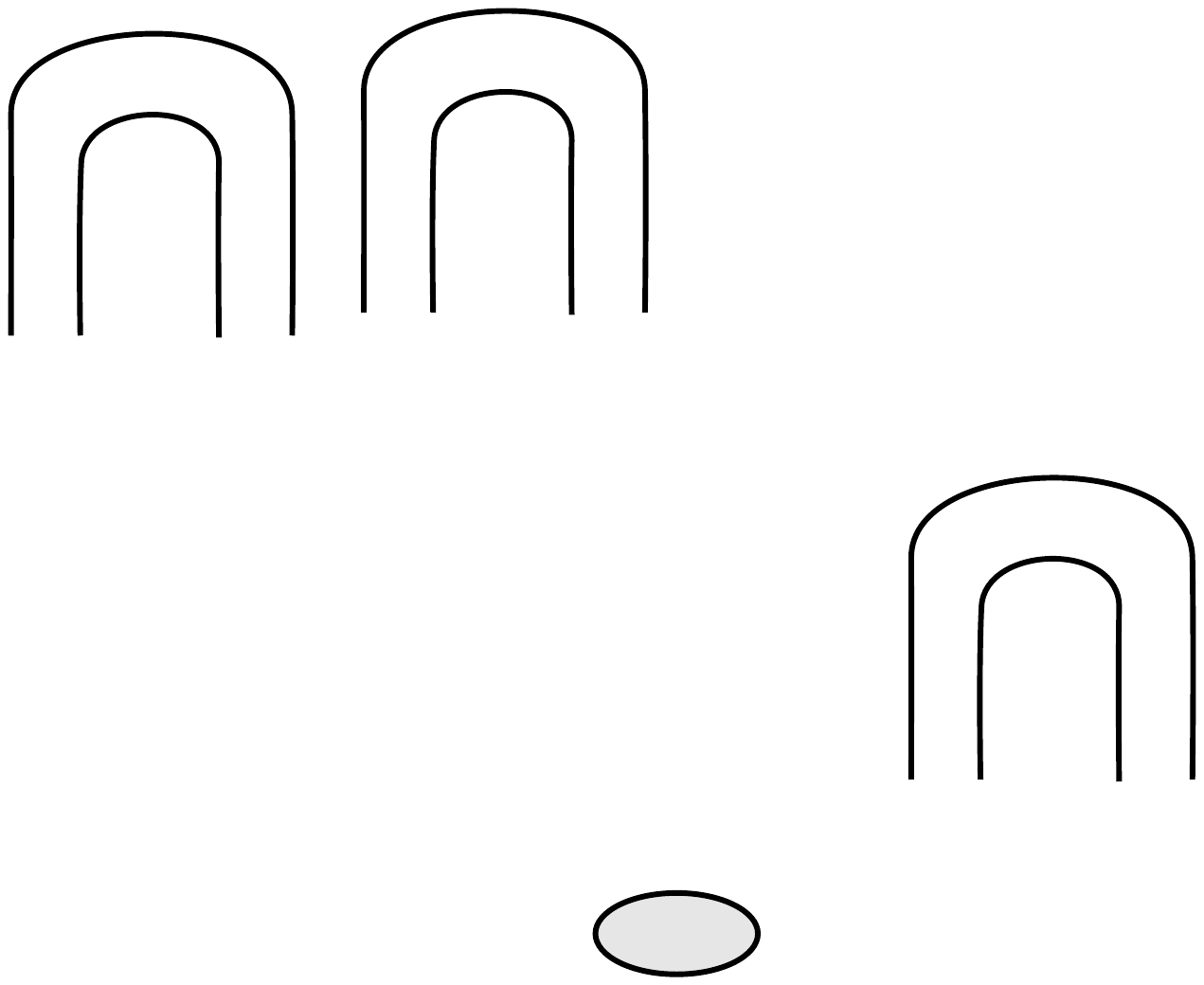}}%
    \put(0.18146344,0.44480575){\makebox(0,0)[lt]{\lineheight{1.25}\smash{\begin{tabular}[t]{l}$\ldots$\end{tabular}}}}%
    \put(0.62238519,0.44971417){\makebox(0,0)[lt]{\lineheight{1.25}\smash{\begin{tabular}[t]{l}$\ldots$\end{tabular}}}}%
    \put(0,0){\includegraphics[width=\unitlength,page=2]{stabproof9.pdf}}%
  \end{picture}%
\endgroup%
 };
	\end{tikzpicture}
 \right)
\stackrel{\eqref{eq:easyslide}}{=}
\mathrm{av}_m\left( 
	\begin{tikzpicture}[baseline={([yshift=-.5ex]current bounding box.center)}]
	\node at (0,0) {\def\svgscale{0.18} 
\begingroup%
  \makeatletter%
  \providecommand\color[2][]{%
    \errmessage{(Inkscape) Color is used for the text in Inkscape, but the package 'color.sty' is not loaded}%
    \renewcommand\color[2][]{}%
  }%
  \providecommand\transparent[1]{%
    \errmessage{(Inkscape) Transparency is used (non-zero) for the text in Inkscape, but the package 'transparent.sty' is not loaded}%
    \renewcommand\transparent[1]{}%
  }%
  \providecommand\rotatebox[2]{#2}%
  \newcommand*\fsize{\dimexpr\f@size pt\relax}%
  \newcommand*\lineheight[1]{\fontsize{\fsize}{#1\fsize}\selectfont}%
  \ifx\svgwidth\undefined%
    \setlength{\unitlength}{841.88976378bp}%
    \ifx\svgscale\undefined%
      \relax%
    \else%
      \setlength{\unitlength}{\unitlength * \real{\svgscale}}%
    \fi%
  \else%
    \setlength{\unitlength}{\svgwidth}%
  \fi%
  \global\let\svgwidth\undefined%
  \global\let\svgscale\undefined%
  \makeatother%
  \begin{picture}(1,0.70707071)%
    \lineheight{1}%
    \setlength\tabcolsep{0pt}%
    \put(0,0){\includegraphics[width=\unitlength,page=1]{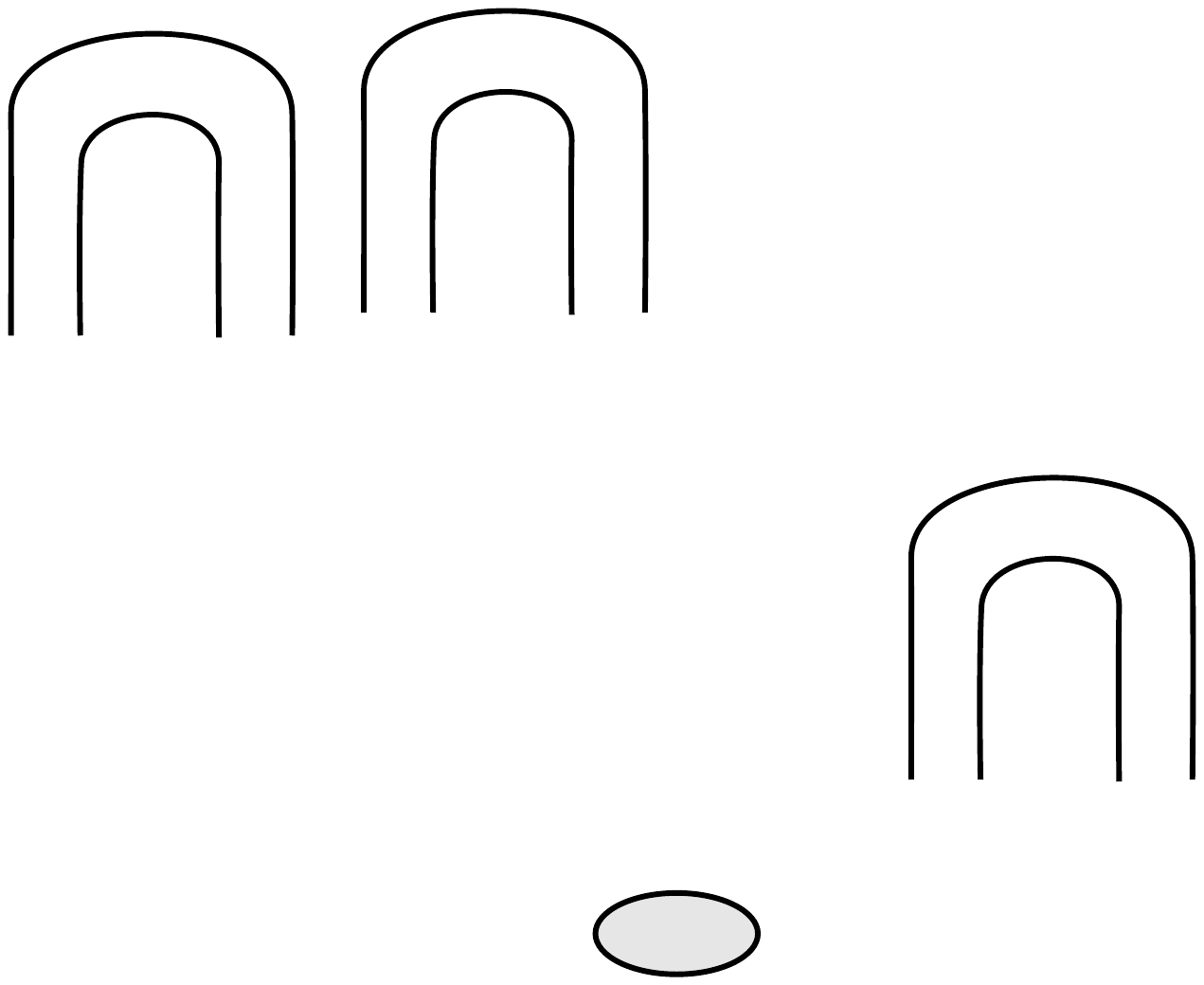}}%
    \put(0.18146344,0.44480575){\makebox(0,0)[lt]{\lineheight{1.25}\smash{\begin{tabular}[t]{l}$\ldots$\end{tabular}}}}%
    \put(0.62238519,0.44971417){\makebox(0,0)[lt]{\lineheight{1.25}\smash{\begin{tabular}[t]{l}$\ldots$\end{tabular}}}}%
    \put(0,0){\includegraphics[width=\unitlength,page=2]{stabproof10.pdf}}%
  \end{picture}%
\endgroup%
 };
	\end{tikzpicture}
 \right)\,.
\end{align}
In the last steps we use the diagrammatic identities \eqref{eq:hardslide} and \eqref{eq:easyslide} to slide the added blue line over the handles until it surrounds the boundary circle. These follow from the invariance under isotopies, two-point moves (R2), three-point moves (R3) and handle slides, see Figure~\ref{fig:stabslide}.

Note that we omitted all red curves from the diagrams in \eqref{eq:stabproproof},  \eqref{eq:hardslide} and \eqref{eq:easyslide} for legibility. This does not restrict generality, as the blue and green curves can slide over them freely due to invariance under the two-point move (R2) and three-point move (R3). Note also that the diagrammatic computations in  \eqref{eq:hardslide} and \eqref{eq:easyslide} have direct generalisations to handles with parallel blue and green circles, such as the first $k$ handles in Figure \ref{fig:standardheegard} and the three handles on the right in \eqref{eq:stabproproof}. This is again due to invariance under isotopies and the two-point and three-point moves. 

To conclude the proof, we multiply both sides of \eqref{eq:stabproproof} by $\dim(\mac)$ and perform an analogous computation, this time by adding a contractible green loop and sliding over the handles until it surrounds the boundary. This proves 
Identity \eqref{eq:stabpropid} and concludes the proof.
\end{Proof}

We now apply Lemma \ref{prop:stabilising} to the trisection diagram $\Sigma_{st}$ of $S^4$ in  \eqref{eq:stabilisation}, whose green and blue curves are in standard position. The associated diagram with a disc removed is  $\Sigma'_{st}$ from  Example \ref{ex:stabilisationsurf_gen}.

\begin{Corollary}\label{cor:stabilising}
Let $(\cala,  \mad, \mac,\calm, \Phi)$ be as in Definition \ref{def:stabilisingcond}.
Then $\calm$ is stabilising with respect to $\Phi$ if and only if $\mathrm{av}_m(\Sigma'_{st}) \neq 0$ for all $m\in I_\calm$.
\end{Corollary}

\begin{figure}
\begin{align}
\nonumber
& 
	\begin{tikzpicture}[baseline={([yshift=-.5ex]current bounding box.center)}]
	\node at (0,0) {\def\svgscale{.12} 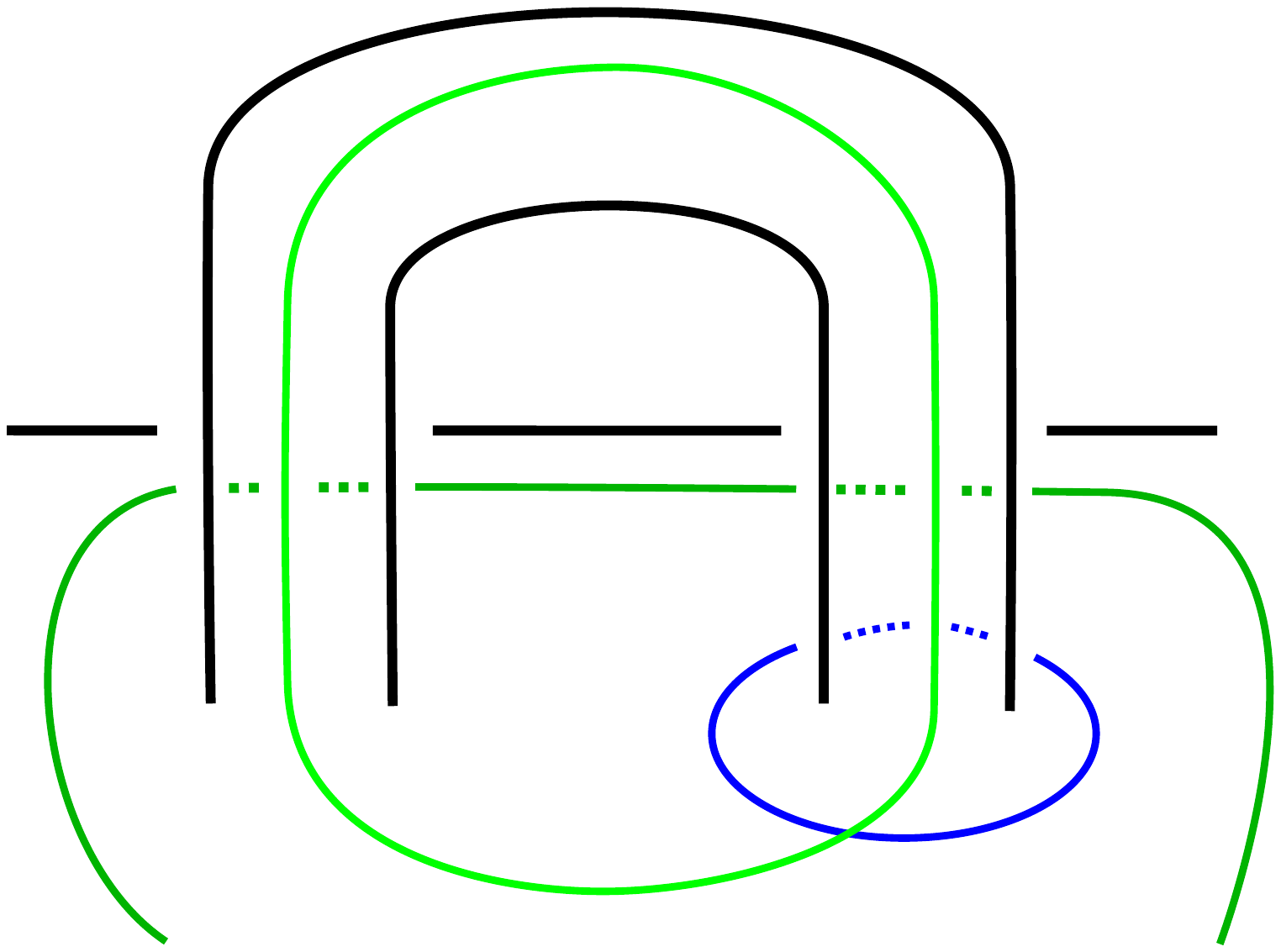 };
	\end{tikzpicture}

\stackrel{\text{R2}}{\rightsquigarrow}

	\begin{tikzpicture}[baseline={([yshift=-.5ex]current bounding box.center)}]
	\node at (0,0) {\def\svgscale{.12} 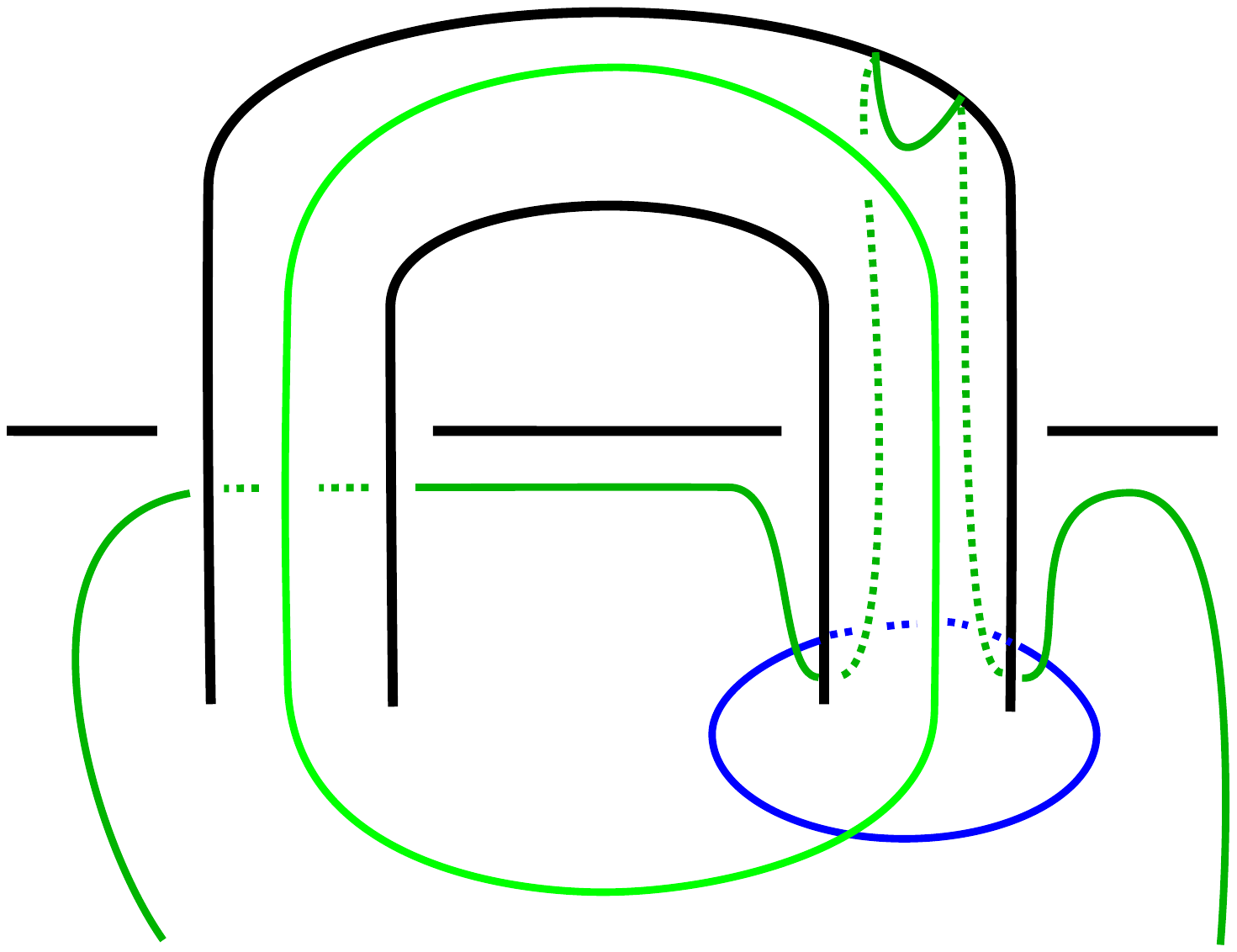 };
	\end{tikzpicture}

\stackrel{\text{handle slide}}{\rightsquigarrow}

	\begin{tikzpicture}[baseline={([yshift=-.5ex]current bounding box.center)}]
	\node at (0,0) {\def\svgscale{.12} 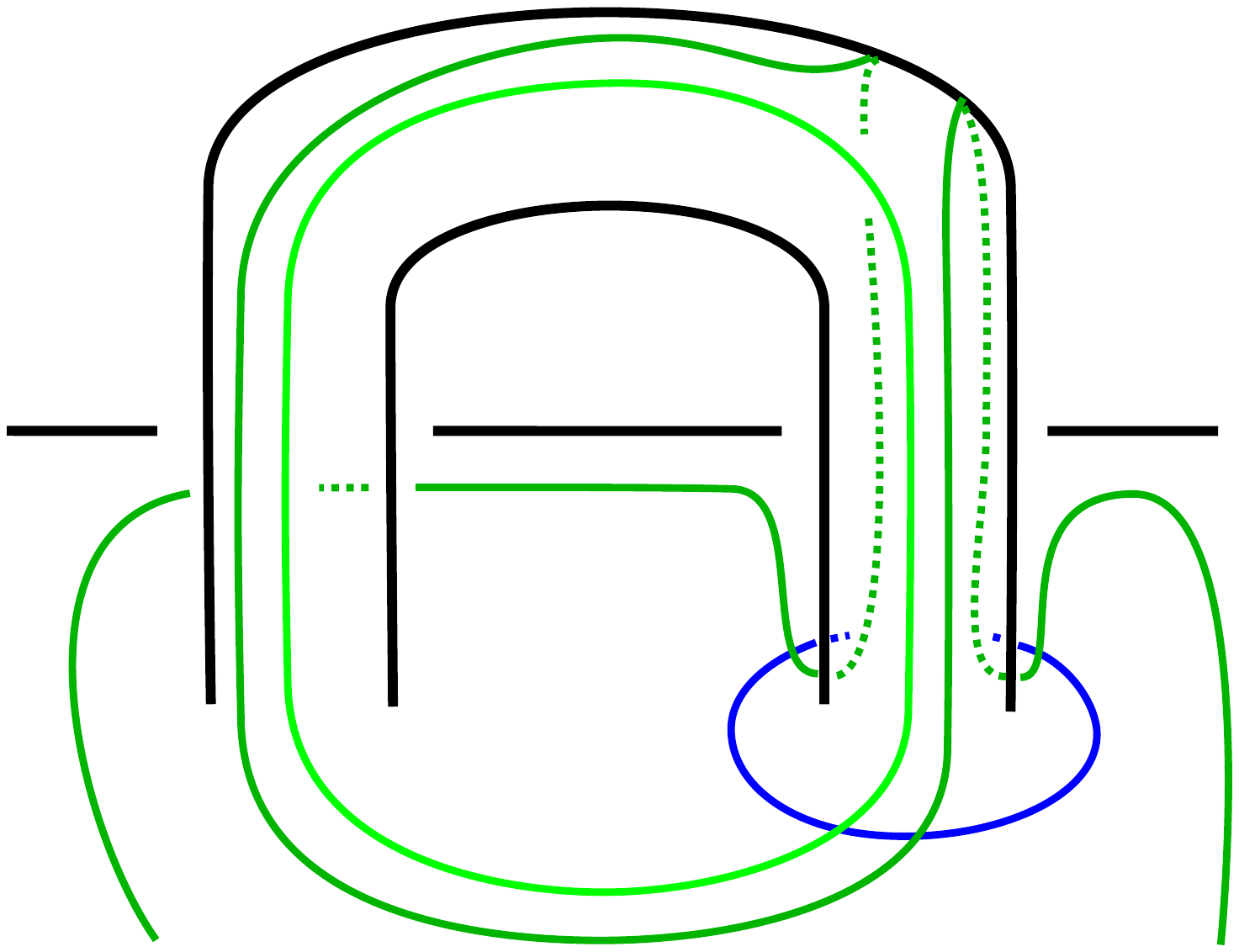 };
	\end{tikzpicture}
\\
\nonumber
\stackrel{\text{isotopy}}{\rightsquigarrow}
& 
	\begin{tikzpicture}[baseline={([yshift=-.5ex]current bounding box.center)}]
	\node at (0,0) {\def\svgscale{.12} 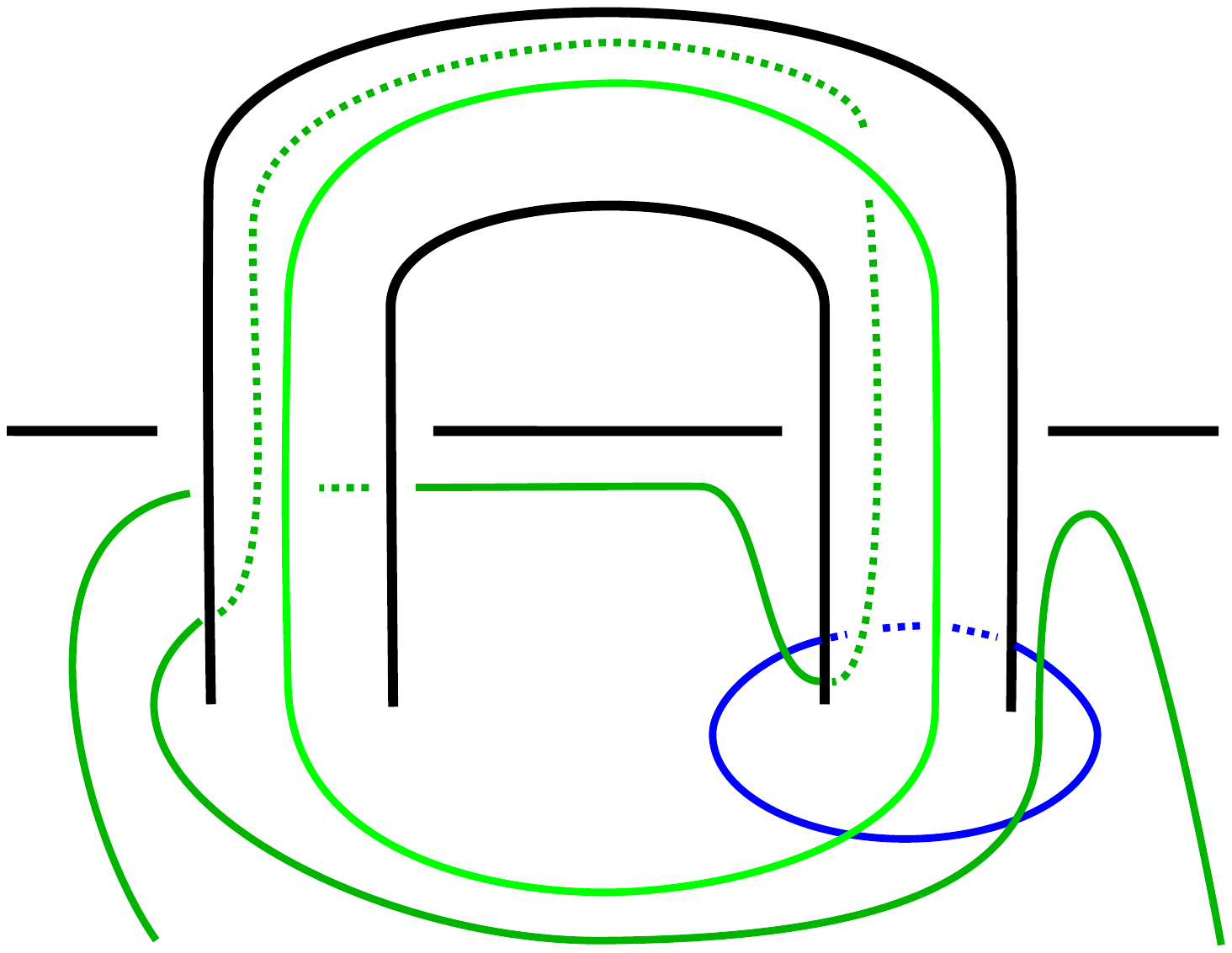 };
	\end{tikzpicture}

\stackrel{\text{R2}}{\rightsquigarrow}

	\begin{tikzpicture}[baseline={([yshift=-.5ex]current bounding box.center)}]
	\node at (0,0) {\def\svgscale{.12} 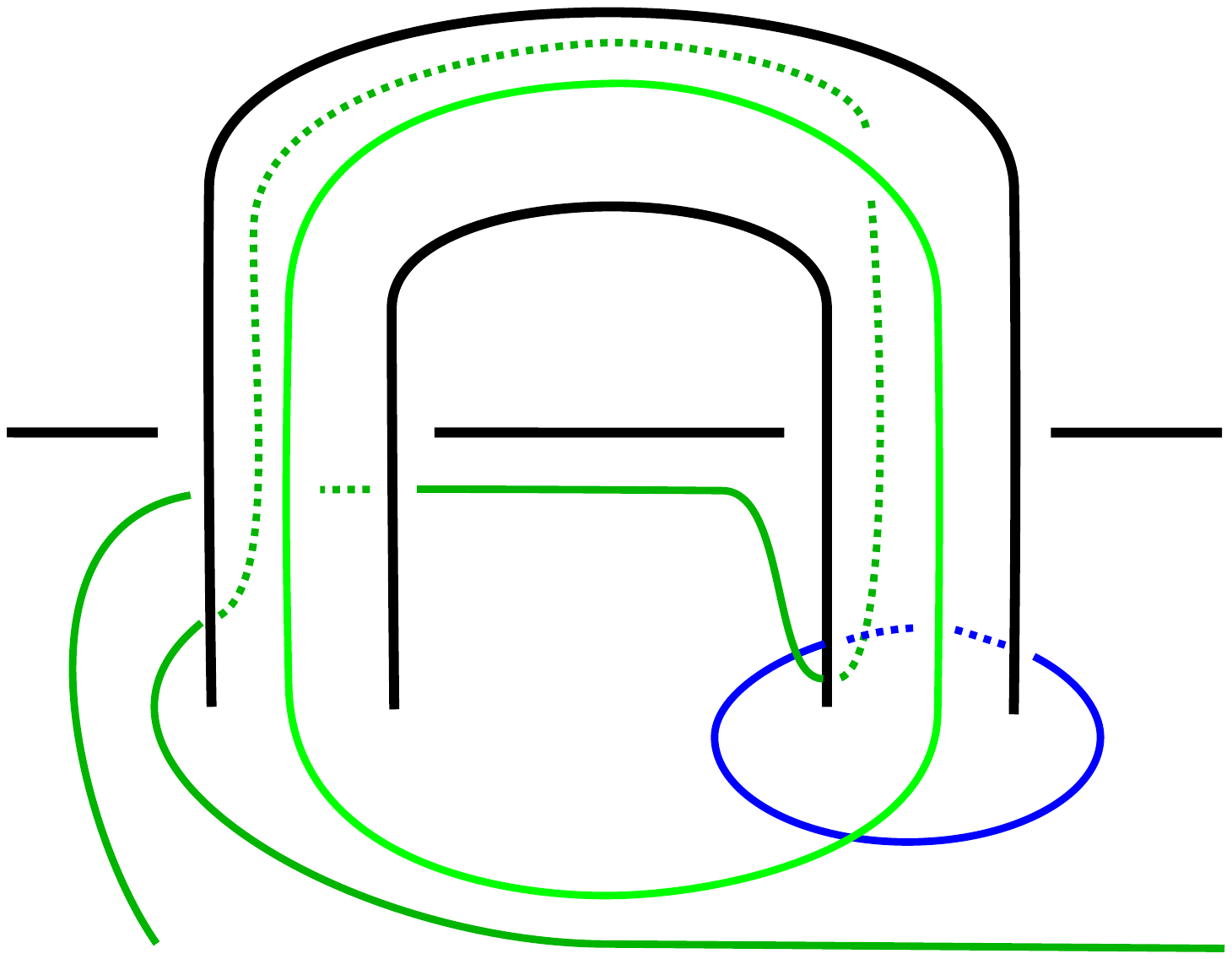 };
	\end{tikzpicture}

\stackrel{\text{isotopy}}{\rightsquigarrow}

	\begin{tikzpicture}[baseline={([yshift=-.5ex]current bounding box.center)}]
	\node at (0,0) {\def\svgscale{.12} 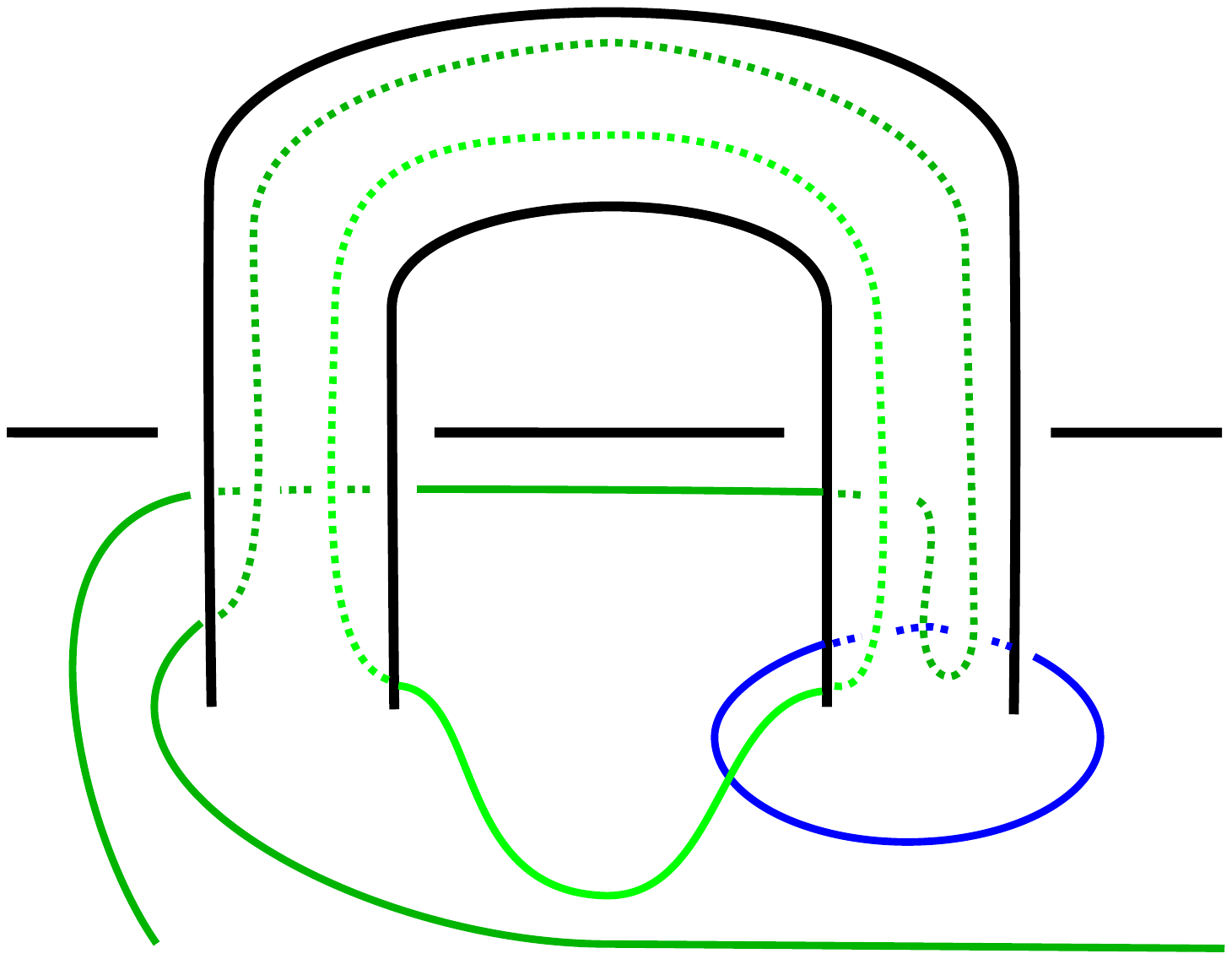 };
	\end{tikzpicture}
\\
\stackrel{\text{handle slide}}{\rightsquigarrow}
& 
	\begin{tikzpicture}[baseline={([yshift=-.5ex]current bounding box.center)}]
	\node at (0,0) {\def\svgscale{.12} 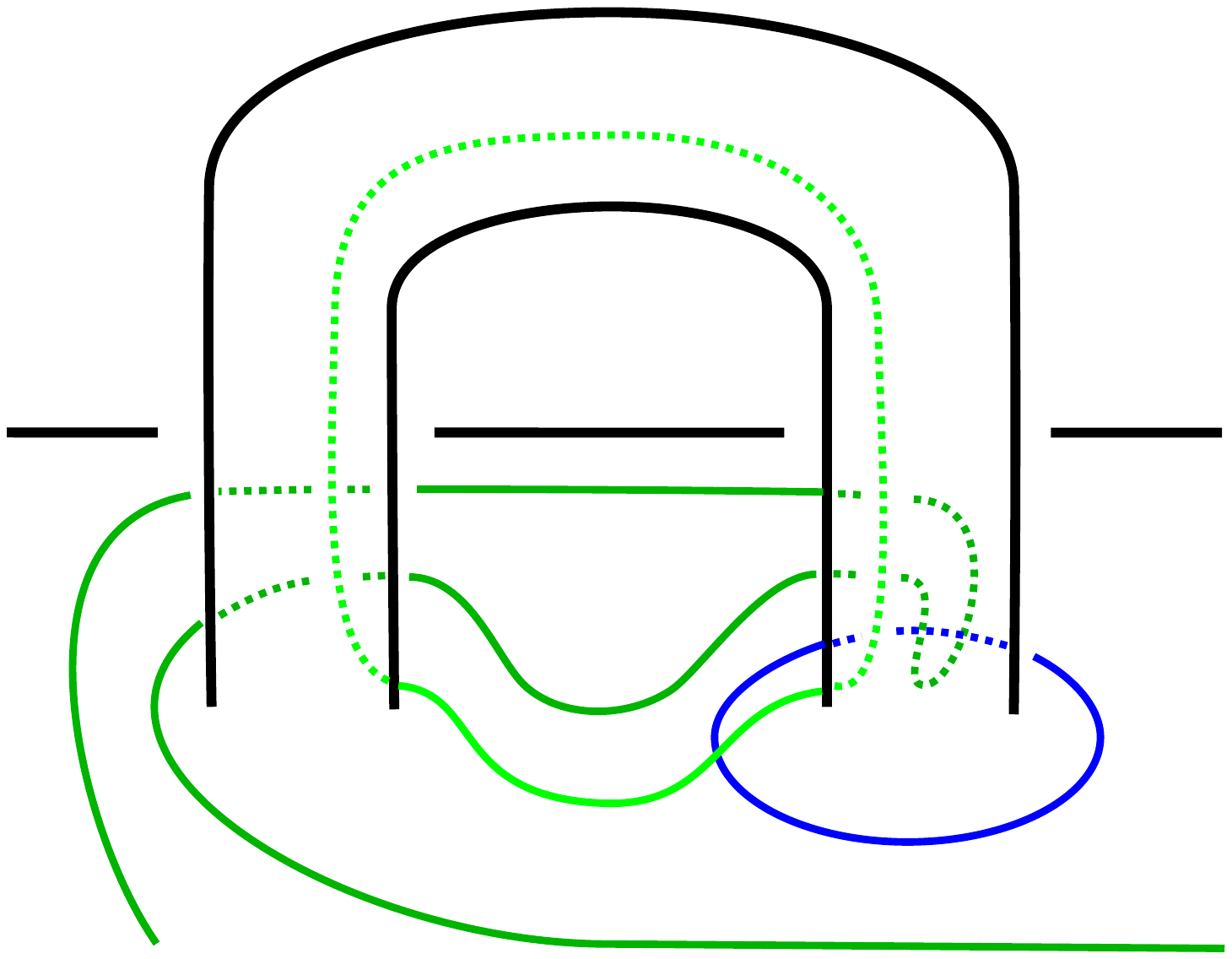 };
	\end{tikzpicture}

\stackrel{\text{R2}}{\rightsquigarrow}

	\begin{tikzpicture}[baseline={([yshift=-.5ex]current bounding box.center)}]
	\node at (0,0) {\def\svgscale{.12} 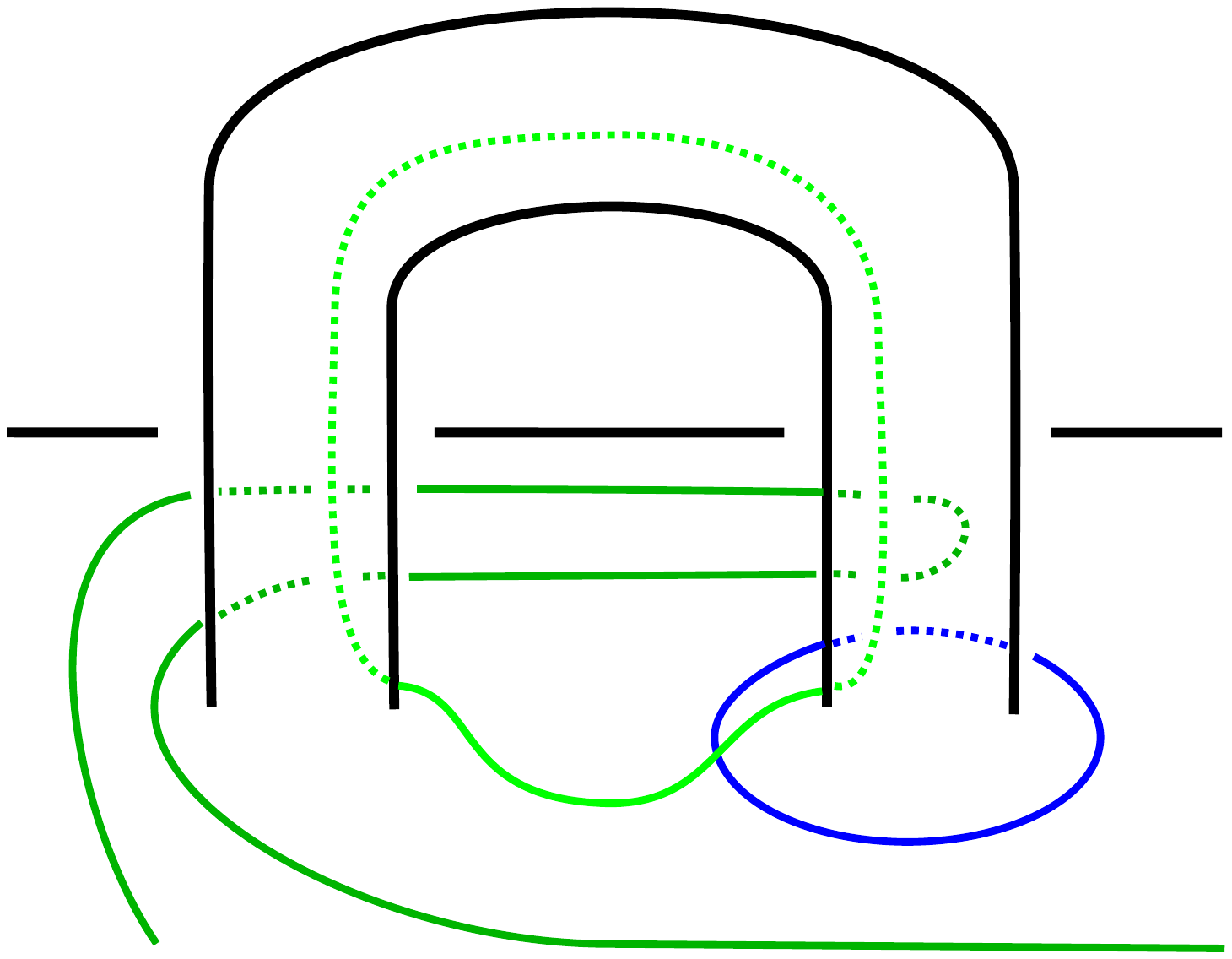 };
	\end{tikzpicture}

\stackrel{\text{isotopy}}{\rightsquigarrow}

	\begin{tikzpicture}[baseline={([yshift=-.5ex]current bounding box.center)}]
	\node at (0,0) {\def\svgscale{.12} 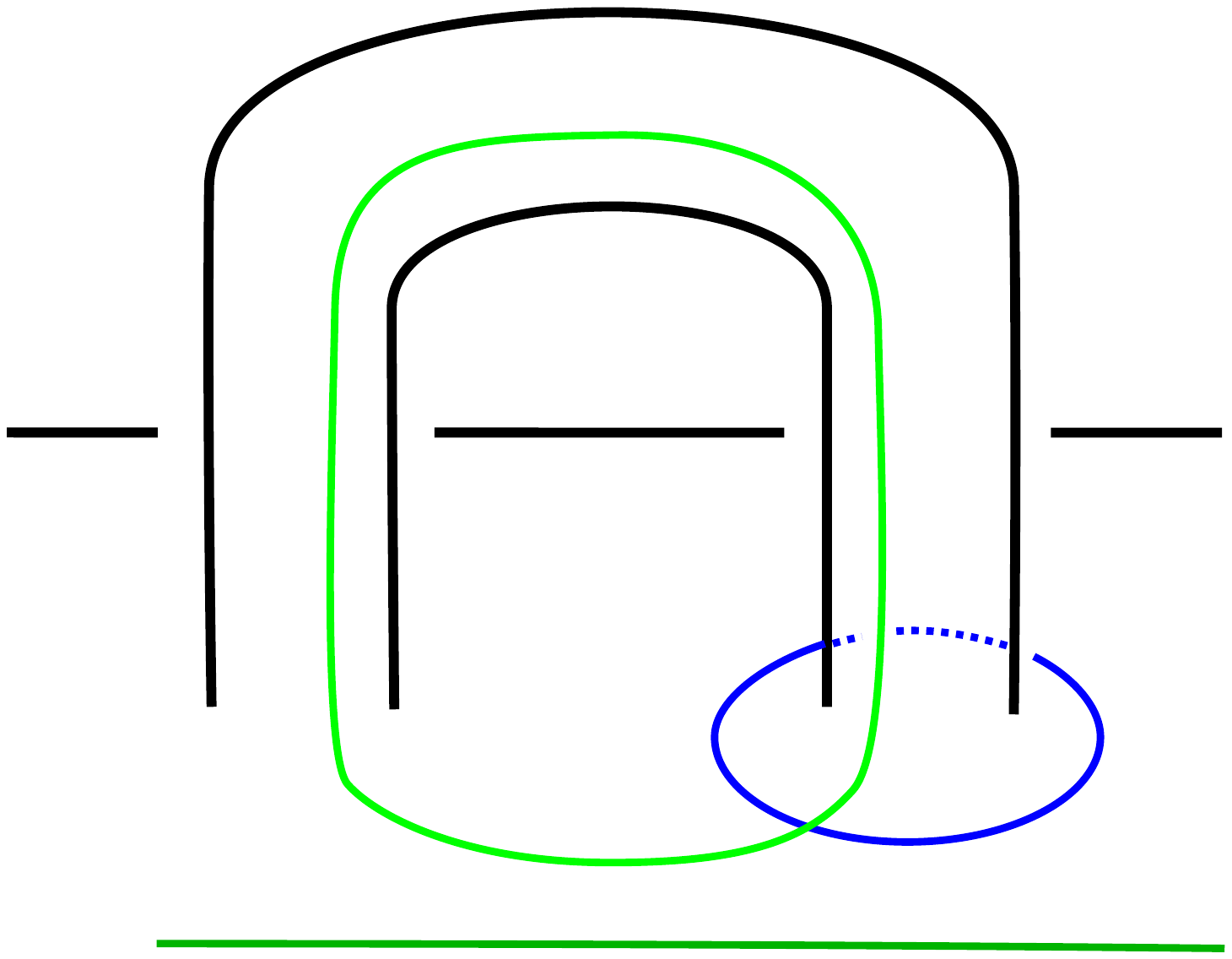 };
	\end{tikzpicture}

\label{eq:hardslide}
\end{align}

\begin{align}
\nonumber
& 
	\begin{tikzpicture}[baseline={([yshift=-.5ex]current bounding box.center)}]
	\node at (0,0) {\def\svgscale{.12} 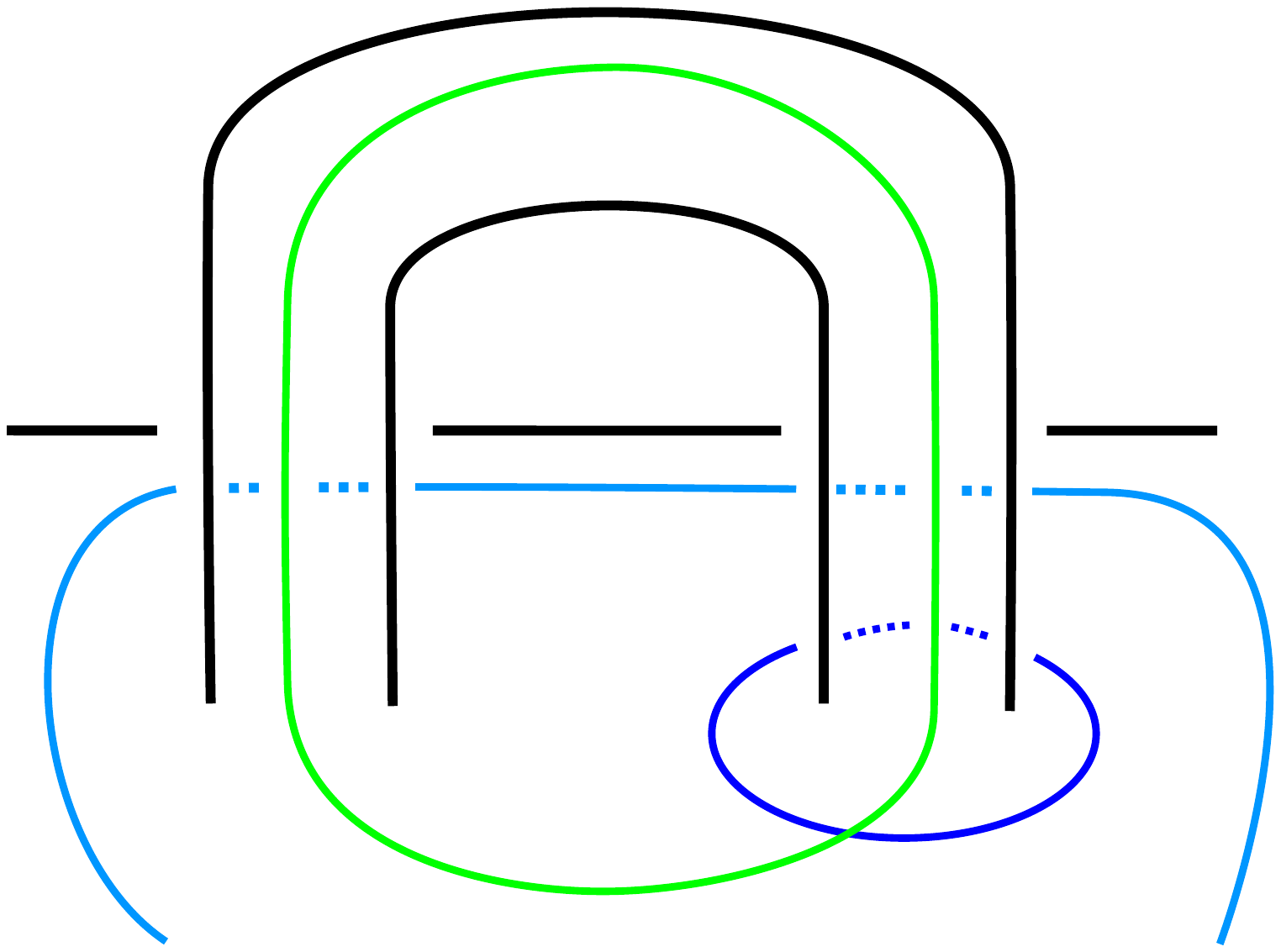 };
	\end{tikzpicture}

\stackrel{\text{handle slide}}{\rightsquigarrow}

	\begin{tikzpicture}[baseline={([yshift=-.5ex]current bounding box.center)}]
	\node at (0,0) {\def\svgscale{.12} 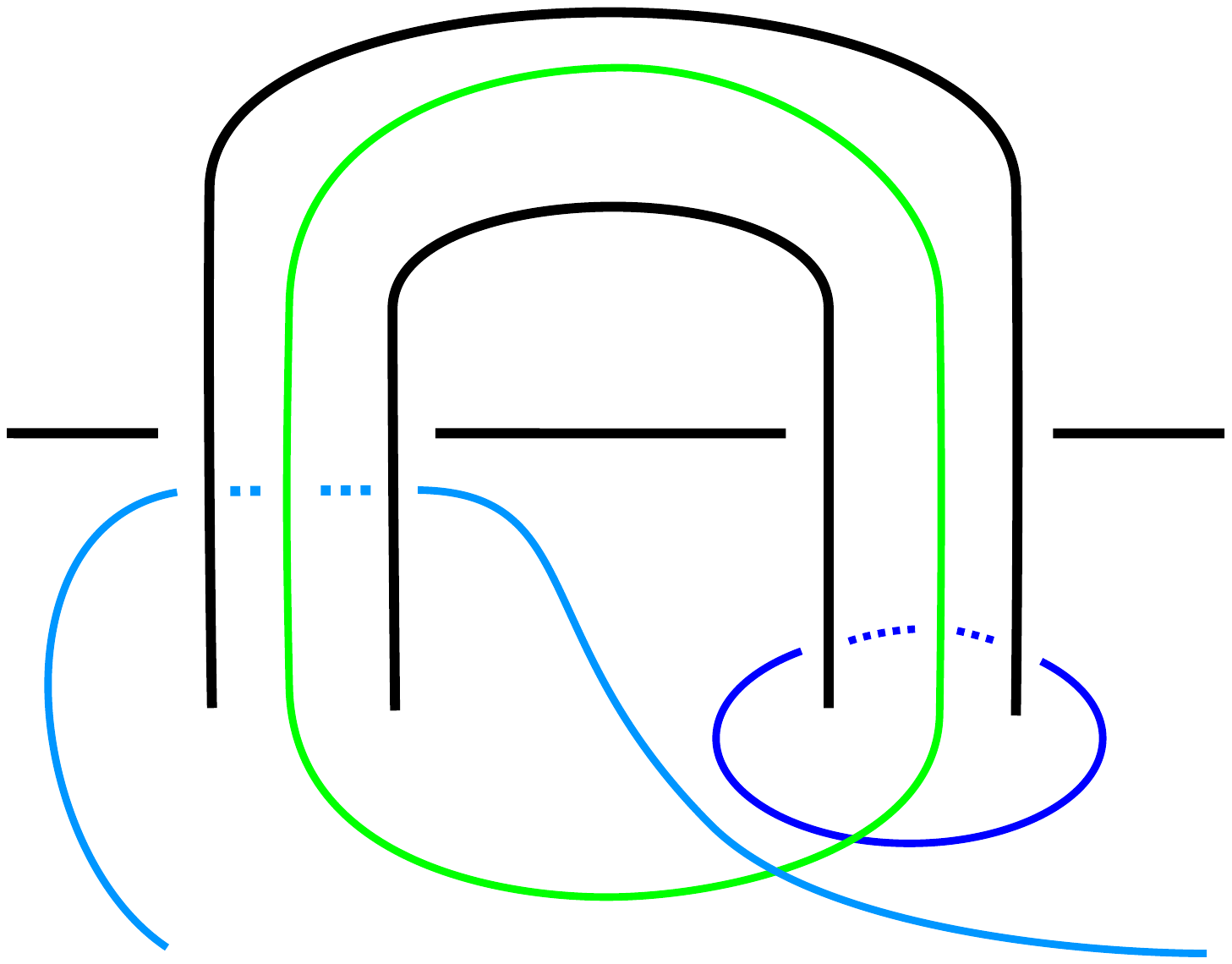 };
	\end{tikzpicture}

\stackrel{\text{isotopy}}{\rightsquigarrow}

	\begin{tikzpicture}[baseline={([yshift=-.5ex]current bounding box.center)}]
	\node at (0,0) {\def\svgscale{.12} 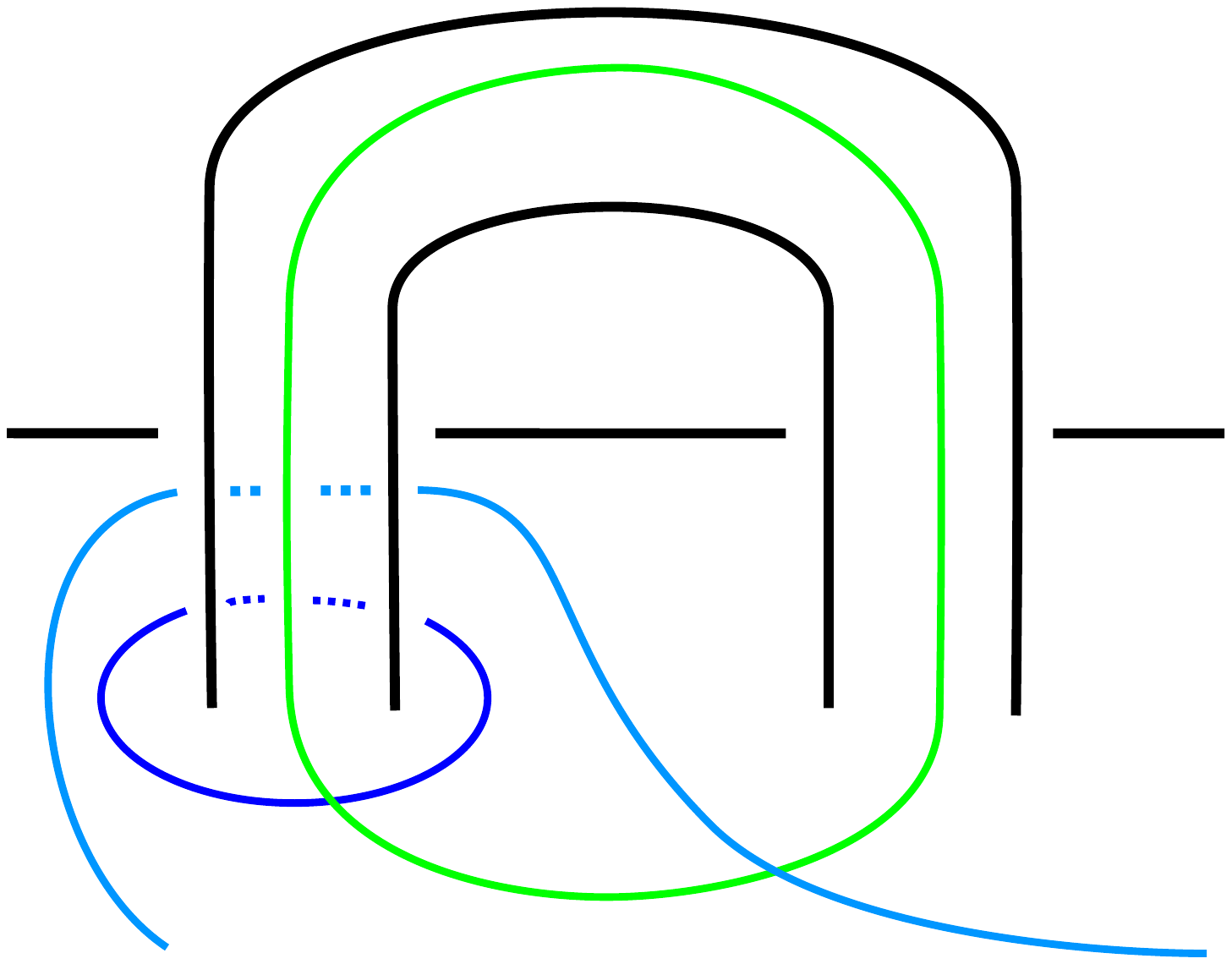 };
	\end{tikzpicture}
\\
\stackrel{\text{handle slide}}{\rightsquigarrow}
& 
	\begin{tikzpicture}[baseline={([yshift=-.5ex]current bounding box.center)}]
	\node at (0,0) {\def\svgscale{.12} 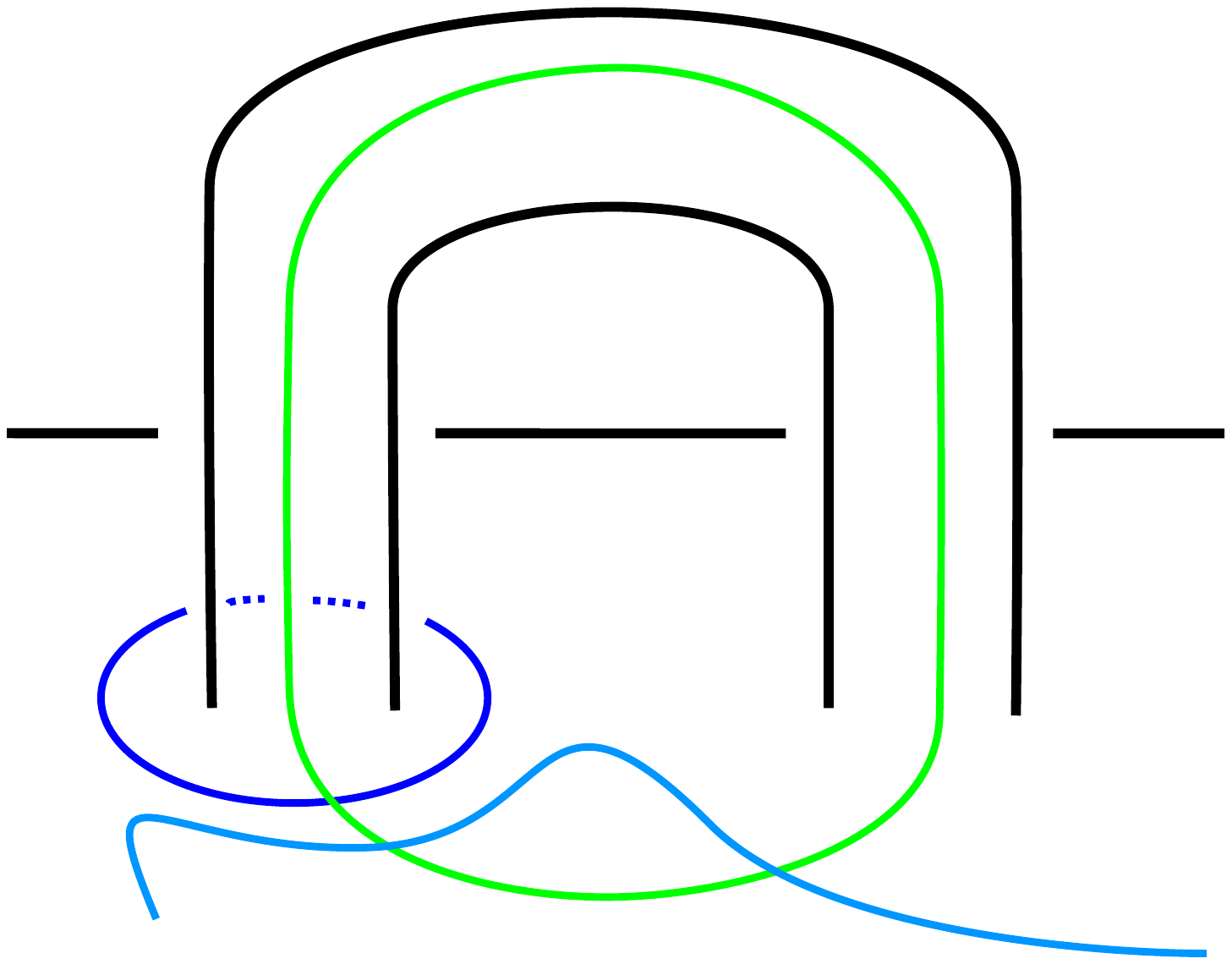 };
	\end{tikzpicture}

\stackrel{\text{R2}}{\rightsquigarrow}

	\begin{tikzpicture}[baseline={([yshift=-.5ex]current bounding box.center)}]
	\node at (0,0) {\def\svgscale{.12} 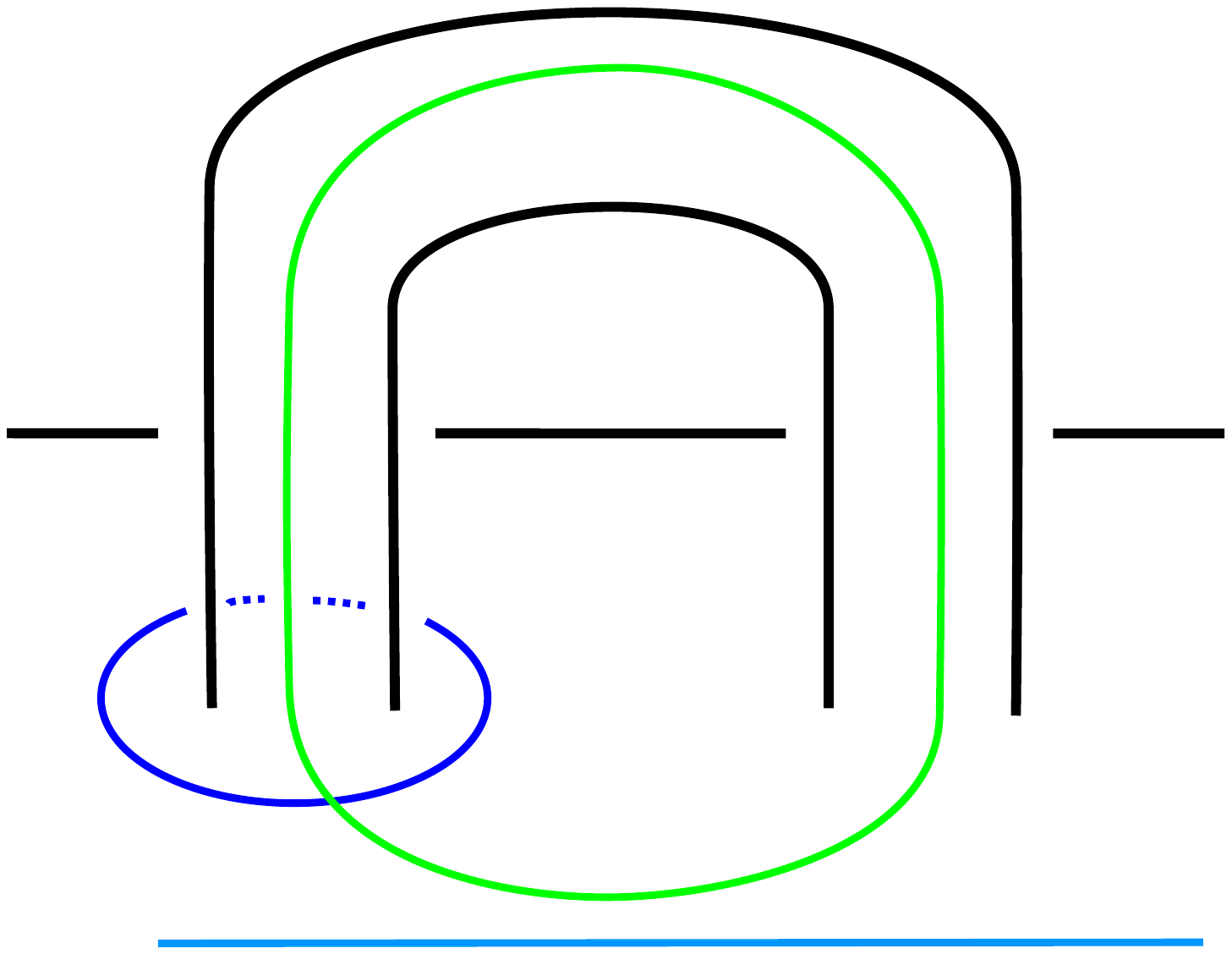 };
	\end{tikzpicture}

\stackrel{\text{isotopy}}{\rightsquigarrow}

	\begin{tikzpicture}[baseline={([yshift=-.5ex]current bounding box.center)}]
	\node at (0,0) {\def\svgscale{.12} 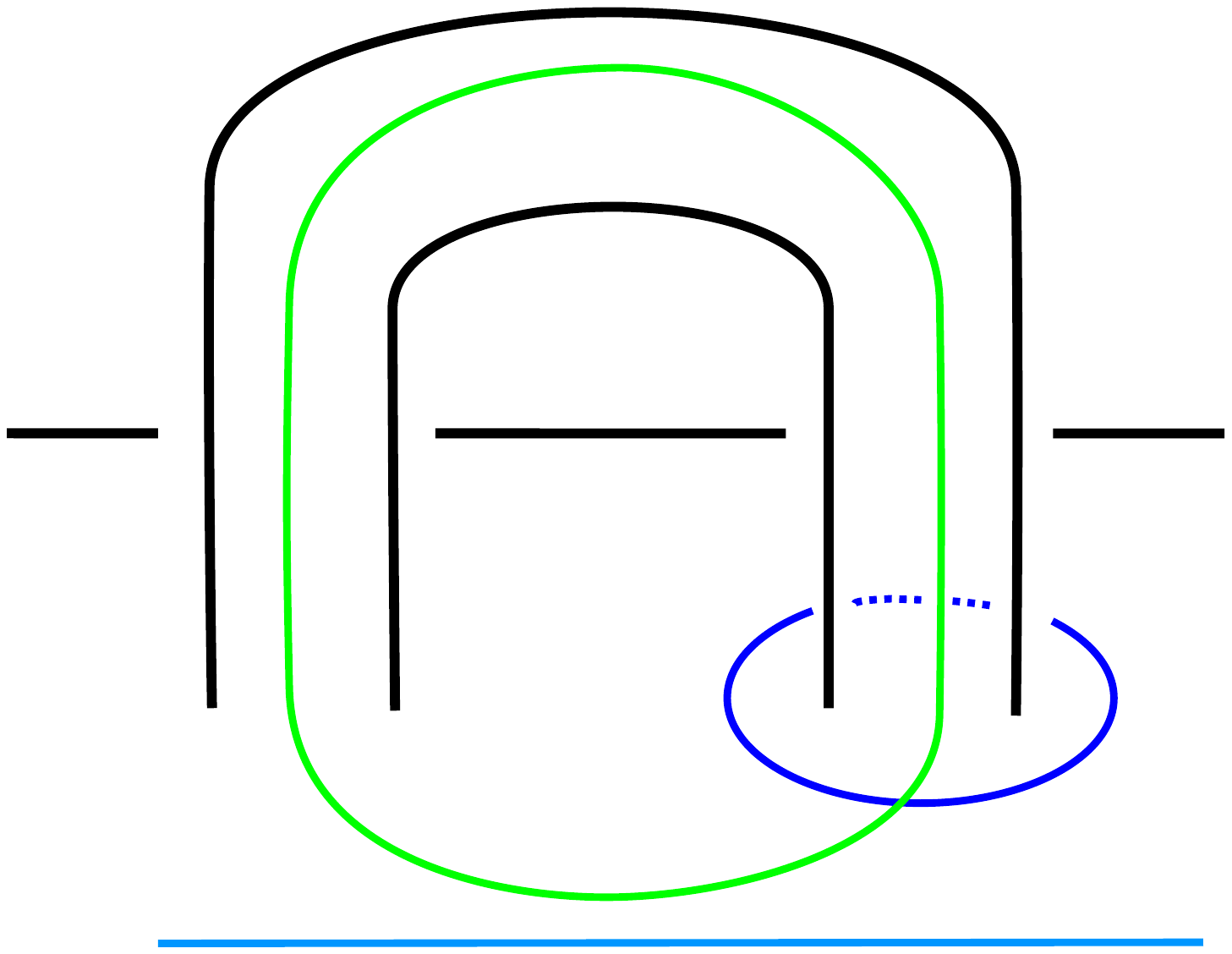 };
	\end{tikzpicture}

\label{eq:easyslide}
\end{align}
\caption{Steps in the computation \eqref{eq:stabproproof}. }
\label{fig:stabslide}
\end{figure}

The proof of Theorem \ref{th:gentrisection} shows that our categorical data is a rather natural framework for trisection diagrams. The three colours of the diagram are associated with the three spherical fusion categories, trisection diagrams can be evaluated in an existing diagrammatic calculus for this data, and invariance under most trisection moves from Theorem \ref{th:intmoves} is already built into this description.

In principle, one can also work with categorical data in which the three spherical fusion categories play a symmetric role.  
Without the condition that the $(\mac,\mad)$-bimodule category structure on $\calm$ is indecomposable, the categorical data for the trisection invariants essentially reduces to a $(\cala,\mad,\mac)$-trimodule category with a trimodule trace - a module category $\calm$ over  three spherical fusion categories $\cala$, $\mad$ and $\mac$ together with coherent isomorphisms that permute their actions and a trace that is compatible with all three module category structures. However, for this data it is not obvious how the averaged evaluation of a trisection diagram can be normalised to ensure invariance under the (de)stabilisation move, and one loses the simple characterisation from Corollary \ref{cor:stabilising}.

\section{Hopf algebraic from categorical trisection invariants}
\label{sec:trisecbimod}
In this section we show how the trisection invariants from \cite{CCC} in Definition \ref{def:cccinv}  arise as a special case of our  trisection invariants from Theorem \ref{th:gentrisection}, if the underlying field is $\C$. 
Roughly speaking, the three Hopf algebras  of a Hopf triplet $(A,B,C)$  in Definition \ref{def:cccinv} define spherical fusion categories 
\begin{align}\label{eq:tripccc}
\mac=C^{*op}\Mod, \qquad \mad=B^{*}\Mod,\qquad \cala=A^{*cop}\Mod.
\end{align}
 The category $\calm=\mathrm{Vect}_\C$ of finite-dimensional vector spaces becomes a $(\mac,\mad)$-bimodule category as in Example \ref{ex:vecmodcat}, and
the Hopf triplet defines a pivotal functor $\Phi:\cala\to\End_{(\mac,\mad)}(\calm)$. 

Conversely, each categorical datum $(\cala,\mad,\mac, \calm, \Phi)$  in Theorem \ref{th:gentrisection} where $\calm=\mathrm{Vect}_\C$ defines a Hopf triplet $(A,B,C)$ given by \eqref{eq:tripccc}, up to equivalence.  We  show that our 4-manifold invariants from Theorem \ref{th:gentrisection}  reduce to  the ones from Definition \ref{def:cccinv} for this data. 

Our proof proceeds in three steps. In Section \ref{subsec:reptheoret} we reformulate the trisection invariant for a Hopf triplet $(A,B,C)$ from Definition \ref{def:cccinv} in terms of representations of  $A^*$, $B^*$, $C^*$. In Section \ref{sec:tripletbimod} we show that a Hopf triplet $(A,B,C)$ defines a $(\mac,\mad)$-bimodule category structure on $\calm=\mathrm{Vect}_\C$ and a pivotal functor $\Phi: \cala\to \End_{(\mac,\mad)}(\calm)$. Conversely, each datum $(\cala,\mad,\mac,\calm,\Phi)$ with $\calm=\mathrm{Vect}_\C$ defines a Hopf triplet.
In Section \ref{subsec:trisecbimod} we  prove that the 4-manifold invariant from Theorem \ref{th:gentrisection} for this data reduces to the trisection invariant from \cite{CCC}.

\subsection{Representation theoretical formulation of the Hopf algebraic invariant}
\label{subsec:reptheoret}

To identify the trisection invariant  from \cite{CCC} 
as a special case of the 4-manifold invariant from Theorem \ref{th:gentrisection},  we formulate it in terms of Hopf algebra representations.

For this we interpret the skew pairings $\tau: A\times B\to \C$ from Definition \ref{def:skewpair} as elements $\tau\in A^*\oo B^*$. 
We use the symbolic notation
$\rho=\low\rho 1\oo \low \rho 2\in A^*\oo B^*$ and $\mu=\low\mu 1\oo\low\mu 2\in B^*\oo A^*$ for elements
of $A^*\oo B^*$ and $B^*\oo A^*$.
For a multiple tensor product  $X=\ldots\oo A^*\oo\ldots\oo B^*\oo \ldots$ of algebras with $A^*$ in the $j$th and $B^*$ in the $k$th component, we write
$\rho^{jk}\in X$ and $\mu^{kj}\in X$ for the elements obtained by inserting $\rho$ and $\mu$ in the $k$th and $j$th components and placing units in all other components.

A direct computation using Hopf algebra duality and Identities \eqref{eq:tauforms}, \eqref{eq:convinv} and \eqref{eq:idd} then reformulates the Hopf doublets and Hopf triplets from Section \ref{sec:hopfdoubtrip} as follows.

\begin{Lemma} \label{lem:hopfdoublettrans} Let $A$, $B$, $C$ be finite-dimensional semisimple complex Hopf algebras.
\begin{compactenum}
\item  An element $\tau\in A^*\oo B^*$ defines a skew pairing  $\tau: A\times B\to \C$ if and only if 
\begin{align}\label{eq:tenstau}
&(\Delta_{A^*}\oo \id)(\tau)=\tau^{13}\cdot \tau^{23} & &(\id\oo\Delta_{B^*})(\tau)=\tau^{13}\cdot \tau^{12}\\
&(\epsilon_{A^*}\oo \id)\tau=1_{B^*}  & &(\id\oo \epsilon_{B^*})\tau=1_{A^*}.\nonumber
\end{align}
In this case, $\tau$ is invertible with inverse  $\tau^\inv=(S_{A^*}\oo \id)(\tau)=(\id\oo S_{B^*})(\tau)$.\\[-1ex]

\item Suppose that  $\tau_{AB}\in A^*\oo B^*$, $\tau_{BC}\in B^*\oo C^*$ and $\tau_{CA}\in C^*\oo A^*$ define skew pairings. Then they define a Hopf triplet $(A,B,C)$ if and only if
\begin{align}\label{eq:tripletYB}
\tau_{BC}^{12}\cdot \tau_{AB}^{31}\cdot \tau^{23}_{CA}=\tau_{CA}^{23}\cdot \tau_{AB}^{31}\cdot \tau_{BC}^{12}\in B^*\oo C^{*op}\oo A^*.
\end{align}
\end{compactenum}
\end{Lemma}

\begin{Remark} As they describe a 2-cocycle, Equations \eqref{eq:tenstau} and \eqref{eq:tripletYB} generalise the conditions on a universal $R$-matrix. In particular,
 Equation \eqref{eq:tripletYB} can be viewed as a generalised QYBE. In terms of the multiplication of the algebra $A^*\oo B^*\oo C^*$ and in the symbolic notation introduced above it reads
\begin{align}\label{eq:tripletYB2}
&(\tau_{AB(1)}\cdot \tau_{CA(2)})\oo (\tau_{BC(1)}\cdot \tau_{AB(2)})\oo (\tau_{CA(1)}\cdot \tau_{BC(2)})\\
&=(\tau_{CA(2)}\cdot \tau_{AB(1)})\oo (\tau_{AB(2)}\cdot \tau_{BC(1)})\oo (\tau_{BC(2)}\cdot \tau_{CA(1)})\in A^*\oo B^*\oo C^*.\nonumber
\end{align}
\end{Remark}

We now consider a Hopf triplet $(A,B,C)$ given by  skew pairings $\tau_{AB}$, $\tau_{BC}$, $\tau_{CA}$. We fix sets of representatives of the isomorphism classes of simple representations of $A^*$, $B^*$, $C^*$
\begin{align}\label{eq:repsabc}
\{\rho_i^A: A^*\to \End_\C(V_i^A)\}_{i\in I^A},\quad \{\rho_j^B: B^*\to \End_\C(V_j^B)\}_{j\in I^B},\quad \{\rho_k^C: C^*\to \End_\C(V_k^C)\}_{k\in I^C}
\end{align}
and set $I=I^A\cup I^B\cup I^C$.
For a trisection diagram $T$ we denote again by  $R,B,G$  its sets of red, blue and green curves and set $\Gamma=R\cup B\cup G$. We assign to the red, blue and green curves the Hopf algebras $A$, $B$, $C$,  as in Section \ref{sec:cccinv}, and denote by  $X_\lambda\in\{A,B,C\}$  the Hopf algebra associated to a curve $\lambda\in \Gamma$, such that $X_\lambda=C$, $X_\lambda=B$ and $X_\lambda=A$ for $\lambda\in G$, $\lambda\in B$ and $\lambda\in R$, respectively.

\begin{Definition}Let $(A,B,C)$ be a Hopf triplet and  $I^A$, $I^B$ and $I^C$ as in \eqref{eq:repsabc}. A \textbf{representation labelling} of a trisection diagram $T$ is a map $h: \Gamma\to I$ with $h(\lambda)\in I^{X_\lambda}$ for all $\lambda\in \Gamma$.
\end{Definition}

\begin{Proposition} \label{prop:repthtrisec}
Let $(A,B,C)$ be a Hopf triplet and   $I^A$, $I^B$ and $I^C$ as in \eqref{eq:repsabc}.
Then the following yields the trisection bracket $\langle T\rangle_{A,B,C}$ from Definition \ref{def:trisecbr}, up to a rescaling of the integrals:
\\[-3ex]
\begin{compactenum}
\item Equip each  curve of a trisection diagram $T$ with an orientation and a basepoint.\\[-2ex]

\item Each  representation labelling $h: \Gamma\to I$ of $T$ defines a vector space $V_{h,T}=\bigotimes_{\lambda\in \Gamma} V^{X_\lambda}_{h(\lambda)}$ and an endomorphism $\phi_{h,T}: V_{h,T}\to V_{h,T}$ constructed as follows:\\[-2ex]

\begin{compactitem}
\item Assign to an intersection point $p\in \lambda\cap \mu$ in $T$
of colours $(X_\lambda,X_\mu)\in\{(A,B), (B,C), (C,A)\}$ and
 with sign $\epsilon_p\in\{\pm 1\}$ as in \eqref{eq:intsign} the element  $\tau_{X_\lambda X_\mu}^{\epsilon_p}\in X_\lambda^*\oo X_\mu^*$.\\[-2ex]

\item Let the elements  $\tau_{X_\lambda X_\mu}^{\pm 1}\in X_\lambda^*\oo X_\mu^*$ act on the factors $V^{X_\lambda}_{h(\lambda)}$ and $V^{X_\mu}_{h(\mu)}$ in $V_{h,T}$,  in the order  given by the basepoints and the orientations of $\lambda$ and $\mu$. \\[-2ex]
\end{compactitem}

\item Multiply the endomorphism $\phi_{h,T}: V_{h,T}\to V_{h,T}$ with a factor  $\dim_\C(V_{h(\lambda)}^{X_\lambda})$ for each curve $\lambda\in \Gamma$, take its trace,  and sum over all representation labellings $h: \Gamma\to I$. 
\end{compactenum}
This gives the following formula
\begin{align}\label{eq:trisecbr}
\langle T\rangle_{A,B,C}=\sum_{h: \Gamma \to I} \prod_{\lambda\in \Gamma} \dim_\C(V^{X_\lambda}_{h(\lambda)}) \cdot \tr(\phi_{h,T}).
\end{align}
\end{Proposition}

\begin{Proof} Let $T$ be a trisection diagram of genus $g$. To show that the right-hand-side of \eqref{eq:trisecbr} coincides with the trisection bracket from Definition \ref{def:trisecbr}, we 
consider the algebra 
$$
H_T=B^{*\oo g}\oo (C^{*op})^{\oo g}\oo (A^{*})^{\oo g},
$$
where the different copies of $B^*$, $C^{*op}$ and $A^{*}$ are assigned to the blue, green and red curves of $T$.  
Inserting the element  $\tau_{X_\lambda X_\mu}^{\epsilon_p}$ for an intersection point $p\in \lambda\cap \mu$ into the factors of $H_T$ for $\lambda$ and $\mu$ and inserting units in the other factors yields an element $t_p\in H_T$. This element corresponds to the pairings $\tau_{X_\lambda X_\mu}^{\epsilon_p}$ in step 5 of Definition \ref{def:trisecbr}.  
Applying the coproducts in step 4 in Definition \ref{def:trisecbr} and combining the pairings in step 5 corresponds to arranging the elements $t_p\in H_T$ on each factor of $H_T$ according to the  ordering   defined by the basepoints and orientations of the associated curve and then multiplying them.  This yields an element 
$$
\phi_T=\prod_{p\in T}  t_p \in  H_T.
$$ 
The trisection bracket from Definition \ref{def:trisecbr} is then obtained by pairing the element $\phi_T\in H_T$ with the element $\lambda_T=\ell_B^{\oo g}\oo \ell_C^{\oo g}\oo \ell_A^{\oo g}\in H^*_T$, where $\ell_A$, $\ell_B$ and $\ell_C$ are normalisable integrals of $A$, $B$, $C$. 

To relate this to the right-hand-side of formula \eqref{eq:trisecbr} we note that the algebra $H_T$ has a Hopf algebra structure, namely as the $g$-fold tensor product $H_T=H^{\oo g}$ of the Hopf algebra 
$$
H=D(D(B^{op}, C^{cop}), A^{op})^{*}.$$
By the remarks after Lemma \ref{lem:hopfconc} we have $H=B^{*}\oo C^{*op}\oo A^*$ \emph{as an algebra}, but the comultiplication is twisted with a 2-cocycle defined by the skew pairings. As the antipode identifies simple representations of $C^*$ and $C^{*op}$, simple representations of $H$ are in bijection with triples of simple representations of $B^*$, $C^*$ and $A^*$. It follows that isomorphism classes of simple representations of $H_T$ are in bijection with representation labellings $h:\Gamma\to I$, and the associated representation  is on  $V_{h,T}$.  

By construction, the endomorphism $\phi_{h,T}: V_{h,T}\to V_{h,T}$ in Proposition \ref{prop:repthtrisec} is the action of the element $\phi_T\in H_T$ in the simple representation $\rho_h: H_T\to \End_\C(V_{h,T})$ defined by the representation labelling $h: \Gamma\to I$. Thus, the right-hand side of \eqref{eq:trisecbr} can be rewritten as the following  sum over a set of representatives of  simple representations of $H_T$
\begin{align}\label{eq:repphelp}
\sum_{\rho_h: H_T\to \End_\C(V_{h,T})} \!\!\!\!\!\dim_{V_{h,T}}\cdot \tr (\rho_h(\phi_T)).
\end{align}
By   \eqref{eq:Haarintdualrep} the element in \eqref{eq:repphelp} is equal to  $\lambda'_T(\phi_T)$, where $\lambda'_T=\ell^{\oo g}\in H^*_T$ and $\ell$ is a normalisable integral of $H^*=D(D(B^{op}, C^{cop}), A^{op})$. Lemma \ref{lem:intdoub} implies that $\ell'=\ell_B\oo \ell_C\oo\ell_A$ is another normalisable integral of $H^*$ for all normalisable integrals $\ell_A\in A$, $\ell_B\in B$ and $\ell_C\in C$. Hence, $\ell$ and $\ell'$ are related by a rescaling of the integrals of $A$, $B$ and $C$, and the right-hand side of \eqref{eq:trisecbr} is a trisection bracket.
\end{Proof}

\subsection{Description of Hopf triplets in terms of bimodule categories}
\label{sec:tripletbimod}

To show that the trisection invariants from \cite{CCC} are a special case of the trisection invariants from Theorem \ref{th:gentrisection} we  relate their algebraic data. 
We first show that every Hopf triplet $(A,B,C)$ defines spherical fusion categories $\cala,\mad,\mac$, a $(\mac,\mad)$-bimodule category structure on $\calm=\mathrm{Vect}_\C$ and 
a pivotal functor $\Phi:\cala\to\End_{(\mac,\mad)}(\calm)$, as in Definition \ref{def:diagevalav}.

For each finite-dimensional complex semisimple Hopf algebra $H$, we equip its representation category $H\Mod$ with the standard pivotal structure induced by the pivotal structure on $\mathrm{Vect}_\C$. This gives  $H\Mod$ the structure of a spherical fusion category such that all traces are given by the traces of $\mathrm{Vect}_\C$.
We use the tensor equivalence $H^{cop}\Mod\cong (H\Mod)^{rev}$ to
 identify right module categories over  $H\Mod$  with left module categories over $H^{cop}\Mod$.

\begin{Proposition} \label{prop:hopfmod}
Let $(A,B,C)$ be a Hopf triplet and 
$$\cala=A^{*cop}\Mod, \qquad \mad=B^*\Mod,\qquad \mac=C^{*op}\Mod.$$ 
\begin{compactenum}
\item The skew pairing $\tau_{BC}$ defines a
$(\mac,\mad)$-bimodule category structure on $\mathrm{Vect}_\C$ given by tensoring on the left, by the associator in $\mathrm{Vect}_\C$ and the coherence isomorphism $q:\lhd(\rhd \times\id_\mad)\Rightarrow \rhd(\id_\mac\times\lhd)$  
from Definition \ref{def:modulecat} given on a
 $B^{*cop}$-module $N$, $C^{*op}$-module $P$ and  vector space $V$ by
\begin{align}\label{eq:bdef}
q_{P,V,N}: N\oo P\oo V\to P\oo N\oo V, \quad n\oo p\oo v\mapsto (\tau_{BC(2)}\rhd p)\oo (\tau_{BC(1)}\rhd n)\oo v. 
\end{align}
\item This defines a  $D(B^{op}, C^{cop})^{*}\Mod$-module category structure on $\mathrm{Vect}_\C$.\\[-2ex]

\item There is a pivotal functor $\Phi:\cala\to \End_{(\mac,\mad)}(\mathrm{Vect}_\C)$ induced by the Hopf algebra homomorphism 
\begin{align}\label{eq:pivfunc}
\phi: D(B^{op}, C^{cop})\to A^*,\quad \phi(b\oo c)(a)=\tau_{AB}(\low a 1 ,b)\tau_{CA}(c,\low a 2).
\end{align}
It assigns to an $A^{*cop}$-module $M$ the  bimodule endofunctor $M\oo -:\mathrm{Vect}_\C\to \mathrm{Vect_\C}$ with coherence isomorphisms given by the maps
\begin{align}\label{eq:bimodfuncstruc}
&s^M_{P,V}: M\oo P\oo V\to P\oo M\oo V,\quad m\oo p\oo v\mapsto (\tau_{CA(1)}^{\inv}\rhd p)\oo(\tau_{CA(2)}^{\inv}\rhd m)\oo v\\
&t^M_{V,N}: N\oo M\oo V\to M\oo N\oo V, \quad n\oo m\oo v\mapsto (\tau_{AB(1)}\rhd m)\oo(\tau_{AB(2)}\rhd n)\oo v\nonumber.
\end{align}
\end{compactenum}
\end{Proposition}

\begin{Proof} 1.~We identify $B^*\Mod$-right module categories with $B^{*cop}\Mod$-left module categories. 
The associators in $\mathrm{Vect}_\C$ define coherence isomorphisms $l$ and $r$ in Definition \ref{def:modulecat} that are natural in all arguments and make the diagrams \eqref{eq:pentagoncdef} to \eqref{eq:triangled} commute.

The isomorphisms $q_{P,V,N}$ from \eqref{eq:bdef} are  natural in all arguments. It remains to show that they  satisfy the coherence axioms \eqref{eq:pentagonbdef} for a bimodule category. As  $l$ and $r$ in \eqref{eq:pentagonbdef} are given by the  associators in $\mathrm{Vect}_\C$, this reduces to the identities
\begin{align}\label{eq:bids}
&Q_{N\oo N',P}=(Q_{N,P}\oo \id)\circ (\id\oo Q_{N',P})   &
&Q_{N,P\oo P'}=(\id\oo Q_{N,P'})\circ (Q_{N,P}\oo\id),
\end{align}
for $C^{*op}$-modules $P$, $P'$ and  $B^{*cop}$-modules $N$, $N'$ and the map 
\begin{align}\label{eq:bnp}
Q_{N,P}=q_{P,\C,N}: N\oo P\to P\oo N, \quad n\oo p\mapsto (\tau_{BC(2)}\rhd p)\oo(\tau_{BC(1)}\rhd n).
\end{align}
A direct computation yields
\begin{align*}
Q_{N\oo N', P}(n\oo n'\oo p)&=(\tau_{BC(2)}\rhd p)\oo (\tau_{BC(1)(2)}\rhd n)\oo (\tau_{BC(1)(1)}\rhd n')\\
&=[\sigma_{23}\circ \sigma_{12}\circ \sigma_{23}((\Delta\oo\id)(\tau_{BC}))]\rhd(p\oo n \oo n')\\
(Q_{N,P}\oo \id)\circ (\id\oo Q_{N',P})(n\oo n'\oo p)&=( (\tau'_{BC(2)}\cdot \tau_{BC(2)})\rhd p)\oo (\tau_{BC(1)}\rhd n)\oo (\tau'_{BC(1)}\rhd n')\\
&=[\sigma_{23}\circ \sigma_{12}\circ \sigma_{23} ( \tau_{BC}^{13}\cdot \tau_{BC}^{12})]\rhd (p\oo n\oo n')
\\
Q_{N, P\oo P'}(n\oo p\oo p')
&=(\tau_{BC(2)(1)}\rhd p)\oo (\tau_{BC(2)(2)}\rhd p')\oo (\tau_{BC(1)}\rhd n)\\
&=[\sigma_{23}\circ \sigma_{12} ( (\id\oo \Delta)(\tau_{BC}))]\rhd(p \oo p'\oo n)\\
(\id\oo Q_{N,P'})\circ (Q_{N,P}\oo\id)(n\oo p\oo p')
&=( \tau_{BC(2)}\rhd p)\oo (\tau'_{BC(2)}\rhd p')\oo ( (\tau'_{BC(1)}\cdot \tau_{BC(1)})\rhd n)
\\
&=[\sigma_{23}\circ \sigma_{12} (\tau_{BC}^{13}\cdot \tau_{BC}^{12})]\rhd ( p\oo p'\oo n),
\end{align*}
where $\sigma_{ij}$ is the linear map that exchanges the $i$th and $j$th factor in a tensor product  and all (co)products  are the ones of $B^*$ and $C^*$.  Identity \eqref{eq:tenstau} then gives \eqref{eq:bids}.

2.~This follows from the tensor equivalences
$$
 (B^{*cop}\Mod)\boxtimes (C^{*op}\Mod)\cong (B^{*cop}\oo C^{*op})\Mod \cong D(B^{op}, C^{cop})^*\Mod,
$$
where the last equivalence follows, because $D(B^{op}, C^{cop})$ is obtained from $B^{op}\oo C^{cop}$ by twisting. More explicitly, we have $D(B^{op}, C^{cop})^*\cong B^{*cop}\oo C^{*op}$ as an algebra.  The formula for the multiplication in Proposition \ref{prop:2cocyc} implies that
the comultiplication of $D(B^{op}, C^{cop})^*$ is given by
\begin{align*}
\Delta_{D(B^{op}, C^{cop})^*} (\beta\oo \gamma)
&=(1\oo \tau_{BC(2)}\oo \tau_{BC(1)}\oo 1)\cdot (\low\beta 2\oo \low\gamma 1\oo\low\beta 1\oo\low\gamma 2)\cdot (1\oo \tau_{BC(2)}^\inv\oo \tau_{BC(1)}^\inv\oo 1).
\end{align*}
Hence, this equivalence is given by the identity functor and the  coherence isomorphisms 
$$
\phi_{V, W}: V\oo W\to V\oo W, \quad v\oo w\mapsto  ((1\oo \tau_{BC(2)})\rhd v)\oo ((\tau_{BC(1)}\oo 1)\rhd w).
$$ 
3.~The map $\phi$ is  a Hopf algebra homomorphism by Definition \ref{def:triplet} of a Hopf triplet.  It induces a tensor functor
$\Phi: A^{*cop}\Mod\to D(B^{op}, C^{cop})^{cop}\Mod$.

By Example \ref{ex:hmodmodule} the fusion category $D(B^{op}, C^{cop})^{cop}\Mod$ is tensor equivalent to
$\End_{H\Mod}(\mathrm{Vect}_\C)$, where $H=D(B^{op}, C^{cop})^*$.  
Using the  formulas in Example \ref{ex:hmodmodule} we then find that the $(\mac,\mad)$-bimodule functor assigned to  an $A^{*cop}$-module $M$ is the functor 
 $M\oo -:\mathrm{Vect}_\C\to\mathrm{Vect}_\C$ with the coherence isomorphisms from \eqref{eq:bimodfuncstruc}. They are determined by the isomorphisms
 \begin{align}\label{eq:cromaps}
 &S_{M,P}=s^M_{P,\C}: M\oo P\to P\oo M, \; m\oo p\mapsto (\tau_{CA(1)}^{\inv}\rhd p)\oo (\tau_{CA(2)}^{\inv}\rhd m)\\&T_{N,M}=t^M_{\C,N}: N\oo M\to M\oo N, \; n\oo m\mapsto (\tau_{AB(1)}\rhd m)\oo (\tau_{AB(2)}\rhd n).\nonumber\end{align}
\end{Proof}

Proposition \ref{prop:hopfmod} shows that any Hopf triplet  defines a triple $(\cala,\mad,\mac)$ of spherical fusion categories, a $(\mac,\mad)$-bimodule category structure on $\mathrm{Vect}_\C$ such that the trace of $\mathrm{Vect}_\C$ is a $(\mac,\mad)$-bimodule trace and a pivotal functor $\Phi:\cala\to \End_{(\mac,\mad)}(\mathrm{Vect}_\C)$. The last two statements follow, because the pivotal structure and thus all traces are given by the pivotal structure of $\mathrm{Vect}_\C$.  We  now show  that the converse statement also holds. 

For this, recall
that a  \textbf{fibre functor} for a fusion category $\mac$  is a tensor functor $F:\mac\to\mathrm{Vect}_\C$. It defines a  Hopf algebra structure on the set $\End(F)$ of endo-natural transformations of $F$. Its algebra structure is given by   the composition of endo-natural transformations and the identity natural transformation.
Its coalgebra structure is induced by the monoidal structure of $F$, and the antipode by the duals of $\mac$.
 
By \cite[Th 5.3.12]{EGNO}  this defines a bijection between pairs $(\mac,F)$ of a fusion category $\mac$ and a fibre functor $F:\mac\to \mathrm{Vect}_\C$, up to tensor equivalence and isomorphism of fibre functors, and finite-dimensional semisimple complex Hopf algebras, up to Hopf algebra isomorphisms. This equivalence  assigns to a finite-dimensional semisimple Hopf algebra $H$ the fusion category $\mac=H\Mod$ together with the forgetful functor $F:H\Mod\to\mathrm{Vect}_\C$.

\begin{Proposition}\label{prop:reducetrip} 
Any datum $(\cala, \mad,\mac,\calm=\mathrm{Vect}_\C,\Phi)$ as in Definition  \ref{def:diagevalav}  defines a Hopf triplet.
\end{Proposition}

\begin{Proof}
By \cite[Ex 7.4.6]{EGNO} for any fusion category $\mac$ the $\mac$-module category structures on $\mathrm{Vect}_\C$ are in bijection with fibre functors on $\mac$. More precisely, there is a bijection between $\mac$-module category structures on $\mathrm{Vect}_\C$  and tensor functors $G_\mac:\mac\to \End_{ex}(\mathrm{Vect}_\C)$, where $\End_{ex}(\mathrm{Vect}_\C)$ is the fusion category of exact $\C$-linear endofunctors of $\mathrm{Vect}_\C$, which yields a fibre functor as follows.

 The functor  $G_\mac$ for a $\mac$-module category structure on $\mathrm{Vect}_\C$ sends $c\in \Ob\mac$ to the endofunctor $c\rhd-$ and a morphism $f:c\to c'$ to the natural transformation $f\rhd-: c\rhd -\Rightarrow c'\rhd -$, see  \cite[Prop 7.3.3]{EGNO}. As there is a tensor equivalence $\End_{ex}(\mathrm{Vect}_\C)\cong \mathrm{Vect}_\C$ that sends an exact $\C$-linear endofunctor $F:\mathrm{Vect}_\C\to\mathrm{Vect}_\C$ to $F(\C)$ and a natural transformation $\kappa: F\Rightarrow F'$ to  $\kappa_\C: F(\C)\to F'(\C)$, this defines a fibre functor for $\mac$. 

Thus, by \cite[Ex 7.4.6]{EGNO} there are finite-dimensional semisimple complex Hopf algebras $C$ and $B$ such that $\mac\cong C^{*op}\Mod$ and $\mad\cong B^*\Mod$  and  such that the $\mac$-module and $\mad$-right module structure on $\mathrm{Vect}_\C$ are given by tensoring on the left. 
Because the trace on $\mathrm{Vect}_\C$ is a $(\mac,\mad)$-bimodule trace, as imposed in Definition \ref{def:diagevalav}, the pivotal structures of $\mac$ and $\mad$ are the standard pivotal structures on $C^{*op}\Mod$ and $B^{*}\Mod$.

The natural isomorphism $q$ from Definition \ref{def:modulecat} for the $(\mac,\mad)$-bimodule category structure on $\mathrm{Vect}_\C$ is  determined  by its component  
$q_{C^{*},\C, B^{*}}: B^{*}\oo C^*\to C^*\oo B^*$, where $C^*$ and $B^*$ are the regular representations of $C^{*op}$ and $B^{*cop}$. By naturality, the latter is determined by 
$$
\tau_{BC}=q_{C^{*},\C, B^{*}}^\inv(1\oo 1)\in B^*\oo C^*
$$
The coherence conditions in the commuting diagrams \eqref{eq:pentagonbdef} then translate into the 
following equations, which give $(B,C)$ the structure of a Hopf doublet by Lemma \ref{lem:hopfdoublettrans} 
\begin{align*}
&(\Delta_{B^*}\oo \id)\tau_{BC}=\tau_{BC}^{13}\cdot \tau_{BC}^{23} & &(\id\oo \Delta_{C^*})\tau_{BC}=\tau_{BC}^{13}\cdot \tau_{BC}^{12}\\
&(\epsilon_{B^*}\oo\id)\tau_{BC}=1_{C^*} & &(\id\oo\epsilon_{C^*})\tau_{BC}=1_{B^*}.
\end{align*}
The first two equations follow directly from the computations in the proof of Proposition \ref{prop:hopfmod}. The second two follow from the conditions on the coherence isomorphisms and the fact that the action on $\C$ is given by the counit.
Also by  Proposition \ref{prop:hopfmod}, this defines a $D(B^{op}, C^{cop})^*\Mod$ module category structure on $\mathrm{Vect}_\C$.

By the proof of Proposition \ref{prop:hopfmod} and Example \ref{ex:hmodmodule}, the  fusion category 
$\End_{(\mac,\mad)}(\mathrm{Vect}_\C)$ is  equivalent to $D(B^{op}, C^{cop})^{cop}\Mod$. As the pivotal structure of $\End_{(\mac,\mad)}(\mathrm{Vect}_\C)$ is induced by the $(\mac,\mad)$-bimodule trace on $\mathrm{Vect}_\C$, see 
Corollary \ref{lem:macmodule} and Equation 
\eqref{eq:pivcond}, it follows that it is the standard pivotal structure induced by the trace on $\mathrm{Vect}_\C$.

The forgetful functor $F':D(B^{op}, C^{cop})^{cop}\Mod\to \mathrm{Vect}_\C$ thus defines a pivotal fibre functor for $\End_{(\mac,\mad)}(\mathrm{Vect_\C})$. Composing it with the pivotal functor $\Phi:\cala \to \End_{(\mac,\mad)}(\mathrm{Vect}_\C)$ yields a pivotal fibre functor for $\cala$. This defines  a finite-dimensional semisimple complex Hopf algebra
$A$ with $\cala\cong A^{*cop}\Mod$, a functor $\Phi': A^{*cop}\Mod\to D(B^{op}, C^{cop})^{cop}\Mod$ and a commuting diagram 
\begin{align*}
\xymatrix{
 \End_{(\mac,\mad)}(\mathrm{Vect}_\C) \ar[r]^{\cong\qquad} & D(B^{op}, C^{cop})^{cop}\Mod \ar[rd]^{F'}\\
& & \mathrm{Vect}_\C\\
\cala \ar[uu]^{\Phi} \ar[r]_{\cong\qquad}& A^{*cop}\Mod  \ar[uu]^{\Phi'}\ar[ru]_F,
}
\end{align*}
where  $F'$ and $F$ are the forgetful functors.
As $F'\Phi'=F$, this defines a Hopf algebra homomorphism
$$
\phi: D(B^{op}, C^{cop})^{cop}\xrightarrow[\cong]{\chi'}\End(F') \xrightarrow{\mu'\mapsto \mu'\Phi'} \End(F)\xrightarrow[\cong]{\chi^\inv} A^{* cop},
$$
where  $\chi': D(B^{op}, C^{cop})^{cop}\to \End(F')$ 
and $\chi: A^{*cop}\to \End(F)$
are the Hopf algebra isomorphism that assigns to an element $b\oo c\in D(B^{op}, C^{cop})^{cop}$ the endo-natural transformation $\mu^{b\oo c}: F'\Rightarrow F'$ with components $\mu^{b\oo c}_M: M\to M$, $m\mapsto (b\oo c)\rhd m$ and to an element $\alpha\in A^{*cop}$ the endo-natural transformation $\mu^\alpha: F\Rightarrow F$ with components $\mu^\alpha_N: N\to N$, $n\mapsto \alpha\rhd n$.
With Definition \ref{def:triplet}, Lemma \ref{lem:hopfconc} and the identities in the proof of Proposition \ref{prop:hopfmod} this  shows that $(A,B,C)$ is a Hopf triplet.
\end{Proof}

\subsection{Hopf algebraic trisection invariants are  special cases of  categorical ones}
\label{subsec:trisecbimod}

Propositions \ref{prop:hopfmod} and \ref{prop:reducetrip} show that Hopf triplets $(A,B,C)$ correspond precisely to those data $(\cala, \mad,\mac, \calm,\Phi)$  in
Definition \ref{def:diagevalav}, where  $\calm=\mathrm{Vect}_\C$. 
We now show that with the identification of the data from Proposition \ref{prop:hopfmod}, the representation theoretical  trisection bracket from Proposition \ref{prop:repthtrisec}  gives the averaged evaluation from Definition \ref{def:diagevalav} for any trisection diagram $T$.
By considering the trisection diagram $\Sigma_{st}$ for $S^4$ from \eqref{eq:stabilisation} we then obtain the following theorem.

\bigskip
\begin{Theorem} \label{th:cccth}For any Hopf triplet $(A,B,C)$ and trisection diagram $T$ 
the averaged evaluation of  $T$ for the data from Proposition \ref{prop:hopfmod}  coincides with its trisection bracket from Proposition \ref{prop:repthtrisec}:
$$
\mathrm{av}(T)=\langle T\rangle_{A,B,C}.
$$
If $\mathrm{av}(\Sigma_{st})=\langle \Sigma_{st}\rangle_{A,B,C}\neq 0$ and $\xi\in\C^\times$ with $\xi^3=\mathrm{av}(\Sigma_{st})=\langle \Sigma_{st}\rangle_{A,B,C}$, the
 4-manifold invariants from Theorem \ref{th:gentrisection} and
Definition \ref{def:cccinv} coincide for all closed, connected 4-manifolds.
\end{Theorem}

\begin{Proof} 
We fix sets of representatives
$I^A$, $I^B$ and $I^C$   of the isomorphism classes of simple representations of  $A^{*}$, $B^{*}$ and $C^{*}$ as in \eqref{eq:repsabc}. This fixes sets of representatives of the simple objects in $\cala$, $\mad$ and $\mac$, respectively. 
We orient the  green, blue and red curves of $T$ and label them with elements of $I_\mac$, $I_\mad$ and $I_\cala$. This defines a labelling map $h: \Gamma\to I$ as in Proposition \ref{prop:repthtrisec} that 
assigns to each red, blue or green curve $\lambda\in \Gamma$ a representation space $V_{h(\lambda)}^{X_\lambda}$ of $X_\lambda=A$, $X_\lambda=B$ or $X_\lambda=C$. 

The result is a labelled surface diagram as in Definition  \ref{def:labsurf}. To evaluate it according to Definition \ref{def:surfdiagev}, we cut $T$ 
along a set of cutting generators to obtain a tricoloured chord diagram. Each oriented chord inherits a labelling by elements of $I^A$, $I^B$ and $I^C$, and  boundary segments are labelled with $\C$. 
By Definition \ref{def:chord}  an endpoint of a chord labelled by a representation $M$  is labelled with vectors in 
\begin{align*}
\Hom_\C(\C, M\oo\C)\cong \Hom_\C(\C, M)\cong M\qquad \Hom_\C(M\oo \C,\C)\cong \Hom_\C(M,\C)\cong M^*.
\end{align*}
if the chord is incoming and outgoing, respectively, at the endpoint.
We choose dual bases of the representation spaces in $I^A$, $I^B$ and $I^C$  and label the chord endpoints with elements of these  bases. 

To evaluate of the resulting tricoloured chord diagram according to Definition \ref{def:tricev} we 
cut the boundary circle of the chord diagram, straighten its boundary to a vertical line
and assign the  natural isomorphisms defined by  \eqref{eq:bdef} and \eqref{eq:bimodfuncstruc}  to the chord crossings
\begin{align}\label{eq:bluered}
&\begin{tikzpicture}[scale=.5]
\draw[color=blue, line width=1pt, ->,>=stealth] (0,0) node[anchor=south]{$N$}--(2,-2);
\draw[color=red, line width=1pt, ->,>=stealth] (2,0) node[anchor=south]{$M$}--(0,-2);
\end{tikzpicture} 
&
&\begin{tikzpicture}[scale=.5]
\draw[color=red, line width=1pt, ->,>=stealth] (0,0) node[anchor=south]{$M$}--(2,-2);
\draw[color=blue, line width=1pt, ->,>=stealth] (2,0) node[anchor=south]{$N$}--(0,-2);
\end{tikzpicture}\\
&t^M_{-,N}: N\oo M\oo -\Rightarrow M\oo N\oo -
&
&t^{M\inv}_{-,N}: M\oo N\oo -\Rightarrow N\oo M\oo -
\nonumber\\
&n\oo m\oo -\mapsto (\tau_{AB(1)}\rhd m)\oo (\tau_{AB(2)}\rhd n)\oo -
&
&m\oo n\oo -\mapsto (\tau^\inv_{AB(2)}\rhd n)\oo (\tau^\inv_{AB(1)}\rhd m)\oo -\nonumber
\end{align}

\begin{align} \label{eq:greenblue}
&\begin{tikzpicture}[scale=.5]
\draw[color=green, line width=1pt, ->,>=stealth] (0,0) node[anchor=south]{$P$}--(2,-2);
\draw[color=blue, line width=1pt, ->,>=stealth] (2,0) node[anchor=south]{$N$}--(0,-2);
\end{tikzpicture}
&
&\begin{tikzpicture}[scale=.5]
\draw[color=blue, line width=1pt, ->,>=stealth] (0,0) node[anchor=south]{$N$}--(2,-2);
\draw[color=green, line width=1pt, ->,>=stealth] (2,0) node[anchor=south]{$P$}--(0,-2);
\end{tikzpicture}\\
&q^\inv_{P,-,N}: P\oo N\oo -\Rightarrow N\oo P\oo -
&
&q_{P,-,N}: N\oo P\oo -\Rightarrow P\oo N\oo -
\nonumber
\\
&p\oo n\oo -\mapsto (\tau_{BC(1)}^{\inv}\rhd n)\oo (\tau_{BC(2)}^{\inv}\rhd p)\oo -
&
&n\oo p\oo -\mapsto (\tau_{BC(2)}\rhd p)\oo (\tau_{BC(1)}\rhd n)\oo -\nonumber
\end{align}
\begin{align} \label{eq:redgreen}
&\begin{tikzpicture}[scale=.5]
\draw[color=red, line width=1pt, ->,>=stealth] (0,0) node[anchor=south]{$M$}--(2,-2);
\draw[color=green, line width=1pt, ->,>=stealth] (2,0) node[anchor=south]{$P$}--(0,-2);
\end{tikzpicture}
&
&\begin{tikzpicture}[scale=.5]
\draw[color=green, line width=1pt, ->,>=stealth] (0,0) node[anchor=south]{$P$}--(2,-2);
\draw[color=red, line width=1pt, ->,>=stealth] (2,0) node[anchor=south]{$M$}--(0,-2);
\end{tikzpicture}
\\
&s^{M}_{P,-}: M\oo P \oo -\Rightarrow P\oo M\oo -
&
&s^{M\inv}_{P,-}: P\oo M \oo -\Rightarrow M\oo P\oo -
\nonumber
\\
&m\oo p\oo -\mapsto (\tau_{CA(1)}^{\inv}\rhd p)\oo (\tau_{CA(2)}^{\inv}\rhd m)\oo -
&
&p\oo m\oo -\mapsto (\tau_{CA(2)}\rhd m)\oo (\tau_{CA(1)}\rhd p)\oo -\nonumber
\end{align}
Here $M\in I_\cala$, $N\in I_\calb$ and $P\in I_\mac$ are simple representations of $A^*$, $B^*$ and $C^*$.

Composing the functors that label the chords in the diagrams corresponds to the tensor products of the associated representation spaces. Composing them with the natural isomorphisms for the crossings corresponds to letting each isomorphism from 
\eqref{eq:bnp} and \eqref{eq:cromaps} act on the associated representation spaces. 

The evaluation of the chord diagram  according to Definition \ref{def:tricev} and Section \ref{subsec:bimodeval} then yields a complex number,  obtained as follows: 
\begin{compactenum}[(i)]
\item For  each horizontal slice that intersects chord crossings, take the  tensor product of all representation spaces for chords intersecting the slice above and below the crossings and combine it with the maps \eqref{eq:bnp} and \eqref{eq:cromaps} for crossings on the slice. \\[-2ex]
\item For a horizontal slice that intersects a vertex on the vertical line, take the tensor product of the representation spaces that intersect  slices above and below the vertex and relate them by the linear map associated to the vertex.\\[-2ex]

\item Compose the linear maps associated to horizontal slices from (i) and (ii) from the top to the bottom of the diagram.
\end{compactenum}
Summing over dual bases assigned to the vertices on the vertical line combines the contributions of different chords in the diagram that correspond to the same curve on the surface. 

By comparing with Proposition \ref{prop:repthtrisec} we find that for each labelling $h: \Gamma\to I$ of the red, blue and green curves in $T$ with elements of $I^A$, $I^B$ and $I^C$, the evaluation of the surface diagram given by $T$ is
$$
\mathrm{ev}_h(T)= \tr(\phi_{h,T}),
$$
where $\phi_{h,T}: V_{h,T}\to V_{h,T}$ is the linear map from Proposition \ref{prop:repthtrisec} and $V_{h,T}=\bigotimes_{\lambda\in \Gamma} V_{h(\lambda)}^{X_\lambda}$ the tensor product of the representation spaces assigned to the curves of $T$. 
The averaged evaluation of $T$ is obtained by multiplying with the dimensions of all representation spaces assigned to curves $\lambda\in\Gamma$ and by summing over all labellings $h: \Gamma\to I$.
Comparing with formula \eqref{eq:trisecbr} gives
$$
\mathrm{av}(T)=\sum_{h:\Gamma\to I} \prod_{\lambda\in\Gamma  }\dim_\C(V^{X_\lambda}_{h(\lambda)})\cdot \tr(\phi_{h,T}) =\langle T\rangle_{A,B,C}.
$$
For the second claim, note that $I_\calm=\{\C\}$ for $\calm=\mathrm{Vect}_\C$.
Recall that the trisection invariant from 
Definition \ref{def:cccinv} and the 4-manifold invariant from  Theorem \ref{th:gentrisection} 
are  obtained by normalising with a third root of 
$\mathrm{av}(\Sigma_{st})$ and of $\mathrm{av}_\C(\Sigma'_{st})$, respectively, where
$\Sigma_{st}$ is the trisection diagram of $S^4$ from
\eqref{eq:stabilisation} and $\Sigma'_{st}$ is the associated surface diagram from Example \ref{ex:stabilisationsurf_gen} with a disc removed and with boundary label $\C$. 
Lemma \ref{lem:gluing}  implies  $\mathrm{av}(\Sigma_{st})=\mathrm{av}_\C(\Sigma'_{st})$,
and
the condition $\mathrm{av}_\C(\Sigma'_{st})\neq 0$ in Lemma \ref{prop:stabilising} reduces to $\langle \Sigma_{st}\rangle_{A,B,C}\neq 0$.
\end{Proof}

Theorem \ref{th:cccth} exhibits the Hopf algebraic trisection invariants from Definition \ref{def:cccinv} as a special case of the 4-manifold invariants in Theorem \ref{th:gentrisection}: the former correspond precisely to those 4-manifold invariants in Theorem \ref{th:gentrisection} where the bimodule category is $\calm=\mathrm{Vect}_\C$.

\section{Four manifold invariants beyond Hopf triplets}
\label{sec:newexamples}

In this section we give examples of the categorical trisection  invariants from Theorem \ref{th:gentrisection} that do not arise from Hopf triplets. 
Such invariants were anticipated already in \cite[Rem 5.10]{CCC}, where it is stated that the Hopf algebraic trisection invariants should have generalisations that  include B\"arenz-Barrett invariants  \cite{BB}
for data that do not arise from Hopf algebra representations. 

Section \ref{subsec:bbinv} shows that the categorical trisection invariants are such a generalisation. We consider a pivotal functor $\phi:\cala\to \mac$ from a spherical fusion category $\cala$ into a modular fusion category $\mac$. We show that this defines a set of algebraic data as in Theorem \ref{th:gentrisection} and that the resulting  trisection invariants  are the B\"arenz-Barrett invariants \cite{BB}, up to a rescaling that depends  on the algebraic data and the Euler characteristic. 

In Section \ref{subsec:DWinvariant} we compute the  trisection invariant from Theorem \ref{th:gentrisection} for the case where $\mac$ and $\mad$ are categories of vector spaces graded by finite groups $C$ and $B$ and the $(\mac,\mad)$-bimodule category $\calm$ is the category of vector spaces graded by a finite $C\times B^{op}$-set $M$. In this case, our trisection invariants define a generalisation of the invariants from \cite{CCC} to a certain triple of \emph{weak} Hopf algebras, as anticipated  in  \cite[Rem 5.10]{CCC}. They admit a simple diagrammatic description, in which the diagrams encode the weak Hopf algebra structure.

\subsection{B\"arenz-Barrett invariants}
\label{subsec:bbinv}

In this section we show that a modular fusion category $\mac$  and a pivotal functor $\phi:\cala\to \mac$ yields the input data for Theorem \ref{th:gentrisection} in such a way that the resulting categorical trisection invariant coincides with the B\"arenz-Barrett  invariant \cite{BB} for $\phi:\cala\to \mac$, up to a rescaling given by the algebraic data and the Euler characteristic.

In the following  let  $\mac$ be a modular fusion category with braiding $c:\oo\Rightarrow\oo^{op}$ and $I_\mac$ a set of representatives of the isomorphism classes of simple objects in $\mac$.  We denote by $\mac^{rev}$ the modular fusion category with the opposite tensor product and braiding, $c^{rev}_{x,y}=c_{y,x}$, and by $\overline \mac$ the modular fusion category with the opposite braiding, $\overline c_{x,y}=c_{y,x}^\inv$.

Recall that in a modular fusion category $\mac$ we have the well-known  cutting strands identity due to Lickorish \cite{Lk} and Roberts \cite{Ro1},  which holds  for any object $x\in \Ob\mac$ 
\begin{align}\label{eq:cuttingid}
\begin{tikzpicture}[scale=.6]
\node at (-.5,0)[anchor=east]{$\sum_{c\in I_\mac} \dim(c)$}; 
\draw[line width=1pt] (2,0) ellipse  (2cm and 1cm);
\node at (4,0) [anchor=west]{$c$};
\draw[draw=none, fill=white] (2,1) circle (.3);
\draw[line width=1pt] (2,2) node[anchor=south]{$x$}--(2,-.8) node[sloped, pos=0.5, allow upside down]{\arrowIn};
\draw[line width=1pt,->,>=stealth] (2,-1.2)--(2,-2) node[anchor=north]{$x$};
\node at (5.5,0){$=$};
\node at (6,0)[anchor=west]{$\dim(\mac)\cdot  \sum_{\iota}$};
\draw[line width=1pt] (10,2) node[anchor=south]{$x$}--(10,.5) node[sloped, pos=0.5, allow upside down]{\arrowIn};
\draw[fill=black] (10,.5) circle (.15) node[anchor=west]{$\iota$};
\draw[fill=black] (10,-.5) circle (.15) node[anchor=west]{$\iota$};
\draw[line width=1pt] (10,-.5)--(10,-2) node[anchor=north]{$x$} node[sloped, pos=0.5, allow upside down]{\arrowIn};
\end{tikzpicture}
\end{align}
Here, $\iota$ runs over 
dual bases of $\Hom_\mac(x,e)$ and $\Hom_\mac(e,x)$, and the orientation of the circle  is irrelevant.

The B\"arenz-Barrett invariant \cite{BB} of a compact smooth oriented connected 4-manifold $X$
generalises the Crane-Yetter invariant \cite{CY,CYK}, Broda's invariant \cite{Br} as well as the refinement of the Broda invariant by Roberts \cite{Ro1,Ro2} and Petit's dichromatic invariant \cite{P}. 

It is defined in terms of a handle decomposition of $X$. This handle decomposition  is encoded in a framed, bicoloured link diagram. Recall that such a diagram can be obtained from a trisection diagram, as explained in Remark \ref{rem:linkdiagram}. It has link components of two colours,  red link components as well as black circles that encircle strands of the red components. The former represent the attaching maps for the 2-handles, the latter the attaching maps for the 1-handles.
We denote by $R$ the set of red components, by $C$ the set of black circles
and define $g:=|R|$ and $k:=|C|$.

The algebraic data for the B\"arenz-Barrett invariant  is a ribbon fusion category $\mac$ with trivial twists on all transparent objects, a spherical fusion category $\cala$ and a pivotal functor $\phi:\cala\to\mac$. In the following for simplicity we restrict attention to the case, where $\mac$ is a modular fusion category and adapt the definition to our notation. 
We fix sets of representatives $I_\cala$ and $I_\mac$ of the isomorphism classes of simple objects in $\cala$ and $\mac$.

\begin{Theorem}\cite[Def 3.6, Th 3.9 and Prop 4.13]{BB}\label{th:bbtheorem}\\
Let $\phi:\cala\to\mac$ be a pivotal functor from a spherical fusion category $\cala$ to a modular fusion category $\mac$. Then the following defines an invariant of closed connected 4-manifolds $X$
    $$
I_{BB}(X)=\dim(\cala)^{k-g}\cdot \dim(\mac)^{-k}\cdot n(\phi)^{-k}\sum_{l: R\to I_\cala}\sum_{l': C\to I_\mac}\prod_{\rho\in R} \dim\, l(\rho)\prod_{\gamma\in C}\dim\, l'(\gamma) \cdot \mathrm{ev}_{\mac}(L),
$$
where
\begin{compactitem}
\item $L$ is a framed link diagram  that specifies a handlebody decomposition of $X$,
\item $n(\phi)=\sum_{a\in I_\cala}  \dim(a) \cdot \dim_\C \Hom_\mac(e, \phi(a))$,
\item the sums run over labellings $l$ and $l'$ of red curves and black circles with elements of $I_\cala$ and  $I_\mac$,
\item $\mathrm{ev}_{\mac}(L)$ denotes the evaluation of the  link diagram $L$ in  $\mac$ whose black circles are labelled by $l'$ and whose red components by the images of $l$ under $\phi$.
\end{compactitem}
\end{Theorem}

To show that the B\"arenz-Barrett invariant for a modular fusion category $\mac$ is a special case of the trisection invariant from Theorem \ref{th:gentrisection}, we consider the following algebraic data. 
We  view the modular fusion category $\mac$ as a $(\mac,\mac^{rev})$-bimodule category or, equivalently, as a $\mac\boxtimes\mac$-module category with 
$$
\rhd=\oo:\mac\times\mac\to\mac\qquad \lhd=\oo^{op}:\mac\times\mac^{rev}\to\mac,
$$
  $l:\rhd(\oo\times\id)\Rightarrow\rhd(\id\times\rhd)$ and $r: \lhd(\lhd\times\id)\Rightarrow \lhd(\id\times\oo^{rev})$ given by the associators in $\mac$ and with the coherence isomorphism $q:\lhd(\rhd\times\id)\Rightarrow\rhd(\id\times\lhd)$ given by the braiding in $\mac$
\begin{align}\label{eq:braidcross}
q_{x,m,y}= a_{x,y,m}  \circ (c_{y,x}\oo 1_m)\circ a_{y,x,m}^\inv: y\oo(x\oo m)\to x\oo (y\oo m).
\end{align}
We then have an induced pivotal  functor $\chi:\mac\to \End_{(\mac,\mac^{rev})}(\mac)$ that assigns
\begin{compactitem}
\item  to an object $x\in \Ob\mac$ the bimodule endofunctor   $F_x=x\oo -:\mac\to\mac$
with the coherence isomorphisms $s^{x}: F_x\rhd\Rightarrow\rhd(\id\times F_x)$ and $t^{x}:\lhd (F_x\times\id)\Rightarrow F_x\lhd$ given by the braiding
\begin{align}\label{eq:funcbraid}
&s^{x}_{y,m} = a_{y,x,m}\circ (c_{x,y}\oo 1_m)\circ a_{x,y,m}^\inv: x\oo (y\oo m)\to y\oo (x\oo m)\\
&t^{x}_{y,m} = a_{x,y,m}\circ (c_{x,y}^\inv\oo 1_m)\circ a_{y,x,m}^\inv: y\oo (x\oo m)\to x\oo (y\oo m).\nonumber
\end{align}
\item to a morphism $\alpha: x\to x'$ the natural transformation $\alpha\oo -:x\oo -\Rightarrow x'\oo -$.
\end{compactitem}
This 
identifies $\mac$ with a fusion subcategory  of  $\End_{(\mac,\mac^{rev})}(\mac)\cong Z(\mac)\cong \mac\boxtimes \overline\mac$, where the second equivalence is due to $\mac$ being modular, see \cite[Prop 8.20.12]{EGNO}. It follows that  any pivotal functor $\phi:\cala\to\mac$ defines a pivotal functor $\Phi=\chi \phi:\cala\to \End_{(\mac,\mac^{rev})}(\mac)$.

The diagrammatic calculus  from Section \ref{sec:diagcalc}, where green and blue lines are labelled by objects of $\calc$, and red lines are labelled by objects of  $\mac \subseteq \End_{(\mac,\mac^{rev})}(\mac)$ reduces to the usual diagrammatic calculus for the modular fusion category $\mac$. More specifically, the choice of the coherence isomorphism \eqref{eq:braidcross} identifies
crossings between green and blue lines in the  diagrammatic calculus from Section \ref{sec:diagcalc} with  overcrossings of the blue  over the green line  in the diagrammatic calculus for  $\mac$. The choice of the coherence isomorphisms in  \eqref{eq:funcbraid} identifies crossings on red lines with overcrossings of the red lines.
\begin{align}\label{eq:gbpic}
&\begin{tikzpicture}[scale=.3]
\begin{scope}[shift={(0,0)}]
\draw[green, line width=1pt] (2,2)--(-2,-2) node[anchor =north]{$x$};
\draw[blue, line width=1pt] (-2,2)--(2,-2) node[anchor =north]{$y$};
\end{scope}
\node at (3,0){$=$};
\begin{scope}[shift={(6,0)}]
\draw[green, line width=1pt] (2,2)--(-2,-2) node[anchor =north]{$x$};
\draw[draw=none, fill=white] (0,0) circle (.5);
\draw[blue, line width=1pt] (-2,2)--(2,-2) node[anchor =north]{$y$};
\end{scope}
\end{tikzpicture}
& 
&\begin{tikzpicture}[scale=.3]
\begin{scope}[shift={(0,0)}]
\draw[blue, line width=1pt] (2,2)--(-2,-2) node[anchor =north]{$x$};
\draw[green, line width=1pt] (-2,2)--(2,-2) node[anchor =north]{$y$};
\end{scope}
\node at (3,0){$=$};
\begin{scope}[shift={(6,0)}]
\draw[green, line width=1pt] (-2,2)--(2,-2) node[anchor =north]{$x$};
\draw[draw=none, fill=white] (0,0) circle (.5);
\draw[blue, line width=1pt] (2,2)--(-2,-2) node[anchor =north]{$y$};
\end{scope}
\end{tikzpicture}\\
&q_{x,-,y}: y\oo(x\oo -)\Rightarrow x\oo (y\oo -) & &q^\inv_{x,-,y}: x\oo(y\oo -)\Rightarrow y\oo (x\oo -)\nonumber
\end{align}
\begin{align}\label{eq:rgpic}
&\begin{tikzpicture}[scale=.3]
\begin{scope}[shift={(0,0)}]
\draw[green, line width=1pt] (2,2)--(-2,-2) node[anchor =north]{$y$};
\draw[red, line width=1pt] (-2,2)--(2,-2) node[anchor =north]{$F_x$};
\end{scope}
\node at (3,0){$=$};
\begin{scope}[shift={(6,0)}]
\draw[green, line width=1pt] (2,2)--(-2,-2) node[anchor =north]{$y$};
\draw[draw=none, fill=white] (0,0) circle (.5);
\draw[red, line width=1pt] (-2,2)--(2,-2) node[anchor =north]{$x$};
\end{scope}
\end{tikzpicture}
& 
&\begin{tikzpicture}[scale=.3]
\begin{scope}[shift={(0,0)}]
\draw[red, line width=1pt] (2,2)--(-2,-2) node[anchor =north]{$F_x$};
\draw[green, line width=1pt] (-2,2)--(2,-2) node[anchor =north]{$y$};
\end{scope}
\node at (3,0){$=$};
\begin{scope}[shift={(6,0)}]
\draw[green, line width=1pt] (-2,2)--(2,-2) node[anchor =north]{$y$};
\draw[draw=none, fill=white] (0,0) circle (.5);
\draw[red, line width=1pt] (2,2)--(-2,-2) node[anchor =north]{$x$};
\end{scope}
\end{tikzpicture}\\
&s^x_{y,-}: x\oo (y\oo -)\Rightarrow y\oo (x\oo -) & &s^{x\inv}_{y,-}: y\oo (x\oo -)\Rightarrow x\oo (y\oo -)\nonumber
\end{align}

\begin{align}\label{eq:rbpic}
&\begin{tikzpicture}[scale=.3]
\begin{scope}[shift={(0,0)}]
\draw[blue, line width=1pt] (2,2)--(-2,-2) node[anchor =north]{$y$};
\draw[red, line width=1pt] (-2,2)--(2,-2) node[anchor =north]{$F_x$};
\end{scope}
\node at (3,0){$=$};
\begin{scope}[shift={(6,0)}]
\draw[blue, line width=1pt] (2,2)--(-2,-2) node[anchor =north]{$y$};
\draw[draw=none, fill=white] (0,0) circle (.5);
\draw[red, line width=1pt] (-2,2)--(2,-2) node[anchor =north]{$x$};
\end{scope}
\end{tikzpicture}
& 
&\begin{tikzpicture}[scale=.3]
\begin{scope}[shift={(0,0)}]
\draw[red, line width=1pt] (2,2)--(-2,-2) node[anchor =north]{$F_x$};
\draw[blue, line width=1pt] (-2,2)--(2,-2) node[anchor =north]{$y$};
\end{scope}
\node at (3,0){$=$};
\begin{scope}[shift={(6,0)}]
\draw[blue, line width=1pt] (-2,2)--(2,-2) node[anchor =north]{$y$};
\draw[draw=none, fill=white] (0,0) circle (.5);
\draw[red, line width=1pt] (2,2)--(-2,-2) node[anchor =north]{$x$};
\end{scope}
\end{tikzpicture}\\
&t^{x\inv}_{y,-}: x\oo (y\oo -)\Rightarrow y\oo (x\oo -) &
&t^{x}_{y,-}: y\oo (x\oo -)\Rightarrow x\oo (y\oo -)\nonumber
\end{align}

The diagrammatic identities \eqref{pic:rm2}  from Section \ref{sec:diagcalc} then depict the invertibility of the braiding in $\mac$ and the diagrammatic identities \eqref{pic:hexagon} the dodecagon identity  in $\mac$.

Note that the choice of the coherence isomorphisms in \eqref{eq:braidcross} and \eqref{eq:funcbraid} is a convention. Any braiding $c_{x,y}$ in $\mac$ in these expressions may be replaced by the opposite braiding $c_{y,x}^\inv$, as long as there is one colour that braids over the other two colours  in  diagrams \eqref{eq:gbpic} to \eqref{eq:rbpic}. This condition  ensures that  Identity \eqref{pic:hexagon} holds.
As different conventions for the coherence isomorphisms  in \eqref{eq:braidcross} and \eqref{eq:funcbraid} amount to a permutation of colours in the associated diagrammatic calculus, we ignore the other choices.

We now relate the trisection invariant from Theorem \ref{th:gentrisection} for this algebraic data to the B\"arenz-Barrett invariant from Theorem \ref{th:bbtheorem}.The following Proposition \ref{prop:bbtri} and  Theorem \ref{th:bbtric} are a generalisation of Theorem 5.7 in \cite{CCC}, where an analogous statement is shown for B\"arenz-Barrett invariants arising from  Hopf algebra representations. Our result is more general, and the proof becomes simpler.
That this datum satisfies the stabilisation condition follows from the proof of Proposition \ref{prop:bbtri}.

\begin{Proposition} \label{prop:bbtri}Let $\mac$ be a modular fusion category and $\phi:\cala\to \mac$ a pivotal functor from a spherical fusion category $\cala$. 
Then the averaged evaluation of a trisection diagram $T$, whose green and blue lines are labelled by objects of $\calc$, and red lines are labelled by images of objects of $\cala$ in $\End_{(\mac,\mac^{rev})}(\mac)$, is computed as follows:
\begin{compactenum}
\item Transform the trisection diagram $T$ into a $(g,k)$-trisection diagram $T'$, where the green and blue curves are in standard position as in Figure \ref{fig:standardheegard}, by applying the moves from Theorem \ref{th:intmoves}.\\[-2ex]
\item Consider the framed link  defined by the red curves of $T'$, with the framing defined by the surface.\\[-2ex]
\item Erase the last $(g-k)$ green and blue curves of $T'$. Replace each of the first  $k$ pairs of green and blue curves in Figure \ref{fig:standardheegard} by a black circle and lift the circle  slightly off the surface such that it surrounds the red strands intersecting this pair. This yields a set $C$ of $k$ black circles.\\[-2ex]
\item Label each red curve of $T'$  by the image  of an object in $I_\cala$ under $\phi$ and each black circle  by an object in $ I_\mac$.  
Evaluate the resulting link diagram $L$ for $\mac$,  multiply by the dimensions of all labels, sum over all such labellings and multiply by $\dim(\mac)$
\end{compactenum}
\begin{align}\label{eq:trisecCY}
\mathrm{av}(T)=\dim(\mac) \sum_{l: R\to I_\cala}\sum_{l': C\to I_\mac} \prod_{\rho\in R} \dim\, l(\rho)\, \prod_{\gamma\in C} \dim\, l'(\gamma)\;\mathrm{ev}_\mac(L).
\end{align}
\end{Proposition}

\begin{figure}
\begin{align}
& &
	\begin{tikzpicture}[baseline={([yshift=-.5ex]current bounding box.center)}]
	\node at (0,0) {\def\svgscale{0.25} 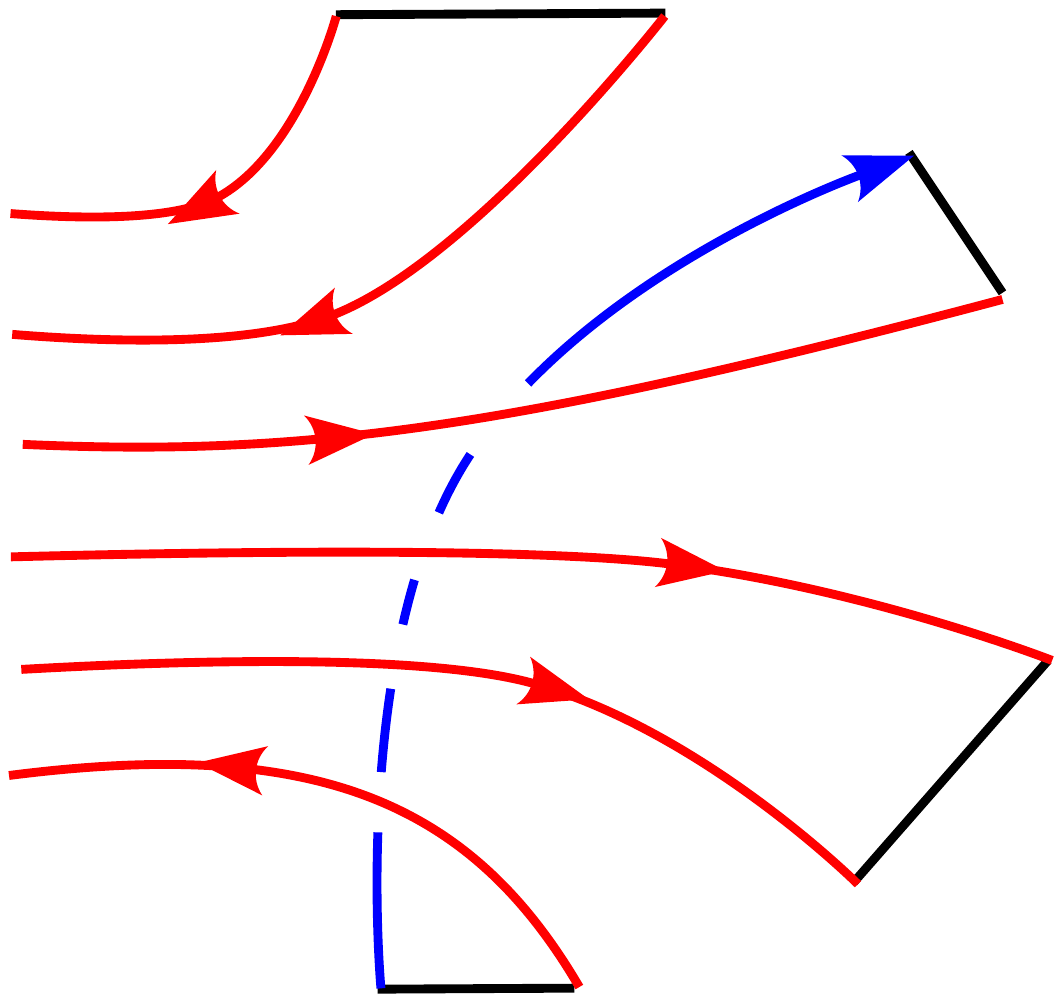 };
	\end{tikzpicture}

&
&=
&
&
	\begin{tikzpicture}[baseline={([yshift=-.5ex]current bounding box.center)}]
	\node at (0,0) {\def\svgscale{0.25} 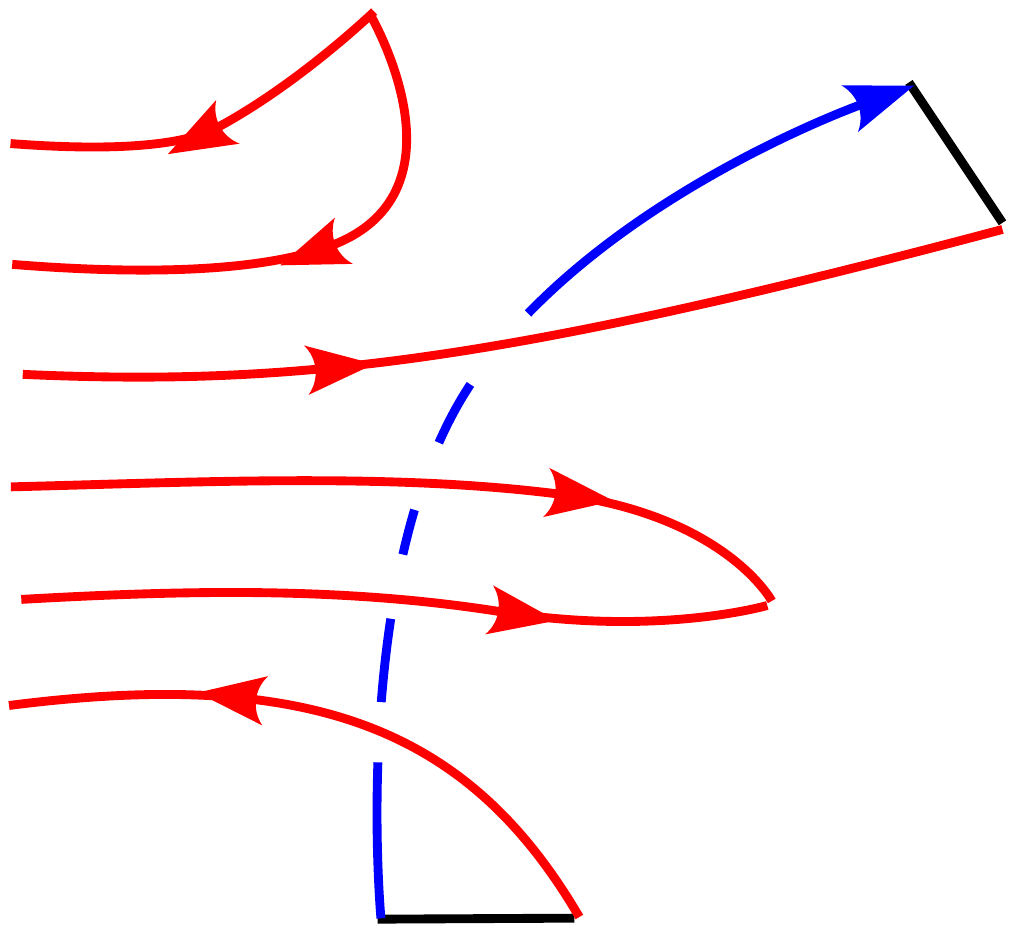 };
	\end{tikzpicture}
\nonumber\\
= & &
	\begin{tikzpicture}[baseline={([yshift=-.5ex]current bounding box.center)}]
	\node at (0,0) {\def\svgscale{0.25} 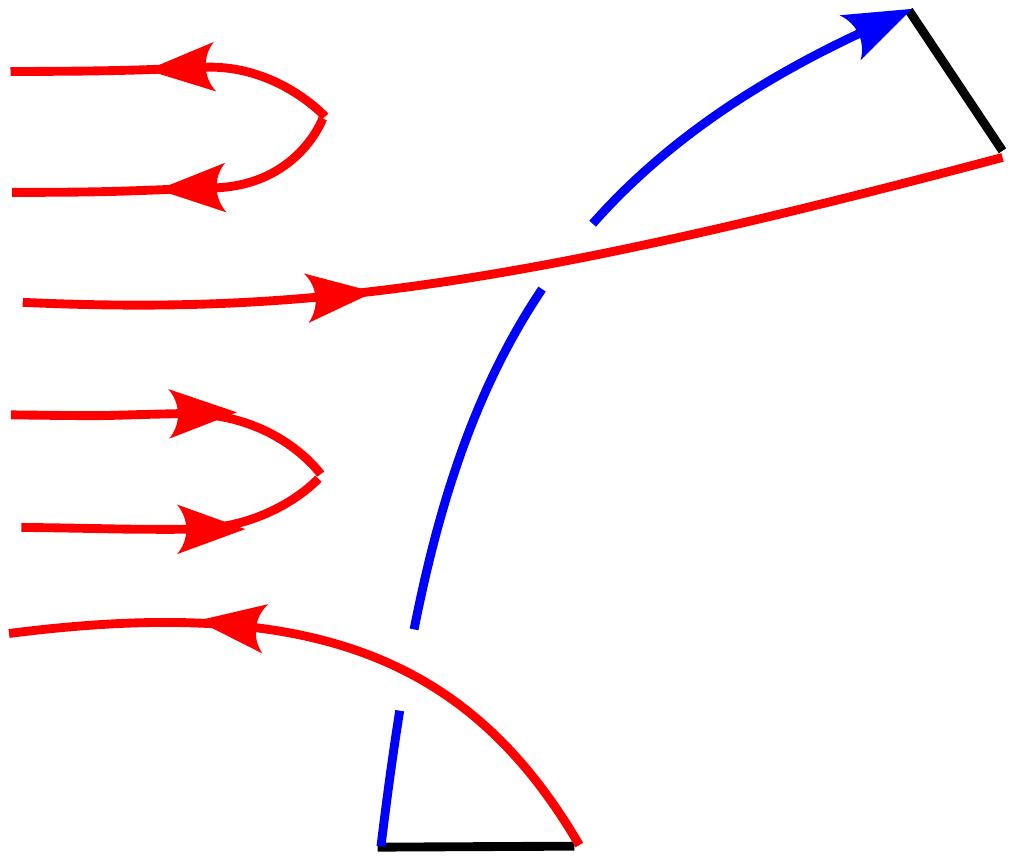 };
	\end{tikzpicture}

&
&=
&
&
	\begin{tikzpicture}[baseline={([yshift=-.5ex]current bounding box.center)}]
	\node at (0,0) {\def\svgscale{0.25} 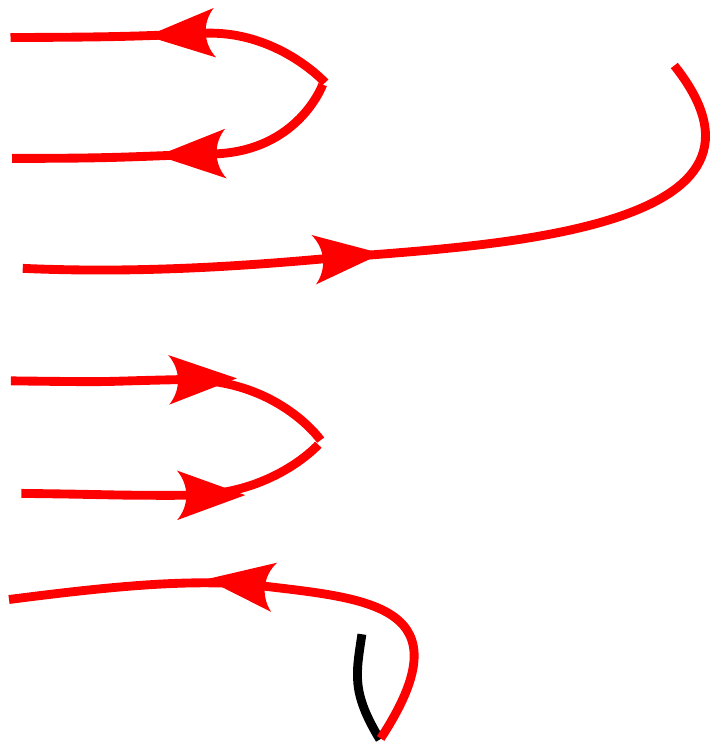 };
	\end{tikzpicture}
\nonumber\\
= & &
	\begin{tikzpicture}[baseline={([yshift=-.5ex]current bounding box.center)}]
	\node at (0,0) {\def\svgscale{0.25} 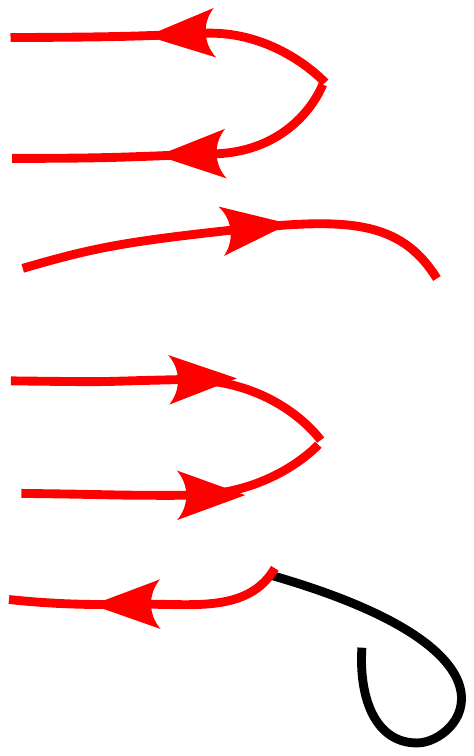 };
	\end{tikzpicture}

&
&=
&
&
	\begin{tikzpicture}[baseline={([yshift=-.5ex]current bounding box.center)}]
	\node at (0,0) {\def\svgscale{0.25} 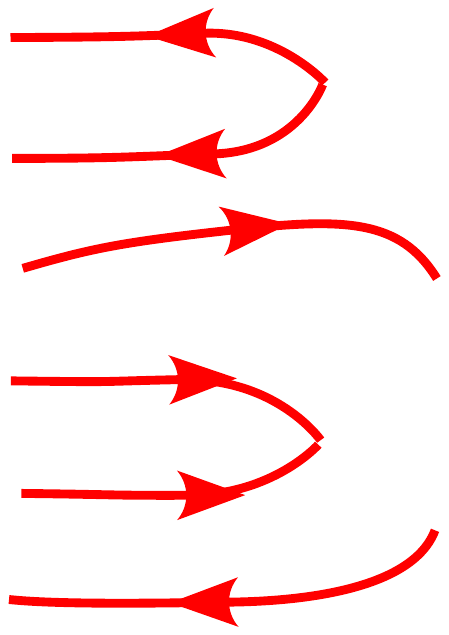 };
	\end{tikzpicture}
\nonumber\\
= & &
	\begin{tikzpicture}[baseline={([yshift=-.5ex]current bounding box.center)}]
	\node at (0,0) {\def\svgscale{0.25} 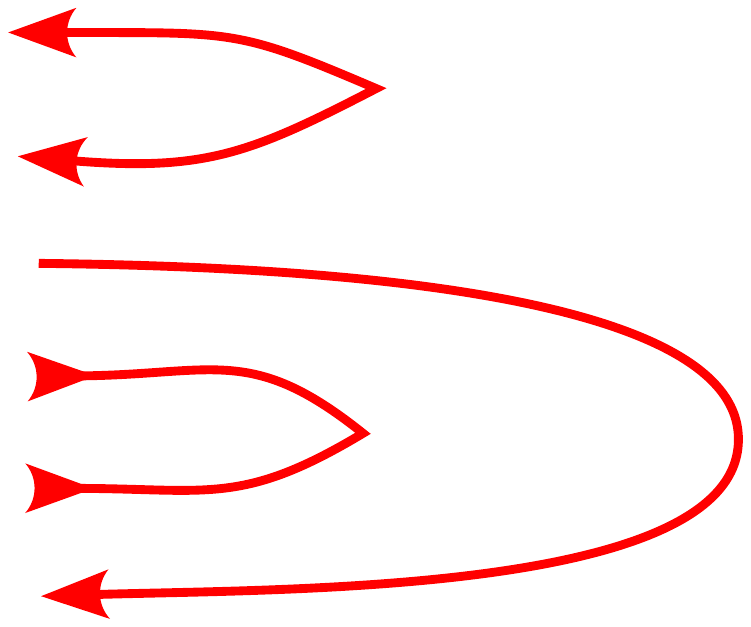 };
	\end{tikzpicture}

\end{align}
\caption{Averaged evaluation of the trisection diagram $T''$ for $m=e$ for  one of  the first $k$ handles.\newline
}
\label{fig:parhandle}
\end{figure}

\begin{figure}
\begin{align}
& &
	\begin{tikzpicture}[baseline={([yshift=-.5ex]current bounding box.center)}]
	\node at (0,0) {\def\svgscale{0.25} 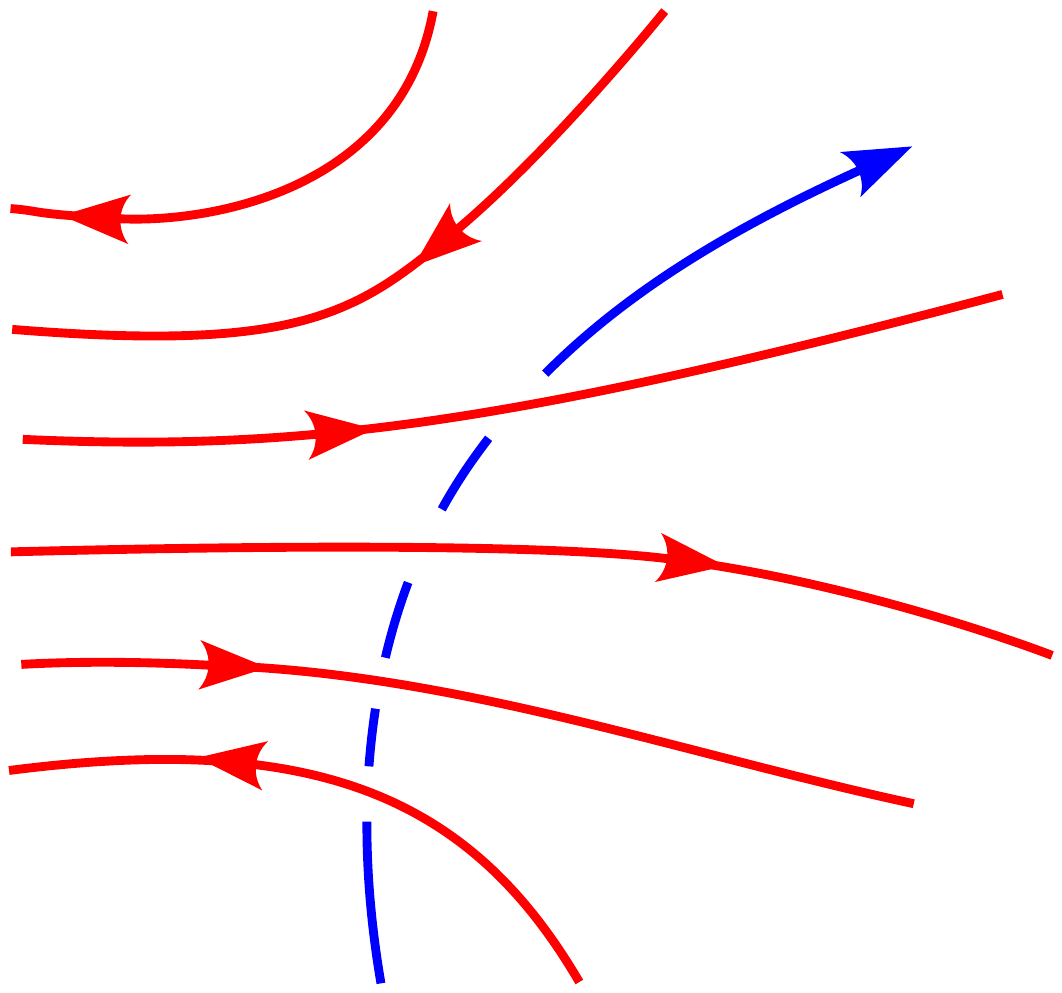 };
	\end{tikzpicture}

&
&=
&
&
	\begin{tikzpicture}[baseline={([yshift=-.5ex]current bounding box.center)}]
	\node at (0,0) {\def\svgscale{0.25} 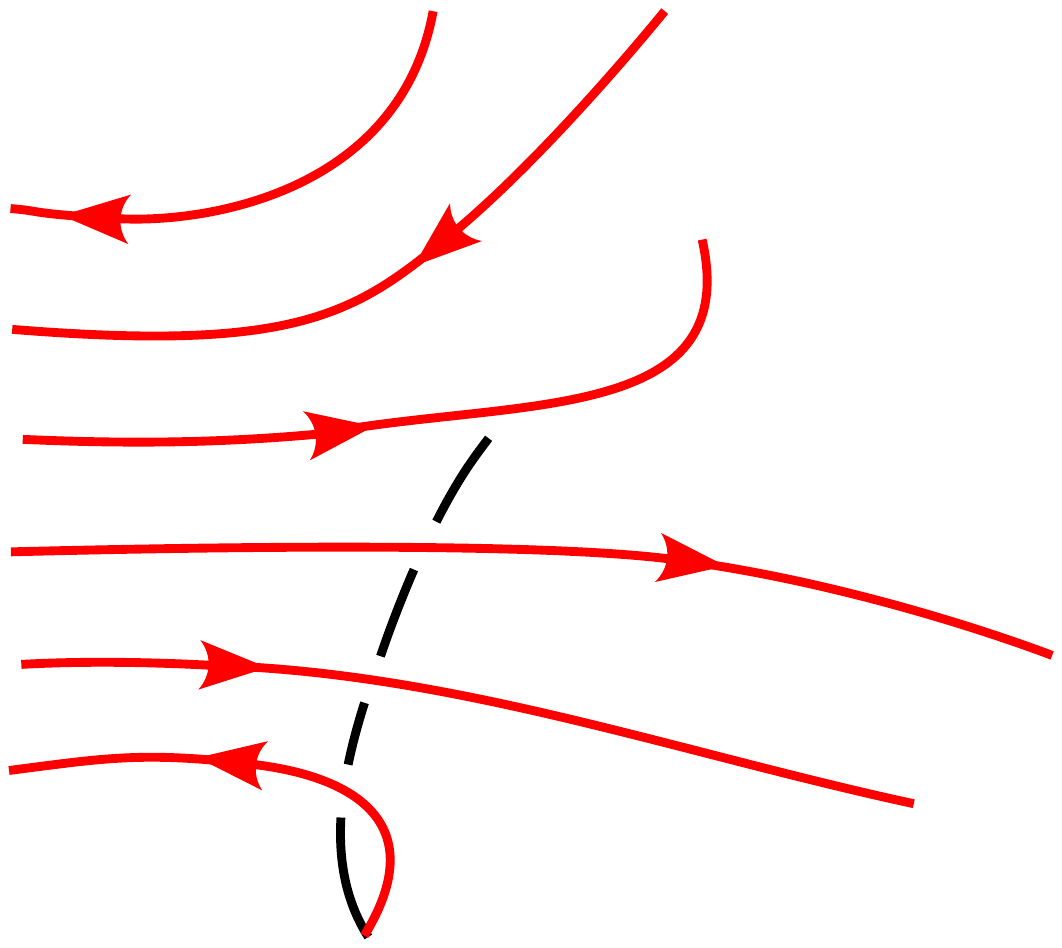 };
	\end{tikzpicture}
\nonumber\\
= & &
	\begin{tikzpicture}[baseline={([yshift=-.5ex]current bounding box.center)}]
	\node at (0,0) {\def\svgscale{0.25} 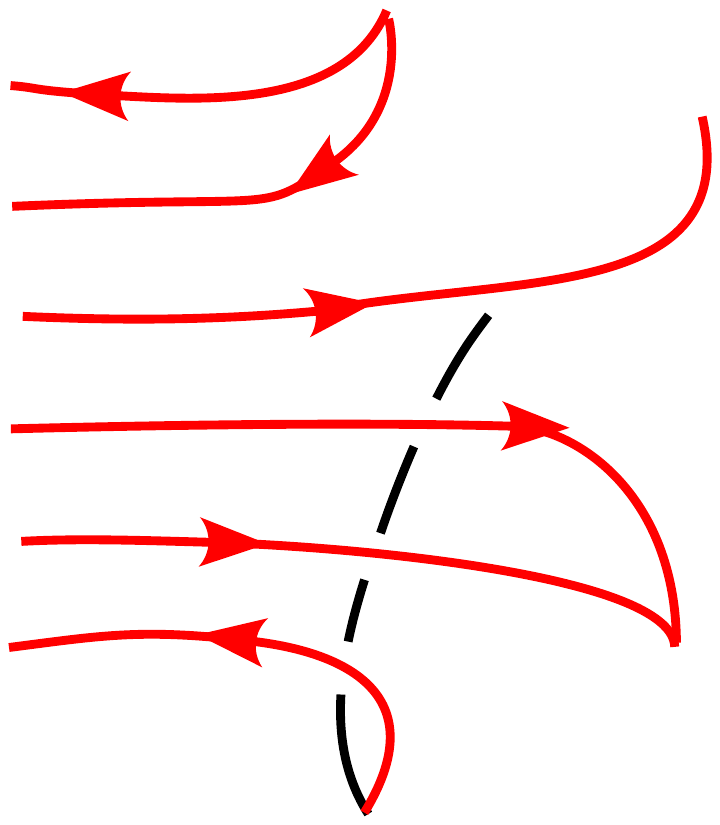 };
	\end{tikzpicture}

&
&=
&
&
	\begin{tikzpicture}[baseline={([yshift=-.5ex]current bounding box.center)}]
	\node at (0,0) {\def\svgscale{0.25} 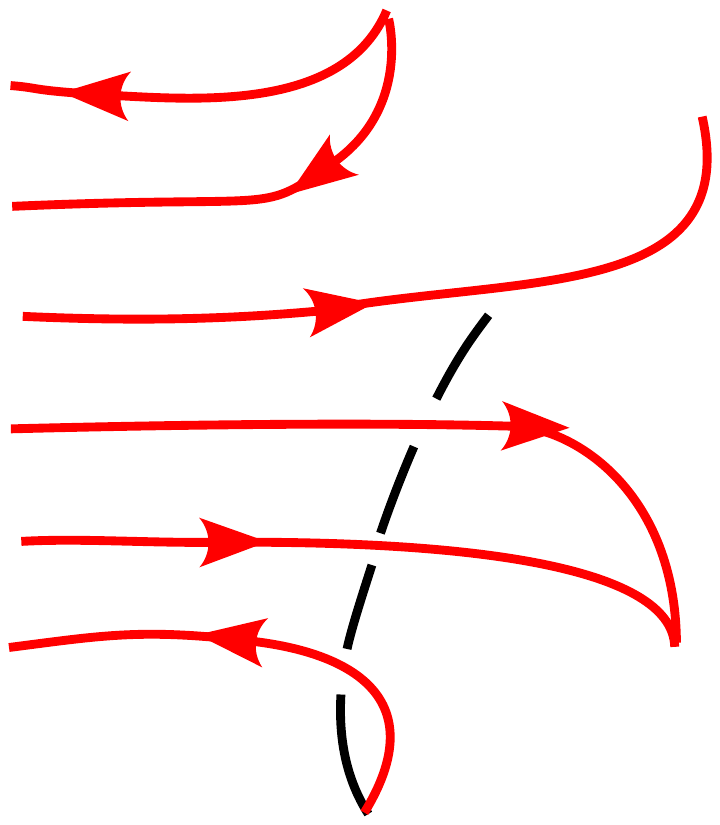 };
	\end{tikzpicture}
\nonumber\\
= & &
	\begin{tikzpicture}[baseline={([yshift=-.5ex]current bounding box.center)}]
	\node at (0,0) {\def\svgscale{0.25} 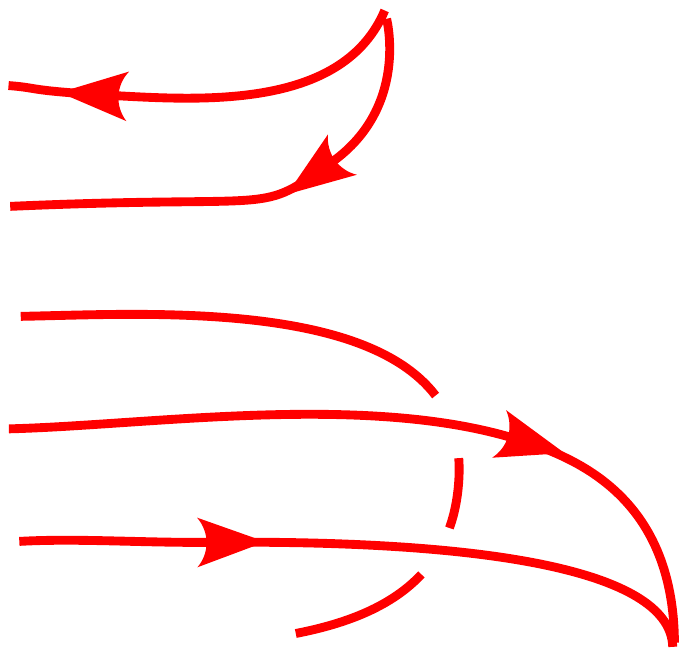 };
	\end{tikzpicture}

&
&=
&
&
	\begin{tikzpicture}[baseline={([yshift=-.5ex]current bounding box.center)}]
	\node at (0,0) {\def\svgscale{0.25} 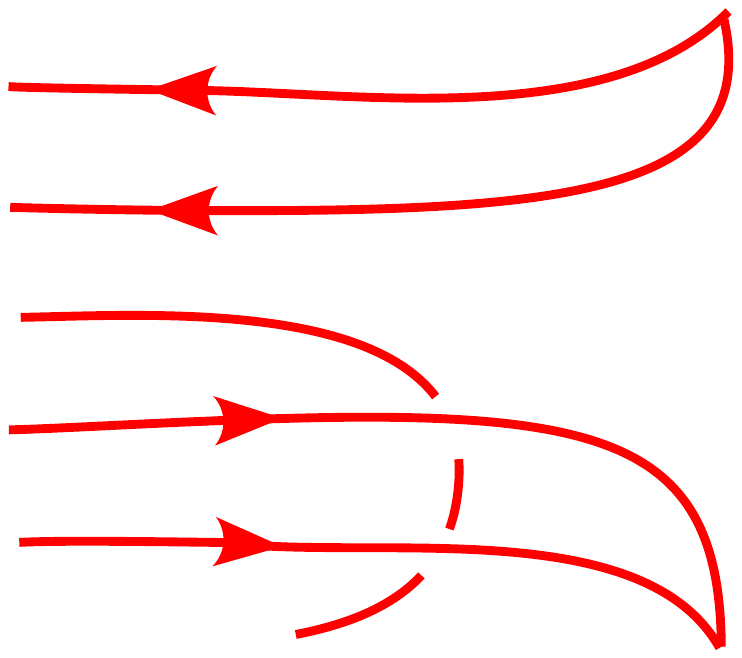 };
	\end{tikzpicture}
\nonumber\\
= & &
	\begin{tikzpicture}[baseline={([yshift=-.5ex]current bounding box.center)}]
	\node at (0,0) {\def\svgscale{0.25} 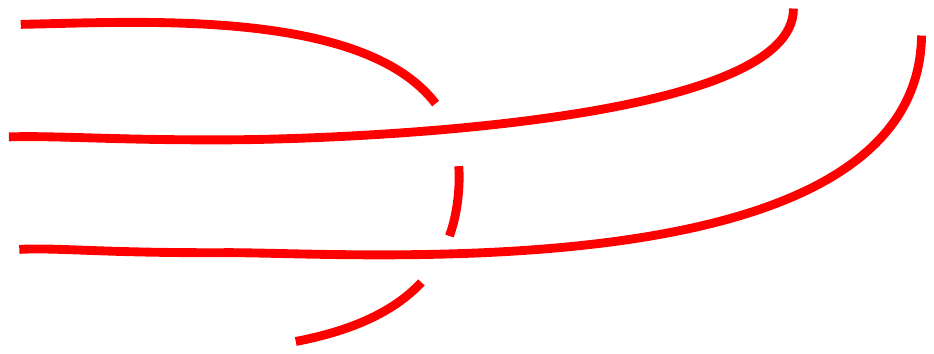 };
	\end{tikzpicture}

\end{align}
\caption{Averaged evaluation of the trisection diagram $T''$ for $m=e$ for  one of  the last $(g-k)$ handles.
}
\label{fig:orthhandle}
\end{figure}

\begin{Proof}
First observe that a trisection diagram $T$ may always be transformed into a trisection diagram $T'$ with the green and blue curves in standard position, see Theorem \ref{th:intmoves} and the discussion that follows. By the proof of Theorem \ref{th:gentrisection} the
averaged evaluation of such a diagram is unchanged by the transformation.

The rest of the proof proceeds as follows. In step 1 we argue that $\mathrm{av}(T')=\dim(\mac)\cdot \mathrm{av}_e(T'')$, where $T''$ is obtained by removing a disc from $T'$ and labelling the boundary with the tensor unit $e$. We then show in steps 2 and 3 that $\mathrm{av}_e(T'')$ is obtained by the procedure described in points $1-4$, up to multiplication by $\dim(\mac)$.

1.~Let $T''$ be  obtained by removing a disc from $T'$ and assigning a boundary label $m\in I_\mac$. By Lemma \ref{prop:stabilising}, there is a constant $C_T\in \C$, independent of   $m\in I_\mac$,  with  $\mathrm{av}_m(T'')=C_T\cdot \dim(m)$ for all $m\in I_\mac$.  Lemma \ref{lem:gluing}  then implies 
$$\mathrm{av}(T')=\sum_{m\in I_\mac} \dim(m)^2\cdot C_T=\dim(\mac)\cdot C_T.$$ 
We can compute $C_T$ by choosing the tensor unit $m=e\in I_\mac$, which yields  $C_T=C_T\cdot \dim(e)=\mathrm{av}_e(T'')$.
To compute the averaged evaluation of $T''$, we cut $T''$ along a set of cutting generators as in Definition \ref{def:surfdiagev}. This yields  a $(4g+1)$-gon $P$, where $4g$ sides are identified pairwise and the remaining side corresponds to the boundary circle of $T''$. Each corner of $P$ is labelled by $m=e$. 
As in Definition \ref{def:surfdiagev}, we assign simple objects in $I_\mac$ to all boundary segments of $P$ that are not labelled by $m=e$ and 
multiply by their dimensions.

2.~We consider one of  the  first $k$ handles with parallel blue and green circles in Figure \ref{fig:standardheegard}. As both, the green and blue curve are labelled by $\mac$, we can apply a handle slide in the surface diagram to transform the green curve into a contractible circle.  The contribution of this green circle to $\mathrm{av}_m(T'')$  is then given by
$$
\sum_{c\in I_\mac} \dim(c)\cdot  \mathrm{tr}(1_c)=\sum_{c\in I_\mac} \dim(c)^2=\dim(\mac).
$$
This allows us to replace the green curves of the  first $k$ handles in $T''$ by a factor $\dim(\mac)^k$. 

It remains to compute the contributions of the blue curves of the first $k$ handles in Figure \ref{fig:standardheegard}.  Each such handle corresponds to  four sides of $P$, which are identified pairwise, as shown in Figure \ref{fig:parhandle}. Their contribution is computed in Figure \ref{fig:parhandle} with the diagrammatic calculus for modular fusion categories.

The result is that the images of the simple objects in $I_\cala$ under $\phi$ that label the red line segments are coupled to the tensor unit in $\mac$ on the two sides of each such  blue circle, as shown in Figure \ref{fig:parhandle}. By the cutting strands identity  \eqref{eq:cuttingid} this amounts to placing a circle around the red strands intersecting a blue curve, labelling it with objects in $I_\mac$, multiplying with their dimensions, summing over all such labellings and then multiplying with $\dim(\mac)^{-k}$. This factor cancels the contribution of the green circle. We may thus replace the contribution of each of the first $k$ pairs of green and blue curves  with the one of a black circle evaluated in $\mac$.
Note that the diagrammatic proof in Figure \ref{fig:parhandle} generalises to any number of segments of red curves intersecting the blue curve. Red and green curves that do not intersect the blue curve are unaffected. 

3.~We consider the last $(g-k)$ handles in Figure \ref{fig:standardheegard}, where the green and blue curves form longitude-meridian pairs. Again, each such handle corresponds to four sides of $P$ that are identified as in Figure \ref{fig:orthhandle}. The contribution of such a handle is evaluated in Figure \ref{fig:orthhandle}. 

The result is that the red segments intersecting the blue and green curves can be disconnected  from the handle 
and braid over each other, in such a way that any curve segment intersecting the blue curve crosses over any  curve segment intersecting the green curve. This amounts to arranging the red curve segments on the surface in such a way that those intersecting the blue curve are above the ones intersecting the green curve in Figure \ref{fig:standardheegard} and then erasing the blue and green curve. Again,  this diagrammatic computation can be performed analogously for any number of red edge segments intersecting the blue and green curves. Segments of red curves that do not intersect them  are unaffected.
\end{Proof}

Proposition \ref{prop:bbtri} allows us to  relate the trisection invariant from  Theorem \ref{th:gentrisection}
for the algebraic data in this section to the B\"arenz-Barrett invariant from Theorem \ref{th:bbtheorem}. For this,  recall from Remark \ref{rem:linkdiagram} that the link diagram describing a handlebody decomposition of a 4-manifold $X$ is obtained from a trisection diagram whose green and blue curves are in standard position by the procedure in Proposition \ref{prop:bbtri}. Formula \eqref{eq:trisecCY} thus coincides with the B\"arenz-Barrett invariant up to a rescaling. By rescaling with the constant for the surface diagram $\Sigma'_{st}$ in Example \ref{ex:stabilisationsurf_gen} we then obtain the following theorem.

\begin{Theorem}\label{th:bbtric} Let  $\phi:\cala\to\mac$ be a pivotal functor from a spherical fusion category $\cala$ to a modular fusion category $\mac$.  For any closed, connected 4-manifold $X$ the trisection invariant $|X|$ from Theorem \ref{th:gentrisection} and the B\"arenz-Barrett invariant $I_{BB}(X)$ from Theorem \ref{th:bbtheorem} are related by 
$$
|X|=\dim(\mac)\cdot [\dim(\mac)\cdot \dim(\cala)^\inv\cdot  n(\phi)]^{2/3-\chi(X)/3}\cdot I_{BB}(X).
$$
\end{Theorem}

\begin{Proof} Let $X$ be a 4-manifold given by a $(g,k)$-trisection diagram $T$ whose green and blue curves are in standard position. 
By Theorem \ref{th:gentrisection} the 4-manifold invariant $|X|$ is given by $|X|=\xi^{-g}\cdot \mathrm{av}(T)$, where 
$\xi^3=C_{st}=\dim(\mac)^\inv\cdot \mathrm{av}(\Sigma_{st})$ and
 $\Sigma_{st}$ is the trisection diagram of $S^4$  from \eqref{eq:stabilisation}. This trisection diagram defines the link diagram in Example \ref{ex:trisecstablink}.
Proposition \ref{prop:bbtri} yields the averaged evaluation 
\begin{align}\label{eq:avstd}
\mathrm{av}(\Sigma_{st})
&=\dim(\mac)^2 \cdot\dim(\cala)^2\cdot n(\phi).
\end{align}
Comparing  Theorem \ref{th:bbtheorem} and Proposition \ref{prop:bbtri} and using formula \eqref{eq:trisecresc} then gives 
\begin{align*}
|X|=\xi^{-g}\cdot  \mathrm{av}(T)=\xi^{-g}\cdot \dim(\mac)^{k+1}\cdot \dim(\cala)^{g-k} \cdot n(\phi)^k\cdot I_{BB}(X)\qquad \xi^3=\dim(\mac)\cdot \dim(\cala)^2\cdot n(\phi),
\end{align*}
and Formula \eqref{eq:euler} for the Euler characteristic $\chi(X)$  implies the result.
\end{Proof}

Let us remark that trisection invariants  in this section have an obvious generalisation. Consider a pivotal functor $\psi:\mad\to\mac$ from a spherical fusion category $\mad$ and the $\mad$-right module category structure given by $
\lhd=\oo^{op}(\id\times\psi): \mac\times\mad\to\mac$.
With the $\mac$-module category structure from above, this  defines a transitive $(\mac,\mad)$-bimodule category structure on $\mac$ and a pivotal functor $\phi:\cala\to \End_{(\mac,\mad)}(\mac)$.

Step 2.~in the proof of Proposition \ref{prop:bbtri} could be carried out analogously for this data, but  would give rise to a factor $\dim(\mad)^k\cdot \dim(\mac)^{-k}$. However, step 3.~would no longer give the evaluation of the link diagram from Theorem \ref{th:bbtheorem}, summed over the simple objects of $\cala$ and $\mac$. 
This indicates that the trisection invariants from Theorem \ref{th:gentrisection} should in general be distinct from 
  B\"arenz-Barrett invariants.

\subsection{4-manifold invariants from  graded vector spaces}
\label{subsec:DWinvariant}

The simplest  categorical data for Theorem \ref{th:gentrisection}  that are not given by Hopf triplets arise from bimodule categories over categories of vector spaces graded by finite groups. In this section, we consider the 4-manifold invariants from Theorem \ref{th:gentrisection} for the bimodule categories from Example \ref{ex:modvectgomega}  over  the spherical fusion categories  from Example \ref{ex:vectgomega}. We assume that all cochains in 
Examples \ref{ex:vectgomega} and \ref{ex:modvectgomega} are trivial.  

We first show in Section \ref{subsec:catdata} that these categorical data define  generalisations of  Hopf triplets, where the three Hopf algebras are replaced by  \emph{weak} Hopf algebras and related  by  pairings that satisfy  conditions analogous to the ones in Definitions \ref{def:skewpair} and \ref{def:triplet}.
In Section \ref{subsec:dwcomputation} we compute the associated 4-manifold invariant from Theorem \ref{th:gentrisection}. We show that it coincides with a generalisation of the trisection invariant   from Definition \ref{def:cccinv} for this  \emph{weak} Hopf triplet. We also show that this 4-manifold invariant is not stronger than the ones in \cite{CCC}, as it is obtained from a trisection invariant for group algebras and their duals by a rescaling.

\subsubsection{Categorical data and weak Hopf algebra structure}
\label{subsec:catdata}

We consider the spherical fusion categories
 $\mac=\mathrm{Vec}_C$ and $\mad=\mathrm{Vec}_B$  from Example \ref{ex:vectgomega} for finite groups $C$ and $B$, equipped with trivial cocycles and the standard spherical structure
$$\mac=\mathrm{Vec}_C=\C^C\Mod=\C[C]^*\Mod\qquad\qquad \mad =\mathrm{Vec}_B=\C^B\Mod=\C[B]^*\Mod.$$ 
A $(\mac,\mad)$-bimodule category $\calm$ can be described equivalently as a $\mac\boxtimes\mad^{rev}$-module category and consequently as a module category over
$\mathcal K=\Vec_{K}$, where $K=C\times B^{op}$, also equipped with the trivial cocycle and standard spherical structure.

By Example \ref{ex:modvectgomega} a finite transitive $C\times B^{op}$-set $M$ defines 
an  indecomposable $(\Vec_V,\Vec_B)$-bimodule category 
$
\calm =\Vec_M
$. 
The functors
$\rhd: \Vec_C\times\Vec_M\to \Vec_M$ and $\lhd: \Vec_M\times \Vec_B\to \Vec_M$ are given  by the $C\times B^{op}$-action
on $M$, the  coherence isomorphisms  by the associators of $\mathrm{Vect}_\C$.

We now show that   the associated  spherical fusion category  $\End_{(\mac,\mad)}(\calm)$ is tensor equivalent to the representation category of a \emph{weak} Hopf algebra.
 For basic definitions, see  Appendix \ref{sec:hopfback} and the references therein. We denote by $\C^M$ the algebra of functions on a finite set $M$, with the pointwise multiplication, and by $\langle M\rangle_\C$ the free vector space generated by $M$. We identify $\C^{M\times M}\cong \C^M\oo \C^M$.

\begin{Lemma} \label{lem:weakhopf}Let $K$ be a finite group and  $M$ a finite transitive $K$-set. 
\begin{enumerate}
\item The following defines a weak Hopf algebra structure on  $\C^{M\times M}\oo \C[K]$, denoted $\C^{M\times M}\rtimes\C[K]$,
\begin{align}\label{eq:weakhopf}
&(\delta_{m}\oo\delta_{n}\oo h)\cdot (\delta_{p}\oo\delta_{q}\oo k)
= \delta_{m}(h\rhd p)\delta_{n}(h\rhd q)\,\delta_{m}\oo\delta_{n}\oo hk\\
&1=\Sigma_{m,n\in M} \delta_{m}\oo\delta_{n}\oo 1\nonumber\\
&\Delta(\delta_{m}\oo\delta_{n}\oo k)=\Sigma_{p\in M} (\delta_{m}\oo\delta_p\oo k)\oo (\delta_p\oo \delta_{n}\oo k)\nonumber\\
&\epsilon(\delta_{m}\oo\delta_{n}\oo k)=\delta_{m}(n)\nonumber\\
&S(\delta_{m}\oo\delta_{n}\oo k)= \delta_{k^\inv\rhd n}\oo\delta_{k^\inv\rhd m}\oo k^\inv.\nonumber
\end{align}

\item Its dual is the  weak Hopf algebra  $\langle M\times M\rangle_\C\oo \C^K$ given by
\begin{align}\label{eq:weakhopf2}
&(m\oo n\oo \delta_h)\cdot (p\oo q\oo\delta_k)=\delta_h(k)\delta_{n}(p) m\oo q\oo \delta_k\\
&1=\Sigma_{m\in M}\Sigma_{k\in K} m\oo m\oo\delta_k\nonumber\\
&\Delta(m\oo n\oo\delta_k)=\Sigma_{x\cdot y=k} m\oo n\oo\delta_x\oo (x^\inv\rhd m)\oo (x^\inv\rhd n)\oo\delta_y\nonumber\\
&\epsilon(m\oo n\oo \delta_k)=\delta_k(1)\nonumber\\
&S(m\oo n\oo\delta_k)=(k^\inv\rhd n)\oo (k^\inv\rhd m)\oo\delta_{k^\inv}.\nonumber
\end{align}
\item The weak Hopf algebras $\C^{M\times M}\rtimes\C[K]$ and $\langle M\times M\rangle_\C\oo \C^K$ have  integrals 
\begin{align}\label{eq:ints}
&\lambda=\sum_{m\in M}\sum_{k\in K} \delta_m\oo\delta_m\oo k &  &\ell=\sum_{m,n\in M} m\oo n\oo\delta_1.
\end{align}
\end{enumerate}
\end{Lemma}

\begin{Proof}
All claims  follow by straightforward, but lengthy computations using Definitions  \ref{def:weakhopf} and \ref{def:intweakhopf}. 
\end{Proof}

Note that for a finite  transitive  $K$-set $M$ with at least two elements  the weak Hopf algebras from Lemma \ref{lem:weakhopf} are in general \emph{not} Hopf algebras. Thus, the source of the weakness is the number of elements of $M$ or, put differently, the fact that the associated $\mathrm{Vec}_K$-module category $\mathrm{Vec}_{M\times M}$ has more than one  isomorphism class of simple objects.  

\begin{Proposition} \label{prop:repweak}Let $K$ be a finite group and $M$ a finite transitive $K$-set.  
Then the spherical fusion category $\End_{\Vec_K}(\Vec_M)$  is tensor equivalent to the representation category  $\C^{M\times M}\rtimes\C[K]\Mod$.
\end{Proposition}

\begin{Proof} 
A $\C^{M\times M}\rtimes\C[K]$-module structure
$\rho: \C^{M\times M}\rtimes \C[K]\to \End_\C(V)$
on a finite-dimensional complex vector space $V$ is given by a direct sum decomposition 
\begin{align}\label{eq:splitvrep}
V=\bigoplus_{m_1,m_2\in M} V_{m_1,m_2}\qquad V_{m_1,m_2}=\rho(\delta_{m_1}\oo\delta_{m_2}\oo 1)V
\end{align}
and isomorphisms
$\psi(m_1,m_2,k)=\rho(\delta_{k\rhd m_1}\oo\delta_{k\rhd m_2}\oo k): V_{m_1,m_2}\to V_{k\rhd m_1, k\rhd m_2}$  with
\begin{align}\label{eq:repid}
\psi(k\rhd m_1,k\rhd m_2,h)\circ \psi(m_1,m_2,k)=\psi( m_1,m_2,hk)
\end{align}
for all $h,k\in K$ and $m_1,m_2\in M$. 
The multiplication law in \eqref{eq:weakhopf} implies that the linear maps $\rho(\delta_{m_1}\oo\delta_{m_2}\oo 1): V\to V$ are commuting projectors and the Identity  \eqref{eq:repid}. 

A morphism of $\C^{M\times M}\rtimes\C[K]$-modules from $V$ to $W$ is a $\C$-linear map $f: V\to W$ 
that satisfies 
$$f(V_{m_1,m_2})\subset W_{m_1,m_2},\qquad   
f\circ \psi^V(m_1,m_2,k)=\psi^W(m_1,m_2,k)\circ f\qquad  \forall m_1,m_2\in M, k\in K.$$
The tensor product of two $\C^{M\times M}\rtimes\C[K]$-modules $V$ and $W$ is 
$V\oo W=\Delta(1)\rhd(V\oo_\C W)$
with the $\C^{M\times M}\rtimes\C[K]$-module structure defined via the comultiplication. This yields
\begin{align}\label{eq:tensprodmod}
(V\oo W)_{m_1,m_2}=\bigoplus_{n\in M} V_{m_1,n}\oo W_{n,m_2}.
\end{align}
The linear isomorphisms $\psi^{V\oo W}(m_1,m_2,k)$ are induced by 
\begin{align}\label{eq:psipr}
\psi^V(m_1,n,k)\oo \psi^W(n,m_2,k): V_{m_1,n}\oo W_{n,m_2}\to V_{k\rhd m_1, k\rhd n}\oo W_{k\rhd n, k\rhd m_2}.
\end{align}
The tensor unit is the $\C^{M\times M}\rtimes\C[K]$-module $\C^M$ with 
\begin{align}\label{eq:trivmod}
\C^M_{m_1,m_2}=\delta_{m_1}(m_2) \C\delta_{m_1}\quad\rho(\delta_{m_1}\oo \delta_{m_2}\oo k)\delta_p=\delta_{m_1}(m_2) \delta_{m_1}(k\rhd p) \delta_{m_1}\quad \forall m_1,m_2,p\in M,\, k\in K.
\end{align}
As $M$ is a transitive $K$-set, a $\Vec_K$-module endofunctor $F:\Vec_M\to\Vec_M$ is determined uniquely on the objects and morphisms by the multiplicities $\nu(m_1,m_2,F)$ in the direct sum decomposition
\begin{align}\label{eq:dirsumfunc}
F(m_2)=\bigoplus_{m_1\in M} m_1^{\oplus\nu(m_1,m_2,F)} \qquad \nu(m_1,m_2,F)=\dim_\C V_{m_1,m_2},\quad V_{m_1,m_2}=\Hom_\calm(m_1, F(m_2)),
\end{align}
which satisfy $\nu(m_1,m_2,F)=\nu(k\rhd m_1,k\rhd m_2, F)$ for all $m_1,m_2\in M$ and $k\in K$.

The coherence isomorphism $s^F: F\rhd\Rightarrow \rhd(\id\times F)$ is determined uniquely by its action on the morphism spaces $\Hom_\calm(m_1, F(m_2))$ or, equivalently, by the linear isomorphisms
\begin{align}\label{eq:psicat}
\psi(m_1,m_2,k): \Hom_\calm(m_1, F(m_2))\to \Hom_\calm(k\rhd m_1, F(k\rhd m_2)),\quad \alpha\mapsto s^{F\inv}_{k,m_2}\circ (k\rhd \alpha).
\end{align}
The coherence condition for $s^F$ in the commutative pentagon \eqref{eq:pentagoncfunc}  is then equivalent to \eqref{eq:repid}.
Thus, the coherence isomorphism $s^F$ of a $\Vec_K$-module endofunctor $F:\Vec_M\to\Vec_M$ defines linear isomorphisms \eqref{eq:psicat} that satisfy 
 \eqref{eq:repid}. 
 
 Conversely, given  direct sum decompositions as in \eqref{eq:splitvrep} and for all $k\in K$ and $m_1,m_2\in M$ linear isomorphisms $\psi(m_1,m_2,k): V_{m_1,m_2}\to V_{k\rhd m_1, k\rhd m_2}$ satisfying
 \eqref{eq:repid}, we
 define a $\Vec_K$-module endofunctor $F:\Vec_M\to\Vec_M$ by \eqref{eq:dirsumfunc}. Its  coherence isomorphism $s^F$ is obtained by inserting for $\alpha$ in \eqref{eq:psicat}  the  inclusions that characterise the direct sums in \eqref{eq:splitvrep} and \eqref{eq:dirsumfunc}.

It is then clear from \eqref{eq:trivmod} and \eqref{eq:dirsumfunc}  that the identity functor $\id:\calm\to\calm$ with the trivial coherence isomorphism corresponds to the tensor unit in $\C^{M\times M}\rtimes\C[K]\Mod$. As $\Vec_K$-module endofunctor $F:\Vec_M\to\Vec_M$ are $\C$-linear and respect direct sums, it follows that the composite of module functors  corresponds to the tensor product of $\C^{M\times M}\rtimes\C[K]$-modules in \eqref{eq:tensprodmod}
and \eqref{eq:psipr}.
\end{Proof}

\begin{Corollary}
\label{lem:funcrep}Any $\Vec_K$-module endofunctor $F:\Vec_M\to\Vec_M$  defines  for all $m_1,m_2\in M$ a  representation of $\mathrm{Stab}(m_1,m_2)=\mathrm{Stab}(m_1)\cap \mathrm{Stab}(m_2)$ given by the linear isomorphisms in \eqref{eq:psicat}
$$\rho: \mathrm{Stab}(m_1,m_2)\to \End_\C (V_{m_1,m_2}), \;k\mapsto \psi(m_1,m_2,k).$$
\end{Corollary}

\begin{Remark} Parts of the proof of Proposition \ref{prop:repweak} are fully analogous to  results of Lischka \cite{L}, see in particular the  proof of  \cite[Lemma 5.7]{L} and \cite[Ex 5.11]{L}. However, we also require the monoidal structure of $\End_{\Vec_k}(\Vec_M)$ and a formulation in terms of  weak Hopf algebras.
\end{Remark}

In the following we will also require the simple representations of the weak Hopf algebra $\C^{M\times M}\rtimes \C[K]$. These can be constructed in analogy to the construction of the simple representations of the Drinfeld double of a finite group by Gould  \cite{G}. The following lemma is an adaptation of  Corollary 7.3  in \cite{MM} and the discussion before it, which in turn is a direct generalisation of the construction in \cite{G}.

\begin{Lemma} \label{lem:replem}Let $K$ be a finite group and $M$ a finite transitive  $K$-set.
\begin{compactenum}
\item Simple representations of  $\C^{M\times M}\rtimes \C[K]$ are given by pairs $(\mathcal O, \psi)$ of a $K$-orbit $\mathcal O\subset M\times M$ and a simple representation $\psi: \mathrm{Stab}(m_1,m_2)\to \End_\C(V)$ for some  $(m_1,m_2)\in \mathcal O$. \\[-2ex]

\item The simple representation $\rho_{(\mathcal O,\psi)}$ for  a pair $(\mathcal O,\psi)$ is on  $\langle\mathcal O\rangle_\C\oo V$ and  obtained by
fixing elements $h_{p,q}\in K$  with $h_{p,q}\rhd (m_1,m_2)=(p,q)$ for all $(p,q)\in \mathcal O$ and $h_{m_1,m_2}=1$. It is given by
\begin{align}
\rho_{(\mathcal O,\psi)}(\delta_m\oo\delta_n\oo k) (h_{p,q}\oo v)=\delta_m(p)\delta_n(q)\; h_{k\rhd p,k\rhd q}\oo \psi( h_{k\rhd p, k\rhd q}^\inv k h_{p,q} )v.
\end{align}
\end{compactenum}
\end{Lemma}

We now focus on the group  $K=C\times B^{op}$ for finite groups $B,C$ and the weak Hopf algebras
\begin{align}\label{eq:weaktrip}
H_A=\langle M\times M\rangle_\C\oo \C^K \qquad
H_B=\C^{M\times M}\rtimes\C[B^{op}]\qquad 
H_C=\C^{M\times M}\rtimes\C[C].
\end{align}
The  canonical pairing between $H_A$ and its dual $\C^{M\times M}\rtimes\C[K]$  induces  pairings of $H_A$ with  $H_B$ and $H_C$.
We also obtain a  pairing between $H_B$ and $H_C$, such that the three pairings define the structure of a Hopf triplet. 
For this, note that  Definition \ref{def:skewpair} of a skew pairing and 
Definition \ref{def:triplet} of a Hopf triplet generalise directly to weak Hopf algebras.

\begin{Lemma} \label{lem:hopftripweak}The weak Hopf algebras  from \eqref{eq:weaktrip}  form  a Hopf triplet $(H_A^{cop}, H_B, H_C^{op,cop})$ with 
\begin{align}\label{eq:skewa}
&\tau_{CA}: H_C^{op,cop}\times H_A^{cop}\to \C, & &\tau_{CA}(\delta_p\oo\delta_q\oo c, m\oo n\oo\delta_h)=\delta_h(c)\delta_p(m)\delta_q(n),
\\
&\tau_{AB}: H_A^{cop}\times H_B\to \C, & &\tau_{AB}(m\oo n\oo\delta_h, \delta_p\oo\delta_q\oo b)=\delta_h(b)\delta_p(m)\delta_q(n)\nonumber\\
&\tau_{BC}: H_B\times H_C^{op,cop}\to \C, & &\tau_{BC}(\delta_p\oo\delta_q\oo b, \delta_m\oo\delta_n\oo c)=\delta_n(p)\delta_n(c\rhd q)\delta_p(m\lhd b). \nonumber
\end{align}
\end{Lemma}

\begin{Proof}The conditions in Definitions \ref{def:skewpair} and \ref{def:triplet} are verified by direct, but lengthy computations. 
\end{Proof}

\begin{Remark}\label{rem:reducem} If $M=\{\bullet\}$ is a one-point set, the Hopf triplet from Lemma \ref{lem:hopftripweak} reduces to
$$
H_A=\C^{C\times B^{op}}=\C[C\times B^{op}]^*\qquad H_B=\C[B^{op}]\qquad H_C=\C[C]
$$
with the skew pairings of $H_A$ with $H_B$ and $H_C$ induced by Hopf algebra duality and the skew pairing between $B$ and $C$ given by $\tau_{BC}(b,c)=1$ for all $b\in B$ and $c\in C$.
\end{Remark}

\subsubsection{Computation of the 4-manifold invariants}
\label{subsec:dwcomputation}

In this section we 
compute the evaluation of trisection diagrams for  the categorical data in Section \ref{subsec:catdata}. We   
show that associated 4-manifold invariant from  Definition \ref{th:gentrisection}  coincides with the trisection invariant from Definition \ref{def:cccinv} for the \emph{weak} Hopf triplet  $(H_A^{cop}, H_B, H_C^{op,cop})$ from Lemma \ref{lem:hopftripweak}.
  
We consider a surface diagram $X$  as in Definition \ref{def:surfdiag}, but restrict attention to the case with  no curve endpoints at the boundary. Then  a labelling $l$  of $X$ by the spherical fusion categories $\cala=\End_{(\mac,\mad)}(\calm)$, $\calb=\C^{B^{op}}\Mod$ and $\mac=\C^{C}\Mod$  as in Definition \ref{def:labsurf} is an assignment of
\begin{compactitem}
\item an element $c_\lambda\in C$ to each green curve $\lambda$ of $X$,
\item an element $b_\lambda\in B$ to each blue curve $\lambda$ of $X$,
\item a representation $\rho_\lambda: \C^{M\times M}\rtimes \C[C\times B^{op}]\to \End_\C(V_\lambda)$ to each red curve $\lambda$ of $X$,
\item an element $m\in M$ to the boundary of $X$, if present.
\end{compactitem}

We denote by $\Gamma$ the set of curves on $X$. 
Recall that the \emph{regions} of $X$ are the connected components of the complement of all curves on $X$. 
The \emph{segments} of a curve $\lambda$  on $X$ are the connected components of the complement of its intersection points. When we speak of regions or labels to the left or right of an oriented curve, we always mean left and right viewed in the direction of its orientation.

\begin{Proposition}\label{prop:evalweak}
The evaluation $\mathrm{ev}_l(X)$ of the labelled surface diagram $(X,l)$ is  obtained as follows:
\begin{compactenum}
\item Label the regions of $X$ by  elements of $M$ such that the labelling matches the boundary label. \\[-2ex]
\item Take the product of the $\delta$-functions assigned to each segment of a blue or green curve 
as follows
\begin{align}\label{eq:chorddelta}
\begin{tikzpicture}[scale=.5]
\draw[line width=1pt, green,->,>=stealth] (0,2)--(0,0) node[anchor=north]{$c$};
\node at (1,1)[anchor=west]{$m_2$};
\node at (-1,1)[anchor=east]{$m_1$};
\node at (0,-1)[anchor=north]{$\delta_{m_1}(c\rhd m_2)$};
 \end{tikzpicture}
 \qquad\qquad
 \begin{tikzpicture}[scale=.5]
\draw[line width=1pt, blue,->,>=stealth] (0,2)--(0,0) node[anchor=north]{$b$};
\node at (1,1)[anchor=west]{$m_2$};
\node at (-1,1)[anchor=east]{$m_1$};
\node at (0,-1)[anchor=north]{$\delta_{m_1}(m_2\lhd b)$};
 \end{tikzpicture}
\end{align}
 \item For each red curve $\lambda$ choose a basepoint and 
assign
\begin{compactitem}
\item to the basepoint the endomorphism $\rho_\lambda(\delta_{m_1}\oo\delta_{m_2}\oo 1)\in \End_\C(V_\lambda)$, where $m_1\in M$ and $m_2\in M$ are the labels of the regions to the right and left of the basepoint,\\[-2ex]

\item to an intersection point $p$ of $\lambda$ with a green curve labelled $c$ or a blue curve labelled $b$, the endomorphisms $\rho_\lambda(1\oo 1\oo c^{-\epsilon_p})$ or $\rho_\lambda(1\oo 1\oo b^{\epsilon_p})$, where $\epsilon_p$ is the sign  from \eqref{eq:intsign}.\\[-2ex]
\end{compactitem}
Compose these endomorphisms in the order of the intersection points on $\lambda$, starting with the endomorphism at the basepoint.  
Multiply by the trace of the resulting endomorphism of $V_\lambda$.\\[-2ex]

\item Sum the resulting number over all assignments of elements of $M$ to the regions of $X$, while keeping the boundary label fixed, if present.
\end{compactenum}
\end{Proposition}

\begin{Proof}
We show that the number obtained by this procedure coincides with the evaluation $\mathrm{ev}_l(X)$ from Definition \ref{def:surfdiagev}.
For this, note   that labelling the regions of $X$ by elements of $M$ and cutting $X$ along a set of cutting generators yields a labelled tricoloured chord diagram with a labelling of the boundary segments by $M$. For each chord endpoint of a blue or green curve the associated morphism space from Definition \ref{def:surfdiagev}  is either 0 or $\C$ and given by the factors in \eqref{eq:chorddelta}. Thus, summing over bases of these morphism spaces corresponds to inserting the factors from  \eqref{eq:chorddelta} in step 2. 

For a chord endpoint of a red curve $\lambda$ with a label $m_1\in M$ on the right and $m_2\in M$ on the left the associated morphism space is $V^\lambda_{m_1,m_2}=\rho_\lambda(\delta_{m_1}\oo\delta_{m_2}\oo 1)V_\lambda$  or its dual
\begin{align*}
\begin{tikzpicture}[scale=.5]
\draw[line width=1pt] (0,3) node[anchor=west]{$m_1$}--(0,0) node[anchor=west]{$m_2$};
\draw[line width=1pt, color=red,](0,1.5)--(-1.5,0) 
node[sloped, pos=0.5, allow upside down]{\arrowOut}
node[pos=.5,anchor=south east]{$\rho_\lambda$};
\node at (0,-1) [anchor=north]{$V^\lambda_{m_1,m_2}$};
\end{tikzpicture}
\qquad
\begin{tikzpicture}[scale=.5]
\draw[line width=1pt] (0,3) node[anchor=west]{$m_2$}--(0,0) node[anchor=west]{$m_1$};
\draw[line width=1pt, color=red,](-1.5, 3)--(0,1.5) 
node[sloped, pos=0.5, allow upside down]{\arrowOut}
node[pos=.5,anchor=north east]{$\rho_\lambda$};
\node at (0,-1) [anchor=north]{$V^{\lambda*}_{m_1, m_2}$};
\end{tikzpicture}
\qquad
\begin{tikzpicture}[scale=.5]
\draw[line width=1pt] (0,3) node[anchor=west]{$m_2$}--(0,0) node[anchor=west]{$m_1$};
\draw[line width=1pt, color=red,](-1.5,0)--(0,1.5) 
node[sloped, pos=0.5, allow upside down]{\arrowOut}
node[pos=.5,anchor=south east]{$\rho_\lambda$};
\node at (0,-1) [anchor=north]{$V^{\lambda *}_{m_1,m_2}$};
\end{tikzpicture}
\qquad
\begin{tikzpicture}[scale=.5]
\draw[line width=1pt] (0,3) node[anchor=west]{$m_1$}--(0,0) node[anchor=west]{$m_2$};
\draw[line width=1pt, color=red,](0,1.5)--(-1.5, 3)
node[sloped, pos=0.5, allow upside down]{\arrowOut}
node[pos=.5,anchor=north east]{$\rho_\lambda$};
\node at (0,-1) [anchor=north]{$V^\lambda_{m_1, m_2}$};
\end{tikzpicture}
\end{align*}
Summing over these morphism spaces  amounts to summing over dual bases of $V^\lambda_{m_1,m_2}$
and $V^{\lambda*}_{m_1,m_2}$.

By \eqref{pic:bimodule}  crossings between blue and green curves are not labelled by any additional data.  Also by \eqref{pic:bimodule} and by construction of the representation in \eqref{eq:psicat} a crossing between a red curve $\lambda$ and a green or blue curve labelled by $c\in C$ or $b\in B$ is labelled as follows
\begin{align}\label{eq:orcross}
\begin{tikzpicture}[scale=.5]
\draw[line width=1pt, color=red] (2,2)node[anchor=south west]{$\rho_\lambda$} --(0,0) node[sloped, pos=0.3, allow upside down]{\arrowOut} ; 
\draw[line width=1pt, color=green] (0,2) node[anchor=south east]{$c$}--(2,0) node[sloped, pos=0.3, allow upside down]{\arrowOut};
\node at (0,1)[anchor=east]{$m_1$};
\node at (1,0)[anchor=north]{$m_2$};
\node at (1,-1) [anchor=north] {$\rho_\lambda(\delta_{m_1}\oo\delta_{m_2}\oo c)$};
\end{tikzpicture}
\quad
\begin{tikzpicture}[scale=.5]
\draw[line width=1pt, color=green] (2,2)node[anchor=south west]{$c$} --(0,0) node[sloped, pos=0.3, allow upside down]{\arrowOut} ; 
\draw[line width=1pt, color=red] (0,2) node[anchor=south east]{$\rho_\lambda$}--(2,0) node[sloped, pos=0.3, allow upside down]{\arrowOut};
\node at (2,1)[anchor=west]{$m_2$};
\node at (1,0)[anchor=north]{$m_1$};
\node at (1,-1) [anchor=north] {$\rho_\lambda(\delta_{m_1}\oo\delta_{m_2}\oo c^\inv)$};
\end{tikzpicture}
\quad
\begin{tikzpicture}[scale=.5]
\draw[line width=1pt, color=red] (2,2)node[anchor=south west]{$\rho_\lambda$} --(0,0) node[sloped, pos=0.3, allow upside down]{\arrowOut} ; 
\draw[line width=1pt, color=blue] (0,2) node[anchor=south east]{$b$}--(2,0) node[sloped, pos=0.3, allow upside down]{\arrowOut};
\node at (0,1)[anchor=east]{$m_1$};
\node at (1,0)[anchor=north]{$m_2$};
\node at (1,-1) [anchor=north] {$\rho_\lambda(\delta_{m_1}\oo\delta_{m_2}\oo b)$};
\end{tikzpicture}
\quad
\begin{tikzpicture}[scale=.5]
\draw[line width=1pt, color=blue] (2,2)node[anchor=south west]{$c$} --(0,0) node[sloped, pos=0.3, allow upside down]{\arrowOut} ; 
\draw[line width=1pt, color=red] (0,2) node[anchor=south east]{$\rho_\lambda$}--(2,0) node[sloped, pos=0.3, allow upside down]{\arrowOut};
\node at (2,1)[anchor=west]{$m_2$};
\node at (1,0)[anchor=north]{$m_1$};
\node at (1,-1) [anchor=north] {$\rho_\lambda(\delta_{m_1}\oo\delta_{m_2}\oo b^\inv)$};
\end{tikzpicture}
\end{align}
The linear endomorphisms at different crossings of a given red chord are composed according to \eqref{eq:repid} and \eqref{eq:psicat}, which yields the product of the associated endomorphisms of $V_\lambda$. 
As the basepoint of the red curve is associated with a factor $\rho_\lambda(\delta_{m_1}\oo \delta_{m_2}\oo 1)$ in the product from step 4, Identity \eqref{eq:weakhopf}  allows us to replace these factors by $\rho_\lambda(1\oo 1\oo c^{\mp\epsilon_p})$ and $\rho_\lambda(1\oo 1\oo b^{\pm \epsilon_p})$, where $\epsilon_p$ is defined as in \eqref{eq:intsign}.
The summation over the morphism spaces at the chord endpoints composes the endomorphisms for  different chords representing a given red curve and takes the trace of the resulting endomorphism of $V_\lambda$.

The summation over the elements of $M$ labelling the different segments at the vertical line of the chord diagram corresponds to the summation over the labelling with simple objects in  $\calm$ in Definition \ref{def:surfdiagev}. Labels associated to regions that do not intersect the vertical line of the chord diagram correspond to the composition of functors labelling the chords at the boundaries of these regions, see in particular formula \eqref{eq:tensprodmod}. For the blue and green curves, this summation is absorbed by the factors from \eqref{eq:chorddelta}. 
\end{Proof}

With the evaluation of a labelled surface diagram from Proposition \ref{prop:evalweak} we can compute its averaged evaluation. This requires a summation over the assignments of elements of $C$ and $B$ to green and blue curves and over the simple representations of the weak Hopf algebra $\C^{M\times M}\rtimes \C[K]$ associated to the red curves. The latter reduces to a condition on the elements of $C$ and $B$ assigned to its crossings.
This condition and  Condition \eqref{eq:chorddelta} are captured in the following definition.

\begin{Definition} \label{def:admisslabel} Let $X$ be a surface diagram with no curve endpoints on the boundary. A labelling
\begin{compactitem}
\item of the green and blue curves of $X$ by elements of $C$ and $B$, 
\item of each region of $X$ by an element of $M$,
\end{compactitem}
 that matches the boundary label, if present, is called \textbf{admissible}, if 
\begin{compactenum}[(i)]
\item for each  segment of a blue or green curve labelled $c\in C$ or $b\in B$  with adjacent labels $m_1,m_2\in M$  one has  $m_1=c\rhd m_2$ or $m_1=m_2\lhd b$, as in \eqref{eq:chorddelta},
\item for each red curve $\lambda$ the signed ordered product $k_\lambda\in K$ of the elements of $C$ and $B$ associated with its intersection points  as in \eqref{eq:orcross} vanishes.
 \end{compactenum}
\end{Definition}

Note that  the definition of $k_\lambda$ in (ii) requires the choice of an orientation and of a basepoint for the red curve on $X$. The signed oriented product is then given by multiplying the elements $c^{-\epsilon_p}$ and $b^{\epsilon_p}$ for each intersection point $p$ of $\lambda$  in the order in which they appear on $\lambda$. As different choices of basepoints yield conjugates of $k_\lambda$ and orientation reversal its inverse,  Condition  (ii) in Definition \ref{def:admisslabel} is independent of the choice of the orientation and the basepoint.

\begin{Proposition} \label{prop:aveval} Let $X$ be a surface diagram with no curve endpoints on the boundary and with boundary label $m\in M$, if $X$ has a boundary.  Then the  averaged evaluation  of  $X$ is
$$
\mathrm{av}_{(m)}(X)=|l_X|\cdot |B|^r\cdot |C|^r,
$$
where $|l_X|$ is the number of  admissible labellings  and  $r$  the number of red curves of $X$.  
\end{Proposition}

\begin{Proof} The averaged evaluation of  $X$ is obtained by summing  its evaluation over all assignments of elements in $C$ and $B$ to the green and blue curves of $X$, respectively, and by summing over all assignments of simple representations of $\C^{M\times M}\rtimes\C[K]$ to its red curves. 
By Proposition \ref{prop:evalweak} its evaluation vanishes, unless Condition (i) is met.  By Proposition \ref{prop:evalweak} crossings between green and blue curves do not contribute to the evaluation. Each red curve $\lambda$ on $X$ contributes a factor
\begin{align*}
\sum_{\rho_\lambda} \dim_\C(\rho_\lambda) \cdot \chi_{\rho_\lambda}(k_\lambda),
\end{align*}
where $\rho_\lambda$ runs over all simple representations of $\C^{M\times M}\rtimes\C[K]$ and $\chi_{\rho_\lambda}=\mathrm{tr}(\rho_\lambda)$ is its character. 

By Lemma \ref{lem:replem} the simple representations of $\C^{M\times M}\rtimes\C[K]$ are given by pairs $(\mathcal O_\lambda, \psi_\lambda)$ of a $K$-orbit $\mathcal O_\lambda\subset M\times M$ and an irreducible representation $\psi_\lambda$ of its stabiliser group. For a fixed labelling of the regions of $X$  by elements of $M$ with the label $(m_1,m_2)\in M\times M$ at the basepoint of $\lambda$, only the orbit $\mathcal O$ with $(m_1,m_1)\in \mathcal O$ contributes. We then  sum over  simple representations $\psi_\lambda$ of  the stabiliser group  $\mathrm{Stab}(m_1,m_2)$. Also by Lemma \ref{lem:replem}, we have $\dim_\C(\rho_\lambda)=|\mathcal O|\cdot \dim_\C (\psi_\lambda)$ and $\chi_{\rho_\lambda}(k_\lambda)=\chi_{\psi_\lambda}(k_\lambda)$ for $k_\lambda\in \mathrm{Stab}(m_1,m_2)$ and $\chi_{\rho_\lambda}(k_\lambda)=0$ else. This gives
 \begin{align*}
\sum_{\rho_\lambda} \dim_\C(\rho_\lambda)\cdot  \chi_{\rho_\lambda}(k_\lambda)&= |\mathcal O|\cdot \sum_{\psi_\lambda} \dim_\C(\psi_\lambda)\cdot  \chi_{\psi_\lambda}(k_\lambda)
=|\mathcal O|\cdot |\mathrm{Stab}(m_1,m_2)|\cdot \delta_1(k_\lambda)\\
&=|K|\cdot \delta_1(k_\lambda)=|B|\cdot |C|\cdot \delta_1(k_\lambda),
\end{align*}
where we used the orthogonality relation for the characters of $\mathrm{Stab}(m_1,m_2)$ and then the orbit-stabiliser formula. This yields Condition (ii) in Definition \ref{def:admisslabel} and an overall factor $|B|^r\cdot |C|^r$.
\end{Proof}

As a direct consequence of  Proposition \ref{prop:aveval} we obtain a combinatorial formula for the 4-manifold invariant from Theorem \ref{th:gentrisection} for the $(\Vec_C,\Vec_B)$-bimodule category $\Vec_M$ and  $\Phi=\id$.

\begin{Corollary} \label{cor:dwtrisec}The  trisection invariant from Theorem \ref{th:gentrisection} for a trisection diagram $T$ is 
$$
|T|=|l_T| \cdot |B|^{-g/3}\cdot |C|^{-g/3},
$$
where $g$ is the genus of $T$ and  $|l_T|$ is the number  of admissible labellings of  $T$. 
\end{Corollary}

\begin{Proof} For a trisection diagram of genus $g$, Proposition \ref{prop:aveval} yields $\mathrm{av}(T)=|l_T|\cdot |B|^g\cdot |C|^g$, as  $r=g$. For the surface diagram $\Sigma'_{st}$ from Example \ref{ex:stabilisationsurf_gen}, we obtain $\mathrm{av}_m(\Sigma'_{st})=|l_{st}|\cdot |B|^3\cdot |C|^3$, where $|l_{st}|$ is the number of admissible labellings of $\Sigma'_{st}$. 
This yields $\xi=\sqrt[3]{|l_{st}|} \cdot |B|\cdot |C|$ in Theorem \ref{th:gentrisection} and 
$$|T|=\xi^{-g} \cdot \mathrm{av}(T)=|l_T|\cdot |l_{st}|^{-g/3}.$$
By inspecting diagram \eqref{eq:chordstablab} one finds that every admissible labelling of $\Sigma'_{st}$ for a given boundary label $m$ is characterised by the conditions  $c_1=b_2=c_3=b_3=1$, $m_1=m\lhd b_1^\inv$, $m_2=c_2\rhd m$ and $m_3=m$. Thus admissible labellings of $\Sigma'_{st}$
are in bijection with elements  of $C\times B$ and $|l_{st}|=|B|\cdot |C|$. 
\end{Proof}

 With Corollary \ref{cor:dwtrisec} we can relate the 4-manifold invariant from Theorem \ref{th:gentrisection} to the trisection invariant from Definition \ref{def:cccinv} for the weak Hopf triplet $(H_A^{cop}, H_B, H_C^{op,cop})$ from Lemma \ref{lem:hopftripweak}.
 For this, note that  
 Definition \ref{def:trisecbr} of the trisection bracket generalises directly to weak Hopf algebras, as it makes use only of the comultiplication, the integrals and the skew pairings.  
 
 To obtain agreement with the 4-manifold invariant from Theorem \ref{th:gentrisection}, we only require a slight modification to Definition \ref{def:cccinv}, namely 
 that the rescaling  $I(T)=\xi^{-g}\cdot \langle T\rangle$ from Definition \ref{def:cccinv} is with a constant $\xi$ that satisfies
  $\xi^3=\langle \Sigma_{st}\rangle \cdot |M|^\inv$, where $\Sigma_{st}$ is the trisection diagram of $S^4$ from \eqref{eq:stabilisation}.
  If we modify the definition of the Hopf algebraic trisection invariant in this way for  weak Hopf triplets, we obtain the following proposition.

\begin{Proposition} \label{th:coincide}For any trisection diagram $X$ the trisection invariant from  Theorem \ref{th:gentrisection} coincides with the trisection invariant from Definition \ref{def:cccinv} for the weak Hopf triplet $(H^{cop}_A,H_B,H^{op,cop}_C)$.
\end{Proposition}

\begin{Proof} We visualise elements of the weak Hopf algebras $\langle M\times M\rangle_\C\oo \C^K$ and $\C^{M\times M}\oo \C[K]$ by oriented red and cyan lines as follows
\begin{align}\label{eq:linelabels}
&\begin{tikzpicture}[scale=.8]
\draw[line width=1pt, red, ->,>=stealth] (0,1.5)--(0,0) node[anchor=north]{$h$};
\node at (.2,.2)[anchor=west]{$n$};
\node at (-.2,.2)[anchor=east]{$m$};
\node at (.2,1.3)[anchor=west]{$h^\inv\rhd n$};
\node at (-.2,1.3)[anchor=east]{$h^\inv\rhd m$};
\end{tikzpicture}
&
&\begin{tikzpicture}[scale=.8]
\draw[line width=1pt, cyan, ->,>=stealth] (0,1.5)--(0,0) node[anchor=north]{$k$};
\node at (-.2,1.3)[anchor=east]{$p$};
\node at (-.2,.2)[anchor=east]{$q$};
\node at (.2,1.3)[anchor=west]{$k^\inv\rhd p$};
\node at (.2,.2)[anchor=west]{$k^\inv\rhd q$};
\end{tikzpicture}
\\
&\qquad\quad m\oo n\oo \delta_{h}
&
&\delta_p\oo\delta_q\oo k.\nonumber
\end{align}
For simplicity we often omit the labels at the start of a red line and on the left side of a cyan line.
The multiplications of the weak Hopf algebras $\langle M\times M\rangle_\C\oo \C^K$ and $\C^{M\times M}\oo \C[K]$ are then given by
\begin{align*}
&\begin{tikzpicture}[scale=1]
\draw[line width=1pt, red, ->,>=stealth] (0,1.5)--(0,0) node[anchor=north]{$h$};
\node at (.2,.2)[anchor=west]{$n$};
\node at (-.2,.2)[anchor=east]{$m$};
\node at (.75,.75){$\cdot$};
\draw[line width=1pt, red, ->,>=stealth] (1.5,1.5)--(1.5,0) node[anchor=north]{$k$};
\node at (1.7,.2)[anchor=west]{$q$};
\node at (1.3,.2)[anchor=east]{$p$};
\node at (2,.75)[anchor=west]{$= \delta_h(k)\delta_n(p)$};
\draw[line width=1pt, red, ->,>=stealth] (5,1.5)--(5,0) node[anchor=north]{$h$};
\node at (5.2,.2)[anchor=west]{$q$};
\node at (4.8,.2)[anchor=east]{$m$};
\end{tikzpicture}
\qquad\qquad
&\begin{tikzpicture}[scale=1]
\draw[line width=1pt, cyan, ->,>=stealth] (0,1.5)--(0,0) node[anchor=north]{$h$};
\node at (-.2,1.3)[anchor=east]{$m$};
\node at (-.2,.2)[anchor=east]{$n$};
\node at (.75,.75){$\cdot$};
\draw[line width=1pt, cyan, ->,>=stealth] (1.5,1.5)--(1.5,0) node[anchor=north]{$k$};
\node at (1.3,1.3)[anchor=east]{$p$};
\node at (1.3,.2)[anchor=east]{$q$};
\node at (2,.75)[anchor=west]{$= \delta_m(h\rhd p)\delta_n(h\rhd q)$};
\draw[line width=1pt, cyan, ->,>=stealth] (6.5,1.5)--(6.5,0) node[anchor=north]{$hk$};
\node at (6.3,1.3)[anchor=east]{$m$};
\node at (6.3,.2)[anchor=east]{$n$};
\end{tikzpicture}
\end{align*}
This states that the multiplications of the labels at the tip of the red and cyan arrows are the  multiplications of $\C^K$ and $\C[K]$, respectively, and enforces that the labels in the regions between the lines match. The comultiplication takes the form
\begin{align}\label{eq:diagcomult}
&\begin{tikzpicture}[scale=1]
\draw[line width=1pt,->,>=stealth, red] (0,1.5)--(0,0) node[anchor=north]{$h$};
\node at (-.2,.2)[anchor=east]{$m$};
\node at (.2,.2)[anchor=west]{$n$};
\node at (.7,.75)[anchor=west]{$\xrightarrow{\Delta}$};
\draw[line width=1pt,->,>=stealth, red] (5,2.25)--(5,1) node[anchor=north]{$y$};
\draw[line width=1pt,->,>=stealth, red] (5,.5) --(5,-.75) node[anchor=north]{$x$};
\node at (1.5,.7) [anchor=west]{$\sum_{x\cdot y=h}$};
\node at (4.8,-.55)[anchor=east]{$m$};
\node at (5.2,-.55)[anchor=west]{$n$};
\node at (4.8,1.2)[anchor=east]{$x^\inv\rhd m$};
\node at (5.2,1.2)[anchor=west]{$x^\inv\rhd n$};
\end{tikzpicture}
\qquad\quad
&\begin{tikzpicture}[scale=1]
\draw[line width=1pt,->,>=stealth, cyan] (0,1.5)--(0,0) node[anchor=north]{$k$};
\node at (-.2,.2)[anchor=east]{$q$};
\node at (-.2,1.3)[anchor=east]{$p$};
\node at (.7,.75)[anchor=west]{$\xrightarrow{\Delta}$};
\draw[line width=1pt,->,>=stealth, cyan] (4,2.25)--(4,1) node[anchor=north]{$k$};
\draw[line width=1pt,->,>=stealth, cyan] (4,.5) --(4,-.75) node[anchor=north]{$k$};
\node at (1.5,.7) [anchor=west]{$\sum_{r\in M}$};
\node at (3.8,-.55)[anchor=east]{$q$};
\node at (3.8,.3)[anchor=east]{$r$};
\node at (3.8,1.2)[anchor=east]{$r$};
\node at (3.8,2.05)[anchor=east]{$p$};
\end{tikzpicture}
\end{align}
Note that this is the comultiplication of $\C^K$ and $\C[K]$ for the labels at the tips of the red and cyan arrows, respectively, and that it ensures that the labels next to these arrows match.
The antipodes  correspond to  reversal of the arrows
\begin{align}\label{eq:diaganti}
&\begin{tikzpicture}[scale=.8]
\draw[line width=1pt, red, ->,>=stealth] (0,0)--(0,1.5) node[anchor=south]{$h^\inv$};
\node at (.2,1.3)[anchor=west]{$h^\inv\rhd n$};
\node at (-.2,1.3)[anchor=east]{$h^\inv\rhd m$};
\end{tikzpicture}
&
&\begin{tikzpicture}[scale=.8]
\draw[line width=1pt, cyan, ->,>=stealth] (0,0)--(0,1.5) node[anchor=south]{$k^\inv$};
\node at (.2,1.3)[anchor=west]{$k^\inv\rhd p$};
\node at (.2,.2)[anchor=west]{$k^\inv\rhd q$};
\end{tikzpicture}
\\
&\qquad S(m\oo n\oo \delta_h)
&
&S(\delta_p\oo \delta_q\oo k), \nonumber
\end{align}
and the integrals are given by
\begin{align}\label{eq:ella}
&\begin{tikzpicture}[scale=.8]
\node at (-1,.5)[anchor=east] {$\ell_A=\sum_{m,n\in M}$};
\draw[line width=1pt, red, ->,>=stealth] (0,1.5)--(0,0) node[anchor=north]{$1$};
\node at (.2,.2)[anchor=west]{$n$};
\node at (-.2,.2)[anchor=east]{$m$};
\end{tikzpicture}
&
&\begin{tikzpicture}[scale=.8]
\node at (-1,.5)[anchor=east] {$\lambda=\sum_{m\in M, k\in K}$};
\draw[line width=1pt, cyan, ->,>=stealth] (0,1.5)--(0,0) node[anchor=north]{$k$};
\node at (-.2,1.3)[anchor=east]{$m$};
\node at (-.2,.2)[anchor=east]{$m$};
\end{tikzpicture}
\end{align}
For the weak Hopf subalgebras $\C^{M\times M}\rtimes\C[C]$ and $\C^{M\times M}\rtimes \C[B^{op}]$ of $\C^{M\times M}\rtimes\C[K]$  we use an analogous diagrammatic calculus, but with green and blue lines, respectively, instead of cyan lines and with the labels at the tip of the arrows restricted to $C$ and $B$. Their  integrals are given by
\begin{align}\label{eq:ellb}
&\begin{tikzpicture}[scale=.8]
\node at (-1,.5)[anchor=east] {$\ell_C=\sum_{m\in M, c\in C}$};
\draw[line width=1pt, green, ->,>=stealth] (0,1.5)--(0,0) node[anchor=north]{$c$};
\node at (-.2,1.3)[anchor=east]{$m$};
\node at (-.2,.2)[anchor=east]{$m$};
\end{tikzpicture}
\qquad\qquad
\begin{tikzpicture}[scale=.8]
\node at (-1,.5)[anchor=east] {$\ell_B=\sum_{m\in M, b\in B}$};
\draw[line width=1pt, blue, ->,>=stealth] (0,1.5)--(0,0) node[anchor=north]{$b$};
\node at (-.2,1.3)[anchor=east]{$m$};
\node at (-.2,.2)[anchor=east]{$m$};
\end{tikzpicture}
\end{align}
The duality pairing between $\langle M\times M\rangle_\C\oo \C^K$ and $\C^{M\times M}\rtimes\C[K]$ and the skew pairing between $H_B$ and $H_C^{op,cop}$ from Lemma \ref{lem:hopftripweak} are then visualised by the following crossing diagrams
\begin{align}\label{eq:crossings}
&\begin{tikzpicture}[scale=.5]
\draw[line width=1pt, red, ->,>=stealth] (2,2)--(-2,-2) node[anchor=north east]{$h$};
\draw[line width=1pt, cyan, ->,>=stealth] (-2,2)--(2,-2) node[anchor=north west]{$k$};
\node at (-2,1.4) [anchor=east]{$p$};
\node at (-2.2,-1.6) [anchor=east]{$m$};
\node at (1.2,-1.6) [anchor=east]{$q$};
\node at (-1.2,-1.6) [anchor=west]{$n$};
\end{tikzpicture}
&
&\begin{tikzpicture}[scale=.5]
\draw[line width=1pt, green, ->,>=stealth] (2,2)--(-2,-2) node[anchor=north east]{$c$};
\draw[line width=1pt, blue, ->,>=stealth] (-2,2)--(2,-2) node[anchor=north west]{$b$};
\node at (-2,1.4) [anchor=east]{$p$};
\node at (-2.2,-1.6) [anchor=east]{$n$};
\node at (1.2,-1.6) [anchor=east]{$q$};
\node at (1.2, 1.6) [anchor=east]{$m$};
\end{tikzpicture}
\\
&\langle m\oo n\oo \delta_h, \delta_p\oo \delta_q\oo k\rangle
&
&\langle \delta_m\oo \delta_n\oo c, \delta_p\oo \delta_q\oo b\rangle\nonumber
\\
&=\delta_m(p)\delta_n(q)\delta_h(k)
&
&=\delta_n(p)\delta_n(c\rhd q)\delta_p(m\lhd b).\nonumber
\end{align}
They state that the labels at the tips of the arrows must coincide for each crossing of a red with a green, blue or cyan line and that the labels in each region of the diagram coincide. 
Due to the properties of the pairings,  the two interpretations of a diagram in which a line crosses two parallel lines of  a different colour coincide. Such a diagram can either be realised by applying a comultiplication to the first line or a multiplication to the pair of lines.  Crossings with opposite signs are obtained by reversing line orientations with the antipode, as in \eqref{eq:diaganti}.

We now consider a trisection diagram $X$ and construct the associated trisection bracket for the weak Hopf triplet $(H_A^{cop}, H_B, H_C^{op,cop})$
 from  Definition \ref{def:trisecbr}. We choose for every curve on $X$ an orientation and a basepoint outside its intersection points. 
We assign to every red curve the  integral $\ell_A$ from \eqref{eq:ella} and to every  blue and green curve the  integral $\ell_B$ and $\ell_C$  from \eqref{eq:ellb}. 

For a curve $\lambda$ with $n_\lambda$ intersection points with curves of different colours, applying
an $(n_\lambda-1)$-fold comultiplication to the associated integral  as in Step 3.~in Definition \ref{def:trisecbr}
  amounts to subdividing $\lambda$ into $n_\lambda$ segments with labels given by  the diagrammatic comultiplication formula \eqref{eq:diagcomult}, applied to the integral. 
We choose these segments in such a way that each intersection point lies in one such segment. 

Applying powers of antipodes according to the sign of the intersections and then the skew pairings as in Steps 4.~and 5.~of Definition \ref{def:trisecbr} amount to imposing identities \eqref{eq:crossings}, possibly with a line orientation reversal according to the signs of the crossings. 

From formula \eqref{eq:ellb} and the second identity in \eqref{eq:diagcomult} it is apparent that the integrals of $H_C$ and $H_B$ assigned to the green and blue curves give rise to a summation over all possible labellings of green curves by elements of $C$ and blue curves by elements of $B$. They also enforce that the labelling by $M$ at the two segments containing the basepoint agree.  

Likewise, formula \eqref{eq:ella},  the first identity in \eqref{eq:diagcomult} and the first identity in \eqref{eq:crossings} show that 
each  integral associated to a red curve of $X$ enforces Condition (ii) in Definition \ref{def:admisslabel}. Condition (i) in Definition \ref{def:admisslabel} is enforced by the labelling of the green and blue lines and by the two identities for the crossings in \eqref{eq:crossings}, which force  the labels in each region of $X$ to agree. 

Hence, the trisection bracket from  Definition \ref{def:trisecbr} for the weak Hopf triplet $(H_A^{cop}, H_B, H_C^{op,cop})$ is the number of admissible labellings of $X$ from Definition \ref{def:admisslabel}. The trisection invariant is  obtained by rescaling it with
$\xi^{-g}$, where $\xi^3=\langle \Sigma_{st}\rangle\cdot |M|^\inv=|B|\cdot |C|$.
By Corollary \ref{cor:dwtrisec} this gives the trisection invariant from Theorem \ref{th:gentrisection}.
\end{Proof}

Theorem \ref{th:coincide} shows that the trisection invariants defined by Hopf triplets from \cite{CCC} indeed have generalisations in terms of weak Hopf algebras, as anticipated in \cite[Rem 5.10]{CCC}. 
Interestingly, all of the three Hopf algebras involved in the triplet become weak, although the underlying spherical fusion categories $\Vec_C$ and $\Vec_B$ admit a description as representation categories of  the \emph{strong} Hopf algebras $\C^C$ and $\C^{B}$. 
However, the resulting trisection invariant is 
 directly related to the one for the strong Hopf triplet in Remark \ref{rem:reducem}.

\begin{Corollary} For any trisection diagram $T$ 
the trisection invariant $|T|$ for  the weak Hopf triplet $(H_A^{cop}, H_B, H_C^{op,cop})$ and
 the trisection invariant $I_{CCC}(T)$  for the strong Hopf triplet from Remark \ref{rem:reducem} are related by a rescaling:
$$
|T|=|M|\cdot I_{CCC}(T).
$$
\end{Corollary}

\begin{Proof} We show that the number $l_T$ of admissible labellings of $T$  is given by
\begin{align}\label{eq:admisslabel}
|l_T|=|M|\cdot |l_T^{B,C}|
\end{align}
where $|l_T^{B,C}|$ is the number of labellings of the  blue and green curves of $T$ by elements of $B$ and $C$ that satisfy Condition (ii) in Definition \ref{def:admisslabel}.
As $|l_T^{B,C}|$ is the number of admissible labellings for the case where $|M|=1$, the claim then follows from Corollary \ref{cor:dwtrisec}.

To prove \eqref{eq:admisslabel}, choose a labelling $l_T^{B,C}$ of the blue and green curves of $T$ such that  Condition (ii) in Definition \ref{def:admisslabel} is satisfied. 
Cutting  $T$ along all red curves yields a surface diagram $D$ on the sphere with $2g$ discs removed, where $g$ is the genus of $T$. Admissible labellings of $T$ that extend $l_T^{B,C}$  are  in bijection with labellings  of the regions of $D$ with elements of $M$  that satisfy  (i) in Definition \ref{def:admisslabel}. 

For each element  $m_r\in M$ assigned to a fixed region $r$ of $D$ , there is at most one such labelling. This follows, because every region of $D$ can be reached from $r$ by crossing green and blue lines, and the labels on the regions to the left and right of a green and blue lines are related by the action of the group elements on the line.

To show that there is at least one such labelling, note that  crossings of green and blue curves do not give rise to any additional conditions on $m_r$. Thus, additional conditions on  $m_r$ can only arise from sequences of regions of $D$ associated to non-contractible paths in $D$. Each such path is a composite of paths that go once around a  boundary circle. Condition (ii) in Definition \ref{def:admisslabel} ensures that this imposes no restrictions on $m_r$. Thus, there are exactly $|M|$ admissible labellings of $T$ that extend $l_{T}^{B,T}$. 
\end{Proof}

\section*{Acknowledgements}

The authors would like to thank the Isaac Newton Institute for Mathematical Sciences, Cambridge, for support and hospitality during the programme The Physics and Mathematics of Boundaries, Impurities, and Defects, where work on this paper was undertaken. This work was supported by EPSRC grant EP/Z000580/1. V.\,M.\,gratefully acknowledges the support of the Austrian Science Fund FWF under the grant ESP~312 ``Defects in low-dimensional TQFTs'', and of the Research Council of Lithuania under the grant S-PD-22-79 ``Lattice systems in topological quantum field theories''.
F.T. was supported by the Heilbronn Institute for Mathematical Research. C.~M.~thanks David Reutter for discussions on the relation of Douglas-Reutter and B\"arenz-Barrett invariants.

\begin{appendix}
\section{Background on Hopf algebras}
\label{sec:hopfback}

\renewcommand{\theDef}{A.\arabic{Def}}

\subsection*{Basic definitions, facts and notations}
\label{sec:hopfnot}

In this section we summarise some basic results on Hopf algebras. For more background, see for instance the books \cite{Mj} by Majid, \cite{Mo} by Montgomery and \cite{R} by Radford.  
We only consider  complex finite-dimensional Hopf algebras.   We denote the  structure maps of a \textbf{Hopf algebra} $H$ by
\begin{align*}
&\text{multiplication:} & &m: H\oo H\to H & &\\
&\text{unit:} & &\eta:\C\to  H, z\mapsto z 1\\
&\text{comultiplication:} & &\Delta: H \to H\oo H, \\
&\text{counit:} & &\epsilon: H\to\C\\
&\text{antipode:} & &S: H\to H
\end{align*}
and recall that they satisfy the defining conditions
\begin{align}\label{eq:alg}
&m\circ (m\oo\id)=m\circ (\id\oo m) & &m\circ (\id\oo\eta)=\id=m\circ (\eta\oo \id)\\
&(\Delta\oo \id)\circ\Delta=(\id\oo\Delta)\circ \Delta
& &(\id\oo \epsilon)\circ \Delta=\id=(\epsilon\oo \id)\circ \Delta\label{eq:coalg}\\
&\Delta\circ m=(m\oo m)\circ \sigma_{23}\circ (\Delta\oo\Delta) & &\epsilon\circ m=\epsilon\oo \epsilon\quad \Delta\circ \eta=\eta\oo \eta \quad \epsilon\circ \eta=\id_\C\label{eq:algcoalg}\\
&m\circ (S\oo\id)\circ\Delta=\eta\circ\epsilon=m\circ (\id\oo S)\circ\Delta\label{eq:antip},
\end{align}
where $\sigma_{23}$ is the map that switches the second and third factor in a tensor product. Conditions \eqref{eq:alg} state that $(H,m,\eta)$ is an associative unital algebra, Conditions \eqref{eq:coalg} that $(H,\Delta, \epsilon)$ is a coassociative counital coalgebra, Conditions \eqref{eq:algcoalg} that $\Delta$, $\epsilon$ are algebra homomorphisms or, equivalently, $m$, $\eta$ are coalgebra homomorphisms, and \eqref{eq:antip} is the defining condition of the antipode.

Recall that the antipode of a Hopf algebra is 
unique
and an anti-algebra and anti-coalgebra map
\begin{align*}
&S\circ m=m^{op}\circ (S\oo S), & &S\circ \eta=\eta,  & &(S\oo S)\circ \Delta=\Delta^{op}\circ S, & &\epsilon\circ S=S.
\end{align*}
Recall also that the antipode of a finite-dimensional complex Hopf algebra $H$ is  a linear isomorphism, see for instance \cite[Th 7.1.14]{R}. It was shown by Larson and Radford \cite{LRa,LRb}, for a textbook reference see \cite[Th 16.1.2]{R}, that it is an involution, $S^2=S\circ S=\id$, if and only if  $H$ is semisimple. 

We use Sweedler notation without summation signs  and write  $\Delta(h)=\low h 1\oo\low h 2$ for  the comultiplication of $H$. We set  $\Delta^{(0)}:=\id$ and $\Delta^{(-1)}:=\epsilon$ and denote an $n$-fold coproduct  by 
$$\Delta^{(n)}(h)=(\Delta\oo \id^{\oo(n-1)})\circ\ldots\circ  (\Delta\oo\id)\circ\Delta(h)=\low h 1\oo\low h 2\oo\ldots\oo \low h {n+1}.$$ 
The coassociativity and counitality axioms and the defining condition on the antipode then read
\begin{align}\label{eq:sweedcond}
&h_{(1)(1)}\oo h_{(1)(2)}\oo h_{(2)}=h_{(1)}\oo h_{(2)(1)}\oo h_{(2)(2)}=:\low h 1\oo \low h 2\oo\low h 3\nonumber\\
&\epsilon(\low h 1)\low h 2=\epsilon(\low h 2)\low h 1=h\nonumber \\
&S(\low h 1)\low h 2=\low h 1 S(\low h 2)=\epsilon(h) 1.
\end{align}
Note that  involutivity of the antipode also implies $
S(\low h 2)\low h 1=\low h 2 S(\low h 1)=\epsilon(h) 1
$ for all $h\in H$. In particular, this holds for any finite-dimensional  semisimple complex Hopf algebra $H$.

For a finite-dimensional complex Hopf algebra $H$ we denote by  $H^{op}$, $H^{cop}$ and $H^{op,cop}$, respectively,  the Hopf algebras with the \textbf{opposite} multiplication, with  the opposite comultiplication and with the opposite multiplication and  comultiplication. Its \textbf{dual Hopf algebra} is denoted $H^*$.

A \textbf{homomorphism of Hopf algebras} from $H$ to $H'$ is a linear map $f:H\to H'$ with
\begin{align}\label{eq:hopfhom}
&m'\circ (f\oo f)=f\circ m & &f\circ \eta=\eta' & &\Delta'\circ f=(f\oo f)\circ\Delta & &\epsilon'\circ f=\epsilon.
\end{align}
Note that these identities imply $S'\circ f=S$, see for instance \cite[Lemma 7.1.3]{R}.

Twisting a Hopf algebra with a 2-cocycle yields a new Hopf algebra with the same comultiplication, counit and unit, but with a changed multiplication and antipode. There is also a dual version, in which the multiplication, unit and counit are unchanged, but the comultiplication and counit are twisted. For an overview, see 
Aljadeff, Etingof, Gelaki and Nikshych \cite{AEGN}.

\begin{Definition} \cite[Def 7.7.1]{R}\label{def:cocyc}\\ 
A {\bf 2-cocycle} for a Hopf algebra $H$ is a bilinear map $\sigma:H\times H\to \C$ that satisfies for all $h,k,l\in H$
$$
\sigma(\low h 1, \low k 1)\,\sigma(\low h 2\low k 2,l)=\sigma(\low k 1,\low l 1)\,\sigma(h, \low k 2\low l 2).
$$
It is called \textbf{normalised}, if $\sigma(h,1)=\sigma(1,h)=\epsilon(h)$ for all $h\in H$. It is called  \textbf{convolution invertible}, if there is a cocycle $\sigma^\inv: H\times H\to \C$ with
$$
\sigma^\inv(\low h 1,\low k 1) \sigma(\low h 2, \low k 2)=\epsilon(h)\epsilon(k)=\sigma^\inv(\low h 2,\low k 2)\sigma(\low h 1,\low k 1)
\qquad\forall h,k\in H.$$
\end{Definition}

\begin{Proposition}\cite[Prop 7.7.3]{R}\label{prop:twist} Let $H$ be a  Hopf algebra and $\sigma: H\times H\to \C$ a convolution invertible 2-cocycle for $H$. Then there is a Hopf algebra $H_\sigma$ with the same unit, counit and comultiplication as $H$ and the multiplication and antipode given by 
\begin{align*}
h\cdot_\sigma k=\sigma(\low h 1,\low k 1) \low h 2\low k 2 \sigma^\inv(\low h 3,\low k 3)\qquad S_\sigma(h)=\sigma(\low h 1, S(\low h 2)) S(\low h 3) \sigma^\inv(S(\low h 4), \low h 5).
\end{align*}
\end{Proposition}

\subsection*{Integrals}\label{subsec:integrals}

In this section we summarise basic facts about integrals of finite-dimensional semisimple complex Hopf algebras. For  more background see for instance the textbooks \cite{Mj,Mo,R}.

\begin{Definition}\label{def:inthopf} Let $H$ be a Hopf algebra.
\begin{compactenum}
\item A \textbf{left integral} in $H$ is an element $\ell\in H$  with $h\cdot \ell=\epsilon(h)\ell$ for all $h\in H$.
\item  A \textbf{right integral} in $H$ is an element $\ell\in H$  with $\ell \cdot h=\epsilon(h)\ell$ for all $h\in H$.
\item An \textbf{integral} in $H$ is an element $\ell\in H$ that is both, a left and a right integral.
\end{compactenum}
A left or right integral $\ell\in H$ is called \textbf{normalisable}, if $\epsilon(\ell)\neq 0$ and \textbf{normalised} if $\epsilon(\ell)=1$.  
\end{Definition}

Left and right integrals coincide for every finite-dimensional semisimple complex Hopf algebra $H$, see for instance \cite[Th 10.3.2]{R}. In this case the integrals in $H$ form a one-dimensional subspace of $H$, and every non-zero integral is normalisable, see \cite[Cor 10.3.3]{R}. This implies in particular that $H$ has a unique normalised integral. Note that this also holds for the dual $H^*$, which is  finite-dimensional semisimple, whenever $H$ is finite-dimensional semisimple, see for instance  \cite[Th 16.1.2]{R}. 

We summarise a number of well-known properties of the integrals of finite-dimensional semisimple complex Hopf algebras, which go back to  Larson and Sweedler \cite{LS}. They are either variants, special cases or direct consequences  of statements in the textbooks \cite{R} and \cite{Mo}, see in particular Exercise 10.1.1, Proposition 10.4.2 and Theorem 16.1.2 in  \cite{R} and Theorem 2.1.3 in \cite{Mo}.

\medskip
\begin{Lemma}\label{lem:haarprops} Let $H$ be a finite-dimensional semisimple complex Hopf algebra with dual $H^*$.
\begin{compactenum}

\item All integrals $\ell\in H$ and $\lambda\in H^*$ satisfy $S(\ell)=\ell$, $S(\lambda)=\lambda$ and for all $h\in H$ and $\alpha\in H^*$
$$
\lambda(\low h 1)\,\low h 2=\lambda(\low h 2)\,\low h 1=\lambda(h) 1\qquad \low\alpha 1(\ell)\,\low\alpha 2=\low\alpha 2(\ell)\,\low\alpha 1=\alpha(\ell) \epsilon.
$$

\item If $\ell\in H$ and $\lambda\in H^*$ are integrals with $\lambda(\ell)=1$, then $\epsilon(\ell)\cdot\lambda(1)=\dim_\C H$
and 
$$
h=\lambda(h\low\ell 1)\, S(\low \ell 2)\qquad \alpha=(\alpha \low\lambda 1)(\ell)\, S(\low\lambda 2)\qquad \forall h\in H,\alpha\in H^*.
$$
\item For any basis $\{b_i\}$ of $H$ and the dual basis $\{\beta^i\}$ of $H^*$ the following are integrals that satisfy $\epsilon(\ell')=\lambda'(1)=\dim_\C H$
$$
\ell'=\sum_i \beta^i(b_{i(1)}) \, b_{i(2)}\in H\qquad\lambda'=\sum_i\beta^i(-\cdot b_i)\in H^*.
$$
\item Coproducts of an integral $\ell\in H$ are invariant under cyclic permutations: for all $n\in\N$
$$
\low \ell 1\oo\low \ell 2\oo\
\ldots\oo \low \ell n=\low\ell n\oo \low\ell 1\oo\ldots\oo \low \ell {n-1}=\ldots=\low \ell 2\oo\ldots\oo \low \ell n\oo \low\ell 1.
$$
\end{compactenum}
\end{Lemma}

We  also require a description of the integrals of a finite-dimensional complex semisimple Hopf algebra $H$ and of its dual $H^*$
in terms of simple representations of $H$. 
For this let $(\rho_i: H\to \End_\C(V_i))_{i\in I}$ be a set of representatives of the isomorphism classes of simple representations of $H$. Suppose  that $0\in I$ and $\rho_0=\epsilon: H\to \C$ is the trivial representation. We then have an algebra isomorphism
\begin{align}\label{eq:algiso}
\phi: H\xrightarrow{\cong}\bigoplus_{i\in I} \End_\C(V_i),\quad h\mapsto (\rho_i(h))_{i\in I},
\end{align}
that can be used to describe the Hopf algebra structures of $H$ and $H^*$ in terms of the simple representations of $H$, see for instance Balsam and Kirillov  \cite[Sec 1]{BK}. In this description the following elements  $\ell\in H$ and $\lambda\in H^*$ are integrals with $\epsilon(\ell)=1$ and $\lambda(\ell)=1$
\begin{align}\label{eq:Haarintdualrep}
\ell=(\delta_{i0} \id_\C)_{i\in I}\qquad\qquad
\lambda=
\sum_{i\in I} \dim_\C(V_i) \,\mathrm{tr}_{V_i}(\rho_i(-)).
\end{align}
That $\ell$ is normalised follows, because the algebra isomorphism \eqref{eq:algiso} identifies the counit with the map $\epsilon:\oplus_{i\in I} \End_\C(V_i)\to \C$, $(f_i)_{i\in I}\mapsto f_0$.
That $\lambda(\ell)=1$  follows directly by inserting $\ell$ into $\lambda$.

\subsection*{Weak Hopf algebras}
\label{subsec:weakhopf}

We summarise the definition of weak Hopf algebras and of their integrals. For more background see the  article by  B\"ohm, Nill and Szlach\'anyi \cite{BNS},
the article  \cite{NV} by Nikshych and Vainerman for a representation theoretical perspective, or Section 7.23 of the textbook \cite{EGNO} for a brief overview. 

A weak Hopf algebra is a generalisation of a Hopf algebra, in which the Hopf algebra axioms are weakened in  a way that is compatible with Hopf algebra duality. 
Just as a Hopf algebra, a weak Hopf algebra is both, an associative unital algebra and a coassociative counital coalgebra. However, the compatibility conditions between the algebra structure and the coalgebra structure are relaxed, and the  antipode axioms change accordingly. 

\begin{Definition}\label{def:weakhopf} \cite[Def 2.1]{BNS}, \cite[Def 7.23.1]{EGNO}  A \textbf{weak Hopf algebra} over $\C$ is a complex vector space $H$ together with an associative unital algebra structure $(H,m,\eta)$ and a coassociative counital coalgebra structure $(H,\Delta,\epsilon)$ such that the following conditions are satisfied.
\begin{compactenum}[(i)]
\item $\Delta(h\cdot k)=\Delta(h)\cdot \Delta(k)$ for all $h,k\in H$,\\[-2ex]
\item $(\Delta\oo\id)\circ \Delta(1)=(\Delta(1)\oo 1)\cdot (1\oo \Delta(1))=(1\oo\Delta(1))\cdot (\Delta(1)\oo 1)$\\
$\epsilon(ghk)=\epsilon(g\low h 1)\epsilon(\low h 2 k)=\epsilon(g\low h 2 )\epsilon(\low h 1 k)$ for all $g,h,k\in H$,\\[-2ex]
\item There is a $\C$-linear map $S: H\to H$ such that for all $h\in H$
\begin{align*}
&m\circ (\id\oo S)\circ \Delta(h)=(\epsilon\oo\id)(\Delta(1)\cdot (h\oo 1))\\
&m\circ (S\oo \id)\circ \Delta(h)=(\id\oo \epsilon)((1\oo h)\cdot \Delta(1))\\
&S(h)=S(\low h 1)\low h 2 S(\low h 3).
\end{align*}
\end{compactenum}
\end{Definition}

It is apparent from the definition that every Hopf algebra is a weak Hopf algebra. It is also apparent that
the dual vector space of a finite-dimensional weak Hopf algebra inherits a  weak Hopf algebra structure, as the axioms in Definition \ref{def:weakhopf} are self-dual. 

There are also notions of  left, right and two-sided  integrals for  weak Hopf algebras, defined analogously to 
Definition \ref{def:inthopf}. The difference is that the counit in Definition \ref{def:inthopf} is replaced by the source and target maps, which are given in terms of the coproduct $\Delta(1)=\low 1 1\oo \low 1 2$
$$
\epsilon_s: H\to H, \; h\mapsto \epsilon(h\low 1 2)\,\low 1 1\qquad\qquad 
\epsilon_t: H\to H,\; h\mapsto \epsilon(\low 1 1 h)\,\low 1 2.
$$

\begin{Definition}\label{def:intweakhopf} \cite[Def 3.1]{BNS}\quad Let $H$ be a weak Hopf algebra.
\begin{compactenum}
    \item A \textbf{left integral} in  $H$ is an element $\ell \in H$ with
    $$
    h\cdot \ell=\epsilon_t(h)\cdot \ell=\epsilon(1_{(1)}h)\, 1_{(2)}\cdot  \ell\qquad \forall h\in H.
    $$
    \item A \textbf{right integral} in  $H$ is an element $\ell \in H$ with
    $$
    \ell\cdot h=\ell \cdot \epsilon_s(h)=\epsilon(h 1_{(2)})\,  \ell \cdot 1_{(1)}\qquad \forall h\in H.
    $$
    \item An \textbf{integral} in  $H$ is an element $\ell\in H$ that is a left and a right integral in $H$.
\end{compactenum}
\end{Definition}

\end{appendix}

\pagebreak

\end{document}